\documentclass[9pt,a4paper]{article}
\usepackage[T1]{fontenc}
\usepackage[francais,english]{babel}
\usepackage[applemac]{inputenc}
\usepackage{pb-diagram}
\usepackage{amscd}
\usepackage[vmargin=3.5cm]{geometry}%
\usepackage{color}%
\usepackage{mdwlist}%
\usepackage{graphicx}%
\usepackage{float}%
\usepackage{makeidx}%
\usepackage{amsmath}%
\usepackage{amssymb}%
\usepackage{amsthm}%
\usepackage{amsfonts}%
\usepackage{mathrsfs}%
\usepackage{bbm}%
\usepackage[all,cmtip]{xy}%
\usepackage{rotate}%
\usepackage{yfonts}%
\usepackage{array}%
\newtheorem{Prop}{Proposition}[section]%
\newtheorem{Theo}[Prop]{Theorem}
\newtheorem{Def}[Prop]{Definition}%
\newtheorem{Cor}[Prop]{Corollary}%
\newtheorem{Lem}[Prop]{Lemma}%
\newtheorem{Conj}[Prop]{Conjecture}%
\newtheorem{Hyp}[Prop]{Assumption}%
\newcommand{\A}{\mathbb A}%

\newcommand{\Afiniq}{\mathbb A^{(\infty)}_{\Q}}
\newcommand{\Afinis}[2]{\mathbb A^{(#1)}_{#2}}
\newcommand{\C}{\mathbb C}%
\newcommand{\fp}{\mathbb F_{p}}%
\newcommand{\Fp}{\mathbb F}%
\newcommand{\E}{\mathbf E}
\newcommand{\F}{\mathbf F}
\newcommand{\G}{\mathbf G}
\newcommand{\gb}{\mathbf g}%
\newcommand{\Gd}{\G_{2}}

\newcommand{\kb}{\mathbf{k}}
\newcommand{\N}{\mathbb N}%
\newcommand{\Q}{\mathbb Q}%
\newcommand{\qp}{\mathbb Q_{p}}%
\newcommand{\qpi}{\mathbb Q_{p,\infty}}%
\newcommand{\ql}{\mathbb Q_{\ell}}%
\newcommand{\Cp}{\mathbb C_{p}}%
\newcommand{\R}{\mathbb R}%
\newcommand{\Z}{\mathbb Z}%
\newcommand{\zp}{\mathbb Z_{p}}%
\newcommand{\zl}{\mathbb Z_{\ell}}%
\newcommand{\Acal}{\mathcal A}%
\newcommand{\Bcal}{\mathcal B}
\newcommand{\Ecal}{\mathcal{E}}%
\newcommand{\Fcal}{\mathcal F}%
\newcommand{\Hcal}{\mathcal H}%
\newcommand{\Lcal}{\mathcal L}%
\newcommand{\Ncal}{\mathcal N}%
\newcommand{\Ocal}{\mathcal O}%
\newcommand{\Pcal}{\mathcal P}%
\newcommand{\Kcal}{\mathcal K}%
\newcommand{\Wcal}{\mathcal W}%
\newcommand{\Xcal}{\mathcal X}%
\newcommand{\Ycal}{\mathcal Y}
\newcommand{\Ecali}{\mathscr E}%
\newcommand{\Icali}{\mathscr I}
\newcommand{\Rcali}{\mathscr R}%
\newcommand{\Xcali}{\mathscr X}%
\newcommand{\Aid}{\mathfrak A}%
\newcommand{\aid}{\mathfrak a}%
\newcommand{\Bid}{\mathfrak B}%
\newcommand{\Cid}{\mathfrak C}%
\newcommand{\cid}{\mathfrak c}%
\newcommand{\pid}{\mathfrak p}%
\newcommand{\mgot}{\mathfrak m}%
\newcommand{\GL}{\operatorname{GL}}%

\newcommand{\SL}{\operatorname{SL}}

\newcommand{\Uni}{\operatorname{U}}%
\newcommand{\GU}{\operatorname{GU}}%

\newcommand{\et}{\operatorname{et}}
\newcommand{\ur}{\operatorname{ur}}
\newcommand{\qppur}{\mathbb{Q}_{p}^{p,\ur}}
\newcommand{\zppur}{\mathbb{Z}_{p}^{p,\ur}}
\newcommand{\zppurhat}{\hat{\mathbb{Z}}_{p}^{p,\ur}}
\newcommand{\diamant}[1]{\langle#1\rangle}%
\newcommand{\somme}[2]{\underset{#1}{\overset{#2}\sum}}%
\newcommand{\produit}[2]{\underset{#1}{\overset{#2}\prod}}%
\newcommand{\produittenseur}[2]{\underset{#1}{\overset{#2}\bigotimes}}%
\newcommand{\union}[2]{\underset{#1}{\overset{#2}\bigcup}}%
\newcommand{\intersection}[2]{\underset{#1}{\overset{#2}\bigcap}}%
\newcommand{\uniondisjointe}[2]{\underset{#1}{\overset{#2}\coprod}}%
\newcommand{\sommedirecte}[2]{\underset{#1}{\overset{#2}\bigoplus}}%
\newcommand{\applicationsimple}[3]{\begin{equation}%
\nonumber%
#1:#2\longrightarrow #3%
\end{equation}}%
\newcommand{\application}[5]{\begin{eqnarray}%
\nonumber%
#1:&#2&\longrightarrow #3\\
\nonumber%
&#4&\longmapsto #5
\end{eqnarray}}%
\newcommand{\suiteexacte}[5]{0\fleche#3\overset{#1}{\fleche}#4\overset{#2}{\fleche}#5\fleche0}
\newcommand{\cardinal}[1]{\mid#1\mid}

\newcommand{\limproj}[1]{\underset{\underset{#1}\longleftarrow}\lim}
\newcommand{\liminj}[1]{\underset{\underset{#1}\longrightarrow}\lim}

\newcommand{\Hom}{\operatorname{Hom}}

\newcommand{\Res}{\operatorname{Res}}

\newcommand{\isom}{\overset{\sim}{\longrightarrow}}
\newcommand{\image}{\operatorname{Im}}
\newcommand{\Ker}{\operatorname{Ker}}

\newcommand{\plonge}{\hookrightarrow}
\newcommand{\matrice}[4]{\begin{pmatrix}#1&#2\\ #3&#4\end{pmatrix}}

\newcommand{\tors}{\operatorname{tors}}%
\newcommand{\rank}{\operatorname{rank}}%
\newcommand{\ord}{\operatorname{ord}}
\newcommand{\mot}{\operatorname{mot}}

\newcommand{\loc}{\operatorname{loc}}
\newcommand{\glob}{\operatorname{glob}}

\newcommand{\tenseur}{\otimes}

\newcommand{\Ltenseur}{\overset{\operatorname{L}}{\tenseur}}
\newcommand{\modulo}{\operatorname{ mod }}
\newcommand{\Spec}{\operatorname{Spec}}
\newcommand{\Spf}{\operatorname{Spf}}
\newcommand{\indicatrice}{\mathbbm{1}}%
\newcommand{\Id}{\operatorname{Id}}%
\newcommand{\Aut}{\operatorname{Aut}}%
\newcommand{\End}{\operatorname{End}}%
\newcommand{\Frac}{\operatorname{Frac}}%
\newcommand{\Tate}{\operatorname{Ta}}%
\newcommand{\Cone}{\operatorname{Cone}}%
\newcommand{\Ind}{\operatorname{Ind}}%
\newcommand{\cInd}{\operatorname{c-Ind}}%
\newcommand{\St}{\operatorname{St}}%
\newcommand{\fleche}{\longrightarrow}%
\newcommand{\croix}{^{\times}}%
\newcommand{\Tor}{\operatorname{Tor}}%
\newcommand{\surjection}{\twoheadrightarrow}%
\newcommand{\rhobar}{\bar{\rho}}%
\newcommand{\chibar}{\bar{\chi}}%

\newcommand{\str}{\operatorname{str}}
\newcommand{\rel}{\operatorname{rel}}
\newcommand{\pr}{\operatorname{pr}}
\newcommand{\cris}{\operatorname{cris}}%
\newcommand{\crys}{\operatorname{crys}}%
\newcommand{\sm}{\operatorname{sm}}%
\newcommand{\dR}{\operatorname{dR}}%
\newcommand{\Fil}{\operatorname{Fil}}
\newcommand{\rig}{\operatorname{rig}}

\newcommand{\Hun}{H^{1}}
\newcommand{\Hunf}{H^{1}_{f}}

\newcommand{\Htilde}{\tilde{H}}
\newcommand{\pitilde}{\tilde{\pi}}

\newcommand{\Htildeun}{\tilde{H}^{1}}
\newcommand{\Htildetildeun}{\tilde{\tilde{H}}^{1}}

\newcommand{\RGamma}{\operatorname{R}\Gamma}%
\newcommand{\BF}{\operatorname{BF}}%
\newcommand{\RGammatilde}{\tilde{\operatorname{R}\Gamma}}%
\newcommand{\RGammaf}{\tilde{\operatorname{R}\Gamma}_{f}}%
\newcommand{\RGammalpha}{\tilde{\operatorname{R}\Gamma}_{\alpha,\alpha}}%
\newcommand{\Det}{\operatorname{{D}et}}%
\newcommand{\Ccont}{C^{\bullet}_{\operatorname{cont}}}%
\newcommand{\Tam}{\operatorname{Tam}}%
\newcommand{\Sel}{\operatorname{Sel}}%

\newcommand{\Tr}{\operatorname{Tr}}%
\newcommand{\tr}{\operatorname{tr}}
\newcommand{\Classe}{\operatorname{Cl}}%
\newcommand{\Reg}{\operatorname{Reg}}%
\newcommand{\Fr}{\operatorname{Fr}}%
\newcommand{\Gal}{\operatorname{Gal}}

\newcommand{\Qbar}{\bar{\Q}}%
\newcommand{\qpbar}{\Qbar_{p}}%
\newcommand{\Kbar}{\bar{K}}%
\newcommand{\Fpbar}{\bar{\mathbb F}}%

\newcommand{\ktilde}{\tilde{k}}%
\newcommand{\Thetatilde}{\tilde{\Theta}}%
\newcommand{\kg}{\kappa}%
\newcommand{\la}{\lambda}
\newcommand{\La}{\Lambda}%
\newcommand{\Lacal}{\Lambda_{\Acal,\infty}}%
\newcommand{\s}{\sigma}%
\newcommand{\p}{\varphi}%
\newcommand{\ptilde}{\tilde{\phi}}%
\newcommand{\Si}{\Sigma}
\newcommand{\hgot}{\mathfrak h}%
\newcommand{\Hecke}{\mathbf{T}}%
\newcommand{\Hs}{\Hecke^{\Sigma}}
\newcommand{\Eul}{\operatorname{Eul}}
\newcommand{\eqdef}{\overset{\operatorname{def}}{=}}

\DeclareFontEncoding{OT2}{}{} 
  \newcommand{\textcyr}[1]{%
    {\fontencoding{OT2}\fontfamily{wncyr}\fontseries{m}\fontshape{n}%
     \selectfont #1}}%
\newcommand{\Sha}{{\mbox{\textcyr{Sh}}}}%
\newcommand{\Nekovar}{Nekov\'a\v{r}}%
\newcommand{\BerZ}{Bernstein-Zelevinsky\ }%
\newcommand{\Iwagood}{Iwasawa-suitable\ }%
\newcommand{\essAIG}{essentially absolutely irreducible and generic\ }
\newcommand{\cl}{\operatorname{cl}}
\numberwithin{equation}{subsection}%
\date{}
\begin{document}%
\title{The Iwasawa Main Conjecture for universal families of modular motives}
\author{Olivier Fouquet \& Xin Wan}%
\maketitle

%
\selectlanguage{english}%

\newcommand{\hord}{\mathfrak h^{\ord}}%
\newcommand{\hdual}{\mathfrak h^{dual}}%
\newcommand{\matricetype}{\begin{pmatrix}\ a&b\\ c&d\end{pmatrix}}%
\newcommand{\Iw}{\operatorname{Iw}}%
\newcommand{\psiwa}{\psi_{\operatorname{Iw}}}%
\newcommand{\Lambdaf}{\mathbf{\Lambda}}
\newcommand{\Hi}{\operatorname{Hi}}
\newcommand{\cyc}{\operatorname{cyc}}
\newcommand{\ab}{\operatorname{ab}}
\newcommand{\can}{\operatorname{can}}%
\newcommand{\Fitt}{\operatorname{Fitt}}%
\newcommand{\Tfiwa}{T_{f,\operatorname{Iw}}}%
\newcommand{\Vfiwa}{V_{f,\operatorname{Iw}}}%
\newcommand{\Afiwa}{A^{*}_{f,\operatorname{Iw}}(1)}%
\newcommand{\Af}{\operatorname{A}}%
\newcommand{\Dunzero}{D_{1,0}}%
\newcommand{\Uun}{U_{1}}%
\newcommand{\Uzero}{U_{0}}%
\newcommand{\Uundual}{U^{1}}%
\newcommand{\Uunun}{U^{1}_{1}}%
\newcommand{\Wdual}{\Wcal^{dual}}
\newcommand{\JunNps}{J_{1,0}(\Ncal, P^{s})}%
\newcommand{\Tatepord}{\Tate_{\pid}^{ord}}%
\newcommand{\Kum}{\operatorname{Kum}}%
\newcommand{\zcid}{z(\cid)}%
\newcommand{\kgtilde}{\tilde{\kappa}}%
\newcommand{\kiwa}{\varkappa}%
\newcommand{\kiwatilde}{\tilde{\varkappa}}%
\newcommand{\Hbar}{\bar{H}}%
\newcommand{\Tred}{T/\mgot T}%
\newcommand{\Riwa}{R_{\operatorname{Iw}}}%
\newcommand{\Oiwa}{\Ocal_{\operatorname{Iw}}}%
\newcommand{\Kiwa}{\Kcal_{\operatorname{Iw}}}%
\newcommand{\Sp}{\Sigma^{p}}%
\newcommand{\Aiwa}{A_{\operatorname{Iw}}}%
\newcommand{\Aiwap}{A_{\operatorname{Iw},\pid}}%
\newcommand{\Siwa}{S_{\operatorname{Iw}}}%
\newcommand{\Liwa}{\Lambda_{\operatorname{Iw}}}%
\newcommand{\liwa}{\lambda_{\operatorname{Iw}}}%
\newcommand{\liwaf}{\lambda(f)_{\operatorname{Iw}}}%
\newcommand{\Viwa}{\mathcal V_{\operatorname{Iw}}}%
\newcommand{\pseudiso}{\overset{\centerdot}{\isom}}%
\newcommand{\pseudisom}{\overset{\approx}{\fleche}}%
\newcommand{\carac}{\operatorname{char}}%
\newcommand{\length}{\operatorname{length}}
\newcommand{\eord}{e^{\ord}}%
\newcommand{\eordm}{e^{\ord}_{\mgot}}%
\newcommand{\hordinfini}{\hord_{\infty}}%
\newcommand{\Mordinfini}{M^{\ord}_{\infty}}%
\newcommand{\hordm}{\hord_{\mgot}}%
\newcommand{\hminm}{\hgot^{min}_{\mgot}}%
\newcommand{\Mordm}{M^{\ord}_{\mgot}}%
\newcommand{\Mtwist}{M^{tw}_{\mgot}}
\newcommand{\Xun}{X_{1}}%
\newcommand{\Xundual}{X^{1}}%
\newcommand{\Xunun}{X^{1}_{1}}%
\newcommand{\Xtw}{X^{tw}}
\newcommand{\Inert}{\mathfrak{In}}%
\newcommand{\Ts}{T_{\Sigma}}
\newcommand{\Tla}{T_{\la}}
\newcommand{\Tlaiwa}{T_{\la,\Iw}}
\newcommand{\Vlaiwa}{V_{\la,\Iw}}
\newcommand{\Tpsi}{T_{\psi}}
\newcommand{\Tpsiwa}{T_{\psi,\Iw}}
\newcommand{\Vpsiwa}{V_{\psi,\Iw}}
\newcommand{\Vla}{V_{\la}}
\newcommand{\rla}{\rho_{\la}}
\newcommand{\Tsx}{T_{\Sigma}^{x_{1}}}
\newcommand{\Vs}{V_{\Sigma}}
\newcommand{\Rs}{R_{\Sigma}}
\newcommand{\Rsr}{R_{\Sigma}(\rhobar)}
\newcommand{\Rsrchi}{R_{\Sigma}^{\chi}(\rhobar)}
\newcommand{\Raid}{R(\aid)}
\newcommand{\Faid}{\Frac(\aid)}%
\newcommand{\Taid}{T(\aid)}
\newcommand{\Vaid}{V(\aid)}
\newcommand{\zaid}{\z(\aid)}
\newcommand{\Hsm}{\mathbf{T}_{\mgot_{\rhobar}}^{\Sigma}}
\newcommand{\Hsmprime}{\mathbf{T}_{\mgot_{\rhobar}}^{\Sigma'}}
\newcommand{\Hsmr}{\mathbf{T}_{\mgot_{\rhobar}}^{\Sigma}}
\newcommand{\Hsmx}{\mathbf{T}_{\mgot_{\rhobar}}^{\Sigma,x_{1}}}
\newcommand{\Hsx}{\mathbf{T}_{\mgot_{\rhobar}}^{\Sigma,x_{1}}}
\newcommand{\Hr}{\mathbf{T}_{\mgot_{\rhobar}}^{\Sigma}}
\newcommand{\Ds}{\Delta_{\Sigma}}
\newcommand{\Tsp}{T_{\psi}}%
\newcommand{\Asp}{A_{\Sp}}%
\newcommand{\Vsp}{V_{\Sp}}%
\newcommand{\SK}{\mathscr{S}}%
\newcommand{\Rord}{R^{\ord}}%
\newcommand{\per}{\operatorname{per}}
\newcommand{\z}{\mathbf{z}}
\newcommand{\vbf}{\mathbf{v}}
\newcommand{\vbfbar}{\bar{\mathbf{v}}}
\newcommand{\vbfd}{{\mathbf{v}^{\vee}}}
\newcommand{\zs}{\mathbf{z}_{\Sigma}}
\newcommand{\zsl}{\mathbf{z}_{\Sigma,\Lambdaf}}
\newcommand{\zsf}{\mathbf{z}(f)_{\Sigma}}
\newcommand{\ztilde}{\tilde{\z}}
\newcommand{\zsprime}{\mathbf{z}_{\Sigma'}}
\newcommand{\ziwa}{\z(f)_{\Iw}}
\newcommand{\zsiwa}{\mathbf{z}_{\Sigma,\Iw}}
\newcommand{\zsfiwa}{\mathbf{z}(f)_{\Sigma,\Iw}}
\newcommand{\zka}{{\mathbf{z}_{\operatorname{Ka}}}}
\newcommand{\Ebarbar}{\bar{\bar{E}}}
\newcommand{\Grsym}{\mathfrak S}
\newcommand{\Gr}{\operatorname{Gr}}
\newcommand{\epsi}{\varepsilon}
\newcommand{\Fun}[2]{F^{{\mathbf{#1}}}_{#2}}
\newcommand{\triv}{\operatorname{triv}}
\newtheorem*{TheoA}{Theorem A}%
\newtheorem*{TheoB}{Theorem B}%
\newtheorem{ConjIMC}{Iwasawa Main Conjecture}[section]%
\newtheorem*{ConjETNC}{Theorem B}
\newcommand{\isocan}{\overset{\can}{\simeq}}
\newcommand{\cusps}{\operatorname{cusps}}
\newcommand{\ad}{\operatorname{ad}}
\newcommand{\Lie}{\operatorname{Lie}}
\makeatletter
\newcommand\@biprod[1]{%
  \vcenter{\hbox{\ooalign{$#1\prod$\cr$#1\coprod$\cr}}}}
\newcommand\biprod{\mathop{\mathpalette\@biprod\relax}\displaylimits}
\makeatother
\newcommand{\Etilde}{\tilde{E}}
\newcommand{\Atilde}{\tilde{A}}
\newcommand{\Btilde}{\tilde{B}}
\newcommand{\Var}{\operatorname{Var}}
\newcommand{\Rep}{\operatorname{Rep}}
\newcommand{\Xbar}{\bar{X}}
\newcommand{\Spa}{\operatorname{Spa}}
\newcommand{\Is}{\operatorname{Is}}
\newcommand{\der}{\operatorname{der}}
\newcommand{\CompCoho}{\Htildeun\left(U^{p},\Ocal\right)_{\mgot_{\rhobar}}}
\newcommand{\CompCoh}{\Htildetildeun\left(U^{p},\Ocal\right)_{\mgot_{\rhobar}}}
\newcommand{\CompCohoC}{\Htildeun_{c}\left(U^{p},\Ocal\right)_{\mgot_{\rhobar}}}
\newcommand{\Uaid}{U_{1}(\aid)}
\newcommand{\Uaidl}{U_{1}(\aid)_{\ell}}
\newcommand{\zt}{\tilde{\z}}
\newcommand{\zaidt}{\tilde{\z}(\aid)}

\abstract{Let $p$ be an odd prime. We prove the cyclotomic Iwasawa Main Conjecture of K.Kato for the motive attached to an eigencuspform $f\in S_{k}(\Gamma_{0}(N))$ with arbitrary reduction type at $p$ under mild assumptions on the residual Galois representation $\rhobar_{f}$. Unde the same hypotheses, we also prove the generalized Iwasawa Main Conjecture for $p$-adic families of modular forms. The Iwasawa Main Conjecture for $f$ is deduced by a limit argument involving fundamental lines from a universal Iwasawa Main Conjecture over the universal deformation space of $\rhobar_{f}$, which itself follows from the cyclotomic Iwasawa Main Conjecture for crystalline eigencuspforms and hence from results on the Iwasawa-Greenberg Main Conjecture for Rankin-Selberg products. The main novel ingredients in our proof are as follows: a new way to study the arithmetic of the Fourier-Jacobi coefficients of Eisenstein series for the group $\Uni(3,1)$, an explicit description of the exponential map in a well-chosen family with prescribed ramification to obtain integral comparisons of various $p$-adic $L$-functions and Selmer modules, Iwasawa theory for the universal zeta elements constructed by K.Nakamura, descent techniques for fundamental lines over the universal regular ring underlying the universal deformation space and ramification properties of the latter over the former.}
\newpage%
\tableofcontents
\section{Introduction}
Let $p>2$ be a prime. Let $f\in S_{k}(\Gamma_{0}(N))$ be a normalized eigencuspform of weight $k\geq2$. Let $\rhobar:G_{\Q}\fleche\GL_{2}(\Fp)$ be a residual representation attached to $f$ with coefficients in a finite extension $\Fp$ of $\fp$. One of the main results of this article is the following.
\begin{Theo}
Assume the following three hypotheses.
\begin{enumerate}
\item\label{NonEisenstein} The $G_{\Q}$-representation $\rhobar$ is absolutely irreducible.
\item\label{Langlands} If $\chibar$ is a character of $G_{\qp}$ and $\chibar_{\cyc}$ is the cyclotomic character modulo $p$, the semisimplification of $\rhobar|G_{\qp}$ is equal neither to $\chibar\oplus\chibar$ nor to $\chibar\oplus\chibar_{\cyc}\chibar$.
\item\label{Steinberg} There exists $\ell\nmid p$ such that $\ell||N$ and such that $\dim_{\Fp}\rhobar^{I_{\ell}}=1$ and $\dim_{\Fp}\rhobar^{G_{\ql}}=0$. 
\end{enumerate}
Then the Iwasawa Main Conjecture (\cite[Conjecture 12.10]{KatoEuler}) is true for for the motive $M(f)$.
\end{Theo}
Assumption \ref{NonEisenstein} means that the modular form is residually non-Eisenstein: the residually Eisenstein case has been treated in recent and forthcoming works of F.Castella, G.Grossi, J.Lee and C.Skinner (\cite{CastellaGrossiLeeSkinner}). Both assumption \ref{NonEisenstein} and assumption \ref{Langlands} are required to analyze the universal deformation space of $\rhobar$ and to appeal to the $p$-adic Langlands correspondence. Assumption \ref{Steinberg}, which is equivalent to the local automorphic representation $\pi(f)_{\ell}$ being the special Steinberg twisted by the unramified character sending $\ell$ to $(-1)\ell^{k/2-1}$ and $\bar{\rho}|_{G_\ell}$ being ramified, plays a technical role in comparing the periods appearing in various $p$-adic $L$-functions. Contrary to most works on the Iwasawa Main Conjecture for modular forms, the eigencuspform $f$ in theorem \ref{TheoIntro} may have arbitrary weight $k\geq2$ and arbitrary ramification at $p$. This is allowed by new techniques which enable us to bypass essential difficulties that prevented previous methods to go beyond the ordinary or low-weight crystalline cases.
\subsection{The Iwasawa Main Conjecture for modular motives}
\subsubsection{Historical introduction}
\paragraph{The Iwasawa Main Conjecture}Let $p>2$ be a prime. Let $\Q_{\infty}/\Q$ be the cyclotomic $\zp$-extension of $\Q$, let $\Gamma_{\Iw}$ be $\Gal(\Q_{\infty}/\Q)$ and let $\Lambda_{\Iw}\eqdef\zp[[\Gamma_{\Iw}]]$ be the classical Iwasawa algebra (\cite{SerreIwasawa}).

Let $M$ be a pure motive over $\Q$. The classical Iwasawa Main Conjecture for $M$ (\cite{KatoHodgeIwasawa,KatoViaBdR}) is a web of conjectures which describe together the $p$-adic variation of the algebraic part of the special values of the $L$-function of $M^{*}(1)\tenseur\chi$ as $\chi$ ranges over the set of Dirichlet characters of order a power of $p$ in terms of $\Lambda_{\Iw}$-adic cohomological invariants of $M$. As such, it can be understood as a conjectural description of the zeroes of a $p$-adic $L$-function (when such an object is known to exist) and hence a $p$-adic variant of the Generalized Riemann Hypothesis, or as a description of the Galois-module structure of certain cohomological invariants of $M$ by means of congruences modulo $p$ between special values of $L(M^{*}(1),\chi,0)$ for various $\chi$.

For the motive $\Q(1)$, the main conjecture was formulated by K.Iwasawa (\cite{IwasawaMainConjecture,MazurWiles}) as a conjectural equality between the ideal generated by the $p$-adic $L$-function $L_{p}$ of Kubota-Leopoldt (\cite{KubotaLeopoldt}) and the characteristic ideal of $\Classe_{\infty}[p^{\infty}]$, the inverse limit on $n$ of the $p$-torsion part of the class group of the ring $\Z[\zeta_{p^{n}}]$ seen as a $\Lambda_{\Iw}$-module. Similar conjectural formulations relating a suitably defined $p$-adic $L$-function and the characteristic ideal of a suitably defined $\Lambda_{\Iw}$-adic Selmer module were found in increasing generality by B.Mazur for ordinary modular abelian varieties (\cite{MazurRational,MazurSwinnertonDyer}), by R.Greenberg for motives with ordinary reduction at $p$ in the sense of \cite{PerrinRiouOrdinaires} (\cite{GreenbergIwasawaRepresentation,GreenbergIwasawaMotives}) and by B.Perrin-Riou for motives with crystalline reduction at $p$ (\cite{PerrinRiouLpadique}).

After the work of S.Bloch and K.Kato (\cite{BlochKato}) and their reformulation in \cite{FontainePerrinRiou}, it has been understood that even for a single motive $M$, the so-called \emph{motivic fundamental line} $\Delta(M)$ of $M$ should play a crucial role in the general statement of conjectures on special values of $L$-functions of motives. Let $M/\Q$ be a motive with coefficients in a number field $L\subset\C$, which we may think of in the context of this introduction as a system of realizations $\{M_{B},M_{\dR},\{M_{\et,\pid}\}_{\pid\in\Spec\Ocal_{L}}\}$ related by comparison isomorphisms and satisfying the Weight Monodromy Conjecture (\cite[3.13.2]{IllusieMonodromie}). If  $M_{\et,\pid}$ is the $\pid$-adic étale realization of $M$ and if $j:\Spec\Q\plonge\Spec\Z[1/p]$ is the natural morphism, we write $\RGamma_{\et}(\Z[1/p],-)$ for $\RGamma_{\et}(\Spec\Q,j_{*}(-))$. The motivic fundamental line of $M$ is a conjectural free $L$-vector space $\Delta(M)$ of dimension 1 defined using motivic cohomology and which satisfies the following properties.
\begin{enumerate}
\item There is a canonical isomorphism $\per_{\C}:\Delta(M)\tenseur_{L}\C\isocan\C$.
\item For any prime $\pid\in\Spec\Ocal_{L}$, there is a canonical isomorphism
\begin{equation}\nonumber
\per_{\pid}\Delta(M)\tenseur_{L}L_{\pid}\isocan\Det^{-1}_{L_{\pid}}\RGamma_{\et}(\Z[1/p],M_{\et,\pid})\tenseur_{L_{\pid}}\Det^{-1}_{L_{\pid}}M_{\et,\pid}(-1)^{+}.
\end{equation}
\end{enumerate}
To $M$ is attached a partial $L$-function $L_{\{p\}}(M,s)$ defined as the formal Euler product
\begin{equation}\nonumber
L_{\{p\}}(M,s)=\produit{\ell\nmid p}{}\frac{1}{\det\left(1-\Fr(\ell)x|M_{\et,\pid}^{I_{\ell}}\right)_{x=\ell^{-s}}}
\end{equation}
which is believed to define a complex $L$-function with an analytic continuation at $s=0$. Denote by $L_{\{p\}}^{*}(M,0)$ the first non-zero term in the Taylor expansion of $L_{\{p\}}(M,s)$ at $s=0$.

Conjectures on special values of $L$-functions then become the following compact statement.
\begin{Conj}\label{ConjTNCintro}
There exists a basis $\z(M)\in\Delta(M)$ called the \emph{zeta element} of $M$ such that 
\begin{equation}\nonumber
\per_{\C}(\z(M)\tenseur1)=L^{*}_{\{p\}}(M^{*}(1),0)
\end{equation}
and such that 
\begin{equation}\nonumber
\Ocal_{L_{\pid}}\cdot\per_{\pid}(\z(M)\tenseur1)=\Det^{-1}_{\Ocal_{L_{\pid}}}\RGamma_{\et}(\Z[1/p],T)\tenseur_{\Ocal_{L_{\pid}}}\Det^{-1}_{\Ocal_{L_{\pid}}}T(-1)^{+}
\end{equation}
for any $G_{\Q}$-stable $\Ocal_{L_{\pid}}$-submodule $T\subset M_{\et,\pid}$.
\end{Conj}
Here we recall that the fact that $\per_{\C}$ and $\per_{\pid}$ are well-defined is in general an integral part of the conjecture.

Suppose now for simplicity of notations that $M/\Q$ is a pure motive with coefficients in $\Q$. Suppose also that $\Xcali$ is $p$-adic family of $G_{\Q}$-representations containing $M_{\et,p}$ (that is to say a rigid analytic $p$-adic space or a $\zp$-scheme whose $\qpbar$-points are Galois representations and such that one such points is $M_{\et,p}$). The advantage of the formulation of conjecture \ref{ConjTNCintro} is that it is amenable to $p$-adic interpolation within the family $\Xcali$. We describe the case of the family of cyclotomic twists of $M$ (\cite{KatoHodgeIwasawa}). Define as in \cite{KatoViaBdR} the \emph{$\Lambda_{\Iw}$-adic fundamental line} $\Delta(M)_{\Iw}$ as the free $\Lambda_{\Iw}$-module of rank 1
\begin{equation}\nonumber
\Delta(M)_{\Iw}=\Det^{-1}_{\Lambda_{\Iw}}\RGamma_{\et}(\Z[1/p],T\tenseur_{\zp}\Lambda_{\Iw})\tenseur_{\Lambda_{\Iw}}\Det^{-1}_{\Lambda_{\Iw}}(T\tenseur_{\zp}{\Lambda_{\Iw}})(-1)^{+}
\end{equation}
for $T$ any $G_{\Q}$-stable $\zp$-module inside the $p$-adic étale realization $M_{\et,p}$ of $M$ (the resulting $\Lambda_{\Iw}$-adic fundamental line is then independent of the choice of $T$ in the sense that there there is a canonical isomorphism between fundamental lines attached to two different such choices). For any finite order character $\chi\in\hat{\Gamma}_{\Iw}$ with values in a finite extension $L_{\chi}/\Q$ and any integer $r\in\Z$, there is then a canonical isomorphism
\begin{equation}\nonumber
\chi_{\cyc}^{r}\chi:\Delta(M)_{\Iw}\tenseur_{\Lambda_{\Iw}}\Ocal_{\chi}\isocan\Det^{-1}_{\Ocal_{\chi}}\RGamma_{\et}(\Z[1/p],T(r)\tenseur\chi)\tenseur_{\Ocal_{\chi}}\Det^{-1}_{\Ocal_{\chi}}(T(r)\tenseur\chi)(-1)^{+}
\end{equation}
where we denote by $\Ocal_{\chi}$ the ring of integers of a finite extension of $\qp$ containing $L_{\chi}$ and where $T(r)\tenseur{\chi}$ is the $\Ocal_{\chi}$-module $T\tenseur_{\zp}\Ocal_{\chi}$ with $G_{\Q}$-action twisted by $\chi_{\cyc}^{r}\chi$.

The Iwasawa Main Conjecture for the motive $M$ is the following statement, which expresses the fact that there should exist a basis of $\Delta(M)_{\Iw}$ which $p$-adically interpolates the zeta elements $\z(M(r)\tenseur\chi)\in\Delta(M(r)\tenseur\chi)$ of the various motives $M(r)\tenseur\chi$ as $r$ ranges over $\Z$ and $\chi$ ranges over the finite order characters of $\Gamma_{\Iw}$.
\begin{Conj}\label{ConjIMCintro}
There exists a basis $\z(M)_{\Iw}\in\Delta(M)_{\Iw}$ called the \emph{$\La_{\Iw}$-adic zeta element} such that for all finite order character $\chi\in\hat{\Gamma}_{\Iw}$ and all integer $r\in\Z$, $\chi_{\cyc}^{r}\chi(\z(M)_{\Iw})$ is the image of the zeta element $\z(M(r)\tenseur\chi)$ of $M(r)\tenseur\chi$ through $\per_{\pid}(-\tenseur1)$. Equivalently, for all such $\chi$ and $r$, $\chi_{\cyc}^{r}\chi(\z(M)_{\Iw})$ is a basis of
\begin{equation}\nonumber
\Det^{-1}_{\Ocal_{\chi}}\RGamma_{\et}(\Z[1/p],T(r)\tenseur\chi)\tenseur_{\Ocal_{\chi}}\Det^{-1}_{\Ocal_{\chi}}(T(r)\tenseur\chi)(-1)^{+}
\end{equation}
and 
\begin{equation}\nonumber
\per_{\C}\left(\per^{-1}_{\pid}\left(\chi_{\cyc}^{r}\chi(\z(M)_{\Iw})\right)\tenseur1\right)=L^{*}_{\{p\}}(M^{*}(1),\chi^{-1},-r).
\end{equation}
\end{Conj}
\paragraph{The case of modular motives}Let $f\in S_{k}(\Gamma_{1}(N))$ be a normalized eingencuspform of weight $k\geq2$ and let $M(f)$ be the pure motive attached to $f$ (\cite{SchollMotivesModular}). Let $E/\qp$ be a finite extension containing the eigenvalues of $f$ and let $\Ocal$ be its ring of integers. We fix $T\subset M(f)_{\et,\pid}$ a $G_{\Q}$-stable $\Ocal$-module and denote by
\begin{equation}\nonumber
\rho_{f}:G_{\Q}\fleche\GL_{2}(\Ocal),\ \rhobar_{f}:G_{\Q}\fleche\GL_{2}(\Fp)
\end{equation}
the corresponding choices of $G_{\Q}$-representations. In \cite{KatoEuler}, the following remarkable result towards conjecture \ref{ConjIMCintro} for the motive $M(f)$ was proved.
\begin{Theo}[\cite{KatoEuler} Theorem 12.4,12.5]\label{TheoKatoIntro}
There exists a non-zero element
\begin{equation}\nonumber
\z(f)_{\Iw}\in\Delta(M(f))_{\Iw}\tenseur_{\Lambda_{\Iw}}\Frac(\Lambda_{\Iw})
\end{equation}
satisfying the following two properties.
\begin{enumerate}
\item Let $1\leq r\leq k-1$ be an integer and $\chi\in\hat{\Gamma}_{\Iw}$ be of finite order. If
\begin{equation}\nonumber
L_{\{p\}}(M(f)^{*}(1),\chi^{-1},-r)\neq0
\end{equation}
then the fundamental line $\Delta(M(f)(r)\tenseur\chi)$ and the period isomorphisms $\per_{\C}$ and $\per_{\pid}$ all exist and satisfy the expected properties. Moreover, $\per_{\pid}^{-1}\left(\chi_{\cyc}^{r}\chi(\z(f)_{\Iw})\right)$ belongs to $\Delta(M(f)(r)\tenseur\chi)$ so that 
\begin{equation}\nonumber
\per_{\C}(\per_{\pid}^{-1}\left(\chi_{\cyc}^{r}\chi(\z(f)_{\Iw})\right)\tenseur1)
\end{equation}
is well-defined and
\begin{equation}\nonumber
\per_{\C}(\per_{\pid}^{-1}\left(\chi_{\cyc}^{r}\chi(\z(f)_{\Iw})\right)\tenseur1)=L_{\{p\}}(M(f)^{*}(1),\chi^{-1},-r).
\end{equation}
\item If $k>2$ or if $M(f)$ has potential good reduction and if the image of $\rho_{f}|G_{\Q(\zeta_{p^{\infty}})}$ contains a subgroup conjugated to $\SL_{2}(\zp)$, then there is an inclusion $\Delta(M(f))_{\Iw}^{-1}\subset\Lambda_{\Iw}\cdot\z(f)_{\Iw}$ of invertible $\Lambda_{\Iw}$-modules inside $\Delta(M(f))_{\Iw}\tenseur_{\Lambda_{\Iw}}\Frac(\Lambda_{\Iw})$. Equivalently, $\z(f)_{\Iw}$ can then be seen as a class in $\Hun_{\et}(\Z[1/p],T\tenseur\Lambda_{\Iw})$ and there is a divisibility of characteristic ideals
\begin{equation}\nonumber
\carac_{\La_{\Iw}}H^{2}_{\et}(\Z[1/p],T\tenseur\La_{\Iw})\mid \carac_{\La_{\Iw}}H^{1}_{\et}(\Z[1/p],T\tenseur\La_{\Iw})/\La_{\Iw}\cdot\z(f)_{\Iw}.
\end{equation}
\end{enumerate}
\end{Theo}
By the Iwasawa Main Conjecture for the modular motive $M(f)$, we henceforth mean the following statement\footnote{Conjecture \ref{ConjIMCweak} is weaker than conjecture \ref{ConjIMCintro}, as it predicts nothing on the behavior of the zeta element at the finitely many $\chi$ where $L(M(f)^{*}(1),\chi^{-1},-r)$ vanishes.}.
\begin{Conj}[\cite{KatoEuler} Conjecture 12.10]\label{ConjIMCweak}
Let $\z(f)_{\Iw}\in\Delta(M(f))_{\Iw}\tenseur_{\Lambda_{\Iw}}\Frac(\Lambda_{\Iw})$ be the basis of theorem \ref{TheoKatoIntro}. Then $\z(f)_{\Iw}$ is a basis of $\Delta(M(f))_{\Iw}$. Equivalently, there is an equality of characteristic ideals
\begin{equation}\nonumber
\carac_{\La_{\Iw}}H^{2}_{\et}(\Z[1/p],T\tenseur\La_{\Iw})=\carac_{\La_{\Iw}}H^{1}_{\et}(\Z[1/p],T\tenseur\La_{\Iw})/\La_{\Iw}\cdot\z(f)_{\Iw}.
\end{equation}
\end{Conj}
K.Kato proved the first statement of theorem \ref{TheoKatoIntro} by introducing a family of distinguished elements $\{\z_{M,N}\in K_{2}(Y(M,N))\}_{M,N}$ where $\z_{M,N}$ is the Steinberg product of two well-defined Siegel units. A delicate reciprocity law links the images of the classes $\z_{M,N}$ through the dual exponential map composed with the Chern class map to products of Eisenstein series and hence to special values of $L$-functions. The second statement then follows from the method of Euler systems, pioneered by F.Thaine, K.Rubin and especially V.Kolyvagin (\cite{Thaine,RubinMain,KolyvaginEuler}) and axiomatized in \cite{PerrinRiouEuler,KatoEulerOriginal,RubinEuler}. In this method, one studies the localization properties of a system of classes $\{\z_{m}\in\Hun_{\et}(\Z[1/p,\zeta_{m}],T)\}_{m\in\N}$ constructed from the $\z_{M,N}$ which mimics the formal behavior of the $L$-function of $M^{*}(1)$ to deduce a crude bound on the length of $H^{2}_{\et}(\Z[1/p],T/p^{n}T)$. Then, one obtains an inclusion of fundamental lines or a divisibility between characteristic ideals as in theorem \ref{TheoKatoIntro} above by $p$-adic interpolation alongside the cyclotomic Iwasawa algebra.

The first proof of the original Iwasawa Main Conjecture for the motive $\Q(1)$ by B.Mazur and A.Wiles (\cite{MazurWiles}) relied on a completely different technique. The starting point of the proof was the observation that the Main Conjecture predicted that a zero of the $p$-adic $L$-function of $M$ (if it exists) should be reflected in the Selmer group. The idea was then to realize this $p$-adic $L$-function as the constant term of the Fourier expansion of a $p$-adic family of automorphic Eisenstein series $\E$ with coefficients in a certain ring $A$. Whenever the $p$-adic $L$-function vanishes modulo a given power $\Pcal^{n}$ of a prime ideal $\Pcal\in\Spec A$, this automorphic form admits a congruence modulo $\Pcal^{n}$ with a $p$-adic family of cuspidal automorphic forms $\F$ with coefficients in $A$. In that situation, the Galois representation attached to $\F$ modulo $\Pcal^{n}$ is isomorphic to the Galois representation attached to $\E$ modulo $\Pcal^{n}$ and is thus reducible, even though the Galois representation attached to $\F$ and with coefficients in $A$ should be irreducible. Ribet's lemma (pioneered in \cite{RibetUnramified}) exploits this reducibility modulo $\Pcal^{n}$ of an irreducible representation to construct a non-trivial extension between Galois representations of lower dimensions, and hence an element in a certain Selmer $A$-module, explaining why the vanishing of a $p$-adic $L$-function has a cohomological counterpart.

In the case of the motive $\Q(1)$, the Eisenstein series $\E$ is a $p$-adic family of genuine Eisenstein series for $\GL_{2}$ and $\F$ is a Hida family of eigencuspforms (\cite{WilesIwasawa}). For modular motives, C.Skinner and E.Urban realized the $p$-adic $L$-function of an eigencuspform $f$ with good ordinary reduction, or more precisely the two-variable $p$-adic $L$-function attached to the automorphic base change of the automorphic representation $\pi(f)$ attached to $f$ to a well-chosen imaginary quadratic number field, as the constant term of a Hida family of Eisenstein series for the reductive group $\Uni(2,2)$ (\cite[Theorem 3.29]{SkinnerUrban}). In this way, they were able to prove the following theorem. (Attentive readers will remark that theorem \ref{TheoSUintro} below is slightly different from \cite[Theorem 3.29]{SkinnerUrban}. First of all, \cite{SkinnerUrban} has the supplementary assumptions $k\equiv 2\modulo p-1$, which was lifted by the second named author in \cite{XinWanHilbert}. In addition, to deduce the following result from \cite{SkinnerUrban}, one would need an argument comparing the canonical period and the Gross period, which might not exist in the relevant literature in a satisfying generality for our purpose. In order to provide a self-contained argument, indications on how to reprove theorem \ref{TheoSUintro} in a way that bypasses this issue are given in section \ref{SubOrdGreenberg}.) 
\begin{Theo}[\cite{SkinnerUrban} Theorem 3.19]\label{TheoSUintro}
Let $f\in S_{k}(\Gamma_{0}(N))$ with $p\nmid N$ be an eigencuspform of weight $k\geq 2$. Assume the following hypotheses on $\rho_{f}$.
\begin{enumerate}
\item The $G_{\Q}$-representation $\rhobar_{f}$ is absolutely irreducible.
\item There exists $\ell||N$, such that $\dim_{\Fp}\rhobar_{f}^{I_{\ell}}=1$ and $\dim_{\Fp}\rhobar_{f}^{G_{\ell}}=0$.
\item The restriction of $\rho_{f}$ to $G_{\qp}$ fits into a short exact sequence of $\Ocal[G_{\qp}]$-modules
\begin{equation}\nonumber
\suiteexacte{}{}{\chi_{1}}{\rho_{f}|G_{\qp}}{\chi_{2}}
\end{equation}
with $\chibar_{1}\neq\chibar_{2}$.
\end{enumerate}
Then there is a divisibility of characteristic ideals
\begin{equation}\nonumber
\carac_{\La_{\Iw}}H^{1}_{\et}(\Z[1/p],T\tenseur\La_{\Iw})/\La_{\Iw}\cdot\z(f)_{\Iw}\mid\carac_{\La_{\Iw}}H^{2}_{\et}(\Z[1/p],T\tenseur\La_{\Iw}).
\end{equation}
Combined with theorem \ref{TheoKatoIntro}\footnote{As is recalled in the proof of lemma \ref{LemEuler}, the assumptions of theorem \ref{TheoSUintro} are enough to deduce theorem \ref{TheoKatoIntro}; see also the discussion after \cite[Theorem 1.4]{SkinnerSteinberg}.}, this yields an equality
\begin{equation}\nonumber
\carac_{\La_{\Iw}}H^{2}_{\et}(\Z[1/p],T\tenseur\La_{\Iw})=\carac_{\La_{\Iw}}H^{1}_{\et}(\Z[1/p],T\tenseur\La_{\Iw})/\La_{\Iw}\cdot\z(f)_{\Iw}.
\end{equation}
and conjecture \ref{ConjIMCweak} is then true for $M(f)$.
\end{Theo}
Together, theorems \ref{TheoKatoIntro} and \ref{TheoSUintro} establish the Iwasawa Main Conjecture for the modular motive $M(f)$ under some conditions on the image of $\rhobar_{f}$ and under the assumption that $\rho_{f}|G_{\qp}$ be reducible, or equivalently that $f$ has ordinary reduction at $p$. In works of the second named author and I.Sprung, the ordinary hypothesis on $\rho_{f}|G_{\qp}$ was relaxed, but only for elliptic curves, that is to say for rational eigencuspforms of weight 2.

In this manuscript, we combine the strength of both approaches to obtain results on the Iwasawa Main Conjecture for eigencuspforms with arbitrary reduction type at $p$ (either crystalline, semistable or potentially semistable of arbitrary even weights from the point of view of $p$-adic Hodge theory, or principal series, Steinberg or supercuspidal from the point of view of automorphic theory).
\subsubsection{The main theorems}
We now state the full version of our main theorems.
\begin{Theo}\label{TheoIntro}
Let $p>2$ be an odd prime and let $\Fp/\fp$ be a finite extension. Let
\begin{equation}\nonumber
\rhobar:G_{\Q}\fleche\GL_{2}(\Fp)
\end{equation}
be an absolutely irreducible and odd Galois representation. Assume that the following assumptions on $\rhobar$ hold.
\begin{enumerate}
\item There is no character $\chibar:G_{\qp}\fleche\Fp\croix$ such that the semisimplification of $\rhobar|G_{\qp}$ is equal to $\chibar\oplus\chibar_{\cyc}\chibar$.
\item There exists $\ell\nmid p$ such that $\rhobar|G_{\ql}$ is a ramified extension
\begin{equation}\nonumber
\suiteexacte{}{}{\mu\chi_{\cyc}^{1-k/2}}{\rhobar|G_{\ql}}{\mu\chi_{\cyc}^{-{k/2}}}
\end{equation}
where $\mu:G_{\ql}\fleche\{\pm1\}$ is the non-trivial unramified quadratic character. 
\end{enumerate}
Let $f\in S_{k}(\Gamma_{0}(N))$ be a normalized eigencuspform of weight $k\geq2$ such that $\rhobar_{f}\simeq\rhobar$. If $\ell||N$, then the Iwasawa Main Conjecture (conjecture \ref{ConjIMCweak}) holds for $M(f)$.
\end{Theo}
We draw the attention of the reader on three features of this result. First, note that the three main hypotheses of theorem \ref{TheoIntro} are requirements solely on $\rhobar_{f}$, or equivalently on the residual representation $\rhobar$. This means that the conclusion of the theorem may be replaced by the statement that the Iwasawa Main Conjecture holds for all eigencuspforms of level $\Gamma_{0}$ which are points of a suitably defined universal deformation ring. Second, theorem \ref{TheoIntro} requires very mild hypotheses on $\rhobar_{f}|G_{\qp}$ and none on $\rho_{f}|G_{\qp}$. This means that we may allow arbitrary bad reduction of $M(f)$ at $p$ or, equivalently but from an automorphic point of view, that $\pi(f)_{p}$ is allowed to be supercuspidal. The generality allowed in the reduction type at $p$ is in contrast with previous works on the topic, which applied only to modular points which were in addition ordinary or crystalline with low weight  at $p$. Given a fixed residual representation $\rhobar$ and a fixed deformation type $\Sigma$ such that any modular point of $\Spec\Rs(\rhobar)[1/p]$ belongs to $S_{k}(\Gamma_{0}(M))$ for some $M\in\N$ exactly divisible by $\ell$, the set of modular points of $\Spec\Rs(\rhobar)[1/p]$ which are ordinary or crystalline with low weight at $p$ has positive codimension in the deformation space $\Spec\Rs(\rhobar)[1/p]$ whereas the set of all modular points is Zariski-dense. In that sense, previous results on the Iwasawa Main Conjecture did not apply to almost all modular points of $\Spec\Rs(\rhobar)$. Theorem \ref{TheoIntro} applies to all of them. Finally, we do not make any assumption on the $\mu$-invariant of $f$. This is in contrast to \cite{EmertonPollackWeston,NakamuraUniversal} which contains results similar to theorem \ref{TheoIntro} under the hypothesis that the $\mu$-invariant of $f$ vanishes, an assumption that is expected to always hold in our setting but which at the moment invariably resisted proof (we describe the role of this hypothesis and the method used to avoid it in the description of the strategy of the proof in section \ref{SubFromTo} below).

In \cite[Conjecture 3.2.2]{KatoViaBdR}, K.Kato formulated what he called the generalized Iwasawa Main Conjecture on the $p$-adic variation of Iwasawa Main Conjectures as a motive varies in a $p$-adic family. Our second main theorem is that this conjecture holds for the universal three-variable power series family of deformations of modular residual representations. In spirit, it says that the Iwasawa Main Conjecture is not only true at each individual modular points of the universal deformation space of a modular residual Galois representation $\rhobar_{f}$ but also that it varies continuously on that space. Because the statement of \cite[Conjecture 3.2.2]{KatoViaBdR} for universal families is itself quite complicated, we refer to section \ref{SubETNCuniv} in the body of the text for the precise meaning of the objects involved in the theorem below.
\begin{Theo}[\cite{KatoViaBdR} Conjecture 3.2.2 for universal $p$-adic families of modular forms]\label{TheoConjUnivWeak}
Let $\rhobar:G_{\Q}\fleche\GL_{2}(\Fp)$ be a residual Galois representation which satisfies the hypotheses of theorem \ref{TheoIntro}. For $\Si$ a finite set of primes, let $\Ts$ be the universal deformation of $\rhobar$ unramified outside $\Si$. We view $\Ts$ as a free module of finite rank over the power-series ring $\Lambdaf\simeq\Ocal[[X_{1},X_{2},X_{3}]]$. Let
\begin{equation}\nonumber
\zs:\Ts(-1)^{+}\fleche\Hun_{\et}(\Z[1/\Sigma],\Ts)
\end{equation}
be the zeta morphism of theorem \ref{TheoNakamura} below. Then $\zs$ induces an isomorphism
\begin{equation}\nonumber
\triv_{\zs}:\Ds(\Ts)\isocan\Lambdaf.
\end{equation}
Equivalently, \cite[Conjecture 3.2.2]{KatoViaBdR} is true for the triple $\left(\Spec\Z[1/\Si],\Lambdaf,\Ts\right)$.
\end{Theo}
In order to help the reader get a better sense of the meaning of this result, we make two remarks. First, if the generalized Iwasawa Main Conjecture is true, then conjecture \ref{ConjIMCweak} is also true except possibly for modular points in a closed subset of positive codimension. In particular, theorem \ref{TheoConjUnivWeak} implies theorem \ref{TheoIntro} for most modular points. Second, when the deformation problem for $\rhobar$ is unobstructed, then the universal deformation ring of $\rhobar$ is a power-series ring and theorem \ref{TheoConjUnivWeak} is the most general statement of the Iwasawa Main Conjecture possible.

We record two corollaries of our main theorem. Write $V$ for $M(f)_{\et,\pid}(k/2)$ and let $T\subset V$ be a $G_{\Q}$-stable $\Ocal$-lattice. For $\ell$ a rational prime, define $\Hunf(G_{\ql},V)$ as in \cite[\S3]{BlochKato} and put
\begin{equation}\nonumber
\Sel_{\Q}(f)=\ker\left(\Hun(G_{\Q},V/T)\fleche\produit{\ell}{}\frac{\Hun(G_{\ql},V/T)}{\image\left(\Hun_{f}(G_{\ql},V)\right)}\right).
\end{equation}
\begin{Cor}\label{CorSelmer}
Let $f\in S_{k}(\Gamma_{0}(N))$ be a normalized eigencuspform satisfying all the hypotheses of theorem \ref{TheoIntro}. Then $\Sel_{\Q}(f)$ is a finite group if and only if $L(f,k/2)\neq0$.
\end{Cor}
Corollary \ref{CorSelmer} is a converse of \cite[Theorem]{KatoEuler}, which established that $\Sel_{\Q}(f)$ is a finite group if $L(f,k/2)\neq0$.
\begin{Cor}
Let $A/\Q$ be an abelian variety of $\GL_{2}$-type of conductor $N$ whose associated weight $2$ cusp form is $f$. Assume that $L(A,1)\neq0$ and that the Galois representation $A[p]$ satisfies all the hypotheses of theorem \ref{TheoIntro}. Then
\begin{equation}\nonumber
v_{p}\left(L(A,1)/\Omega_f\right)=v_{p}\left(\cardinal{\Sha(A/\Q)[p^{\infty}]}\produit{q|N}{}\Tam_{q}(A/\Q)\right).
\end{equation}
Equivalently, the $p$-part of the Birch and Swinnerton-Dyer Conjecture for $A$ holds.
\end{Cor}
Note that the period appearing in the above corollary is the one for the associated modular form. If $A$ is an elliptic curve over $\mathbb{Q}$, then this period is the same as the Neron period up to the Manin constant. In an unpublished work of B.Edixhoven, it is shown that the prime divisors of the Manin constant are equal to $2,3,5$ or $7$. We also note that it would be very interesting to compare and combine the results of this manuscript with the remarkable results obtained in \cite{ColmezWang}.
\subsection{Outline of the proof}
The proof of theorems \ref{TheoIntro} and \ref{TheoConjUnivWeak} rely extensively on techniques of $p$-adic interpolation on various $p$-adic families of Galois representations, including the universal deformation ring $\Rs(\rhobar)$ of the residual representation $\rhobar$.

Under our hypotheses, $\Rs(\rhobar_{f})$ is known to be a reduced, complete intersection, local ring of relative dimension 3 over $\zp$ isomorphic to a $p$-adic Hecke algebra $\Hsmr$. Moreover, $\Spec\Hsmr[1/p]$ contains a Zariski-dense subset $\Xcali^{\sm}$ which contains the set of classical points corresponding to classical eigencuspforms satisfying strong supplementary properties at $p$: for $x\in\Xcali^{\sm}$, $\rho_{x}|G_{\qp}$ is either ordinary or crystalline and short (depending on whether $\rhobar|G_{\qp}$ is reducible or irreducible). As it is a complete intersection local ring of relative dimension 3, $\Rs(\rhobar)\simeq\Hsmr$ is a finite, free $\Lambdaf$-module for $\Lambdaf$ a power-series ring in three variables. It turns out that the existence of $\Xcali^{\sm}$ suggests a preferential choice of such a $\Lambdaf$-structure, that is to say an especially convenient explicit description of the morphism $\Lambdaf\plonge\Hsmr$. The map
\begin{equation}\nonumber
\xymatrix{
\Spec\Hsmr\ar[d]\\
\Spec\Lambdaf
}
\end{equation}
corresponding to this $\Lambdaf$-structure may well be ramified in general, and may be so even at points corresponding to classical eigencuspforms, but it is not ramified at the points of $\Spec\Lambdaf$ below the points in $\Xcali^{\sm}$.

The general strategy of the proof of conjecture \ref{ConjIMCweak} for a modular point $\rho_{f}\in\Spec\Hsmr[1/p]$ is then as follows. First, we prove the two-variable Greenberg-Iwasawa Main Conjecture for the Rankin-Selberg $L$-function attached to a point in $\Xcali^{\sm}$ by realizing the $p$-adic $L$-function involved in this conjecture as the constant term of a Hida family of Eisenstein series for the reductive group $\Uni(3,1)$. Second, we show that this two-variable main conjecture for $x\in\Xcali^{\sm}$ implies conjecture \ref{ConjIMCweak} for the same $x$ (this requires a delicate $p$-integrality property of $p$-adic periods, as the $p$-adic period appearing in the interpolation property of the Greenberg-Iwasawa Main Conjecture is not the same as the one involved in conjecture \ref{ConjIMCweak} or, to phrase the problem differently, since special values of classical points on the cyclotomic line are not interpolated by the Greenberg Rankin-Selberg $p$-adic $L$-function). Third, we use an Euler system argument tho show that conjecture \ref{ConjIMCweak} holds at every point except those in a set of large codimension containing the points where $\Hsmr$ is ramified over $\Lambdaf$ (at this stage, the set of points for which conjecture \ref{ConjIMCweak} has not yet been established may thus very well contain the point corresponding to our $\rho_{f}$ of interest). Finally, a limit argument relying on the interpolation properties of fundamental lines and zeta elements on the universal deformation ring allows us to conclude that it holds for $\rho_{f}$ as well.

We now explain each step in more details.

\subsubsection{Classical points in $\Xcali^{\sm}$}\label{SubIntroXsm}For simplicity of exposition, we assume in this introduction that  $\rhobar|G_{\qp}$ is absolutely irreducible and that $x\in\Xcali^{\sm}$ is attached to a classical eigencuspform $f$ such that $\rho_{f}|G_{\qp}$ is a crystalline representation in the image of the Fontaine-Laffaille functor (\cite{FontaineLaffaille}). The general strategy of the proof here is similar to that of previous works of the second named author in the case of supersingular elliptic curves (\cite{XinWanIMC}): the Iwasawa Main Conjecture for the motive $M(f)$ is deduced from the two-variable Greenberg-Iwasawa Main Conjecture for Rankin-Selberg $L$-functions as studied in \cite{XinWanRankinSelberg}.

Choose a quadratic imaginary field ${\Kcal}/\mathbb{Q}$ such that $p$ splits as $v_0\bar{v}_0$. Let ${\Kcal}_\infty/\Kcal$ be the $\mathbb{Z}_p^2$-extension of ${\Kcal}$. Let $\Gamma_{\Kcal}$ be the Galois group $\Gal({\Kcal}_\infty/{\Kcal})$ and let $\La_{\Kcal}$ be the Iwasawa algebra $\zp[[\Gamma_{\Kcal}]]$. For simplicity of notations, we assume in this introduction that the relevant objects are defined over the two-variable power-series ring $\La_{\Kcal}$ even though in fact an extension of the coefficients to the ring of integers of the maximal unramified extension of a finite extension of $\qp$ might be required. Write $\Gamma^\pm$ for the rank one over $\mathbb{Z}_p$ submodules of $\Gamma_{\Kcal}$ on which the action of the complex conjugation $c$ is given by $\pm1$. If $g$ is a modular form with complex multiplication by ${\Kcal}$, then it is ordinary at $p$. If in addition the weight of $g$ is greater than the weight of $f$, then the Galois representation attached to $f\tenseur g$ satisfies Panchishkin's condition (\cite{GreenbergMotives}) so the Iwasawa theory of the Rankin-Selberg product $f\otimes g$ is strongly similar to the theory of Iwasawa theory for ordinary forms. Accordingly, the statement and proof of the Iwasawa Main Conjectures for the seemingly more complicated object $f\tenseur g$ is more accessible than their counterparts for $f$. Moreover, the Rankin-Selberg $L$-function $L(f\tenseur g,s)$ is related to the original $L$-function of $f$. In that case, we deduce results on the Iwasawa theory of the original form $f$ from the Iwasawa theory of Rankin-Selberg $L$-functions developed in \cite{XinWanRankinSelberg}.

We first recall the setting of the two-variable Iwasawa Main Conjecture for Rankin-Selberg products. Let $\Kcal\subset_{f} F\subset\Kcal_{\infty}$ be a finite subextension and let $v$ be a finite place of $\Ocal_{F}$. The Greenberg local condition $\Hun(G_{F_{v}},A)\subset\Hun(G_{F_{v}},A)$ at $v$ is defined to be
\begin{equation}\nonumber
\image\left(\Hun(G_{F_{v}}/I_{v},V^{I_{v}})\fleche\Hun(G_{F_{v}}/I_{v},A^{I_{v}})\right)\subset\Hun(G_{F},A)
\end{equation}
if $v\nmid p$, to be $\Hun_{\Gr}(G_{F_{v}},A)$ if $v|v_{0}$ and to be 0 if $v|\bar{v}_{0}$. The {Greenberg Rankin-Selberg Selmer group} of $f$ is the $\La_{\Kcal}$-module
\begin{equation}\nonumber
\Sel^{\Gr}_{\Kcal}(f)\eqdef\liminj{\Kcal\subset_{f}F\subset\Kcal_{\infty}}\ker\left(\Hun(G_{F,\Sigma},A)\fleche\sommedirecte{v\in\Sigma}{}\Hun(G_{F_{v}},A)/\Hun_{\Gr}(G_{F_{v}},A)\right).
\end{equation}
We write
\begin{equation}\nonumber
X^{\Gr}_{\Kcal}(f)\eqdef\Hom\left(\Sel^{\Gr}_{\Kcal}(f),\qp/\zp\right)
\end{equation}
for the Pontrjagin dual of $\Sel^{\Gr}_{\Kcal}(f)$. On the analytic side, the object corresponding to $X^{\Gr}_{\Kcal}(f)$ is the Greenberg Rankin-Selberg $p$-adic $L$-function $\Lcal^{\Gr}_{\Kcal}(f)\in\Frac\left(\La_{\Kcal}\right)$ constructed in \cite{EischenWan} which interpolates certain special values of the Rankin-Selberg $L$-functions $L(f\tenseur\xi,s)$ attached to certain characters of $\Gamma_{\Kcal}$ (see definition \ref{DefGreenbergRankinSelberg} and proposition \ref{EischenWan} for details). The two-variable Greenberg-Iwasawa Rankin-Selberg Main Conjecture is then the following statement.
\begin{Conj}[Greenberg Main Conjecture]\label{ConjGreenbergIntro}
The $\La_{\Kcal}$-module $X^{\Gr}_{\Kcal}(f)$ is torsion, $\Lcal^{\Gr}_{\Kcal}(f)$ belongs to $\La_{\Kcal}$ and there is an equality of ideals
\begin{equation}\nonumber
\carac_{\La_{\Kcal}}\left(X^{\Gr}_{\Kcal}(f)\right)=\left(\Lcal^{\Gr}_{\Kcal}(f)\right).
\end{equation}
\end{Conj}
As we mentioned already, conjecture \ref{ConjGreenbergIntro} is (relatively) more tractable than conjecture \ref{ConjIMCweak} when the Galois representation attached to $f\tenseur g$ satisfies Panchishkin's condition, that is to say when the weight of $g$ is greater than the weight of $f$. As conjecture \ref{ConjIMCweak} corresponds to the choice of $g$ trivial, it \emph{cannot} follow obviously from conjecture \ref{ConjGreenbergIntro}. Nevertheless, when $f$ is of weight 2, it is proved in \cite{XinWanRankinSelberg,CLW} under some hypotheses that the inclusion of ideals
\begin{equation}\nonumber
\carac_{\La_{\Kcal}}\left(X^{\Gr}_{\Kcal}(f)\right)\subset\left(\Lcal^{\Gr}_{\Kcal}(f)\right)
\end{equation}
holds after inverting $p$. In \cite{XinWanIMC}, the second named author then developed some local $\pm$-theory in a style similar to S.Kobayashi and B.D.Kim and used the explicit reciprocity law for Beilinson-Flach element and Poitou-Tate exact sequence to deal with the remaining power of $p$ and to deduce conjecture \ref{ConjIMCweak} from conjecture \ref{ConjGreenbergIntro}.

Our main result regarding that topic is the following generalization to arbitrary even weight (theorem \ref{TheoGreenberg} in the body of the manuscript).
\begin{Theo}\label{TheoGreenbergIntro}
Let $f\in S_{k}(\Gamma_{0}(N))$ be an eigencuspform of even weight $k$ satisfying the following assumptions.
\begin{enumerate}
\item $\rhobar_{f}|G_{\Kcal}$ is absolutely irreducible.
\item The local representation $\rho_{f}|G_{\qp}$ is crystalline. Moreover assume either $\bar{\rho}_f|_{G_{\qp}}$ is absolutely irreducible or $f$ is ordinary at $p$.
\item There exists $q||N$ (in particular $q\nmid p$) which is not split in $\Kcal$.
\item If $q|N$ is not split in $\Kcal$, then $q||N$.
\suspend{enumerate}
Then the following inclusion of ideals of $\Ocal^{\ur}[[\Gamma_{\Kcal}]]$ holds
$$\carac_{\mathcal{O}^{{\ur}}[[\Gamma_\mathcal{K}]]}(X_{\Kcal}^{\Gr}(f)
\otimes_{\mathcal{O}}\mathcal{O}^{{\ur}})\subseteq(\Lcal_{\Kcal}^{\Gr}(f))$$
up to height-one primes which are pullbacks of primes in $\mathcal{O}[[\Gamma^+]]$. Assume in addition that the following assumption holds.
\resume{enumerate}
\item If $\ell|N$ is not split in $\mathcal{K}$, then $\ell$ is ramified in $\Kcal$ and $\pi(f)_\ell$ is a special Steinberg representation twisted by $\chi_{{\ur}}$ for $\chi_{{\ur}}$ the unramified character sending $\ell$ to $(-1)\ell^{\frac{k}{2}-1}$.
\end{enumerate}
Then
$$\carac_{\mathcal{O}^{{\ur}}[[\Gamma_\mathcal{K}]]}(X_{\Kcal}^{\Gr}(f)
\otimes_{\mathcal{O}}\mathcal{O}^{{\ur}})\subseteq(\Lcal_{\Kcal}^{\Gr}(f))$$
holds.
\end{Theo}

To prove theorem \ref{TheoGreenbergIntro} for modular forms of general even weight, the difficulty is twofold. First of all, the proof of the Greenberg-Iwasawa Main Conjecture in \cite{XinWanIMC} when $f$ has weight 2 relies on a computation of the Fourier-Jacobi expansion of an Eisenstein series. When $f$ has higher weight, one needs to compute the Fourier-Jacobi expansion for vector valued Eisenstein series, which seems formidable. Secondly we do not have in general an explicit enough local theory as in the $\pm$ case (except in the special case when $a_p=0$, and for elliptic curves over $\mathbb{Q}$ but $a_p\not=0$ by I.Sprung), while the work \cite{XinWanIMC} used such a theory in a crucial way.

To prove this theorem, we use the full strength of the joint work \cite{EischenWan} of the second named author and E.Eischen. In this work, vector-valued Klingen Eisenstein families on $\Uni(3,1)$ are constructed from pullbacks of nearly holomorphic Siegel Eisenstein series on $\Uni(3,3)$. Combining the earlier works \cite{XinWanIMC,XinWanRankinSelberg} on the $p$-adic property for Fourier-Jacobi coefficients with the general theory of T.Ikeda (\cite{IkedaFourierJacobi}) shows that the Fourier-Jacobi coefficient of nearly holomorphic Siegel Eisenstein series can be expressed as finite sums of products of Eisenstein series and theta functions on the Jacobi group containing $\Uni(2,2)$ (see lemma \ref{LemCoreIkeda} below\footnote{We thank T.Ikeda for showing us a simple argument using lowest weight representations giving the proof of this lemma.}). In the scalar valued case the Siegel Eisenstein series on $\Uni(3,3)$ is holomorphic, and the local Fourier-Jacobi coefficient at the Archimedean place can be explicitly expressed as \emph{the} product of a Siegel section on $\Uni(2,2)$ and a Schwartz function. In general, it is difficult to compute Archimedean local Fourier-Jacobi coefficients explicitly. We fix one Archimedean weight and vary the $p$-adic nebentypus in families. Instead of computing the Archimedean Fourier-Jacobi integral, we can use a conceptual argument to factor out a finite sum of Archimedean integrals out of it and prove the factor is non-zero. After this, we can apply Ichino's triple product formula to prove that a certain Fourier-Jacobi coefficient is co-prime to the $p$-adic $L$-function we study. This also allows us to remove the square-free conductor assumption for $f$ in the previous works of the second named author.

With theorem \ref{TheoGreenbergIntro} in hand, we prove the following theorem, which establishes conjecture \ref{ConjIMCweak} at points in $\Xcali^{\sm}$.
\begin{Theo}\label{TheoCrysIntro}
Let $f\in S_{k}(\Gamma_{0}(N))$ be a normalized eigencuspform of weight $k\geq2$ satisfying the following hypotheses. 
\begin{enumerate}
\item The $G_{\qp}$-representation $\rho_{f}|G_{\qp}$ is crystalline (equivalently $p\nmid N$) and short.
\item The local residual representation $\rhobar_{f}|G_{\qp}$ is absolutely irreducible.
\item There exists $\ell||N$ such that $\dim_{\Fp}\rhobar^{I_{\ell}}=1$ and $\dim_{\Fp}\rhobar^{G_{\ql}}=0$. 
\end{enumerate}
Then conjecture \ref{ConjIMCweak} holds for $M(f)$.
\end{Theo}
The proof of theorem \ref{TheoCrysIntro} is quite different from the one in \cite{XinWanIMC}. We avoid the search for a two-variable integral local theory analogous to the $\pm$ theory. Our main innovation is to study unramified local Iwasawa theory along an appropriately chosen $\zp$-line in $\Spec\Lambda_{\Kcal}$. Suppose first that the image of $\rho_f$ contains $\SL_{2}(\zp)$, in which case we know an upper bound on the Selmer group by the works of Kato. After carefully studying the control theorem, proving the full equality amounts to the computation of the cardinality (which is finite) of the Selmer group for $T_f(-\frac{k-2}{2})$ twisted by some appropriately defined generic finite order character of $\Gamma_{\Iw}$. In order to do so, we consider deformations of the corresponding Galois representation along a one-variable family which corresponds to the $\mathbb{Z}_p$-extension of ${\Kcal}$ that is totally ramified at $\bar{v}_0$ and unramified at $v_0$. In this family, we do have a nice integral local theory at $v_0$, which we develop in Section \ref{Unramified Iwasawa Theory} using $p$-adic Hodge theory. Roughly speaking, the integral local subspace of interest is generated by the image of the canonical differential class through the exponential map, the crucial point being the boundedness of the denominators involved in the explicit description of the exponential map (proposition \ref{boundedness} and corollary \ref{specialization}). This local subspaces replaces the $\pm$-theory and is enough to prove the sought-for bound along this one dimensional family (it seems to us this one-dimensional line is the only choice of $p$-adic family containing our original point of interest for which the argument works). We note that a deformation-theoretic argument is necessary in this approach: we can not only look at one point on the cyclotomic line as we are unable to rule out the possibility that the Greenberg type $p$-adic $L$-function is identically $0$ on the cyclotomic line.

The shortcoming of the approach as described is that it does not prove anything without the assumption that the image of $\rho_f$ contains $\SL_{2}(\zp)$ since in this case Kato only proved the upper bound for Selmer groups after inverting $p$. Nevertheless, we can still prove the Main Conjecture after inverting $p$. Our idea is to use the analytic Iwasawa theory of J.Pottharst (\cite{JayAnalytic}) and Iwasawa theory for $(\varphi,\Gamma)$-modules (upgraded to a two variable setting), in the context of \Nekovar-Selmer complexes (\cite{SelmerComplexes}). It turns out that Pottharst's trianguline-ordinary theory and its more flexible analytic Iwasawa theory setting allow similar proofs as in the ordinary case. In the two-variable setting there are subtleties to take care of -- for example there is a finite set of height one primes where the regulator map vanishes. Also we need to compare differently constructed analytic $p$-adic $L$-functions: although the two $p$-adic $L$-functions in question agree on all arithmetic points, this does not suffice to uniquely determine them as analytic functions due to the growth conditions.

We prove the theorem first under the assumption that the Satake parameters of $f$ at $p$ are distinct. This assumption on the Satake parameter is conjecturally automatic. We finally remove this assumption in Section \ref{irregular} by a trick of comparing with the three variable $p$-adic $L$-function for Bianchi modular forms constructed by A.Betina and C.Williams (\cite{BetinaWilliams}).
\paragraph{Remark:}When we do not know that the image of $\rho_{f}$ contains a subgroup conjugated to $\SL_{2}(\zp)$, our argument actually shows one divisibility for powers of $p$. More precisely, it shows that the strict Selmer group is bounded below by the index of $\z(f)_{\Iw}$ in the integral Iwasawa cohomology in the case when the zeta element is indeed inside integral Iwasawa cohomology (if not then this statement is empty).

\subsubsection{From $\Xcali^{\sm}$ to all classical points}\label{SubFromTo}We continue our outline of the proof of theorem \ref{TheoIntro}. The second main step is to state a version of the Iwasawa Main Conjecture with coefficients in $\Rs(\rhobar_{f})$. The fundamental idea, which has been known since \cite{KatoViaBdR}, is to express the Iwasawa Main Conjecture as a description of the image of the fundamental line through a specific zeta morphism. Defining a fundamental line over the full universal deformation space which specializes to the classical fundamental lines at classical points is well-known (\cite{FouquetX}), the difficulty therefore lies in the definition of the zeta morphism. Such a zeta morphism has been recently constructed over the universal deformation space by K.Nakamura (\cite{NakamuraUniversal}) using completed cohomology, the Bernstein derivative functor and the results of V.Paskunas on the explicit description of the so-called Montréal functor of P.Colmez (\cite{ColmezFoncteur,PaskunasMontreal}). The key property of this zeta morphism is that it specializes to the zeta morphism constructed in \cite{KatoEuler} at classical points of the universal deformation space. For our purposes, we need a slight strengthening of these results which relate the zeta morphism defined over the full universal deformation space to the ones defined only on irreducible components thereof but which keep track of the Euler factors at primes of bad reduction. This requires a technical result on newvectors in $p$-adic families of automorphic representations of $\GL_{2}(\ql)$ ($\ell\nmid p$) which is certainly well-known to experts but for which we provide a proof in appendix \ref{AppendixCompleted} for the convenience of the reader.

Let $\Hsmr$ be the $p$-adic Hecke algebra isomorphic to the universal deformation ring $\Rs(\rhobar_{f})$, let $\Raid$ be a quotient of $\Hsmr$ by a minimal prime ideal $\aid$, and let $\la(f):\Raid\fleche\Oiwa$ be a modular point of $\Raid$. We write $\Ts$ for the universal Galois representation with coefficients in $\Hsmr$, $\Taid$ for $\Ts\tenseur_{\Hsmr}\Raid$ and $T(f)_{\Iw}$ for $\Taid\tenseur_{\Raid,\la(f)}\Oiwa$. Let $\Ds(\Ts),\Delta(\Taid)$ and $\Delta(f)_{\Iw}$ the fundamental lines of the Galois representations $\Ts,\Taid$ and $T(f)_{\Iw}$ respectively. Finally, let $Q(\Hsmr)$ be the total quotient ring of the reduced ring $\Hsmr$. The existence of the universal zeta morphism $\zs$ of K.Nakamura can be reformulated as the existence of an isomorphism
\begin{equation}\label{EqFunIsomIntro}
\triv_{\zs}:\Ds\isom\frac{x}{y}\Hsmr
\end{equation}
between the fundamental line $\Ds$ and an invertible module $\frac{x}{y}\Hsmr$ inside $Q(\Hsmr)$ which specializes to the isomorphism 
\begin{equation}\nonumber
\triv_{\zsfiwa}:\Ds(f)\isom\frac{x(f)}{y(f)}\Oiwa
\end{equation}
constructed in \cite[Theorem 12.5]{KatoEuler} (recalled as theorem \ref{TheoKato} below) at each classical point. The universal Iwasawa Main Conjecture is then the statement that $\triv_{\zs}$ is an isomorphism 
\begin{equation}\nonumber
\triv_{\zs}:\Ds\isom\Hsmr\subset Q(\Hsmr)
\end{equation}
or equivalently that the invertible module $\frac{x}{y}\Hsmr$ of \eqref{EqFunIsomIntro} is $\Hsmr$ itself. Likewise, the existence of the universal morphism $\z(\aid)$ with coefficients in $\Raid$ can be reformulated as the existence of an isomorphism 
\begin{equation}\nonumber
\triv_{\z(\aid)}:\Delta(\Taid)\isom\frac{x(\aid)}{y(\aid)}\Raid
\end{equation}
for some invertible ideal $\frac{x(\aid)}{y(\aid)}\Raid$ inside $\Frac(\Raid)$ and the Iwasawa Main Conjecture with coefficients in $\Raid$ is the statement that $\triv_{\z(\aid)}$ is an isomorphism 
\begin{equation}\nonumber
\triv_{\z(\aid)}:\Delta(\Taid)\isom\Raid.
\end{equation}
A crucial step in establishing theorem \ref{TheoIntro} is the proof that, independently of the truth of the various Iwasawa Main Conjectures, there exists a pair of commutative diagrams
\begin{equation}\label{EqDiagFunIntro}
\xymatrix{
\Ds(\Ts)\ar[d]_{-\tenseur_{\Hsmr}\Raid}\ar[rr]^{\triv_{\zs}}&&\frac{x}{y}\Hsmr\ar[d]^{\lambda}\\
\Delta(\Taid)\ar[rr]^{\triv_{\z(\aid)}}&&\frac{x(\aid)}{y(\aid)}\Raid
}
\end{equation}
and
\begin{equation}\label{EqDiagFunIntro2}
\xymatrix{
\Delta(\Taid)\ar[rr]^{\triv_{\z(\aid)}}\ar[d]_{-\tenseur_{\Raid,\la_{f}}\Oiwa}&&\frac{x(\aid)}{y(\aid)}\Raid\ar[d]^{\lambda_{f}}\\
\Delta(f)_{\Iw}\ar[rr]^{\triv_{\z(f)_{\Iw}}}&&\frac{x(f)}{y(f)}\Oiwa.
}
\end{equation}
A way to understand this statement is to interpret it as a perfect control theorem for fundamental lines. In particular, we note that such results seem to be inaccessible for characteristic ideals of Selmer modules, not only because such objects are not defined when the ring of coefficients is not known to be normal (as is the case for $\Hsmr$ and $\Raid$) but also because they would require the knowledge that no non-zero pseudo-null submodules may occur in these Selmer modules, which in turn is expected only when the Selmer modules is related to a $p$-adic $L$-function, and thus only under strong supplementary hypotheses on the image of the decomposition group at $p$ through the Galois representation. Another candidate in the formulation of the Iwasawa Main Conjectures which is defined more generally is the Fitting ideal of the Selmer module. The formation of Fitting ideals commutes with arbitrary change of coefficients but unfortunately Fitting ideals are not principal ideals in general, whereas the proof of the commutativity of the diagrams \eqref{EqDiagFunIntro} and \eqref{EqDiagFunIntro2} crucially relies on the properties of invertible modules. Consequently, it seems to us that this part of our proof relies on the use of fundamental lines, that is to say the images of perfect complexes through the determinant functor of \cite{MumfordKnudsen}: only through this choice do we have invertible ideals whose formation commutes with arbitrary base change of coefficients. In that respect, it might be instructive to compare our results with for instance \cite[Theorem 1.4]{NakamuraUniversal}. From the perspective of this manuscript, the proof of that result amounts to the compatibility of the fundamental line with specialization at the maximal ideal under the hypothesis that the $\mu$-invariant does not vanish. The vanishing of the $\mu$-invariant appears however to be a mysterious and hard problem. Working consistently with fundamental lines and relying on the commutativity of the diagrams \eqref{EqDiagFunIntro} and \eqref{EqDiagFunIntro2} (and their generalizations to other specializations maps beside $\Hsmr\fleche\Raid$ and $\Raid\fleche\Oiwa$) allows us to remove the hypothesis that the $\mu$-invariant vanishes; indeed this was one our main motivation for working within this framework.

Ideally, the proof of the Iwasawa Main Conjecture would then proceed in the following way: using the fact that the image of $\triv_{\z(f)_{\Iw}}$ is as expected by the Iwasawa Main Conjecture (namely $\Oiwa$) for suitably chosen points in $\Xcali^{\sm}$, we would prove by descent that similarly the images of $\triv_{\zs}$ and $\triv_{\z(\aid)}$ are respectively $\Hsmr$ and $\Raid$ as expected, or equivalently that the Iwasawa Main Conjecture holds over these spaces. Then, specializing this time to a classical point outside $\Xcali^{\sm}$, we would conclude that it holds for all classical points. This strategy, however, does not quite work as described. The problem is that its first step is to deduce a universal Iwasawa Main Conjecture from the knowledge of Iwasawa Main Conjectures at given points (here the set of points considered is $\Xcali^{\sm}$ but the objection below applies to any set of classical points). The rings $\Hsmr$ and $\Raid$, however, are not known to be factorial. As an invertible module over a non-factorial ring can very well be non-trivial even if it becomes trivial at all classical points, we may not prove a universal main conjecture in this way. To bypass this difficulty, we view $\Ts$ as a finite, free module of $\Lambdaf$ and formulate the Iwasawa Main Conjecture over this ring. As $\Lambdaf$ is regular, the descent argument sketched above can be carried out with coefficients in $\Lambdaf$. However, a new problem now arises in that proving the Iwasawa Main Conjecture at a point $x$ now means proving it at point a $x\in\Spec\Lambdaf$ and not anymore at a point of $\Spec\Hsmr$. The former is significantly harder than the latter since the contributions to the Iwasawa Main Conjecture at $x$ of various two-dimensional Galois representations corresponding to the various points of $\Hsmr$ over $x$ could be entangled in a complicated way. Using the fact that points in $\Xcali^{\sm}$ are unramified over $\Lambdaf$ and the precise structure of $\Hsmr$ as $\Lambdaf$-algebra, we show first that the Iwasawa Main Conjecture for points of $\Xcali^{\sm}$ entails the Iwasawa Main Conjecture with coefficients in $\Lambdaf$ (theorem \ref{TheoConjUnivWeak}), and then that the Iwasawa Main Conjecture with coefficients in $\Lambdaf$ entails the Iwasawa Main Conjecture at all classical points except those which are ramified over $\Lambdaf$.

This is not quite the end of the proof still, as our original classical point of interest $\rho_{f}$ might very well be ramified over $\Lambdaf$, and in that case it is far from obvious that theorem \ref{TheoIntro} for $M(f)$ can be deduced from \ref{TheoConjUnivWeak}. Using a limit argument, we show that there are many points nearby $\rho_{f}$ which satisfy the Iwasawa Main Conjecture and which are unramified over $\Lambdaf$. In that situation, we can show that even if the Iwasawa Main Conjecture is false over the full deformation space $\Hsmr$, it must be true for $M(f)$. The precise working of the proof crucially relies on the fact that both the fundamental lines and the zeta morphisms admit a version with and without Euler factors, and that both commute with specialization at a classical point, even if this classical point is ramified over $\Lambdaf$.

\subsubsection*{Notations}All rings are assumed to be commutative (and unital). The residual field at a prime ideal $\pid\in\Spec A$ is denoted by $\kg(\pid)$. If $A$ is reduced, its total quotient ring is denoted by $Q(R)$. If $A$ is a domain, its fraction field is denoted by $\Frac(A)$. If $A$ is a local ring, we denote by $\mgot_{A}$ its maximal ideal. If $F$ is a field, we denote by $G_{F}$ the Galois group of a separable closure of $F$. If $F$ is a number field and $\Sigma$ is a finite set of finite places of $\Q$, we denote by $F_{\Sigma}$ the maximal Galois extension of $F$ unramified outside $\Sigma\cup\{v|\infty\}$ and by $G_{F,\Sigma}$ the Galois group $\Gal(F_{\Sigma}/F)$. The ring of integer of $F$ is written $\Ocal_{F}$. If $v$ is a finite place of $F$, we denote by $\Ocal_{F,v}$ the unit ball of $F_{v}$, by $\varpi_{v}$ a fixed choice of uniformizing parameter of $\Ocal_{F,v}$ and by $k_{v}$ is the residual field of $\Ocal_{F,v}$. The reciprocity law of local class field theory is normalized so that $\varpi_{v}$ is sent to (a choice of lift of) the geometric Frobenius morphism $\Fr(v)$. Let $E/\qp$ be a finite extension with unit ball $\Ocal$ and residual field $\Fp$. We fix a uniformizing parameter $\varpi$ of $E$. 

For $G$ a topological group, a $G$-representation $(T,\rho,A)$ is an $A$-module $T$ free of finite rank together with a continuous morphism
\begin{equation}\nonumber
\rho:G\fleche\Aut_{A}(T).
\end{equation}

Let $\Q(\zeta_{p^{\infty}})$ be the extension of $\Q$ generated by all roots of unity of order a power of $p$. Then $\Q(\zeta_{p^{\infty}})/\Q$ is an abelian extension and $\Gamma\eqdef\Gal\left(\Q(\zeta_{p^{\infty}})/\Q\right)$ is isomorphic to $\zp\croix$ through the cyclotomic character $\chi_{\cyc}$. We say that a character $\chi:\Gamma\fleche\C_{p}\croix$ is classical if there exist a finite order character $\epsi$ and an integer $n\in\Z$ such that $\chi=\epsi\chi_{\cyc}^{n}$ (so the classical characters are $\Gamma$ are the de Rham characters of $G_{\qp}$ seen as characters of $\Gamma$). We write $\Gamma_{\Iw}$ for the largest subquotient of $\Gamma$ isomorphic to $\zp$. Define $\Q_{n}/\Q$ to be the subfield of $\Q(\zeta_{p^{n+1}})$ with Galois group $G_{n}$ over $\Q$ equal to $\Z/p^{n}\Z$. 

\section{Generalities on $p$-adic Hecke algebras}
\subsection{$p$-adic Hecke algebras}
We write $U^{p}$ (resp. $U_{p}$) for a compact open subgroup of $\GL_{2}(\Afinis{p\infty}{\Q})$ (resp. of $\GL_{2}(\qp)$) and we denote by $U\subset\GL_{2}(\Afiniq)$ a compact open subgroup of the form $U_{p}U^{p}$. If $\ell\nmid p$ is a finite prime, we write $U_{\ell}$ for the compact open subgroup of $\GL_{2}(\Q_{\ell})$ which is the local component of $U$. We denote by $\Sigma(U)$ the finite set of finite primes whose elements are $p$ and the rational primes $\ell\nmid p$ such that $U_{\ell}$ is not a maximal open compact subgroup. If $U=U_{p}U^{p}$, $\Sigma(U)$ depends only on $U^{p}$.

For integers $M$ and $N$, we denote by $U(M,N)\subset\Afinis{\infty}{\Q}$ the compact open subgroup such that 
\begin{equation}\nonumber
U(M,N)_{\ell}=\left\{\matricetype\in\GL_{2}\left(\zl\right)\mid a-1,b\equiv0\modulo \ell^{v_{\ell}(M)}\textrm{ and }c,d-1\equiv0\modulo \ell^{v_{\ell}(N)}\right\}
\end{equation}
for all $\ell$. Denote by $U(N)$ the compact open subgroup $U(N,N)$ or equivalently such that
\begin{equation}\nonumber
U(N)_{\ell}=\left\{\matricetype\in\GL_{2}\left(\zl\right)\mid \matricetype\equiv\matrice{1}{0}{0}{1}\modulo \ell^{v_{\ell}(N)}\right\}
\end{equation}
for all $\ell$ and by $U_{1}(N)$ the compact open subgroup $U(1,N)$ or equivalently such that
\begin{equation}\nonumber
U_{1}(N)_{\ell}=\left\{\matricetype\in\GL_{2}\left(\zl\right)\mid \matricetype\equiv\matrice{*}{*}{0}{1}\modulo \ell^{v_{\ell}(N)}\right\}
\end{equation}
for all $\ell$. By construction, the sets $\Sigma(U_{p}U(N))$ and $\Sigma(U_{p}U_{1}(N))$ are both equal to $\{\ell|Np\}$.

For $U\subset\GL_{2}(\Afiniq)$ a compact open subgroup as above, we denote by $Y(U)$ the affine modular curve of level $U$, that is to say the Shimura curve whose complex points are given by
\begin{equation}\nonumber
Y(U)(\C)=\GL_{2}(\Q)\backslash\left(\C-\R\times\GL_{2}(\Afiniq)\right)/U.
\end{equation}
For simplicity, we write $Y(M,N)$ for $Y(U(M,N))$, $Y(N)$ for $Y(N,N)$ and $Y_{1}(N)$ for $Y(U_{1}(N))$. Note that if $N$ and $m$ are integers dividing $M$, then there is a covering of $\Q(\zeta_{M})$-schemes from $Y(m,M)$ onto $Y_{1}(N)\tenseur_{\Q}\Q(\zeta_{m})$.

We write $H^{i}_{\et}\left(U,-\right)$ for the $i$-th étale cohomology group functor $H^{i}_{\et}\left(Y(U)\times_{\Q}\Qbar,-\right)$ and $H^{i}_{c}\left(U,-\right)$ for the $i$-th étale cohomology group with compact support functor $H^{i}_{c}\left(Y(U)\times_{\Q}\Qbar,-\right)$.

For $\Sigma$ a finite set of finite primes containing $\Sigma(U)$, let $\Hecke^{\Sigma}(U)$ be the $p$-adic Hecke algebra of endomorphisms of $\Hun_{\et}\left(U,\Ocal\right)$ generated as $\Ocal$-algebra by the Hecke operators $T(\ell)$ and $S(\ell)$ for $\ell\notin\Sigma$.

Fix $\Sigma$ as above and consider
\begin{equation}\nonumber
\rhobar:G_{\Q,\Sigma}\fleche\GL_{2}(\Fp)
\end{equation}
an absolutely irreducible and modular (hence odd) Galois representation unramified outside $\Sigma$.  We say that a compact open subgroup $U\subset\GL_{2}\left(\Afiniq\right)$ is allowable for $\rhobar$ (or allowable, for short) if there exists a (necessarily unique) maximal ideal $\mgot_{\rhobar}\in\Spec\Hs(U)$ such that 
\begin{equation}\nonumber
\begin{cases}
T(\ell)\modulo\mgot_{\rhobar}=\tr\rhobar(\Fr(\ell))\\
S(\ell)\modulo\mgot_{\rhobar}=\det\rhobar(\Fr(\ell))
\end{cases}
\end{equation}
holds for all $\ell\notin\Sigma$.

Fix $U'\subset U$ two allowable compact open subgroups. Denote respectively by $\mgot'_{\rhobar}$ and $\mgot_{\rhobar}$ the maximal ideals of $\Hecke^{\Sigma'}(U')$ and $\Hecke^{\Sigma}(U)$ attached to $U'$ and $U$. Then the covering $Y(U')\fleche Y(U)$ induces by restriction a map $\Hecke^{\Sigma'}(U')_{\mgot'_{\rhobar}}\fleche\Hs(U)_{\mgot_{\rhobar}}$ which is always a surjection and which is an isomorphism if $U$ is sufficiently small. In the following, we always assume that $U$ is allowable and sufficiently small in that sense. We write $\Hsmr$ for the common isomorphism class of Hecke algebras $\{\Hs(U)_{\mgot_{\rhobar}}\}_{\Sigma,U}/\simeq$ and $\mgot_{\rhobar}$ for the maximal ideal of $\Hsmr$. Even though we suppress this dependence from the notation, note that $\Hsmr$ does depend on the choice of $U$.

Henceforth, the set $\Sigma$ is fixed once and for all. We write $\Sp$ for $\Sigma\backslash\{p\}$.

A specialization of $\Hsm$ is by definition a local ring morphism $\lambda:\Hsm\fleche S$. According to $\cite{CarayolRepresentationsGaloisiennes}$, there exists a $\Hsmr$-module $\Ts$ free of rank 2 endowed with a continuous $G_{\Q,\Sigma}$-action through the $G_{\Q,\Sigma}$-representation
\begin{equation}\nonumber
\rho_{\Sigma}:G_{\Q,\Sigma}\fleche\Aut_{\Hsmr}\left(\Ts\right)\simeq\GL_{2}(\Hsm)
\end{equation}
uniquely characterized up to isomorphism by the requirement that $\tr(\rho_{\Sigma}(\Fr(\ell)))=T(\ell)$ for all $\ell\notin\Sigma$. To a specialization $\lambda:\Hsm\fleche S$ is attached the $G_{\Q,\Sigma}$-representation $(T_{\lambda},\rho_{\lambda},S)$ defined by $\rho_{\lambda}=\lambda\circ\rho_{\Sigma}$ and characterized uniquely (up to isomorphism) by the fact that $\tr(\rho_{\lambda}(\Fr(\ell)))=\lambda(T_{\ell})$ for all $\ell\notin\Sigma$. If $\aid\in\Spec^{\min}\Hr$ is a minimal prime ideal of $\Hr$, let $\Raid$ be the quotient domain $\Hr/\aid$, let $\Faid$ be te field of fractions of $\Raid$, let $\Taid$ be the $G_{\Q,\Sigma}$-representation $\Ts\tenseur_{\Hr}\Raid$ and let finally $\Vaid$ be the $G_{\Q,\Sigma}$-representation $\Ts\tenseur_{\Hr}\Faid$. If $\la:\Hr\fleche A$ is a specialization with values in a domain, let $\Vla$ be the representation $\Ts\tenseur_{\Hr,\la}\Frac(A)$.

\subsection{Commutative algebra properties of $\Hsmr$}
\subsubsection{Classical primes of the $p$-adic Hecke algebra}
A prime ideal $\pid\subset\Spec\Hs$ with residual field $\kb(\pid)\subset\Qbar_{p}$ is said to be classical if there exists a (necessarily unique) normalized eigencuspform $f\in S_{k}(U)$ of weight $k\geq2$ whose system of eigenvalues
\begin{equation}\nonumber
\lambda_{f}:\Hs\fleche\Qbar_{p}
\end{equation}
coincides with the morphism $\Hs\fleche\kb(\pid)$. It is said to be classical up to a twist if there exists a normalized eigencuspform $f$ of weight $k\geq2$ and a classical character $\chi\in\hat{\Gamma}_{\Iw}$ such that the morphism $\Hs\fleche\kb(\pid)$ is the system of eigenvalues of $f\tenseur\chi$. A specialization $\la:\Hs\fleche S$ is classical (resp. up to a twist) if it factors through a $\Hs/\pid$ for $\pid$ classical (resp. up to a twist). We write $f_{\lambda}$ or $f_{\pid}$ for the eigencuspform attached to $\lambda$. By definition, the $G_{\Q,\Sigma}$ representations $\rho_{\lambda}$ and $\rho_{f_{\lambda}}$ are isomorphic. Let $\pid$ be a prime of $\Hs[1/p]$ such that the morphism $\lambda:\Hs[1/p]\fleche\kb(\pid)$ is classical up to a twist. Then $\pid,\lambda$ or $\rho_{\lambda}$ are said to be crystalline if $D_{\cris}\left(\rho_{\lambda}|G_{\qp}\right)$ is a $\kb(\pid)$-vector space of dimension 2. They are said to be distinguished, crystalline if the Frobenius morphism $\p$ acting on $D_{\cris}\left(\rho_{\lambda}|G_{\qp}\right)$ has two distinct eigenvalues. They are said to be crystalline and short if the weight $k$ of $f_{\lambda}$ satisfies $2\leq k\leq p$, or equivalently if $\rho|G_{\qp}$ is the image by the Fontaine-Laffaille functor of a rank 2 Fontaine-Laffaille module with non-trivial graded piece of the filtration in degree $0$ and $k-1$ (\cite[Théorème 8.4]{FontaineLaffaille}). According to \cite{SaitoMonodromie}, if $\la$ is crystalline, then there exist a compact open subgroup $U\subset\GL_{2}(\Afiniq)$ with $U_{p}$ a maximal compact open subgroup of $\GL_{2}(\qp)$, $f\in S_{k}(U)$ a newform, $\eta:G_{\Q}\fleche\qpbar\croix$ a character unramified at $p$ and $i\in\Z$ such that $\la$ is the system of eigenvalues of $f\tenseur\eta\chi_{\cyc}^{i}$.

All the terminology above (classical point, classical prime, crystalline point, crystalline prime...) is extended in the obvious way to $\Hsmr$ and its quotients.
\subsubsection{Identification with universal deformation rings}
We make the following assumption on $\mgot_{\rhobar}$.
\begin{Hyp}\label{HypTW}
The $G_{\Q,\Sigma}$-representation $\rhobar$ satisfies the following properties.
\begin{enumerate}
\item\label{HypIrr} If $p^{*}=(-1)^{(p-1)/2}p$, then $\rhobar|G_{\Q(\sqrt{p^{*}})}$ is absolutely irreducible.
\item\label{HypLoc} If $\rhobar|G_{\qp}$ is an extension
\begin{equation}\nonumber
\suiteexacte{}{}{\chi_{1}}{\rhobar|G_{\qp}}{\chi_{2}},
\end{equation}
then $\chi_{1}^{-1}\chi_{2}\notin\{1,\bar{\chi}_{\cyc}\}$.
\end{enumerate}
\end{Hyp}
If $\p$ is an endomorphism of the $\Fp[G_{\qp}]$-module of $\rhobar|G_{\qp}$, then $\p$ is scalar unless $(\rhobar|G_{\qp})^{ss}\simeq\chibar\oplus\chibar$. Under assumption \ref{HypTW} then, the endomorphisms of $\rhobar|G_{\qp}$ are scalar.

Let $D_{\rhobar}$ be the deformation functor from the category of finite, local, $W(\Fp)$-algebras with residue field $\Fp$ to the category of sets which attaches to $A$ the set of isomorphism casses of pairs $\{(T,\rho,A),\iota)\}$ such that $(T,\rho,A)$ is a $G_{\Q,\Sigma}$-representation and $\iota:\rho\tenseur_{A}\Fp\simeq\rhobar$ is an isomorphism of $\Fp[G_{\Q,\Sigma}]$-module. It follows from statement \ref{HypIrr} of assumption \ref{HypTW} that $D_{\rhobar}$ is representable by a complete, local, noetherian ring $\Rs(\rhobar)$. We denote by $(\Ts^{u},\rho_{\Sigma}^{u},\Rs(\rhobar))$ the corresponding universal $G_{\Q,\Sigma}$-representation. By definition, there is a map $\Rsr\fleche\Hsmr$ which is surjective as its image contains the image of $\tr(\rho_{\Sigma}^{u}(\Fr(\ell)))$, which is $\tr(\rho_{\Sigma}(\Fr(\ell)))=T_{\ell}$ for all $\ell\notin\Sigma$.

The following lemma is well-known.
\begin{Lem}\label{LemFontaineLaffaille}
Suppose $\rhobar|G_{\qp}$ is irreducible. Then there exists a crystalline and short point $\lambda_{f}:\Hs\fleche\Ocal$ such that $\Rs(\rhobar)\simeq\Rs(\rhobar_{f})$.
\end{Lem}
\begin{proof}
If $\chi:G_{\Q,\Sigma}\fleche W(\Fp)\croix$ is a character, then $\rho\mapsto\rho\tenseur\chi$ induces an isomorphism of functors $D_{\rhobar}\fleche D_{\rhobar\tenseur\chibar}$. Hence, it is enough to prove that there exists such a character $\chi$ and an eigencuspform $f$ such that $\rhobar_{f}\tenseur\chibar\simeq\rhobar$ and such that $\rho_{f}$ is crystalline and short. 

Let us denote by $I^{w}\subset I_{p}$ the wild ramification subgroup of $G_{\qp}$ and by $I^{t}\eqdef I_{p}/I^{w}\simeq \limproj{r}\ \Fp_{p^{r}}\croix$  the tame ramification quotient. Recall that a character $\psi:I^{t}\fleche \Fp\croix$ is fundamental if it is the composition of the natural map
\begin{equation}\nonumber
I^{t}\simeq\limproj{r}\ \Fp_{p^{r}}\croix\fleche\Fp\croix_{p^{n}}
\end{equation} 
with one of the $n$ embeddings $\Fp\croix_{p^{n}}\plonge \Fp\croix$ and that it is of level $n$ if $\Fp_{p^{n}}$ is the smallest subfield through which $\psi$ factors. Denote by $\psi_{i}$ the two fundamental characters of level 2 of $I^{t}$. As $I^{w}$ is a pro-$p$-group and $I^{t}$ is commutative, $\rhobar^{ss}|I^{w}$ is trivial and $\rhobar^{ss}|I_{p}$ is the sum of two characters $\chibar_{1},\chibar_{2}$ of $I^{t}$. Since $G_{\qp}/I_{p}$ acts by conjugation on $I^{t}$ through the monodromy relation, $\chibar_{1}$ and $\chibar_{2}$ are permuted modulo $p$ by elevation to the $p$-th power. Hence, the $\chibar_{i}$ are of level 1 or 2. This permutation is the identity if and only if the $\chibar_{i}$ are both of level 1, if and only if $\rhobar|G_{\qp}$ is reducible. As we have assumed that $\rhobar|G_{\qp}$ is irreducible, the characters $\chibar_{i}$ are both of level 2. Then there is a unique pair $(a,b)$ satisfying $0\leq a< b\leq p-1$ such that
\begin{equation}\nonumber
\rhobar|I^{t}\simeq\matrice{\psi_{1}^{a}\psi_{2}^{b}}{0}{0}{\psi_{1}^{b}\psi_{2}^{a}}\simeq\matrice{\psi_{2}^{b-a}}{0}{0}{\psi_{1}^{b-a}}\tenseur\bar{\chi}_{\cyc}^{a}
\end{equation}
Accordingly, $\rhobar$ has Serre-weight $k=1+pa+b$ (\cite[(2.2.4)]{SerreConjecture}).  Put $\rhobar'=\rhobar\tenseur\bar{\chi}_{\cyc}^{-a}$. Then $\rhobar'|I^{t}$ may be written
\begin{equation}\nonumber
\rhobar'|I^{t}\simeq\matrice{\psi_{2}^{b-a}}{0}{0}{\psi_{1}^{b-a}}
\end{equation}
and so has Serre-weight $k'=1+b-a$. Then $2\leq k'\leq p$. According to the proof of the weight part of Serre's Conjecture (\cite{EdixhovenWeight}), there exists an eigencuspform $f\in S_{k'}$ with level prime to $p$ such that $\rhobar_{f}$ is isomorphic to $\rhobar'$. Then $\lambda_{f\tenseur\chi_{\cyc}^{a}}$ is a crystalline point of $\Hsmr[1/p]$ and so $\Rs(\rhobar)$ is isomorphic to $\Rs(\rhobar_{f})$.

\end{proof}

\begin{Lem}\label{LemPresentation}
There exists an integer $n\geq4$ and elements $(y_{4},\cdots,y_{n})$ such that the ring $\Rs(\rhobar)$ admits a presentation 
\begin{equation}\nonumber
\Rs(\rhobar)\simeq\Ocal[[X_{1},\cdots,X_{n}]]/(y_{4},\cdots,y_{n}).
\end{equation}
\end{Lem}
\begin{proof}
See \cite[Lemma 4.3]{BockleDensity}.
\end{proof}
\begin{Prop}\label{PropTW}
The universal deformation ring $R_{\Sigma}(\rhobar)$ is a flat $\Ocal$-algebra which is a complete intersection ring of Krull dimension 4.\end{Prop}
\begin{proof}
Let $k$ be the Serre-weight of $\rhobar$. According to \cite{EdixhovenWeight}, there exists a classical prime $\pid\in\Spec\Hsm[1/p]$ attached to a classical eigencuspform $f$ of weight $k$. Denote by $\chi:G_{\mathbb{Q}}\fleche\Ocal\croix$ the character $\det\rho_{f}$. Let $R_{\Sigma}^{\chi}(\rhobar)$ be the deformation ring of $\rhobar$ parametrizing deformations $\rho$ of $\rhobar$ with coefficients in complete noetherian local $\Ocal$-algebras, unramified outside $\Sigma$ and such that $\det\rho$ is equal to $\chi$ (in particular, $\rho_{f}$ corresponds to a point in $\Rs^{\chi}(\rhobar)$). Let $\rho$ be an arbitrary deformation of $\rhobar$ corresponding to a point of $\Rs(\rhobar)$ with coefficients in $\Ocal\croix$. Then $\left(\det\rho\right)^{-1}\chi$ has values in $\Ocal\croix$ and $\left(\det\rho\right)^{-1}\chi\equiv \left(\det\rhobar\right)^{-1}\left(\det\rhobar_{f}\right)^{-1}\equiv1\modulo\varpi$. So $\left(\det\rhobar\right)^{-1}\chi$ has values in $1+\varpi\Ocal$ and corresponds to a point of the universal deformation ring $\Rs(\indicatrice)$ parametrizing deformations $\psi:G_{\Q,\Sigma}\fleche\Ocal\croix$ of the trivial character $\indicatrice$ of $G_{\Q,\Sigma}$. As $p$ is odd, the multiplicative group $1+\varpi\Ocal$ is uniquely 2-divisible so characters with values in $1+\varpi\Ocal$ admit canonical square roots. Let $\psi_{\rho}$ be the canonical square root of $\left(\det\rho\right)^{-1}\chi$. Then $\rho\tenseur\psi_{\rho}\equiv\rhobar\modulo\varpi$ and $\det\left(\rho\tenseur\psi_{\rho}\right)=\chi$ so $\rho\tenseur\psi_{\rho}$ corresponds to a point of $\Rs^{\chi}(\rhobar)$. The map $\rho\mapsto\rho\tenseur\psi_{\rho}$ thus induces an isomorphism $\Rs(\rhobar)\simeq\Rs^{\chi}(\rhobar)\hat{\tenseur}_{\Ocal}\Rs(\indicatrice)$. As $\Rs(\indicatrice)$ is isomorphic to a power-series ring in one-variable over a complete intersection $\Ocal$-algebra of relative dimension zero, $R_{\Sigma}(\rhobar)$ is a flat $\Ocal$-algebra which is a complete intersection ring of relative dimension 3 if and only if $R^{\chi}_{\Sigma}(\rhobar)$ is a flat $\Ocal$-algebra which is a complete intersection ring of relative dimension 2. This we now show.

Assume first that $\rhobar|G_{\qp}$ is reducible. Twisting by a character if necessary, we may then assume that $(\rhobar|G_{\qp})^{ss}=\chibar_{1}\oplus\chibar_{2}$ with $\chibar_{1}(I_{p})=\{1\}$. Let $R^{\ord,\chi}_{\Sigma}(\rhobar)$ be the quotient of $\Rs^{\chi}(\rhobar)$ parametrizing deformations $\rho$ which in addition to being points of $R^{\chi}_{\Sigma}(\rhobar)$ are such that there exists a short exact sequence of non-zero $G_{\qp}$-representations 
\begin{equation}\nonumber
\suiteexacte{}{}{\chi_{1}}{\rho|G_{\qp}}{\chi_{2}}
\end{equation}
with $\chi_{1}(I_{p})=\{1\}$. According to the main results of \cite{DiamondHecke,FujiwaraDeformation} (following \cite{WilesFermat,TaylorWiles}), the ring $R^{\ord,\chi}_{\Sigma}(\rhobar)$ is isomorphic to a suitable Hecke algebra and hence flat of relative dimension 0 over $\Ocal$. Suppose now that $\rhobar|G_{\qp}$ is irreducible. Thanks to lemma \ref{LemFontaineLaffaille}, we may then assume that the eigencuspform $f$ above is crystalline of weight $2\leq k\leq p$. Denote by $R^{\crys,\chi}_{\Sigma}(\rhobar)$  the quotient of $\Rs^{\chi}(\rhobar)$ parametrizing deformations $\rho$ which in addition to being points of $R^{\chi}_{\Sigma}(\rhobar)$ are such that $\rho|G_{\qp}$ is crystalline and short. By the modularity result of \cite{DiamondFlachGuo}, the ring $R^{\crys,\chi}_{\Sigma}(\rhobar)$ is isomorphic to a suitable Hecke algebra and hence flat of relative dimension 0 over $\Ocal$. Both when $\rhobar|G_{\qp}$ is reducible and $*=\ord$ and when $\rhobar|G_{\qp}$ is irreducible and $*=\crys$ then, the ring $R^{*,\chi}_{\Sigma}(\rhobar)$ is flat of relative dimension 0 over $\Ocal$.

Next we describe the kernel of the map $R^{\chi}_{\Sigma}(\rhobar)\fleche R^{*,\chi}_{\Sigma}(\rhobar)$. Assume first that $\rhobar|G_{\qp}$ is reducible. Assumption \ref{HypLoc} of \ref{HypTW} and \cite[Lemma 2.2]{SkinnerWilesDur} imply that the universal framed ordinary deformation ring $R^{\square,\ord}(\rhobar|G_{\qp})$ is a regular ring of relative dimension 3. According to \cite[Corollary 7.4]{BockleDemuskin}, $R^{\square,\ord}(\rhobar|G_{\qp})$ is a quotient of $R^{\square}(\rhobar|G_{\qp})$ by a length two regular sequence. As a quotient of a ring by a regular sequence is regular (if and) only if the ring itself was regular, this entails that $R^{\square}(\rhobar|G_{\qp})$ is a regular ring of relative dimension 5 over $\Ocal$ and that kernel of the map $R^{\square}(\rhobar|G_{\qp})\fleche R^{\square,\ord}(\rhobar|G_{\qp})$ is generated by a subset of a system of parameters of cardinal 2. Now we assume that $\rhobar|G_{\qp}$ is irreducible. Then \cite[Proposition 2.2 and Corollary 2.3]{DiamondFlachGuo} (and previous results of \cite{RamakrishnaFlat}) imply that the universal deformation ring $R(\rhobar|G_{\qp})$ is a regular ring of relative dimension 5 and that $R^{\crys}\left(\rhobar|G_{\qp}\right)$ is a regular quotient of relative dimension 2. Let $R^{\chi}(\rhobar|G_{\qp})$ and $R^{\crys,\chi}(\rhobar|G_{\qp})$ respectively denote the quotient of the universal deformation ring parametrizing deformations $\rho$ of $\rhobar|G_{\qp}$ with coefficients in complete noetherian local $\Ocal$-algebras and such that $\det\rho$ is equal to $\chi$ and the quotient parametrizing such deformations which are in addition crystalline and short. By local class field theory, $R^{\chi}(\rhobar|G_{\qp})$ is then a regular ring of relative dimension 3 and $R^{\crys,\chi}(\rhobar|G_{\qp})$ is a regular ring of relative dimension 1. The kernel of $R^{\chi}(\rhobar|G_{\qp})\fleche R^{\crys,\chi}(\rhobar|G_{\qp})$ is thus generated by a regular sequence of length 2.

Denoting by $R^{\chi}(\rhobar|G_{\qp})$ and $R^{*,\chi}(\rhobar|G_{\qp})$ the quotient of the framed or usual universal deformation ring parametrizing deformations $\rho$ of $\rhobar|G_{\qp}$ with coefficients in complete noetherian local $\Ocal$-algebras and such that $\det\rho$ is equal to $\chi$ (resp. which in addition are of type $*$), we consequently see that in both the reducible and irreducible case, there is a commutative diagram
\begin{equation}\nonumber
\xymatrix{
R^{\chi}(\rhobar|G_{\qp})\ar[r]\ar[d]&R^{*,\chi}(\rhobar|G_{\qp})\ar[d]\\
R^{\chi}_{\Sigma}(\rhobar)\ar[r]&R^{*,\chi}_{\Sigma}(\rhobar)
}
\end{equation}
where the vertical map are induced by restriction from $G_{\Q,\Sigma}$ to $G_{\qp}$ and where the kernel of the upper horizontal arrow is generated by a regular sequence of length 2. This implies that $R^{*,\chi}_{\Sigma}(\rhobar)$ is a quotient of $\Rs^{\chi}(\rhobar)$ by a an ideal generated by at most two elements, and thus that $R^{*,\chi}_{\Sigma}(\rhobar)$ is a quotient of $\Rs(\rhobar)$ by an ideal $(x_{1},x_{2},x_{3})$ generated by at most three elements. This entails in particular that $\Rs(\rhobar)$ is of Krull dimension at most 4.

According to lemma \ref{LemPresentation}, the ring $\Rs(\rhobar)$ is of dimension at least 4 so we conclude that it is flat of relative dimension 3 over $\Ocal$. The same lemma shows that the zero-dimensional ring $\Rs^{*,\chi}(\rhobar)/(\varpi)$ admits a presentation
\begin{equation}\nonumber
\Rs^{*,\chi}(\rhobar)/(\varpi)\simeq\Ocal[[X_{1},\dots,X_{n}]]/(\varpi,x_{1},x_{2},x_{3},y_{4},\dots,y_{n})
\end{equation}
As $\Rs^{*,\chi}(\rhobar)/(\varpi)$ is of Krull dimension zero, this implies that $(\varpi,x_{1},x_{2},x_{3},y_{4},\dots,y_{n})$ is a regular sequence in $\Ocal[[X_{1},\dots,X_{n}]]$. In particular, $\Rs^{*,\chi}(\rhobar)/(\varpi)$, $\Rs^{*,\chi}(\rhobar)$ and $\Rs(\rhobar)$ are complete intersection rings of dimension 0, 1 and 4 respectively. It follows in addition that $(x_{1},x_{2},x_{3})$ is a regular sequence in $\Rs(\rhobar)$.
\end{proof}
We note that the proof above establishes the existence of suitable quotients $\Rs^{*,\chi}(\rhobar)$ even without twisting $\rhobar$: in that case, these rings are classical Hecke algebras acting on spaces of eigencuspforms which are ordinary or crystalline and short up to a twist. This entails in particular the following corollary.
\begin{Cor}\label{CorSmooth}
There exists a Zariski-dense, open subset $\Xcali^{\sm}\subset\Spec\Rsr$ containing all points ordinary up to a twist if $\rhobar|G_{\qp}$ is reducible and all points crystalline and short up to a twist if $\rhobar|G_{\qp}$ is irreducible such that the map
\begin{equation}\nonumber
\xymatrix{
\Spec\Rsrchi\ar[d]\\
\Spec\Ocal
}
\end{equation}
is formally smooth at $x\in\Xcali^{\sm}$. In particular, the ring $\Rsr$ is reduced.
\end{Cor}
\begin{proof}
Let $(x_{2},x_{3})$ the regular sequence of the proof of proposition \ref{PropTW} and let $\Lambdaf$ be the power-series ring $\Ocal[[X_{2},X_{3}]]$. We consider the commutative diagram
\begin{equation}\nonumber
\xymatrix{
\Lambdaf\ar[rr]^(0.5){X_{i}\mapsto x_{i}}\ar[d]_{\pi}&&R^{\chi}_{\Sigma}(\rhobar)\ar[d]^{\pi'}\\
\Ocal\ar[rr]&&R^{*,\chi}_{\Sigma}(\rhobar)
}
\end{equation}
where the vertical maps $\pi$ and $\pi'$ are the quotient maps modulo $(X_{2},X_{3})$ and $(x_{2},x_{3})$ respectively. Let $\pid\in\Spec\Rsrchi$ be a point above $\pi^{*}:\Spec\Ocal\fleche\Spec\Lambdaf$. As $\Rsrchi_{\pid}$ is flat over $\Rsrchi$, the regular sequence $(x_{2},x_{3})$ of $\Rsrchi$ remains a regular sequence of $\Rsrchi_{\pid}$. The quotient $\Rsrchi_{\pid}/(x_{2},x_{3})$ is by construction the localization at a minimal prime ideal of the reduced, classical Hecke algebra of appropriate level and weight after extension of scalars to the fraction field $E$ of $\Ocal$. Hence, it is a separable field extension and $\Spec\Rsrchi_{\pid}/(x_{2},x_{3})\fleche\Spec E$ is an étale morphism. Consequently, $\Spec\Rsrchi\fleche\Spec\Lambdaf$ is unramified at $\pid$. As it is also finite and flat by the proof of \ref{PropTW}, the locus $\Xcali^{\sm}$ of étale points of $\Spec\Rsrchi\fleche\Spec\Lambdaf$ contains $\pid$ and is thus non-empty. Let $U\subset\Spec\Lambdaf$ be the complement of the support of $\Omega^{1}_{\Rsrchi/\Lambdaf}$ regarded as $\Lambdaf$-module. By the above, $U$ is non-empty, formally smooth over $\Spec\Ocal$ and $\Xcali^{\sm}$ is formally smooth over $U$. Hence $\Xcali^{\sm}$ is formally smooth over $\Spec\Ocal$. By construction, $U$ is open and non-empty, hence Zariski-dense in $\Spec\Lambdaf$. So $\Spec\Lambdaf\backslash U$ is of codimension at least 1. As $\Spec\Rsrchi\fleche\Spec\Lambdaf$ is finite, $\Spec\Rsrchi\backslash\Xcali^{\sm}$ is also of codimension at least 1. As $\Rsrchi$ is Cohen-Macaulay, it is equidimensional. Hence $\Xcali^{\sm}$ is Zariski-dense in each irreducible component. Let $\{\aid_{i}|i\in I\}\subset\Spec\Rsrchi$ be the finite set of minimal prime ideals of $\Rsrchi$ and write $\Xcali_{i}\eqdef\Spec\Rsrchi/\aid_{i}$. Then the sets $\Xcali^{\sm}\cap\Xcali_{i}$ are non-empty and disjoint. As $\Xcali^{\sm}\cap\Xcali_{i}$ is étale over $U$, its generic degree is equal to its degree at $\pid$. Hence, each $\Xcali_{i}$ contains a point $\pid_{i}\in\Xcali^{\sm}\cap\Xcali_{i}$ above $\pi^{*}$ and the $\pid_{i}$ are pairwise distinct.
 
As $\Rsrchi$ is a Cohen-Macaulay ring, the set $\{\aid_{i}|i\in I\}$ is also the set of associated primes of $\Rsrchi$. By the properties of $\Xcali^{\sm}$ just established, for each $i\in I$, there exists an element $a_{i}\in\Rsrchi$ such that $a_{j}\mod\aid_{i}$ vanishes if and only if $i\neq j$ and such that $\Spec\Rsrchi_{a_{i}}$ is smooth over $E$. Then $a=\somme{i\in I}{}a_{i}$ does not belong to any minimal prime $\aid_{i}$, hence does not belong to an associated prime of $\Rsrchi$. Hence $a$ is not a zero-divisor and there is an embedding $\Rsrchi\plonge\Rsrchi_{a}$. As $\Spec\Rsrchi_{a}\simeq\uniondisjointe{i\in I}{}\Spec\Rsrchi_{a_{i}}$ and as $\Spec\Rsrchi_{a_{i}}$ is smooth, $\Rsrchi$ is reduced.
\end{proof}
The following proposition is well-known.
\begin{Prop}\label{PropDense}
The set of primes classical up to a twist is Zariski-dense in $\Spec\Rsr$. In particular, $\Rsr$ and $\Hsmr$ are isomorphic. 
\end{Prop}
\begin{proof}
The first assertion follows from corollary \ref{CorSmooth} as in \cite[Proof of theorem 3.7]{BockleDensity}. As the morphism $\Spec\Hsmr\fleche\Spec\Rsr$ is dominant, the kernel of the natural surjection $\Rsr\fleche\Hsm$ is included in the nilradical of $\Rsr$, which is a reduced ring according to corollary \ref{CorSmooth}.
\end{proof}
\begin{Cor}\label{CorCohen}
Let $\Lambdaf$ be the power-series ring $\Ocal[[X_{1},X_{2},X_{3}]]$. Then there is a length 3 regular sequence $(x_{1},x_{2},x_{3})$ inside $\Hsm\simeq\Rs(\rhobar)$ such that the assignment $X_{i}\mapsto x_{i}$ endows $\Hsm$ with a structure of $\Lambdaf$-algebra for which it is finite and free as $\Lambdaf$-module and verifying in addition the following property. Classical points of $\Hsm$ which induce the same morphisms after restriction to $\Lambdaf$ through this specified embedding $\Lambdaf\plonge\Hsm$ are either all crystalline and short up to a twist (resp. ordinary up to a twist) or none of them are. 
\end{Cor}
\begin{proof}
As ${\Rs(\indicatrice)}$ is a power-series ring with coefficients in complete intersection ring of relative dimension zero over $\Ocal$, there exists elements $x\in\in\mgot_{\Rs(\indicatrice)}\backslash\mgot^{2}_{\Rs(\indicatrice)}$ such that $(x,\varpi)$ is a system of parameters of $\Rs(\indicatrice)$. Let $x_{1}\in\mgot_{\Rs(\indicatrice)}\backslash\mgot^{2}_{\Rs(\indicatrice)}$ be such a regular element which we view as an element of $\Rsr$ through the isomorphism $\Rs(\rhobar)\simeq\Rs^{\chi}(\rhobar)\hat{\tenseur}_{\Ocal}\Rs(\indicatrice)$.

According to the proof of proposition \ref{PropTW}, the kernel of $\Rs(\rhobar)/(x_{1})\surjection\Rs^{*,\chi}(\rhobar)$ contains a regular sequence $(x_{2},x_{3})$ which is the image of the regular sequence generating the kernel of $R(\rhobar|G_{\qp})/(x_{1})\surjection R^{*,\chi}(\rhobar|G_{\qp})$. To a crystalline and short (resp. ordinary) point $\rho_{f}$ is thus attached a $\Lambdaf$-structure on $\Rs(\rhobar)$ defined by $X_{i}\mapsto x_{i}$ (note that $x_{2}$ and $x_{3}$ depend on the choice of $\rho_{f}$). Given such a $\Lambdaf$-structure on $\Rs(\rhobar)$, if two characteristic zero specializations $\psi_{i}:\Hsm\fleche S$ with values in a discrete valuation ring coincide after restriction to $\Lambdaf$, then they have the same value on $x_{2}$ and $x_{3}$. So they both factor through the quotient $\Rs^{*,\chi}(\rhobar)$ attached to this choice of $\Lambdaf$-structure or neither of them does.

\end{proof}
Note that the $\Lambdaf$-algebra structure described in the previous proof depends on an initial choice of a classical point $\rho_{f}$. In what follows, we omit references to this dependence.
\begin{Lem}\label{LemFontaineLaffailleTwist}
Let $\lambda_{f}:\Hsmr\fleche\Ocal$ be a classical point of even weight $k$.  Let $\omega:G_{\Q,\Sigma}\fleche\fp\croix\plonge\zp\croix$ be the character giving the action on roots of unity of order $p$.

If $\rhobar|G_{\qp}$ is irreducible, then there exists a crystalline and short point $\lambda_{g}:\Hsmr\fleche\Ocal$ of even weight $k'$ such that either $T(f)\left(-\frac{k-2}{2}\right)$ and $T(g)\left(-\frac{k'-2}{2}\right)$ are congruent modulo $\varpi$ or $T(f)\left(-\frac{k-2}{2}\right)$ and $T(g)\left(-\frac{k'-2}{2}\right)\tenseur\omega^{\frac{p-1}{2}}$ are congruent modulo $\varpi$.

If $\rhobar|G_{\qp}$ is reducible, then there exists a classical point $\lambda_{g}:\Hsmr\fleche\Ocal$ of even weight $k'$ with $\rho_{g}|G_{\qp}$ reducible and such that either $T(f)\left(-\frac{k-2}{2}\right)$ and $T(g)\left(-\frac{k'-2}{2}\right)$ are congruent modulo $\varpi$ or $T(f)\left(-\frac{k-2}{2}\right)$ and $T(g)\left(-\frac{k'-2}{2}\right)\tenseur\omega^{\frac{p-1}{2}}$ are congruent modulo $\varpi$.
\end{Lem}
\begin{proof}
First assume that $\rhobar|G_{\qp}$ is irreducible. According to the lemma \ref{LemFontaineLaffaille} and its proof, there exist a crystalline point $\lambda_{g}:\Hsmr\fleche\Ocal$ of weight $2\leq k'\leq p$ and an integer $a$ such that $\rhobar_{g}\tenseur\chibar_{\cyc}^{a-\frac{k-2}{2}}$ and $\rhobar\tenseur\chibar_{\cyc}^{-\frac{k-2}{2}}$ are isomorphic. Hence  $\rhobar_{g}\tenseur\chibar_{\cyc}^{-\frac{k'-2}{2}}$ is isomorphic to $\left(\rhobar\tenseur\chibar_{\cyc}^{-\frac{k-2}{2}}\right)\tenseur\chibar^{m}$ for $m={-a+\frac{k-k'}{2}}$. As
\begin{equation}\nonumber
\chibar_{\cyc}^{1+2m}=\det\left(\left(\rhobar\tenseur\chibar_{\cyc}^{-\frac{k-2}{2}}\right)\tenseur\chibar^{m}\right)=\det\left(\rhobar_{g}\tenseur\chibar_{\cyc}^{-\frac{k'-2}{2}}\right)=\chibar_{\cyc},
\end{equation}
the integer $m$ is congruent to $0$ or to $(p-1)/2$ modulo $p$. If $m=0$, then $T_{g}\left(-\frac{k'-2}{2}\right)$ is congruent to $T(f)\left(-\frac{k-2}{2}\right)$. If $m=(p-1)/2$, then $T_{g}\left(-\frac{k'-2}{2}\right)\tenseur\omega^{m}$ is congruent to $T(-\frac{k-2}{2})$.

Now assume that $\rhobar|G_{\qp}$ is reducible. Then there exists a point $\lambda_{g}:\Rs^{\chi,\ord}(\rhobar)\fleche\Ocal$ which is classical of a given weight $k'$. As in the first part of the proof, there exists an integer $m$ such that $T(g)\left(-\frac{k'-2}{2}\right)\tenseur\omega^{m}$ is congruent to $T(f)(-\frac{k-2}{2})$. Then, as above, $m$ must be equal to 0 or to $(p-1)/2$ modulo $p-1$.

In both cases, taking determinants (or examining again the proof of \ref{LemFontaineLaffaille}) shows that $k'$ and $k$ are congruent modulo 2 and hence that $k'$ is even.
\end{proof}
\section{Fundamental lines and zeta elements over deformation rings}
\subsection{Classical Iwasawa theory for modular motives}
\subsubsection{Equivariant modular motives and period maps}
Let
\begin{equation}\nonumber
f(z)=\somme{n=1}{\infty}a_{n}q^{n}\in S_{k}(U)
\end{equation}
be an eigencuspform with coefficients in a number field $L\overset{\iota}{\plonge}\C$. Let $M$ be the pure motive over $\Q$ and with coefficients in $L$ attached to $f$ (\cite{SchollMotivesModular}). We assume that $L$ embeds in $E$.

Recall that $\Q_{n}/\Q$ is the subfield of $\Q(\zeta_{p^{n+1}})$ with Galois group $G_{n}$ over $\Q$ equal to $\Z/p^{n}\Z$. If $\chi\in\hat{G}_{n}$ is a character, we denote by $L_{\chi}$ the extension of $\Q$ generated by $L$ and the image of $\chi$. Let $\pid|p$ be a prime ideal of $\Ocal_{L_{\chi}}$ and let $E_{\chi}$ be the corresponding finite extension of $E$. We consider $h^{0}(\Spec\Q_{n})$ the Artin motive attached to the regular representation of $G_{n}$ viewed as a pure motive over $\Q$ with coefficients in $L_{\chi}$ and denote by $h^{0}(\Spec\Q_{n})_{\chi}$ the direct summand of $h^{0}(\Spec\Q_{n})$ on which $G_{n}$ acts through $\chi$.

Let $M_{\chi}$ be the motive $M\times_{\Q}h^{0}(\Spec\Q_{n})_{\chi}$. Be definition, a realization of $M_{\chi}$ is the corresponding realization of $M$ with scalars extended from $L$ to $L_{\chi}$ together with an action of $G_{n}$ through $\chi$. We denote by $V_{\C}$ (resp. $V_{\chi,\C}$) the Betti realization of $M$ (resp. $M_{\chi}$), by $V_{\dR}$ (resp. $V_{\chi,\dR}$) the de Rham realization of $M$ (resp. $M_{\chi}$) and by $V$ (resp. $V_{\chi}$) the $\pid$-adic étale realization of $M$ (resp. $M_{\chi}$). 
\begin{Def}
For $S$ a set of finite primes containing $\{p\}$, the $S$-partial $L$-function $L_{S}(M,\chi,s)$ is the holomorphic complex function satisfying 
\begin{equation}\nonumber
L_{S}(M,\chi,s)\eqdef\produit{\ell\notin S}{}\frac{1}{1-a_{\ell}\chi(\Fr(\ell))\ell^{-s}+\epsi(\ell)\chi(\Fr(\ell))\ell^{k-1-2s}}
\end{equation}
for all $s\in\C$ with $\Re s>>0$. This is also the $S$-partial $L$-function of the motive $M_{\chi}$.
\end{Def}
Note that the $S$-partial $L$-function $L_{S}(M,{\chi},s)$ depends on the choice of $\iota$ though we suppress this dependence from the notation. We fix $S$ a finite set of primes containing $p$, an integer $1\leq r\leq k-1$, an integer $n\in\N$ and a character $\chi\in\hat{G}_{n}$.

The Betti-de Rham comparison isomorphism of $\C[\Gal(\C/\R)]$-modules
\begin{equation}\nonumber
V_{\dR}(r)\tenseur_{\Q}\C\isocan V_{B}(r)\tenseur_{\Q}\C 
\end{equation}
induces a complex period map
\begin{equation}\label{EqPerC}
\per_{\C}:\Fil^{0}V_{\chi,\dR}(r)\tenseur_{L_{\chi}}\C\fleche V_{\chi,\C}(r-1)^{+}\tenseur_{L_{\chi}}\C.
\end{equation}
Under our hypothesis on $r$, $\Fil^{0}V_{\dR}(r)$ is equal to the $L$-vector subspace $S_{k}(U)(f)$ on which $\Hecke_{k}(U)$ acts through $\la(f)$ and \eqref{EqPerC} is an isomorphism (\cite{DeligneFonctionsL}). Besides, the composition of localization at $p$ with the dual exponential map $\exp^{*}$ of \cite{BlochKato}
\begin{equation}\nonumber
\exp^{*}:\Hun(G_{\qp(\zeta_{p^{n}})},V(r))\fleche D_{\dR}^{0}(V(r))
\end{equation}
induces an inverse $p$-adic period map of $E_{\chi}$-vector spaces
\begin{equation}\label{EqPerP}
\per^{-1}_{p}:\Hun_{\et}(\Z[1/S],V_{\chi}(r))\fleche\Fil^{0}V_{\chi,\dR}(r)\tenseur_{L}E_{\chi}
\end{equation}
which is equivariant under the action of $G_{n}$ on both sides. According to \cite{Rohrlich} and \cite[Theorems 12.4,14.5]{KatoEuler}, for all $n\in\N$ and all $\chi\in\hat{G}_{n}$ except possibly finitely many, $\Hun_{\et}(\Z[1/S],V_{\chi}(r))$ is of dimension 1 and \eqref{EqPerP} is an isomorphism. We say that $M_{\chi}(r)$ is strictly critical when this holds (\cite[Section 3.2.6]{KatoViaBdR}).

Suppose that $M_{\chi}$ is strictly critical. The determinant functor applied to the $p$-adic period map \eqref{EqPerP} then yields an isomorphism of free $E_{\chi}$-vector spaces of rank 1
\begin{equation}\nonumber
\xymatrix{
\Det_{E_{\chi}}\Hun_{\et}(\Z[1/S],V_{\chi}(r))\ar[d]_{\isocan}^{\per^{-1}_{p}}\\
\Fil^{0}V_{\chi,\dR}(r)\tenseur_{L}E_{\chi}.
}
\end{equation}
Taking the tensor product with the determinant of $V_{\chi}(r-1)^{+}$ yields an identification  
\begin{equation}\label{EqPerDetP}
\xymatrix{
\Det_{E_{\chi}}\Hun_{\et}(\Z[1/S],V(r))\tenseur_{E_{\chi}}\Det^{-1}_{E_{\chi}}V_{\chi}(r-1)^{+}\ar[d]_{\isocan}^{\per^{-1}_{p}}\\
\Det_{E_{\chi}}\left(\Fil^{0}V_{\chi,\dR}(r)\tenseur_{L}E_{\chi}\right)\tenseur_{E_{\chi}}\Det^{-1}_{E_{\chi}}V_{\chi}(r-1)^{+}.
}
\end{equation}
Similarly, the determinant functor applied to the complex period map \eqref{EqPerC} induces an identification
\begin{equation}\label{EqPerDetC}
\xymatrix{
\left[\Det^{}_{\C}\Fil^{0}V_{\chi,\dR}(r)\tenseur_{L_{\chi}}\C\right]\tenseur\left[\Det^{-1}_{\C}V_{\chi,\C}(r-1)^{+}\tenseur_{L_{\chi}}\C\right]\ar[d]_{\isocan}^{\per_{\C}}\\
\C.
}
\end{equation}  
\begin{Def}
The \emph{motivic fundamental line} $\left(\Delta_{\mot}(M_{\chi}(r)),\per_{p},\per_{\C}\right)$ of the strictly critical motive $M_{\chi}(r)$ is the one-dimensional $L_{\chi}$-vector space
\begin{equation}\label{EqQsub}
\Delta_{\mot}(M_{\chi}(r))\eqdef\Det_{L_{\chi}}\Fil^{0}V_{\chi,\dR}(r)\tenseur_{L_{\chi}}\Det^{-1}_{L_{\chi}}V_{\chi,\C}(r-1)^{+}
\end{equation}
together with the two isomorphisms 
\begin{equation}\nonumber
\per_{p}:\Delta_{\mot}(M_{\chi}(r))\tenseur_{L_{\chi}}E_{\chi}\isom\Det_{E_{\chi}}\Hun_{\et}(\Z[1/S],V_{\chi}(r))\tenseur_{E_{\chi}}\Det^{-1}_{E_{\chi}}(V_{\chi}(r-1))^{+}
\end{equation}
\begin{equation}\nonumber
\per_{\C}:\Delta_{\mot}(M_{\chi}(r))\tenseur_{L_{\chi}}\C\isom\C.
\end{equation}
\end{Def}
Note that motivic fundamental line of $M_{\chi}(r)$ is an $L$-rational subspace both of the target of \eqref{EqPerDetP} and of the source of \eqref{EqPerDetC}. Consequently, to any element $\z$ of the source of \eqref{EqPerDetP} whose image through $\per^{-1}_{p}$ lands in the motivic fundamental line $\Delta_{\mot}(M_{\chi}(r))$ is attached a complex number $\per_{\C}(\per^{-1}_{p}(\z)\tenseur1)$. Let $\per_{\C}(\per^{-1}_{p}(-)\tenseur1)$ be the map defined by this composition on the sub-$L_{\chi}$-subspace $\Delta_{\mot}(M_{\chi}(r))_{\pid}$ equal to the image of $\Delta_{\mot}(M_{\chi}(r))$ through $\per_{p}$. Note also that the motivic fundamental line depends on the choice of $S$.
\paragraph{Remark:}For $n\in\N$, let $M_{n}$ be the motive $M\times h^{0}(\Spec\Q_{n})$ and denote by $V_{n,\C}$, $V_{n,\dR}$ and $V_{n}$ its various realizations. It is conjectured (\cite{KatoHodgeIwasawa}) that there are motivic cohomology groups
\begin{equation}\nonumber
H^{1}_{f}(M_{n}(r)),H^{1}_{f}(M_{n}^{*}(1-r))
\end{equation}
which are finitely generated $L[G_{n}]$-modules of finite projective dimension and an equivariant motivic fundamental line $\Delta_{\mot}(M_{n}(r))$
defined by
\begin{align}\nonumber
\left(\Det_{L[G_{n}]}H^{1}_{f}(M_{n}(r))\right)\tenseur_{L[G_{n}]}\left(\Det^{-1}_{L[G_{n}]}H^{1}_{f}(M_{n}^{*}(1-r))\right)\tenseur_{L[G_{n}]}\\\nonumber
\left(\Det_{L[G_n]}\Fil^{0}V_{n,\dR}(r)\right)\tenseur_{L[G_n]}\left(\Det^{-1}_{L[G_n]}V_{n,\C}(r-1)^{+}\right)
\end{align}
with a canonical isomorphism
\begin{equation}\nonumber
\Delta_{\mot}(M_{n}(r))\tenseur_{L[G_{n}]}L_{\chi}\isocan\Delta_{\mot}(M_{\chi}(r))
\end{equation}
induced by isomorphisms
\begin{equation}\nonumber
H^{1}_{f}(M_{n}(r))\tenseur_{L[G_{n}],\chi}L_{\chi}\simeq0,\ H^{1}_{f}(M_{n}^{*}(1-r))\tenseur_{L[G_{n}],\chi}L_{\chi}\simeq0
\end{equation}
whenever $M_{\chi}$ is strictly critical.
\subsubsection{Zeta morphism for modular motives}
As in the previous subsection, $S$ is a finite set of finite primes containing $p$, $n$ is an integer and $\chi\in\hat{G}_{n}$ is such that $M_{\chi}(r)$ is strictly critical. By definition, the $E_{\chi}$-vector spaces $V_{\chi}(r-1)^{+}$ and $\Hun_{\et}(\Z[1/S],V_{\chi}(r))$ are then isomorphic.

Let
\begin{equation}\nonumber
Z:V_{\chi}(r-1)^{+}\isom\Hun_{\et}(\Z[1/S],V_{\chi}(r))
\end{equation}
be an isomorphism between them. Applying the functor $\Det$ then yields an identification
\begin{equation}\nonumber
\Det_{E_{\chi}}(Z):\Det_{E_{\chi}}\Hun_{\et}(\Z[1/S],V_{\chi}(r))\tenseur_{E_{\chi}}\Det^{-1}_{E_{\chi}}V_{\chi}(r-1)^{+}\simeq E_{\chi}.
\end{equation}
Let $\Delta_{L,Z}(M_{\chi})$ be the $L_{\chi}$-vector space pre-image of $L_{\chi}\subset E_{\chi}$ through this isomorphism.

Let $T(f)_{\Iw}$ be $T(f)\tenseur_{\Ocal}\Oiwa$ for $T(f)$ a $G_{\Q}$-stable lattice inside $V$. Suppose there is a non-zero morphism
\begin{equation}\nonumber
Z_{\Iw}:T(f)_{\Iw}(-1)^{+}\fleche\Hun(\Z[1/S],T(f)_{\Iw})\simeq\Hun(\Z[1/p],T(f)_{\Iw}).
\end{equation}
For all $n\in\N$ and all $\chi\in\hat{G}_{n}$ except possibly finitely many, $Z_{\Iw}$ induces a non-zero morphism 
\begin{equation}\nonumber
Z_{\chi}:V_{\chi}(r-1)^{+}\fleche\Hun(\Z[1/S],V_{\chi}(r)).
\end{equation}
Outside of this finite set and of the finite set of characters $\chi$ such that $M_{\chi}$ is not strictly critical, the source and target of $Z_{\chi}$ are $E_{\chi}$-vector spaces of dimension 1 and $Z_{\chi}$ is a non-zero map between them, hence an isomorphism.
\begin{Def}
A morphism 
\begin{equation}\nonumber
Z:V_{\chi}(r-1)^{+}\isom\Hun_{\et}(\Z[1/S],V_{\chi}(r))
\end{equation}
is \emph{motivic} if $\Delta_{L,Z}(M_{\chi})$ is equal to $\per_{p}\left(\Delta_{\mot}(M_{\chi})\right)$. If a morphism $Z$ is motivic, we say it is the \emph{$S$-partial zeta morphism} of $M_{\chi}$ if $\per_{\C}\left(\per_{p}^{-1}\left(\Det_{E_{\chi}}(Z)^{-1}(1)\right)\tenseur1\right)$ is equal to $L_{S}(M^{*}(1),\chi^{-1},-r)\in\C$. A morphism 
\begin{equation}\nonumber
Z_{\Iw}:T(f)_{\Iw}(-1)^{+}\fleche\Hun(\Z[1/S],T(f)_{\Iw}).
\end{equation}
is the \emph{$S$-partial zeta morphism} of $T(f)_{\Iw}$ if $Z_{\chi}$ is the $S$-partial zeta morphism of $V_{\chi}$ for all $\chi$ such that $M_{\chi}$ is strictly critical. 
\end{Def}
Note that the composition $\per_{\C}\left(\per_{p}^{-1}\left(\Det_{E_{\chi}}(Z)^{-1}(1)\right)\tenseur1\right)$ makes sense as the pre-image of $\Det_{E_{\chi}}(Z)^{-1}(1)$ through $\per_{p}$ belongs to $\Delta_{\mot}(M_{\chi}(r))$ if (and only if) $Z$ is motivic. A zeta morphism is unique, if it exists.

As recalled in the following theorem, one of the major contribution of K.Kato towards the Equivariant Tamagawa Number Conjecture for motives attached to classical eigencuspforms is that they admit zeta morphisms (\cite[Theorem 12.5]{KatoEuler}).
\begin{Theo}\label{TheoKato}
Let $\lambda_{f}:\Hsmr\fleche\Ocal$ be a classical point. For all finite set of primes $S\supset\{p\}$, there exists a zeta morphism
\begin{equation}\nonumber
Z_{S}:T(f)_{\Iw}(-1)^{+}\fleche\Hun_{\et}(\Z[1/S],T(f)_{\Iw}).
\end{equation}
We write $\zsfiwa$ for the $\Sigma$-partial zeta morphism and $\z(f)_{\Iw}$ for the $\{p\}$-partial zeta morphism.
\end{Theo}
\subsection{Zeta elements over deformation rings}
\subsubsection{Review of the results of Nakamura}
In this section, we review the work of K.Nakamura (\cite{NakamuraUniversal}) on the construction of zeta elements over deformation rings. Recently P.Colmez and S.Wang announced the construction of the same object by a rather different (and sophisticated) approach (\cite{ColmezWang}). We strengthen slightly assumption \ref{HypTW}.
\begin{Hyp}\label{HypNakamura}
If $(\rhobar|G_{\qp})^{ss}=\chibar_{1}\oplus\chibar_{2}$, then $\chibar_{1}^{-1}\chi_{2}\neq\bar{\chi}_{\cyc}$.
\end{Hyp}
Compared to assumption \ref{HypTW}, assumption \ref{HypNakamura} adds the slight restriction that $\rhobar|G_{\qp}$ cannot be an extension of the trivial character by the cyclotomic character (up to twist).
\begin{Theo}[K.Nakamura]\label{TheoNakamura}
Under assumptions \ref{HypTW} and \ref{HypNakamura}, there exists a zeta morphism
\begin{equation}\nonumber
\zs:\Ts(-1)^{+}\fleche\Hun_{\et}(\Z[1/\Sigma],\Ts)
\end{equation}
such that for all classical point $\lambda_{f}$ the following diagram commutes 
\begin{equation}\nonumber
\xymatrix{
\Ts(-1)^{+}\ar[d]\ar[rr]^{\zs}&&\Hun_{\et}(\Z[1/\Sigma],\Ts)\ar[d]\\
T(f)_{\Iw}(-1)^{+}\ar[rr]^{\zsfiwa}&&\Hun_{\et}(\Z[1/\Sigma],T(f)_{\Iw})
}
\end{equation}
where $\zsfiwa$ is the zeta morphism of theorem \ref{TheoKato}.
\end{Theo}
\begin{proof}
This is \cite[Theorem 1.1]{NakamuraUniversal}.
\end{proof}
\paragraph{Remark:}In \cite{NakamuraUniversal}, the supplementary assumption $p\geq5$ is made. In this remark, we briefly explain why we can dispense from this assumption. The assumption $p\geq5$ is required in the proofs of theorem 4.7 and theorem B.24 of \textit{loc. cit.}. There, it is made in order to appeal respectively to \cite[Corollary 6.5]{PaskunasDuke} and to \cite[Proposition 6.3]{PaskunasMontreal}, which are stated in \cite{PaskunasMontreal,PaskunasDuke} only under the hypothesis $p\geq5$. However, the relevant parts of \cite[Corollary 6.5]{PaskunasDuke} and \cite[Proposition 6.3]{PaskunasMontreal} have been generalized to all $p$ in \cite[Proposition 2.13,2.33]{PaskunasANT}.\footnote{We thank K.Nakamura and S-N.Tung for explaining this to us.}

\paragraph{}

If $\lambda:\Hsmr\fleche S$ is a specialization, we denote by
\begin{equation}\nonumber
\zs(\lambda):T_{\lambda}(-1)^{+}\fleche\Hun_{\et}(\Z[1/\Sigma],T_{\lambda})
\end{equation}
the morphism induced by the surjection $\Ts(-1)^{+}\fleche T_{\lambda}(-1)^{+}$ and the morphism
\begin{equation}\nonumber
\lambda_{*}:\Hun_{\et}(\Z[1/\Sigma],T)\fleche\Hun_{\et}(\Z[1/\Sigma],T_{\lambda}).
\end{equation}
If $\lambda_{f}:\Hr\fleche\Ocal$ is classical up to a twist, we write $\zsf$ for $\zs(\lambda_{f})$. Suppose $\lambda$ has values in $S$ and that $x_{1}\in\Hsmr$ belongs to the kernel of $\lambda$. To $\lambda$ is then attached a specialization $\lambda_{\Iw}:\Hsmr\fleche S_{\Iw}$ by identifying $\Hsmr$ with $\left(\Hsmr/x_{1}\right)[[\Gamma_{\Iw}]]$ and setting $\lambda_{\Iw}$ equal to the identity on $\Gamma_{\Iw}$. We write $\zs(\lambda)_{\Iw}$ for the zeta element attached to $\lambda_{\Iw}$. Theorem \ref{TheoNakamura} ensures that $\zsfiwa$ is equal to $\zs(\lambda_{f})_{\Iw}$ if $\lambda_{f}:\Hsmr\fleche\Ocal$ is a classical point.
\subsubsection{Iwasawa-suitable specializations}
Henceforth, we consistently assume that $\mgot_{\rhobar}$ satisfies the following supplementary assumption.
\begin{Hyp}\label{HypEnL}
There exists $\ell||N$, $\ell\not=p$ such that $\dim_{\Fp}\rhobar_{f}^{I_{\ell}}=1$.
\end{Hyp}
Let $\lambda:\Hsmr\fleche A$ be a specialization with values in a domain $A$. Recall that $(V_{\la},\rho_{\la},\Frac(A))$ is the $G_{\Q,\Sigma}$-representation $\Tla\tenseur_{A}\Frac(A)$.
\begin{Def}
For $\ell\nmid p$ a prime, the Euler factor $\Eul_{\ell}(\Tla)$ of $\Tla$ is $\det\left(1-\Fr(\ell)|\Vla^{I_{\ell}}\right)$.
\end{Def}
Let $\aid\in\Spec^{\min}\Hr$ be a minimal primer ideal through which $\la$ factors.
\begin{Def}
A specialization $\lambda:\Hsmr\fleche A$ with values in a domain is \emph{relatively pure at $\ell\nmid p$} if the $A$-rank of $H^{0}(G,\Vla)$ is equal to the $\Raid$-rank of $H^{0}(G,\Vaid)$ for $G$ equal either to $G_{\ql}$ or to $I_{\ell}$. It is \emph{relatively pure} if it is relatively pure at $\ell$ for all $\ell\nmid p$.
\end{Def}
Note that in particular if $\la$ is a relatively pure specialization at $\ell$, then 1 is not an eigenvalue of $\Fr(\ell)$ acting on $V^{I_{\ell}}$.

The following proposition follows from the proof of the Weight-Monodromy Conjecture for modular motives (\cite{SaitoMonodromie}).
\begin{Prop}\label{PropMonodromie}
A classical specialization is relatively pure.
\end{Prop}
\begin{proof}
See for instance \cite[Lemma 3.9]{FouquetDihedral}.
\end{proof}
We say that a specialization $\psi:\Hsmr\fleche A$ contains a specialization $\la:\Hsmr\fleche A$ if there exists a morphism of local rings $\la':S\fleche A$ such that the diagram
\begin{equation}\nonumber
\xymatrix{
\Hsmr\ar[d]_{\psi}\ar[r]^{\la}& A\\
S\ar[ru]_{\la'}&
}
\end{equation}
is commutative. 
\begin{Lem}\label{LemEulerCommute}
Let $\psi:\Hr\fleche S$ be a specialization with values in a domain and containing a specialization $\la$ with values in a domain $A$ which is relatively pure (at $\ell$). Then $\psi$ is relatively pure (at $\ell$), the formation of $\Eul_{\ell}(-)$ commutes with all relatively pure specialization contained in $\psi$ and for all such specialization $\Eul_{\ell}(-)$ does not vanish. This applies in particular to the quotient map $\psi:\Hr\fleche\Raid$ for all $\aid\in\Spec^{\min}\Hr$. 
\end{Lem}
\begin{proof}
Let $\ell\nmid p$ be a prime. By construction, $\rank_{S}H^{0}(G_{\ql},V_{\psi})\leq \rank_{A}H^{0}(G_{\ql},\Vla)$ so $H^{0}(G_{\ql},V_{\psi})=0$. Likewise, $\rank_{A}H^{0}(I_{\ell},\Vaid)\leq\rank_{S}H^{0}(I_{\ell},V_{\psi})\leq\rank_{A} H^{0}(I_{\ell},\Vla)$ so $\rank_{S}H^{0}(I_{\ell},V_{\psi})=\rank_{A}H^{0}(I_{\ell},\Vaid)$. So $\psi$ is relatively pure. It is then enough to show that 
\begin{equation}\label{EqEqEqEuler}
\la\left(\Eul_{\ell}(\Taid)\right)=\Eul_{\ell}(\Tla)
\end{equation}
and that 
\begin{equation}\nonumber
\Eul_{\ell}(\Tla)\neq0
\end{equation}
to conclude. Both members of \eqref{EqEqEqEuler} are equal to $\det\left(1-\la(\Fr(\ell))|\Vla^{I_{\ell}}\right)$. The determinant $\det(1-\Fr(\ell)|\Vla^{I_{\ell}})$ vanishes only if $\Vla^{G_{\ql}}\neq\{0\}$. A specialization for which this holds cannot be relatively pure by definition.

The last assertion follows from the previous ones since $\psi:\Hr\fleche\Raid$ contains a classical, hence relatively pure, specialization by proposition \ref{PropDense}.
\end{proof}
\begin{Def}
A specialization $\lambda:\Hsmr/(x_{1})\fleche A$ with values in a domain $A$ is \emph{\Iwagood} if it contains a specialization $\psi:\Hsmr/(x_{1})\fleche\Ocal$ which satisfies the following two properties.
\begin{enumerate}
\item $\Tpsi$ and $\Tpsi^{*}(1)$ are relatively pure.
\item The morphism $\zs(\psi)_{\Iw}$ is non-zero.
\end{enumerate}
We say that a specialization $\lambda:A\fleche S$ is \Iwagood if there exists an Iwasawa-suitable specialization $\lambda'$ such that the diagram
\begin{equation}\nonumber
\xymatrix{
\Hsmr\ar[r]^(0.6){\lambda'}\ar[d]&S\\
A\ar[ru]_{\lambda}&
}
\end{equation}
commutes.
\end{Def}
Note that an \Iwagood specialization is itself relatively pure.
\begin{Lem}\label{LemEuler}
Let $\la:\Hr\fleche\Ocal$ be an \Iwagood specialization. Let $V_{\la,\Iw}$ be the $G_{\Q,\Sigma}$-representation $\Tlaiwa\tenseur_{\Oiwa}\Frac(\Oiwa)$. Then $\Hun_{\et}(\Z[1/p],\Tlaiwa)$ is a free $\Oiwa$-module of rank 1 and the complex 
\begin{equation}\nonumber
\Cone\left(V_{\lambda,\Iw}(-1)^{+}\overset{\z(\lambda)_{\Iw}}{\fleche}\RGamma_{\et}\left(\Z[1/p],V_{\lambda,\Iw}\right)\right)
\end{equation}
is acyclic.
\end{Lem}
\begin{proof}
By assumption \ref{HypTW}, $H^{0}_{\et}(\Z[1/\Sigma],T_{\lambda,\Iw}/x)$ and $H^{0}_{\et}(\Z[1/\Sigma],T_{\lambda,\Iw}/(x,y))$ vanish for all regular sequence $(x,y)\in\Oiwa^{2}$. The isomorphisms
\begin{equation}\nonumber
\RGamma_{\et}\left(\Z[1/\Sigma],T_{\la,\Iw}\right)\Ltenseur_{\Oiwa}\Oiwa/x\simeq\RGamma_{\et}\left(\Z[1/\Sigma],T_{\la,\Iw}/x\right)
\end{equation}
and
\begin{equation}\nonumber
\RGamma_{\et}\left(\Z[1/\Sigma],T_{\la,\Iw}/x\right)\Ltenseur_{\Oiwa/x}\Oiwa/(x,y)\simeq\RGamma_{\et}\left(\Z[1/\Sigma],T_{\la,\Iw}/(x,y)\right)
\end{equation}
then imply first that $H^{1}_{\et}(\Z[1/\Sigma],T_{\lambda,\Iw})$ is a torsion-free $\Oiwa$-module, then that there is an embedding of $H^{1}_{\et}(\Z[1/\Sigma],T_{\lambda,\Iw})/x$ into $H^{1}_{\et}(\Z[1/\Sigma],T_{\lambda,\Iw}/x)$ and hence that $y$ is a regular element in $H^{1}_{\et}(\Z[1/\Sigma],T_{\lambda,\Iw})/x$. Hence $H^{1}_{\et}(\Z[1/\Sigma],T_{\lambda,\Iw})$ is of depth 2 as $\Oiwa$-module and hence an $\Oiwa$-module free of rank 1 by the theorem of Auslander-Buchsbaum and Serre.

Let $\gamma$ be a generator of the free, rank 1 $\Oiwa$-module $T_{\lambda,\Iw}(-1)^{+}$. According to \cite[Theorem 1.1 (ii)]{NakamuraUniversal}, $\zs(\lambda)_{\Iw}(\gamma)$ is then the first class of an Euler system, which is necessarily non-zero since $\la$ is \Iwagood by assumption. As we have just seen, this entails that this Euler system is not $\Oiwa$-torsion. Assumptions \ref{HypTW} and \ref{HypEnL} imply that the image of $\rhobar$ acts irreducibly on $\Fpbar_{p}^{2}$ and has order divisible by $p$, so contains a subgroup conjugated to $\SL_{2}(\Fp_{q})$ for some $q=p^{n}$, so contains a unipotent element $\bar{\s}\neq\Id$. Let $\s$ be a lifting of $\bar{\s}$ to $\rho_{\lambda}(G_{\Q,\Sigma})$. Then the kernel of $\s-1$ is strictly included inside $T_{\lambda}$ and its cokernel is of dimension 1 after tensor-product with the residual field $k$. Hence this kernel is free of rank 1. The representation $T_{\lambda}$ thus satisfies assumptions $(\operatorname{i_{{str}}})$, $\operatorname{(ii_{{str}}})$ and $(\operatorname{iv}_{\pid})$ of \cite[Theorem 0.8]{KatoEulerOriginal}. As $T_{\lambda}$ is not abelian, the conclusion of \cite[Proposition 8.7]{KatoEulerOriginal} also holds. As obtaining this proposition was the sole use of hypothesis $(\operatorname{iii})$ of \cite[Theorem 0.8]{KatoEulerOriginal}, this theorem also holds for $T_{\lambda,\Iw}$. Consequently $H^{2}_{\et}(\Z[1/p],T_{\lambda,\Iw})$ is $\Oiwa$-torsion. The cohomology of the complex
\begin{equation}\nonumber
C=\Cone\left(V_{\lambda,\Iw}(-1)^{+}\overset{\zs(\lambda)_{\Iw}}{\fleche}\RGamma_{\et}\left(\Z[1/p],V_{\lambda,\Iw}\right)\right)
\end{equation}
is concentrated in degree 1 and 2,
\begin{equation}\nonumber
H^{1}(C)=\Hun_{\et}(\Z[1/p],T_{\lambda,\Iw})/\zs(\lambda)_{\Iw}(\gamma)
\end{equation}
and
\begin{equation}\nonumber
H^{2}(C)=H^{2}_{\et}(\Z[1/p],T_{\lambda,\Iw}).
\end{equation}
The complex $C\Ltenseur_{\Oiwa}\Frac(\Oiwa)$ is thus acyclic. 
\end{proof}

In the following propositions, we extend the results of Nakamura to the case of primitive zeta morphisms over irreducible components of the universal deformation space. We remark that it seems impossible to construct primitive zeta morphisms over the whole universal deformation space. If $p\nmid n$, we denote by $\s_{n}$ the arithmetic Frobenius morphism in $\Gal(\Q(\zeta_{p^{\infty}})/\Q)$ or $\Gal(\Q_{\infty}/\Q)$.
\begin{Prop}\label{PropNakamura}
Let $\aid\in\Spec\Hr$ be a minimal prime ideal. Let $c,d$ be integers such that $\ell\nmid cd$ if $\ell\in\Sigma\cup\{2,3\}$. Under assumptions \ref{HypTW} and \ref{HypNakamura}, there exists a zeta morphism
\begin{equation}\nonumber
_{c,d}\zaid:\Taid(-1)^{+}\fleche\Hun_{\et}(\Z[1/p],\Taid)
\end{equation}
such that the following diagram commutes 
\begin{equation}\nonumber
\xymatrix{
\Taid(-1)^{+}\ar[d]\ar[rr]^{_{c,d}\zaid}&&\Hun_{\et}(\Z[1/p],\Taid)\ar[d]\\
T(f)_{\Iw}(-1)^{+}\ar[rr]^{_{c,d}\ziwa}&&\Hun_{\et}(\Z[1/p],T(f)_{\Iw})
}
\end{equation}
for all classical point $\lambda_{f}:\Hr\fleche\Ocal$ factoring through $\Raid$. Here, $_{c,d}\ziwa$ designates any of the morphism appearing in the statement of \cite[Theorem 12.6]{KatoEuler}.

\end{Prop}
\begin{proof}
We use the notations and results of section \ref{AppendixCompleted}.

Let $N(\aid)$ be the tame Artin conductor of the representation $\rho(\aid)$, or equivalently by relative purity the tame conductor of any classical point factoring through $\Raid$. Let $U_{1}(\aid)_{\ell}$ be the compact open subgroup
\begin{equation}\nonumber
U_{1}(\aid)_{\ell}=\left\{\matricetype\in\GL_{2}(\Z_{\ell})|\matricetype\equiv\matrice{*}{*}{0}{1}\modulo \ell^{v_{\ell}(N(\aid))}\right\}.
\end{equation}
and let $U_{1}(\aid)$ be the compact open subgroup
\begin{equation}\nonumber
U_{1}(\aid)\eqdef\produit{\ell\in\Sp}{}U_{1}(\aid)_{\ell}\subset\produit{\ell\in\Sp}{}\GL_{2}(\ql).
\end{equation}
We also view $U_{1}(\aid)$ as a compact open subgroup of $\GL_{2}(\A_{\Q}^{(p\infty)})$, in which case $\Uaidl$ is taken to be compact open maximal if $\ell\notin\Sigma$.

Let $U_{p}\subset \GL_{2}(\qp)$ be a compact open subgroup. In \cite[Sections 2.2, 5.1]{KatoEuler} and \cite[Sections 2.3.6, 2.4.2]{FukayaKatoSharifi}, elements 
\begin{equation}\nonumber
_{c,d}\z_{U_{p}}\in K_{2}(Y(\Uaid U_{p}))
\end{equation}
are constructed (they are denoted there $_{c,d}z_{1,Mp^{r},p^{r}}(u,v)$). These elements are compatible for the norm map 
\begin{equation}\nonumber
K_{2}(Y(\Uaid U'_{p}))\fleche K_{2}(Y(\Uaid U_{p}))
\end{equation}
attached to the covering $Y(\Uaid U'_{p}))\fleche Y(\Uaid U_{p}))$ if $U'_{p}\subset U_{p}$ (\cite[Sections 2.4.8]{FukayaKatoSharifi}). After localizing at $\mgot_{\rhobar}$ and taking the Chern class map, they thus yield morphisms
\begin{equation}\label{EqKatoSharifi}
_{c,d}\zt:\Htilde_{\et}^{1}(\Uaid,\Ocal)_{\mgot_{\rhobar}}(-1)\fleche\Hun_{\et}(\Z[1/\Sigma],\Htilde_{\et}^{1}(\Uaid,\Ocal)_{\mgot_{\rhobar}})
\end{equation}
from the inverse limit on the level of cohomology groups in the tower of modular curves with full level structure at $p$ and tame level $\Uaid$ to the first cohomology group of $\Spec\Z[1/\Sigma]$ with coefficients in this completed cohomology group.

By theorem \ref{TheoEmerton}, we have the following description of the completed cohomology with full level at all primes in $\Sigma$ as $\Hr[G_{\Q,\Sigma}\times\GL_{2}(\qp)\times\produit{\ell\in\Sp}{}\GL_{2}(\ql)]$-module
\begin{equation}\nonumber
\Htilde_{\et}^{1}(\Ocal)_{\mgot_{\rhobar}}\simeq\rho_{\Sigma}\tenseur\pitilde_{p}(\rho_{\Sigma}|G_{\qp})\hat{\tenseur}\produittenseur{\ell\in\Sp}{}\pitilde_{\ell}(\rho_{\Sigma}|G_{\ql}).
\end{equation}
Taking $\Raid$-quotient, $\Uaid$-coinvariants and applying proposition \ref{PropEssentialVectors}, we obtain an isomorphism
\begin{equation}\nonumber
\Htilde_{\et}^{1}(\Uaid,\Ocal)_{\mgot_{\rhobar}}\tenseur_{\Hr}\Raid\simeq\rho(\aid)\tenseur\pitilde_{p}(\rho(\aid)|G_{\qp}).
\end{equation}
Hence \eqref{EqKatoSharifi} induces a morphism of $\Hr[\GL_{2}(\qp)]$-morphism
\begin{equation}\nonumber
_{c,d}\zaidt:\Taid(-1)^{+}\tenseur\pitilde_{p}(\rho(\aid)|G_{\qp})\fleche\Hun_{\et}(\Z[1/\Sigma],\Taid)\tenseur\pitilde_{p}(\rho(\aid)|G_{\qp}).
\end{equation}
Using the description of the (dual of) the so-called Montréal functor (\cite[Section IV.2.3]{ColmezFoncteur}) given in \cite{PaskunasMontreal}, it is shown in \cite[Appendix B]{NakamuraUniversal} that there is an isomorphism
\begin{equation}\nonumber
\End(\pitilde_{p}(\rho_{\Sigma}|G_{\qp}))\simeq \Hr
\end{equation}
where endomorphisms are taken in the category of $\Ocal[\GL_{2}(\qp)]$-modules which are also compact as $\Ocal[[\GL_{2}(\zp)]]$-modules and Pontrjagin dual to a locally admissible $\Ocal[\GL_{2}(\qp)]$-module.\footnote{Here again, we note that the supplementary assumption $p\geq5$ appears in \cite{NakamuraUniversal}. We can dispense from it by modifying the proof exactly in the same way as explained in the remark following the statement of theorem \ref{TheoNakamura}.} The base-change properties of Colmez's Montréal functor and the same proof as in \cite[Appendix B]{NakamuraUniversal} then show that 
\begin{equation}\nonumber
\End(\pitilde_{p}(\rho(\aid)|G_{\qp}))\simeq \Raid.
\end{equation}
There is thus a bijection between
\begin{equation}\nonumber
\Hom_{\Raid[\GL_{2}(\qp)]}\left(\Taid(-1)^{+}\tenseur\pitilde_{p}(\rho(\aid)|G_{\qp}),\Hun_{\et}(\Z[1/\Sigma],\Taid)\tenseur\pitilde_{p}(\rho(\aid)|G_{\qp})\right)
\end{equation}
and
\begin{equation}\nonumber
\Hom_{\Raid[\GL_{2}(\qp)]}\left(\Taid(-1)^{+},\Hun_{\et}(\Z[1/\Sigma],\Taid)\right).
\end{equation}
Let
\begin{equation}\nonumber
_{c,d}\zaid:\Taid(-1)^{+}\fleche\Hun_{\et}(\Z[1/\Sigma],\Taid)
\end{equation}
be the morphism image of $_{c,d}\zaidt$ through this bijection. Because $\Taid\simeq(\Taid/x_{1})[[\Gamma_{\Iw}]]$, the morphism $_{c,d}\zaid$ can be viewed as having values in $\Hun_{\et}(\Z[1/p],\Taid)$. Let $\la(f):\Hr\fleche\Ocal$ be a classical point factoring through $\Raid$. Then $\pitilde_{p}(\rho(\aid)|G_{\qp})\tenseur_{\Raid,\la(f)}\Ocal$ is isomorphic to $\pitilde_{p}(\rho_{f}|G_{\qp})$. Hence, the morphism
\begin{equation}\nonumber
_{c,d}\zaid\tenseur_{\la}1:T(f)_{\Iw}(-1)^{+}\fleche\Hun_{\et}(\Z[1/p],T(f)_{\Iw})
\end{equation}
induced by $_{c,d}\zaid$ coincides with the morphism $_{c,d}\z(f)_{\Iw}$ (which would be written $_{c,d}z_{1,M,p^{\infty}}$ in \cite[Section 3.1.2]{FukayaKatoSharifi}).
\end{proof}
Let $\mu\in\Oiwa$ be the element $(c-\s_{c})(d-\s_{d})$. Let
\begin{equation}\label{EqDefZaid}
\zaid:\Taid(-1)^{+}\fleche\Hun_{\et}(\Z[1/p],\Taid)\tenseur_{\Raid}\Raid[1/\mu]
\end{equation}
be the morphism $(c-\s_{c})^{-1}(d-\s_{d})^{-1}_{c,d}\zaid$. For $\la:\Hr\fleche A$ factoring through $\Raid$, we denote by $_{c,d}\z(\la)_{\Iw}$ the morphism
\begin{equation}\nonumber
_{c,d}\z(\la)_{\Iw}:\Tlaiwa(-1)^{+}\fleche\Hun_{\et}(\Z[1/p],\Tlaiwa).
\end{equation}
The following proposition, which follows closely the strategy of proof of Nakamura, establishes the integrality of zeta morphisms over irreducible components of $\Hr$. Its results are crucial in the proof of proposition \ref{PropDivETNC}, and hence in the proof of our main result.
\begin{Prop}\label{PropIntegral}
Let $\aid\in\Spec\Hr$ be a minimal prime ideal. Under assumptions \ref{HypTW} and \ref{HypNakamura},  the zeta morphism \eqref{EqDefZaid} defines a morphism
\begin{equation}\nonumber
\zaid:\Taid(-1)^{+}\fleche\Hun_{\et}(\Z[1/p],\Taid)
\end{equation}
such that there is an equality
\begin{equation}\label{EqZaidZs}
\zaid=\left(\produit{\ell\in\Sp}{}\Eul_{\ell}(\Taid^{*}(1))\right)^{-1}\zs(\aid)
\end{equation}
of morphisms $\Taid(-1)^{+}\fleche\Hun_{\et}(\Z[1/p],\Taid)$ and such that the following diagram commutes 
\begin{equation}\nonumber
\xymatrix{
\Taid(-1)^{+}\ar[d]\ar[rr]^{\zaid}&&\Hun_{\et}(\Z[1/p],\Taid)\ar[d]\\
T(f)_{\Iw}(-1)^{+}\ar[rr]^{\ziwa}&&\Hun_{\et}(\Z[1/p],T(f)_{\Iw})
}
\end{equation}
for all classical point $\lambda_{f}:\Hr\fleche\Ocal$ factoring through $\Raid$. Here, $\ziwa$ is the zeta morphism of theorem \ref{TheoKato}. In particular, if $\la:\Hr\fleche S$ is an \Iwagood specialization factoring through $\Raid$, then there is an equality
\begin{equation}\nonumber
\z(\la)_{\Iw}=\left(\produit{\ell\in\Sp}{}\Eul_{\ell}(\Tlaiwa^{*}(1))\right)^{-1}\zs(\la_{\Iw})
\end{equation}
of morphisms $\Tlaiwa(-1)^{+}\fleche\Hun_{\et}(\Z[1/p],\Tlaiwa)$.
\end{Prop}
\begin{proof}
We first show that the morphism
\begin{equation}\nonumber
\zaid:\Taid(-1)^{+}\fleche\Hun_{\et}\left(\Z[1/p],T(\aid)\right)\tenseur_{\Raid}\Frac(\Raid)
\end{equation}
defined by \eqref{EqDefZaid} has values in $\Hun_{\et}(\Z[1/p],\Taid)$. 

For $U_{p}\subset\GL_{2}(\qp)$ a compact open subgroup, we define
\begin{equation}\nonumber
\z_{U_{p}}=(c-\s_{c})^{-1}(d-\s_{d})^{-1}_{c,d}\z_{U_{p}}.
\end{equation}
Let $\Hs(U_{1}(\aid)U_{p})_{\mgot_{\rhobar},\Iw}$ be the Hecke algebra $\Hs(U_{1}(\aid)U_{p})_{\mgot_{\rhobar}}[[\Gamma_{\Iw}]]$. By the previous propositions, we may view $\z_{U_{p}}$ as a morphism of $\Oiwa$-modules
\begin{equation}\nonumber
\Hun_{\et}(Y(U_{1}(\aid)U_{p}),\Ocal)_{\mgot_{\rhobar},\Iw}(-1)\fleche\Hun_{\et}(\Z[1/\Sigma],\Hun_{\et}(Y(U_{1}(\aid)U_{p}),\Ocal)_{\mgot_{\rhobar},\Iw})\tenseur_{\Oiwa}\Oiwa[1/\mu]
\end{equation}
equivariant under the action of $\Hs(U_{1}(\aid)U_{p})_{\mgot_{\rhobar},\Iw}$ on both sides. For simplicity of notation, we write $M_{U,\rhobar,\Iw}$ for $\Hun_{\et}(Y(U_{1}(\aid)U_{p}),\Ocal)_{\mgot_{\rhobar},\Iw}$. As before, the morphisms $\z_{U_{p}}$ are compatible with the inclusion $U'_{p}\subset U_{p}$ in the sense that the diagram
\begin{equation}\nonumber
\xymatrix{
M_{U_{p}',\rhobar,\Iw}(-1)^{+}\ar[r]^(0.4){\z_{U'_{p}}}\ar[d]&\Hun_{\et}(\Z[1/p],M_{U_{p}',\rhobar,\Iw}[1/\mu])\ar[d]\\
M_{U_{p},\rhobar,\Iw}(-1)^{+}\ar[r]^(0.4){\z_{U_{p}}}&\Hun_{\et}(\Z[1/p],M_{U_{p},\rhobar,\Iw}[1/\mu])
}
\end{equation} 
is commutative. In order to show that $\zaid$ has values in $\Hun_{\et}(\Z[1/p],\Taid)$, it is thus enough to show that $\z_{U_{p}}$ has values in $\Hun_{\et}(\Z[1/p],M_{U_{p},\rhobar,\Iw})$ for $U_{p}$ small enough. As
\begin{equation}\nonumber
M_{U,\rhobar}\tenseur_{\Ocal}E=\sommedirecte{f}{}\rho_{f}\tenseur\pi(f)^{U_{1}(\aid)U_{p}}
\end{equation}
where the sum is over all $f\in S_{2}(U_{1}(\aid)U_{p})$ attached to classical points of $\Hs(\Uaid U_{p})_{\mgot_{\rhobar}}$, it is enough to show that $\z_{U_{p}}$ has values in $\Hun_{\et}(\Z[1/p],M_{U_{p},\rhobar,\Iw})$ after projection to the eigenspace $M_{U_{p},\rhobar}[f]$ corresponding to the choice of one such $f$. As $\pi(f)^{U_{1}(\aid)U_{p}}\neq0$, it is an irreducible $\Hs(\Uaid U_{p})_{\mgot_{\rhobar}}[1/p]$-module (\cite[Proposition 4.3]{BushnellHenniart}). Hence, it is enough to show that the image of any non-zero vector of $M_{U_{p},\rhobar}[f]$ belongs to $\Hun_{\et}(\Z[1/p],M_{U_{p},\rhobar}[f])$. If $v$ is a vector in $M_{U_{p},\rhobar}^{U_{1}(\aid')U'_{p}}[f]$ such that $f$ is new of level $U_{1}(\aid')U'_{p}$, this is established in \cite[Theorem 12.5,12.6]{KatoEuler}.

If $\la_{f}:\Hr\fleche\Ocal$ is a classical point, then $_{c,d}\zs(f)_{\Iw}$ satisfies
\begin{equation}\nonumber
_{c,d}\zs(f)_{\Iw}=(c-\s_{c})(d-\s_{d})\left(\produit{\ell\in\Sp}{}\Eul_{\ell}(T(f)_{\Iw}^{*}(1))\right)\z(f)_{\Iw}
\end{equation}
 by construction. The remaining assertions of the proposition thus follows by specializations from theorem \ref{TheoNakamura}, lemma \ref{LemEulerCommute} and the density of classical points. 
\end{proof}
\begin{Prop}\label{PropZetaDense}
Let $\lambda:\Hr\fleche A$ be \Iwagood. Let $\Vlaiwa$ be $\Tlaiwa\tenseur_{A}\Frac(A_{\Iw})$. Then the complex
\begin{equation}\nonumber
\Cone\left(V_{\lambda,\Iw}(-1)^{+}\overset{\z(\lambda)_{\Iw}}{\fleche}\RGamma_{\et}\left(\Z[1/p],V_{\lambda,\Iw}\right)\right)
\end{equation}
is acyclic. The set of specializations $\lambda:\Hsmr/(x_{1})\fleche\Ocal$ which are not \Iwagood is of codimension at least 1.
\end{Prop}
\begin{proof}
Let $\pid$ be the kernel of an \Iwagood specialization with values in $\Ocal$ contained in $\la$. Let $\ell\nmid p$ be a prime. If $H^{0}(I_{\ell},\Tla\tenseur_{A}A_{\pid})$ vanishes, then so does $H^{0}(I_{\ell},\Tla\tenseur_{A}A_{\pid}/\pid)$ by relative purity at $\ell$ and
\begin{equation}\nonumber
H^{0}(I_{\ell},\Tla\tenseur_{A}A_{\pid})\tenseur_{A_{\pid}}\kg(\pid)\simeq H^{0}(I_{\ell},\Tla\tenseur_{A}\kg(\pid)).
\end{equation}
If $H^{0}(I_{\ell},\Tla\tenseur_{A}A_{\pid})$ has rank 2, then $\Tla$ is unramified at $\ell$, and so is $\Tla\tenseur_{A}\kg(\pid)$, so that again
\begin{equation}\nonumber
H^{0}(I_{\ell},\Tla\tenseur_{A}A_{\pid})\tenseur_{A_{\pid}}\kg(\pid)\simeq H^{0}(I_{\ell},\Tla\tenseur_{A}\kg(\pid)).
\end{equation}
If finally $H^{0}(I_{\ell},\Tla\tenseur_{A}A_{\pid})$ is a rank 1 $A_{\pid}$-module, then $H^{0}(I_{\ell},\Tla\tenseur_{A}A_{\pid})\tenseur_{A_{\pid}}\kg(\pid)$ and $H^{0}(I_{\ell},\Tla\tenseur_{A}A_{\pid})\tenseur_{A_{\pid}}\Frac(A_{\pid})$ have the same ranks so $H^{0}(I_{\ell},\Tla\tenseur_{A}A_{\pid})$ is free of rank 1 by Nakayama's lemma. There is thus an isomorphism of complexes
\begin{equation}\nonumber
[H^{0}(I_{\ell},\Tla\tenseur_{A}A_{\pid})\overset{1-\Fr(\ell)}{\fleche}H^{0}(I_{\ell},\Tla\tenseur_{A}A_{\pid})]\Ltenseur_{A_{\pid}}\kg(\pid)\simeq[H^{0}(I_{\ell},\Tla\tenseur_{A}\kg(\pid))\overset{1-\Fr(\ell)}{\fleche}H^{0}(I_{\ell},\Tla\tenseur_{A}\kg(\pid))]
\end{equation}
for all $\ell\nmid p$. This family of isomorphisms induces an isomorphism
\begin{equation}\nonumber
\RGamma_{\et}\left(\Z[1/p],\Tla\tenseur_{A}A_{\pid}\right)\Ltenseur_{A_{\pid}}\kg(\pid)\simeq\RGamma_{\et}\left(\Z[1/p],\Tla\tenseur_{A}\kg(\pid)\right).
\end{equation}
In particular, $H^{2}(\Z[1/p],\Tla)$ is $A$-torsion and $\Hun(\Z[1/p],\Tla)$ has $A$-rank 1. This implies that the complex
\begin{equation}\nonumber
\Cone\left(V_{\lambda,\Iw}(-1)^{+}\overset{\z(\lambda)_{\Iw}}{\fleche}\RGamma_{\et}\left(\Z[1/p],V_{\lambda,\Iw}\right)\right)
\end{equation}
is acyclic.

We show that the set of specializations which are not \Iwagood has large codimension. \cite[Lemma 3.9]{FouquetDihedral} implies that $\lambda:\Hr/(x_{1})\fleche\Ocal$ is not relatively pure at $\ell$ only if $\la$ factors through $\Raid$ such that $\Vaid$ is generically special Steinberg and if the $\ell$-adic monodromy of $\Vla$ is trivial. The set of such specializations has codimension at least 1 in each irreducible component of $\Hr$, and the same is then true of the union of these sets. To conclude, it is thus enough to show that the set of spécializations $\lambda:\Hsmr/(x_{1})\fleche\Ocal$ such that $\zs(\lambda)$ is non-zero has codimension at least 1.

Let $\lambda({f}):\Hsmr\fleche\Ocal$ be a classical point of $\Hr$. The isomorphism
\begin{equation}\nonumber
\RGamma_{\et}\left(\Z[1/\Sigma],\Ts\right)\Ltenseur_{\Hsmr,\lambda(f)_{\Iw}}\Ocal_{\Iw}\simeq\RGamma_{\et}(\Z[1/\Sigma],T(f)_{\Iw})
\end{equation}
and the fact that both complexes are acyclic outside degree 1 and 2 shows that there is an isomorphism of $\Oiwa$-modules
\begin{equation}\nonumber
H^{2}_{\et}(\Z[1/\Sigma],\Ts)\tenseur_{\Hsmr,\lambda(f)_{\Iw}}\Oiwa\simeq H^{2}_{\et}(\Z[1/\Sigma],T(f)_{\Iw}).
\end{equation}
In particular, $H^{2}_{\et}(\Z[1/\Sigma],\Ts)$ is $\Hsmr$-torsion. Let $a\in\Hsmr$ a regular element such that $aH^{2}_{\et}(\Z[1/\Sigma],\Ts)=0$. Let $A$ be the localization of $\Hsmr$ at $\{a^{n}|n\in\N\}$. The complex
\begin{equation}\nonumber
\RGamma_{\et}(\Z[1/\Sigma],\Ts)\Ltenseur_{\Hsmr}A\simeq \RGamma_{\et}(\Z[1/\Sigma],\Ts\tenseur_{\Hsmr}A)
\end{equation}
is then a perfect complex of $A$-modules acyclic outside degree 1. As the cohomology of
\begin{equation}\nonumber
\RGamma_{\et}(\Z[1/\Sigma],\Ts)\Ltenseur_{\Hsmr}A/\mgot_{A}\simeq\RGamma_{\et}(\Z[1/\Sigma],\Ts\tenseur_{\Hsmr}A/\mgot_{A})
\end{equation}
vanishes in degree 0 and 2 (in degree 0 by assumption \ref{HypTW}, in degree 2 by Nakayama's lemma), the $A$-module $\Hun_{\et}(\Z[1/\Sigma],\Ts\tenseur_{\Hsmr}A)$
is a free $A$-module, necessarily of rank 1. Let $u$ be the morphism
\begin{equation}\nonumber
u:(\Ts\tenseur_{\Hsmr}A)(-1)^{+}\fleche\Hun_{\et}(\Z[1/\Sigma],\Ts\tenseur_{\Hsmr}A)
\end{equation}
which sends a generator of the free $A$-module $(\Ts\tenseur_{\Hsmr}A)(-1)^{+}$ to a generator of the free $A$-module $\Hun_{\et}(\Z[1/\Sigma],\Ts\tenseur_{\Hsmr}A)$. Then the morphism
\begin{equation}\nonumber
\zs\tenseur1:(\Ts\tenseur_{\Hsmr}A)(-1)^{+}\fleche\Hun_{\et}(\Z[1/\Sigma],\Ts\tenseur_{\Hsmr}A)
\end{equation}
may be written $\zs\tenseur1=\alpha u$ for some $\alpha\in A$. Let $\psi:\Hsmr\fleche\Ocal$ be such that $\psi(a)\neq0$ and let $\psi':A\fleche E$ be the morphism through which $\psi$ factors by the universal property of localization. Then 
\begin{equation}\nonumber
\RGamma_{\et}(\Z[1/\Sigma],\Ts\tenseur_{\Hsmr}A)\Ltenseur_{A,\psi'}E\simeq\RGamma_{\et}(\Z[1/\Sigma],T_{\psi}\tenseur_{\Ocal}E)
\end{equation}
so $H^{2}_{\et}(\Z[1/\Sigma],T_{\psi}\tenseur_{\Ocal}E)$ vanishes, $\Hun_{\et}(\Z[1/\Sigma],T_{\psi}\tenseur_{\Ocal}E)$ is free of rank 1 and
\begin{equation}\nonumber
\Hun_{\et}(\Z[1/\Sigma],\Ts\tenseur_{\Hsmr}A)\tenseur_{A,\psi'}E\simeq\Hun_{\et}(\Z[1/\Sigma],T_{\psi}\tenseur_{\Ocal}E).
\end{equation}
This means that $\psi'_{*}(\zs\tenseur1)=\zs(\psi)$ is non-zero. The set of specializations  $\lambda:\Hsmr\fleche\Ocal$ such that  $\zs(\lambda)$ vanishes is thus of codimension at least 1. This implies that the set of specializations $\lambda:\Hsmr/(x_{1})\fleche\Ocal$ such that $\z(\lambda)_{\Iw}$ vanishes is also of codimension at least 1.
\end{proof}

\subsection{Fundamental lines}
Let $A$ be a quotient of $\Hsm$ or $\Lambdaf$. In particular, $A$ could be $\Hsm$, $\Raid$, $\Lambdaf$, $\Ocal_{\Iw}$, or $\Ocal$. Let $T$ be $\Ts\tenseur_{\Hsm}A$ or $\Ts\tenseur_{\Lambdaf}A$ accordingly. 


To $T$ is attached an étale sheaf over $\Spec\Z[1/\Sigma]$. The following lemma is well-known.
\begin{Lem}\label{LemPerfect}
The complexes $\RGamma_{\et}(\Z[1/\Sigma],T)$ and $\RGamma(G_{\ql},T)$, for all $\ell$, are perfect and commute with arbitrary derived base-change of coefficients.
\end{Lem}
\begin{proof}
See for example \cite[(4.2.9) Proposition]{SelmerComplexes}.
\end{proof}
\begin{Def}\label{DefFun}
The $\Sigma$-partial fundamental line $\Delta_{\Sigma}(T)$ is the graded invertible $A$-module
\begin{equation}\nonumber
\Delta_{\Sigma}(T)\eqdef\Det^{-1}_{A}\RGamma_{\et}\left(\Z[1/\Sigma],T\right)\tenseur_{A}\Det^{-1}_{A}T(-1)^{+}.
\end{equation}
Suppose in addition that $A$ is a domain. For $\ell\in\Sigma$ prime to $p$, the graded invertible module $\Xcali_{\ell}(T)$ is 
\begin{equation}\nonumber
\Xcali_{\ell}(T)=\begin{cases}
\Det_{A}\RGamma(G_{\ql}/I_{\ell},T^{I_{\ell}})&\textrm{ if $\rank_{A}T^{I_{\ell}}\neq1$,}\\
\Det_{A}[A\overset{1-\alpha_{\ell}}{\fleche} A]&\textrm{ if $\rank_{A}T^{I_{\ell}}=1$.}
\end{cases}
\end{equation}
Here, the complex $[A\overset{1-\alpha_{\ell}}{\fleche}A]$ is placed in degree $0,1$. The fundamental line $\Delta_{A}(T)$ is the graded invertible $A$-module
\begin{equation}\nonumber
\Delta(T)\eqdef\Det^{-1}_{A}\RGamma_{\et}\left(\Z[1/\Sigma],T\right)\tenseur_{A}\Det^{-1}_{A}T(-1)^{+}\tenseur_{A}\produittenseur{\ell\in\Sigma}{}\left(\Xcali_{\ell}(T)\tenseur_{A}\Det^{-1}_{A}\RGamma(G_{\ql},T)\right).
\end{equation}
\end{Def}
Note that the previous definition makes sense in all cases thanks to lemma \ref{LemPerfect}. It follows from the same lemma that any morphism $\psi:A\fleche A'$ induces a canonical isomorphism 
\begin{equation}\nonumber
\Delta_{A}(T)\tenseur_{A}A'\isocan\Delta_{A'}(T').
\end{equation}
\begin{Def}
A specialization $\la:\Lambdaf/(x_{1})\fleche A$ with values in a domain $A$ is \emph{\Iwagood} if it contains a specialization $\psi:\Lambdaf/(x_{1})\fleche\Ocal$ such that all specialization $\psi':\Hr/(x_{1})\fleche\Ocal$ above $\psi$ are \Iwagood. We say that a specialization $\lambda:A\fleche S$ is \Iwagood if there exists an Iwasawa-suitable specialization $\lambda'$ such that the diagram
\begin{equation}\nonumber
\xymatrix{
\Lambdaf\ar[r]^(0.6){\lambda'}\ar[d]&S\\
A\ar[ru]_{\lambda}&
}
\end{equation}
commutes.

\end{Def} 
\begin{Prop}\label{PropMonodromy}
Let $A$ be a domain. The fundamental line $\Delta_{A}(T)$ is canonically isomorphic to
\begin{equation}\nonumber
\Det^{-1}_{A}\RGamma_{\et}\left(\Z[1/p],T\right)\tenseur_{A}\Det^{-1}_{A}T(-1)^{+}
\end{equation}
when both are defined. Suppose that $\psi:A\fleche A'$ is an \Iwagood specialization which is contained in $\la$. Then the natural map
\begin{equation}\nonumber
\Delta_{A}(T)\tenseur_{A,\psi}A'{\fleche}\Delta_{A'}(T')
\end{equation}
is an isomorphism.
\end{Prop}
\begin{proof}
If $\Det^{-1}_{A}\RGamma_{\et}\left(\Z[1/p],T\right)$ is defined, then $\RGamma_{\et}\left(\Z[1/p],T\right)$ is a perfect complex of $A$-modules, so $T^{I_{\ell}}$ is a perfect complex of $A$-modules for al $\ell\in\Sigma$. In that case, taking a projective resolution of of $T^{I_{\ell}}$ induces a canonical isomorphism between $\Det_{A}[T^{I_{\ell}}\overset{1-\Fr(\ell)}{\fleche}T^{I_{\ell}}]$ and $\Xcali_{\ell}(T)$. The properties of the $\Det$ functor and the collection of these canonical isomorphisms for all $\ell\nmid p$ induce a canonical isomorphism between 
\begin{equation}\nonumber
\Det^{-1}_{A}\RGamma_{\et}\left(\Z[1/p],T\right)\tenseur_{A}\Det^{-1}_{A}T(-1)^{+}
\end{equation}
first with
\begin{equation}\nonumber
\Det^{-1}_{A}\Cone\left(\RGamma_{\et}(\Z[1/\Sigma],T)\oplus\sommedirecte{\ell\in\Sigma}{}\RGamma(G_{\ql},T)\fleche\sommedirecte{\ell\in\Sigma}{}\RGamma(G_{\ql}/I_{\ell},T^{I_{\ell}})\right)[-1]\tenseur_{A}\Det^{-1}_{A}T(-1)^{+}.
\end{equation}
then with
\begin{equation}\nonumber
\Det^{-1}_{A}\RGamma_{\et}\left(\Z[1/\Sigma],T\right)\tenseur_{A}\Det^{-1}_{A}T(-1)^{+}\tenseur_{A}\produittenseur{\ell\in\Sigma}{}\left(\Xcali_{\ell}(T)\tenseur_{A}\Det^{-1}_{A}\RGamma(G_{\ql},T)\right).
\end{equation}
Let $\psi$ be as in the lemma. Then $H^{0}(I_{\ell},T)$ and $H^{0}(I_{\ell},T')$ have the same rank and $\psi(\alpha_{\ell})=\alpha'_{\ell}$ for all $\ell\nmid p$. Hence $\Xcali_{\ell}(T)\tenseur_{A,\psi}A'$ is canonically isomorphic to $\Xcali_{\ell}(T')$. As all other constituents of $\Delta_{A}(T)$ commute with $-\tenseur_{A,\psi}A'$ for general reasons, the collection of the canonical isomorphismes $\Xcali_{\ell}(T)\tenseur_{A,\psi}A'\isocan\Xcali_{\ell}(T')$ induce a canonical isomorphism $\Delta_{A}(T)\tenseur_{A,\psi}A'\isocan\Delta_{A'}(T')$.
\end{proof}
It follows from lemma \ref{LemEuler} that the fundamental lines of \Iwagood specializations come with canonical isomorphisms to the trivial graded invertible module.
\begin{Def}\label{DefTriv}
Let $\la:\Hr\fleche A$ be an \Iwagood specialization with zeta morphism
$$\z(\la)_{\Iw}:\Vlaiwa(-1)^{+}\isom\Hun_{\et}(\Z[1/p],\Vlaiwa).$$
The trivialization map
$$\triv_{\z(\la)_{\Iw}}:\Delta_{\Aiwa}(\Tlaiwa)\fleche\Frac(\Aiwa)$$
is the composition
\begin{equation}\nonumber
\Delta_{\Aiwa}(\Tlaiwa)\plonge\Delta_{\Frac(\Aiwa)}(V_{\la,\Iw})\isocan\Det_{\Frac(\Aiwa)}\Cone(\z(\lambda)_{\Iw})[-1]\isocan \Frac(\Aiwa).
\end{equation}
where the first canonical isomorphism is
\begin{equation}\nonumber
\RGamma(\Z[1/p],V_{\la,\Iw})\simeq\Hun_{\et}\left(\Z[1/p],V_{\la,\Iw}\right)[-1]
\end{equation}
and where the second follows from the acyclicity of the complex
\begin{equation}\nonumber
\Cone\left(V_{\lambda,\Iw}(-1)^{+}\overset{\z(\lambda)_{\Iw}}{\fleche}\RGamma_{\et}\left(\Z[1/p],V_{\lambda,\Iw}\right)\right).
\end{equation}
If  more generally $\Lambda$ is $\Hr/x_{1}$ or $\Lambdaf/x_{1}$ and if $\la:\Lambda\fleche S$ is an Iwasawa suitable specialization with values in a reduced ring and if $V_{\la,\Iw}$ denotes $T_{\la,\Iw}\tenseur_{S_{\Iw}}Q(\Siwa)$, the trivialization map
$$\triv_{\zs(\liwa)}:\Delta_{\Si}(\Tlaiwa)\fleche Q(\Siwa)$$
is the composition
\begin{equation}\label{EqTrivSigma}
\Delta_{\Si}(\Tlaiwa)\plonge\Delta_{\Si}(V_{\la,\Iw})\isocan\Det_{\Si}\Cone(\zs(\liwa))[-1]\isocan Q(\Siwa).
\end{equation}
\end{Def}
Note that if $\la:\Lambda\fleche S$ is \Iwagood, then the isomorphism 
\begin{equation}\nonumber
\RGamma_{\et}(\Z[1/\Si],\Tlaiwa)\Ltenseur_{\Siwa,\psi}\Oiwa
\end{equation}
for $\psi$ an \Iwagood specialization shows that $H^{2}_{\et}(\Z[1/\Si],\Tlaiwa)$ is $\Siwa$-torsion so that the last isomorphism in \eqref{EqTrivSigma} is well defined.
\begin{Prop}\label{PropDetCom}
Let $\la:\Hr\fleche A$ be an \Iwagood specialization and let $\psi:A\fleche S$ be an \Iwagood specialization contained in $\la$. The canonical isomorphism
\begin{equation}\nonumber
\Delta_{\Aiwa}(\Tlaiwa)\tenseur_{\Aiwa,\psi}\Siwa{\fleche}\Delta_{\Siwa}(\Tpsiwa)
\end{equation}
of proposition \ref{PropMonodromy} fits into a commutative diagram
\begin{equation}\nonumber
\xymatrix{
\Delta_{\Aiwa}(\Tlaiwa)\ar[d]_{-\tenseur_{\Aiwa,\psi}\Siwa}\ar[rr]^{\triv_{\z(\la)_{\Iw}}}&&\frac{x}{y}\Aiwa\ar[d]^{\psi}\\
\Delta_{\Siwa}(\Tpsiwa)\ar[rr]^(0.55){\triv_{\z(\psi)_{\Iw}}}&&\frac{x'}{y'}\Siwa
}.
\end{equation}
In particular, the morphism $\psi$ extends to a morphism $\frac{x}{y}\Aiwa\fleche\frac{x'}{y'}\Siwa$.
\end{Prop}
\begin{proof}
Let $\pid\subset\Spec\Aiwa$ be the kernel of $\psi:\Aiwa\fleche\Siwa$. We have seen in the proof of proposition  \ref{PropZetaDense} that $\RGamma_{\et}(\Z[1/p],\Tlaiwa\tenseur_{\Aiwa}\Aiwap)$ is a perfect complex which commutes with $-\Ltenseur_{\Aiwap}\kg(\pid)$ and that the complex
\begin{equation}\nonumber
\Cone\left((\Tlaiwa\tenseur_{\Aiwa}\Aiwap)(-1)^{+}\overset{\z(\la)_{\Iw}}{\fleche}\RGamma_{\et}(\Z[1/p],\Tlaiwa\tenseur_{\Aiwa}\Aiwap)\right)
\end{equation}
is acyclic. By functoriality of $\Det$, there is a commutative diagram
\begin{equation}\nonumber
\xymatrix{
\Delta_{\Aiwa}(\Tlaiwa)\ar@{^{(}->}[r]\ar[d]_{-\tenseur_{\Aiwa}\Siwa}&\Delta_{\Aiwap}(\Tlaiwa\tenseur_{\Aiwa}\Aiwap)\ar[r]^(0.7){\sim}\ar[d]_{-\tenseur_{A_{\pid}}\kg(\pid)}&\Aiwap\ar[d]^{-\tenseur_{A_{\pid}}\kg(\pid)}\\
\Delta_{\Siwa}(\Tpsiwa)\ar@{^{(}->}[r]&\Delta_{\kg(\pid)}(\Vpsiwa)\ar[r]^(0.6){\sim}&\kg(\pid)
}
\end{equation}
whose vertical arrows are induced by $\psi$. In particular, if $\triv_{\z(\la)_{\Iw}}\left(\Delta_{\Aiwa}(\Tlaiwa)\right)$ is generated by $x/y$, then $y$ may be chosen so that it does not belong to $\pid$ and $\psi(x)/\psi(y)$ generates $\triv_{\z(\psi)_{\Iw}}\left(\Delta_{\Siwa}(\Tpsiwa)\right)$.
\end{proof}

The preceding lemma admits a variant which applies to $\Sigma$-partial fundamental lines.
\begin{Prop}\label{PropDetComSigma}
Let $\la:\Hr\fleche A$ be an \Iwagood specialization and let $\psi:A\fleche S$ be an \Iwagood specialization contained in $\la$. The canonical isomorphism
\begin{equation}\nonumber
\Delta_{\Si}(\Tlaiwa)\tenseur_{\Aiwa,\psi}\Siwa{\isom}\Delta_{\Si}(\Tpsiwa)
\end{equation}
fits into a commutative diagram
\begin{equation}\nonumber
\xymatrix{
\Delta_{\Si}(\Tlaiwa)\ar[d]_{-\tenseur_{\Aiwa,\psi}\Siwa}\ar[rr]^{\triv_{\zs(\liwa)}}&&\frac{x}{y}\Aiwa\ar[d]^{\psi}\\
\Delta_{\Si}(\Tpsiwa)\ar[rr]^(0.55){\triv_{\zs(\psiwa)}}&&\frac{x'}{y'}\Siwa
}.
\end{equation}
In particular, the morphism $\psi$ extends to a morphism $\frac{x}{y}\Aiwa\fleche\frac{x'}{y'}\Siwa$.
\end{Prop}
\begin{proof}
As $\RGamma_{\et}(\Z[1/\Si],-)$ commutes with arbitrary change of coefficients and sends perfect complexes to perfect complexes, the proof is the same as that of proposition \ref{PropDetCom} only easier.
\end{proof}

\subsection{The Equivariant Tamagawa Number Conjectures}\label{SubETNCuniv}
The following conjecture is Kato's statement of the Iwasawa Main Conjecture of classical Iwasawa theory (\cite{KatoHodgeIwasawa}).
\begin{Conj}\label{ConjIMC}
Let $\lambda(f)_{\Iw}:\Hsmr\fleche\Oiwa$ be the specialization attached to a classical point $\lambda(f)$. Then the trivialization morphism
\begin{equation}\nonumber
\triv_{\z(f)_{\Iw}}:\Delta_{\Oiwa}(T(f)_{\Iw})\plonge\Frac(\Oiwa)
\end{equation}
given by the zeta morphism
\begin{equation}\nonumber
\z(f)_{\Iw}:T(f)_{\Iw}(-1)^{+}\fleche\Hun_{\et}(\Z[1/p],T(f)_{\Iw})
\end{equation}
induces an isomorphism
\begin{equation}\nonumber
\triv_{\z(f)_{\Iw}}:\Delta_{\Oiwa}(T(f)_{\Iw})\isocan\Oiwa.
\end{equation}
\end{Conj}
In \cite{KatoViaBdR}, the preceding conjecture was extended to arbitrary $\Lambda$-adic families of motives. In our situation of interest, we obtain the following conjecture.
\begin{Conj}\label{ConjETNC}
Let $\lambda:\Hsmr\fleche S$ be an \Iwagood specialization. Then the trivialization morphism
\begin{equation}\nonumber
\triv_{\zs(\lambda_{\Iw})}:\Delta_{\Sigma}(\Tlaiwa)\plonge Q(\Siwa)
\end{equation}
given by the zeta morphism
\begin{equation}\nonumber
\zs(\liwa):\Tlaiwa(-1)^{+}\fleche\Hun_{\et}(\Z[1/\Sigma],\Tlaiwa)
\end{equation}
induces an isomorphism
\begin{equation}\nonumber
\triv_{\zs(\lambda)}:\Delta_{\Sigma}(T_{\lambda})\isocan\Siwa.
\end{equation}
Assume moreover that $S$ is a domain. Then the trivialization morphism
\begin{equation}\nonumber
\triv_{\z(\lambda)_{\Iw}}:\Delta_{\Siwa}(\Tlaiwa)\plonge\Frac(\Siwa)
\end{equation}
given by the zeta morphism
\begin{equation}\nonumber
\z(\lambda)_{\Iw}:\Tlaiwa(-1)^{+}\fleche\Hun_{\et}(\Z[1/p],\Tlaiwa)
\end{equation}
induces an isomorphism
\begin{equation}\nonumber
\triv_{\z(\lambda)_{\Iw}}:\Delta_{\Siwa}(\Tlaiwa)\isocan\Siwa.
\end{equation}
\end{Conj}
Note in particular the following special cases of interests.
\begin{Conj}\label{ConjUniv}
The trivialization morphism
\begin{equation}\nonumber
\triv_{\zs}:\Delta_{\Sigma}(\Ts)\plonge Q(\Hsmr)
\end{equation}
given by the zeta morphism
\begin{equation}\nonumber
\zs:\Ts(-1)^{+}\fleche\Hun_{\et}(\Z[1/\Sigma],\Ts)
\end{equation}
induces an isomorphism
\begin{equation}\nonumber
\triv_{\zs}:\Delta_{\Sigma}(\Ts)\isocan\Hr.
\end{equation}
\end{Conj}
In view of proposition \ref{PropDense}, conjecture \ref{ConjUniv} is the most general statement on the $p$-adic variation of special values of $L$-functions that can be formulated for modular motives.
\begin{Conj}\label{ConjUnivDomain}
Let $\z(\aid)$ be the morphism $\z(\lambda)$ for $\lambda:\Hsmr\fleche\Raid$. The trivialization morphism
\begin{equation}\nonumber
\triv_{\z(\aid)}:\Delta_{\Raid}(T_{\lambda})\plonge\Faid)
\end{equation}
given by the zeta morphism
\begin{equation}\nonumber
\z(\aid):T(\aid)(-1)^{+}\fleche\Hun_{\et}(\Z[1/p],T(\aid))
\end{equation}
induces an isomorphism
\begin{equation}\nonumber
\triv_{\z(\aid)}:\Delta_{\Raid}(\Taid)\isocan\Raid.
\end{equation}
\end{Conj}
We record the following compatibility property between the conjectures stated so far.
\begin{Prop}\label{PropConjCompatible}
Let $\lambda:\Hsmr\fleche S$ be an \Iwagood specialization (resp. with values in a domain). Assume that conjecture \ref{ConjUniv} is true (resp. conjecture \ref{ConjUnivDomain} for $\Raid$ an irreducible component through which $\lambda$ factors). Then conjecture \ref{ConjETNC} is true for $\lambda$ (resp. the second half of conjecture \ref{ConjETNC} is true for $\lambda$).
\end{Prop}
\begin{proof}
According to proposition \ref{PropDetComSigma}, if $\triv_{\zs}:\Ds(\Ts)\isocan\Hsmr$, then $\triv_{\zs(\lambda)}:\Ds(T_{\lambda})\simeq\Lambda$. Hence, the first half of conjecture \ref{ConjETNC} is true. Next, we establish the second half of this conjecture when $\lambda$ is the quotient map $\Hsmr\fleche\Raid$. First, we record the commutative diagram
\begin{equation}\nonumber
\xymatrix{
\Ds(\Ts)\ar[d]_{-\tenseur_{\Hsmr}\Raid}\ar[rr]^{\triv_{\zs}}&&\Hsmr\ar[d]^{\lambda}\\
\Ds(\Taid)\ar[rr]^{\triv_{\zs(\aid)}}&&\Raid.
}
\end{equation}
Let $\Xcali$ be $\Delta_{\Raid}(\Taid)\tenseur_{\Raid}\Ds(\Taid)^{-1}$. By definition of $\RGamma_{\et}(\Z[1/\Sigma],-)$ and $\RGamma_{\et}(\Z[1/p],-)$, we have
\begin{equation}\nonumber
\Xcali\simeq\produittenseur{\ell\in\Sigma}{}\Xcali_{\ell}(\Taid)\tenseur\Det^{-1}_{\Raid}\RGamma(G_{\ql},\Taid).
\end{equation}
Hence 
\begin{equation}\nonumber
\Xcali\tenseur_{\Raid}\Faid\simeq\produittenseur{\ell\in\Sigma}{}\Xcali_{\ell}(\Vaid^{*}(1))^{-1}
\end{equation}
by Tate's local duality.
Hence $\triv_{\zs(\aid)}$ sends the $\Raid$-submodule $\Delta_{\Raid}(\Taid)$ of $\Ds(\Taid)\tenseur_{\Raid}\Faid$ to $\produit{\ell}{}\Eul_{\ell}(\Vaid^{*}(1))^{-1}\Raid$. As 
\begin{equation}\nonumber
\z(\aid)=\left(\produit{\ell}{}\Eul_{\ell}(\Vaid^{*}(1))^{-1}\right)\zs(\aid),
\end{equation}
there is an isomorphism
\begin{equation}\nonumber
\triv_{\z(\aid)}:\Delta_{\Raid}(\Taid)\simeq\Raid.
\end{equation}
Hence conjecture \ref{ConjETNC} is true in this case.

Now let $\lambda:\Hsmr\fleche S$ be an arbitrary \Iwagood specialization with values in a domain. Then $\z(\lambda)_{\Iw}$ is the image of $\z(\aid)$ through $\lambda$ and $\Delta_{\Raid}(\Taid)\tenseur_{\Raid}\Siwa\simeq\Delta_{\Siwa}(\Tlaiwa)$ so there is an isomorphism
\begin{equation}\nonumber
\triv_{\z(\lambda)_{\Iw}}:\Delta_{\Siwa}(\Tlaiwa)\simeq\Siwa
\end{equation}
and conjecture \ref{ConjETNC} is thus true for $\lambda$.
\end{proof}
The method of Euler system (\cite{PerrinRiouEuler,RubinEuler,KatoEulerOriginal}) yields the following partial result towards conjecture \ref{ConjETNC} when $\la:\Hr\fleche\Ocal$ is an \Iwagood specialization.
\begin{Prop}\label{PropDivETNC}
Let $\la:\Hr\fleche\Ocal$ be an \Iwagood specialization. Then the image of
\begin{equation}\nonumber
\triv_{\z(\la)_{\Iw}}:\Delta_{\Oiwa}(\Tlaiwa)^{-1}\plonge\Frac(\Oiwa)
\end{equation}
is included inside $\Oiwa$.
\end{Prop}
\begin{proof}
This proof relies crucially on the fact that $\z(\la)_{\Iw}$ has values in $\Hun_{\et}(\Z[1/p],\Tlaiwa)$ (rather than $\Hun_{\et}(\Z[1/p],\Tlaiwa)\tenseur_{\Oiwa}\Frac(\Oiwa)$) when $\la:\Hr\fleche\Ocal$ is an \Iwagood specialization, as was proved in proposition \ref{PropIntegral}. As recalled in the proof of lemma \ref{LemEuler}, the representation $\Tlaiwa$ then satisfies all hypotheses of \cite[Theorem 0.8]{KatoEulerOriginal} except possibly $(\operatorname{iii})$ whose role is replaced here by the statement of \cite[Proposition 8.7] {KatoEulerOriginal}. We thus have 
\begin{equation}\nonumber
\carac_{\Oiwa}H^{2}_{\et}(\Z[1/p],\Tlaiwa)\mid\carac_{\Oiwa}H^{1}_{\et}(\Z[1/p],\Tlaiwa)/\image\z(\la)_{\Iw}.
\end{equation}
After localization at a height-one prime, the structure theorem of modules of discrete valuation ring shows that this divisibility is equivalent to the statement that the image of
\begin{equation}\nonumber
\triv_{\z(\la)_{\Iw}}:\Delta_{\Oiwa}(\Tlaiwa)^{-1}\plonge\Frac(\Oiwa)
\end{equation}
is included inside $\Oiwa$.
\end{proof}
We are interested in the following weaker version of conjecture \ref{ConjUniv} in which we view $\Ts$ as a $\Lambdaf$-module.
\begin{Conj}\label{ConjUnivWeak}
The trivialization morphism
\begin{equation}\nonumber
\triv_{\zs}:\Ds(\Ts)\plonge\Frac(\Lambdaf)
\end{equation}
given by the zeta morphism
\begin{equation}\nonumber
\zs:\Ts(-1)^{+}\fleche\Hun_{\et}(\Z[1/\Sigma],\Ts)
\end{equation}
induces an isomorphism
\begin{equation}\nonumber
\triv_{\zs}:\Ds(\Ts)\isocan\Lambdaf.
\end{equation}
\end{Conj}
\section{The residually irreducible, crystalline case}
In this section, we prove conjecture \ref{ConjIMC} for the $G_{\Q}$-representation $(T(f),\rho_{f},\Ocal)$ when $f\in S_{k}(\Gamma_{0}(N))$ is an eigencuspform of even weight $k$ with trivial central character, with coefficients in a number field $F\subset E$, such that $\rho_{f}|G_{\qp}$ is a crystalline and short $G_{\qp}$-representation and such that $\rhobar_{f}|G_{\qp}$ is an irreducible $G_{\qp}$-representation. We write $(T_{f},\rho,\Ocal)$ for the $G_{\Q}$-representation $\left(T(f)\left(k/2\right),\rho_{f}\left(k/2\right),\Ocal\right)$ and let $V_{f}$ be $T_{f}\otimes_{\Ocal}E$.

We record once and for all the following assumptions, that will be crucial in all the section.
\begin{Hyp}\label{HypCrys}
The $G_{\qp}$-representation $\rho|G_{\qp}$ is crystalline and short with irreducible residual representation.
\end{Hyp}
We write $\qpi$ for the cyclotomic $\zp$-extension of $\qp$, $\qppur$ for the unramified $\mathbb{Z}_p$-extension of $\qp$ and $\zppur$ for the unit ball of $\qppur$. Let $\mathbb{Q}_{p,\infty}^{p,{\ur}}$ be the unramified $\mathbb{Z}_p$-extension of $\mathbb{Q}_{p,\infty}$. As before, we let $\Gamma_{\Iw}$ be $\mathrm{Gal}(\mathbb{Q}_\infty/\mathbb{Q})$ and fix one of its topological generator $\gamma$. We write $\Gamma_{p,\infty}^{p,\ur}\simeq\zp^{2}$ for the Galois group $\mathrm{Gal}(\mathbb{Q}_{p,\infty}^{p,{\ur}}/\qp)$. Let ${\Kcal}$ be a quadratic imaginary extension of $\mathbb{Q}$ where $p$ splits as $v_0\bar{v}_0$. Let ${\Kcal}_\infty$ be the $\mathbb{Z}_p^2$-extension of ${\Kcal}$. We write $\Gamma_{\Kcal}$ for the Galois group $\mathrm{Gal}({\Kcal}_\infty/{\Kcal})$. Define $\Lambda_{\Kcal}=\Ocal[[\Gamma_{\Kcal}]]$. Let ${\Kcal}^{v_0}$ and ${\Kcal}^{\bar{v}_0}$ be respectively the $\mathbb{Z}_p$-extension of ${\Kcal}$ unramified outside $v_0$ and unramified outside $\bar{v}_0$. Define $\Gamma_{v_0}=\mathrm{Gal}({\Kcal}^{v_0}/{\Kcal})$ and $\Gamma_{\bar{v}_0}=\mathrm{Gal}({\Kcal}^{\bar{v}_0}/{\Kcal})$. We fix topological generators $\gamma_{v_0}\in\Gamma_{v_0}$ and $\gamma_{\bar{v}_0}\in\Gamma_{\bar{v}_0}$.
\subsection{Review of $p$-adic Hodge theory}
We introduce necessary material in Iwasawa cohomology (including over affinoid rings) and recall standard notions of $p$-adic Hodge theory (\cite{FontaineDeRham,FontaineCorps}). In this section, $(T,\rho,\Ocal)$ and $(V,\rho,E)$ denote $G_{\qp}$-representations.
\subsubsection{Iwasawa Cohomology Groups}
The classical Iwasawa cohomology $H^1_{\mathrm{cl,Iw}}(\mathbb{Q}_{p,\infty}, T)$ is the inverse limit with respect to the co-restriction map.
$$\limproj{\qp\subseteq \mathbb{Q}_{p,n}\subset \mathbb{Q}_{p,\infty}}H^1(G_{\mathbb{Q}_{p,n}}, T).$$
More generally, if $K\subset K_{n}\subset K_{\infty}$ is a tower of abelian separable extensions, we define $H^{i}_{\cl,\Iw}(K_{\infty},T)$ to be the $i$-th cohomology group of the complex $\RGamma_{\cl,\Iw}(K_{\infty},T)$ defined to be the image in the derived category of 
\begin{equation}\nonumber
\limproj{{{n}}}\ \Ccont(G_{K},T\tenseur_{\Ocal}\Ocal[\Gal(K_{n}/K)])
\end{equation}
(see \cite[Section 8]{SelmerComplexes} for details). The Iwasawa cohomology complex $\RGamma_{\cl,\Iw}(K_{\infty},V)$ is the image in the derived category of
\begin{equation}\nonumber
\limproj{{{n}}}\ \Ccont(G_{K},T\tenseur_{\Ocal}\Ocal[\Gal(K_{n}/K)])[1/p]
\end{equation}
for $T$ any $G_{K}$-stable $\Ocal$-lattice inside $V$. 
\begin{Lem}\label{LemCohoLocale}
Assume that $\rhobar$ is irreducible. For $A$ equal to $\Fp,\Ocal$, $\Oiwa$ or $\Ocal[[\Gamma_{1}]]$, the complex $\RGamma(G_{\qp},T\tenseur_{\Ocal}A)$ is a perfect complex of $A$-modules whose cohomology is concentrated in degree 1. Moreover, the $A$-module $H^{1}(G_{\qp},T\tenseur_{\Ocal}A)$ is free of rank 2. In particular, $$H^1_{\mathrm{cl,Iw}}(\mathbb{Q}_{p,\infty}, T)\simeq\Oiwa^{2}$$
and
$$H^1_{\mathrm{cl,Iw}}(\mathbb{Q}_{p,\infty}^{p,{\ur}},T)\simeq \Ocal[[\Gamma_{p,\infty}^{p,\ur}]]^{2}.$$
\end{Lem}
\begin{proof}
As $\RGamma(G_{\qp},-\tenseur_{\Ocal}A)$ sends perfect complexes to perfect complexes and commutes with $-\Ltenseur_{A}S$ for arbitrary $S$, the statements for $A=\Ocal,\Oiwa$ and $\Ocal[[\Gamma_{p,\infty}^{p,\ur}]]$ follow from the statement for $A=\Fp$ by Nakayama's lemma. In that case, $H^0(G_{\qp}, T\tenseur_{\Ocal}\Fp)=H^2(G_{\qp}, T\tenseur_{\Ocal}\Fp)=0$ by the assumption and Tate local duality. Hence, the local Euler characteristic formula
\begin{equation}\nonumber
\frac{\cardinal{H^0(G_{\qp}, T\tenseur_{\Ocal}\Fp)}\cdot \cardinal{ H^2(G_{\qp}, T\tenseur_{\Ocal}\Fp)}}{\cardinal{ \Hun(G_{\qp}, T\tenseur_{\Ocal}\Fp)}}=\cardinal{\Fp}^{-2}
\end{equation}
reduces to $\cardinal{ \Hun(G_{\qp}, T\tenseur_{\Ocal}\Fp)}=\cardinal{\Fp}^2$ and hence to $\Hun(G_{\qp}, T\tenseur_{\Ocal}\Fp)\simeq {\Fp}^2$. 
\end{proof}

\subsubsection{$(\varphi,\Gamma)$-modules}

Let $\mathbb{C}_p$ be the $p$-adic completion of $\bar{\mathbb{Q}}_p$ and let $\mathcal{O}_{\mathbb{C}_p}$ be its unit ball, that is to say the set of elements with $p$-adic norm less than or equal to $1$. Let $\tilde{E}^+$ be the perfect characteristic $p$ ring
\begin{equation}\nonumber
\tilde{E}^+\eqdef\left\{x=(x_{n})_{n\in\N}|x_{n}\in\Ocal_{\Cp}/p\Ocal_{\Cp},\ \forall n\in\N,\ x_{n+1}^{p}=x_{n}\right\}.
\end{equation}
If $x=(x_{n})_{n\in\N}$ belongs to $\tilde{E}^{+}$ and if for all $n\in \N$, $\hat{x}_{n}$ is a lifting of $x_{n}$ to $\Ocal_{\Cp}$, then for all $n\in\N$ the limit
\begin{equation}\nonumber
x^{(n)}=\lim_{k}(\hat{x}_{n+k})^{p^{k}}\in\Ocal_{\Cp}
\end{equation}
does not depend on the choices of the $\hat{x}_{n+k}$. We identify $\tilde{E}^{+}$ with
\begin{equation}\nonumber
\Ocal_{\Cp}^{\flat}\eqdef\left\{x=(x^{(n)})_{n\in\N}|x^{(n)}\in\Ocal_{\Cp},\ \forall n\in\N,\ (x^{(n+1)})^{p}=x^{(n)}\right\}.
\end{equation}
through the map $x\longmapsto(x^{(n)})_{n\in\N}$. If $v_{p}$ is the valuation on $\Cp$ normalized so that $v_{p}(p)=1$, then the valuation
\begin{equation}\nonumber
v(x)=v_{p}(x^{(0)})=\lim_{n}{p^{n}}v_{p}(\hat{x}_{n})
\end{equation}
makes $\tilde{E}^+$ a complete valuation ring (in particular a domain). Let $\tilde{E}$ be the fraction field of $\tilde{E}^+$. 

Fix once and for all non-trivial $p^n$-th roots of unity $\zeta_{p^n}$ with $\zeta_{p^{n+1}}^p=\zeta_{p^n}$ and write $\epsi$ for the element $(1,\zeta_{p},\zeta_{p^{2}},\cdots)\in\Ocal_{\Cp}^{\flat}$. Then
\begin{equation}\nonumber
v(\epsi-1)=\lim_{n}{p^{n}}v_{p}(\zeta_{p^{n}}-1)=\lim_{n}{p^{n}}\frac{1}{(p-1)p^{n-1}}=\frac{p}{p-1}>1.
\end{equation}

We write $\Atilde^{+},\Atilde$ and $[\cdot]$ respectively for the ring of Witt vectors of $\tilde{E}^+$ or $\tilde{E}$ and for the Teichm\"uller lift from $\tilde{E}^+$ to $\Atilde^{+}$ or from $\tilde{E}$ to $\Atilde$. Then $x\in\Atilde^{+}$ may be written in a unique way
\begin{equation}\nonumber
x=\somme{m=0}{\infty}p^{m}[x_{m}]
\end{equation}
with $x_{m}\in\Ocal_{\Cp}^{\flat}$. Let $\p$ be the Frobenius map on $\Etilde^{+}$. We also denote by $\p$ the functorial extension of $\p$ to $\Atilde^{+}$, that is to say the map
\begin{equation}\nonumber
\p\left(x\right)=\somme{m=0}{\infty}p^{m}[x^{p}_{m}].
\end{equation}
J-M.Fontaine observed that the map
\begin{equation}\nonumber
\theta:\Atilde^+\fleche \mathcal{O}_{\mathbb{C}_p}
\end{equation} 
defined either by
\begin{equation}\label{EqThetaUn}
\theta\left((x_{n})_{n\in\N}\right)=\somme{n=0}{\infty}p^{n}x_{n}^{(n)}
\end{equation} 
or by
\begin{equation}\label{EqThetaDeux}
\theta\left(\somme{m=0}{\infty}p^{m}[b_{m}]\right)=\somme{m=0}{\infty}p^{m}b^{(0)}_{m}
\end{equation}
is a surjective ring homomorphism with principal kernel. The ring $B_{\dR}^{+}$ is the completion 
\begin{equation}\nonumber
B_{\dR}^{+}\eqdef\limproj{n\geq 0}\ \Atilde^+[1/p]/(\ker(\theta))^n
\end{equation}
of $\Atilde^{+}$ with respect to $\ker\theta$. Let $\pi\in\Atilde^{+}$ be $[\epsi]-1$. Then $\theta(\pi)=0$. We define
\begin{equation}\nonumber
t=\log(\epsi)=\somme{n\geq 1}{}\frac{(-1)^{n-1}}{n}([\varepsilon]-1)^n\in B_{{\dR}}^+.
\end{equation}
Then $B_{\dR}^{+}$ a discrete valuation ring with maximal ideal $(t)$ and residue field $\mathbb{C}_p$. 
\begin{Lem}
Let $[a]$ be the lift of $a\in\fp$ in $\zp$. Put
\begin{equation}\nonumber
q=\p^{-1}\left(\somme{a\in\mathbb{F}_p}{}[\varepsilon]^{[a]}\right)\in \Atilde^{+}.
\end{equation}
Then $q\in\ker\theta$ but $q/\pi\notin\ker\theta$.
\end{Lem}
\begin{proof}
By \eqref{EqThetaUn}
\begin{equation}\nonumber
\theta(q)=\somme{a\in\fp}{}\zeta_{p}^{a}=0.
\end{equation}
Choose $x\in\Etilde^{+}$ with $x^{(0)}=-p$. Let $\xi\in \Atilde^{+}$ be the element $[x]+p$. Then
\begin{equation}\nonumber
\theta(\xi)=-p+p=0
\end{equation}
by \eqref{EqThetaDeux} and $\xi$ is a generator of the principal ideal $\ker\theta$ by \cite[Proposition 2.4]{FontaineBarsotti}. If $\alpha\in \Atilde^{+}$ satisfies $\theta(\alpha/\pi)=0$, then there exists $\lambda\in\Atilde^{+}$ such that $\alpha=\xi^2\lambda$. So the valuation of the image of $\alpha$ in the residual field $\Cp$ is at least 2. The valuation of the class of $q$ in $\Cp$ is 1.
\end{proof}
\begin{Def}
Let $A_{\crys}^0$ be the divided power envelop of $\Atilde^+$ with respect to $\ker{\theta}$, that is to say the set obtained by adding all elements $a^m/m!$ for all $a\in\ker{\theta}$. Let $A_{\crys}$ and $B_{\crys}$ be the rings $A_{\crys}\eqdef\limproj{n}\ A_{\crys}^0/p^n A_{\crys}^0$ and $B_{\crys}^+{\eqdef}A_{\crys}[1/p]$.
\end{Def}
For $n\in\N$, we write $\gamma_n(x)=\frac{x^n}{n!}$. Let $\Fil^r A_{\crys}$ be $A_{\crys}\cap\Fil^r B_{{\dR}}$ and let $\Fil_p^r A_{\crys}$ be $\{x\in\Fil^r A_{\crys}|\varphi x\in p^r A\}$. Then we have the following lemma.
\begin{Lem}\label{2.16}
Let $a$ be the largest integer such that $(p-1)a<r$. Then for every $x\in\Fil^rA_{\crys}$, $p^a\cdot a!x$ belongs to $\Fil_p^r A_{\crys}$. Moreover, $\Fil^r_pA_{\crys}$ is the associated sub-$\Atilde^+$-module of $A_{\crys}$ generated by ${q}^j\gamma_b(p^{-1}t^{p-1})$ for all integers $j,b$ such that $j+(p-1)b\geq r$.
\end{Lem}
\begin{proof}
This is \cite[Proposition 6.24]{FontaineOuyang}.
\end{proof}
For $0\leq r\leq s\leq\infty$, $r,s\in\mathbb{Q}$, let $\Atilde^{[r,s]}$ be the $p$-adic completion of
\begin{equation}\nonumber
\Atilde^+\left[\frac{p}{[\varepsilon-1]^r}, \frac{[\varepsilon-1]^s}{p}\right]
\end{equation}
(see \cite[Corollaire 2.2]{BergerEquaDiff}). The Frobenius morphism $\p$ on $\Atilde^+$ extends to a map $\p:\Atilde^{[r,s]}\fleche \Atilde^{[pr,ps]}$. For $0\leq r_{1}\leq r_{2}\leq s_{2}\leq s_{1}\leq+\infty$, the natural inclusion
\begin{equation}\nonumber
\Atilde^+\left[\frac{p}{[\varepsilon-1]^{r_{1}}}, \frac{[\varepsilon-1]^{s_{1}}}{p}\right]\plonge\Atilde^+\left[\frac{p}{[\varepsilon-1]^{r_{2}}}, \frac{[\varepsilon-1]^{s_{2}}}{p}\right]
\end{equation}
extends to an injective morphism $\Atilde^{[r_{1},s_{1}]}\plonge\Atilde^{[r_{2},s_{2}]}$. Let $\tilde{B}^{[r,s]}$ be $\Atilde^{[r,s]}[1/p]$ and let $\tilde{B}^{\dag,r}_{\mathrm{rig}}$ be
\begin{equation}\nonumber
\tilde{B}^{\dag,r}_{\mathrm{rig}}\eqdef\intersection{r\leq s\leq\infty}{}\tilde{B}^{[r,s]}.
\end{equation}
The Frobenius map $\p$ on $\Atilde^{[r,s]}$ extends to bijective maps 
\begin{equation}\nonumber
\p:\tilde{B}^{[r,s]}\fleche\tilde{B}^{[pr,ps]}, \p:\tilde{B}^{\dag,r}_{\mathrm{rig}}\fleche \tilde{B}^{\dag,pr}_{\mathrm{rig}}.
\end{equation}
Let $\tilde{B}^\dag_{\mathrm{rig}}$ be
\begin{equation}\nonumber
\tilde{B}^\dag_{\mathrm{rig}}{\eqdef}\union{r\geq0}{}\tilde{B}^{\dag,r}_{\mathrm{rig}}.
\end{equation}
The ring $\tilde{B}^\dag_{\mathrm{rig}}$ comes naturally with an action of the Frobenius morphism $\p$. For all $n\in\N$, the Frobenius morphism induces bijection
\begin{equation}\nonumber
\p^{-n}:\Atilde^{[p^{n-1}(p-1),p^{n-1}(p-1)]}\isom\Atilde^{[(p-1)/p,(p-1)/p]}.
\end{equation}
According to \cite[Proposition 2.11]{BergerEquaDiff}, there is a natural injection
\begin{equation}\nonumber
\tilde{B}^{[\frac{p-1}{p}, \frac{p-1}{p}]}\plonge B_{{\dR}}^+
\end{equation}
and consequently injections
\begin{equation}\label{n}
\iota_n:\tilde{B}^{\dag,p^{n-1}(p-1)}\rightarrow \tilde{B}_{\mathrm{rig}}^{\dag,\frac{p-1}{p}}\hookrightarrow \tilde{B}^{[\frac{p-1}{p},\frac{p-1}{p}]}\hookrightarrow B_{{\dR}}^+
\end{equation}
for each $n\geq 0$. 

Let $K/\qp$ be a finite extension and let $K_{0}$ be its largest unramified subextnsion.  Let $H_{K}\eqdef\Gal(\Kbar/K(\zeta_{p^{\infty}}))$ be the absolute Galois group of the cyclotomic extension of $K$. For $0\leq r\leq\infty$, we write $\tilde{B}^{\dagger,r}_{\rig,K}$ and $\tilde{B}^{\dagger}_{\rig,K}$ for the $H_{K}$-invariants of $\tilde{B}^{\dagger,r}_{\rig}$ and $\tilde{B}^{\dagger}_{\rig}$ respectively. Let $B^{\dagger,r}_{\rig}$ and $B^{\dagger}_{\rig}$ be the closure of $K_{0}[\pi,\pi^{-1}]$ inside $\tilde{B}^{\dag,r}_{\mathrm{rig}}$ and $\tilde{B}^\dag_{\mathrm{rig}}$ respectively. There exists $r(K)\in\R_{+}$ such that for all $r\geq r(K)$, there exists a unique étale extension $B^{\dagger,r}_{\rig,K}$ of ${B}^{\dag,r}_{\mathrm{rig}}$ such that the natural map $B^{\dagger,r}_{\rig,K}\tenseur_{B_{\rig}^{\dagger,r}}\tilde{B}_{\rig}^{\dagger,r}\fleche\tilde{B}^{\dagger,r}_{\rig,K}$ is an isomorphism. The Robba ring $B^{\dagger}_{\rig,K}$ is defined to be
\begin{equation}\nonumber
B^{\dagger}_{\rig,K}\eqdef\liminj{r\geq r(K)}B^{\dagger,r}_{\rig,K}.
\end{equation}
If $K/\qp$ is unramified, then $B^{\dagger,r}_{\rig,K}$ is isomorphic to the ring of power-series converging in the annulus $[p^{-1/r},1[$
\begin{equation}\nonumber
\left\{f(z)=\sum_{n\in\mathbb{Z}}a_{n}z^n|a_n\in K,\forall \rho\in[p^{-1/r},1[,\ \lim_{n\rightarrow\pm\infty}|a_{n}|\rho^{n}=0\right\}
\end{equation}
through the evaluation morphism at $\pi$ (\cite[Proposition 2.31]{BergerEquaDiff}) and so $B^{\dagger}_{\rig,K}$ is isomorphic to the ring $\Rcali$ of power-series which converge in some non-empty annulus $]\rho_{0},1[$
\begin{equation}\nonumber
\left\{f(z)=\sum_{n\in\mathbb{Z}}a_{n}z^n|a_n\in K,\forall \rho\in]\rho_{0},1[,\ \lim_{n\rightarrow\pm\infty}|a_{n}|\rho^{n}=0\right\}.
\end{equation}
In that case, the action of $\p$ on $B^{\dagger,r}_{\rig,K}$ induces the action $\p(z)=(1+z)^{p}-1$ on power-series. Let $\Ecali$ be the field of Laurent series $\somme{n\in\Z}{}a_{n}z^{n}$ with coefficients in $K$, such that $(v_{p}(a_{n}))_{n\in\Z}$ is bounded below and such that $\lim_{n\rightarrow-\infty}v_{p}(a_{n})=+\infty$ together with the action of $\p$ defined above. Then $\Ecali$ is an extension of $\p(\Ecali)$ of degree $p$. For $f\in\Ecali$, let $\psi(f)$ be $p^{-1}\p^{-1}\left(\Tr_{\Ecali/\p(\Ecali)}f\right)$. Then
\begin{equation}\nonumber
\psi=\Ecali\fleche\Ecali
\end{equation}
is a left-inverse of $\p$ which extends by continuity to $\psi:\Rcali\fleche\Rcali$. Any $x\in\Rcali$ may then be written
\begin{equation}\nonumber
x=\somme{i=0}{p-1}(1+z)^{i}\p(x_{i}).
\end{equation}
The map $\psi$ then satisfies
\begin{equation}\nonumber
\psi(x)=x_{0}.
\end{equation}
Moreover, if $K_n$ denotes the extension $K(\zeta_{p^n})$, then the maps $\iota_n$ defined in (\ref{n}) above satisfy
\begin{equation}\nonumber
\iota_n(B^{\dag,p^{n-1}(p-1)}_{\mathrm{rig},K})\rightarrow K_n[[t]].
\end{equation}
\begin{Def}
We say a $B^\dag_{\mathrm{rig},K}$-module $D$ is a $(\varphi,\Gamma_K)$-module of rank $d$ over $B^\dag_{\mathrm{rig},K}$ if
\begin{itemize}
\item $D$ is a finite free $B^\dag_{\mathrm{rig},K}$-module of rank $d$;
\item $D$ is equipped with a $\varphi$-semilinear map $\varphi: D\rightarrow D$ such that
$$\varphi^*(D): B_{\mathrm{rig},K}^\dag\otimes_{\varphi,B^\dag_{\mathrm{rig},K}}D\rightarrow D: a\otimes x\rightarrow a\varphi(x)$$
is an isomorphism;
\item $D$ is equipped with a continuous semilinear action of $\Gamma_K$ which commutes with $\varphi$.
\end{itemize}
\end{Def}
In practice, we only need the notion of $(\varphi,\Gamma_K)$-module for $K$ an extension of $\qp$ included inside $E$ (enlarging $E$ if necessary). In the following, we always make this hypothesis.

Write $r_n=p^{n-1}(p-1)$. According to \cite{CherbonnierColmez}, any $(\varphi,\Gamma_{K})$-module $D$ is overconvergent, \textit{i.e.} there exists $n(D)>0$ and a unique finite free $B^{\dag,r_n(D)}_{\mathrm{rig},K}$-submodule $D^{n(D)}\subseteq D$ of rank $d$ with $B^\dag_{\mathrm{rig}, K}\otimes_{B^{\dag,r_n(D)}_{\mathrm{rig}, K}} D^{n(D)}=D$.

For any $n>n(D)$ we define
$$D_{\mathrm{dif}}^+(D){\eqdef}K_n[[t]]\otimes_{\iota_n, B_{\mathrm{rig},K}^{\dag,r_n}}D^{(n)}$$
and
$$D_{\mathrm{dif}}(D){\eqdef}K_n((t))\otimes_{\iota_n, B_{\mathrm{rig},K}^{\dag,r_n}}D^{(n)}.$$
We also define
$$D_{{\dR}}^K(D)=D_{\mathrm{dif}}(D)^{\Gamma_K=1}, D_{\mathrm{crys}}^K(D)=D[1/t]^{\Gamma_{\Kcal}=1}.$$
The filtration on $D_{{\dR}}^K(D)$ is given by
$$\Fil^iD_{{\dR}}^K(D)=D_{{\dR}}^K(D)\cap t^i D_{\mathrm{dif}}^+(D), i\in\mathbb{Z}.$$
We define a $(\varphi,\Gamma)$-module $D$ to be crystalline (de Rham) if the rank $D$ is equal to the $\mathbb{Z}_p$-rank of $D_{\mathrm{crys}}$ ($D_{{\dR}}$). If $V$ is a representation over some finite extension $L$ of $\qp$, we make all these definitions by regarding it as a $\qp$-representation.

\subsubsection{Bloch-Kato exponential maps}
Recall the fundamental exact sequence in $p$-adic Hodge theory (\cite{BlochKato})
\begin{equation}\label{fundamental}
0\fleche \qp\fleche B_{\crys}^{\varphi=1}\oplus B_{{\dR}}^+\overset{\iota}{\fleche} B_{{\dR}}\fleche 0.
\end{equation}
Tensoring with the Galois representation $V$ and taking $G_{K}$-cohomology yields an exact sequence
\begin{align}\label{EqLongExp}
0\fleche H^{0}(G_{K},V)\fleche D_{\crys}(V)\overset{1-\p,\iota}{\fleche}D_{\crys}(V)\oplus D_{\dR}(V)/\Fil^{0}D_{\dR}(V)\overset{\delta}{\fleche}\Hun(G_{K},V)\\\nonumber
\overset{\delta^{*}}{\fleche} D_{\crys}(V^{*}(1))^{*}\oplus\Fil^{0}D_{\dR}(V)\fleche D_{\crys}(V^{*}(1))^{*}\fleche H^{2}(G_{K},V)\fleche 0
\end{align}
in which $(-)^{*}$ denotes $\Hom(-,E)$. The exponential and dual exponential maps $\exp$ and $\exp^{*}$ of \cite{BlochKato,KatoViaBdR} are respectively the maps
\begin{equation}\nonumber
\exp:\frac{D_{{\dR}}(V)}{\Fil^{0}D_{\dR}(V)}\fleche H^1(G_{K},V)
\end{equation}
and
\begin{equation}\nonumber
\exp^{*}:\Hun(G_{K},V)\fleche\Fil^{0}D_{\dR}(V)
\end{equation}
deduced from the previous exact sequence by considering only the second component of the maps $\delta$ and $\delta^{*}$. The exact sequence \eqref{EqLongExp} and consequently the maps $\exp$ and $\exp^*$ are generalized to de Rham $(\varphi,\Gamma_{K})$-modules over $B^\dag_{\mathrm{rig},K}$ (\cite{NakamuraDeRham}).

As in \cite{BlochKato}, we define the subspace $\Hun_{f}(G_{K},V)$ of $\Hun(G_{K},V)$ to be the image of $\delta$, which is also the kernel of the map
\begin{equation}\nonumber
\Hun(G_{K},V)\fleche\Hun(G_{K},V\tenseur_{\qp}B_{\crys}).
\end{equation}
By definition, $\Hun_{f}(G_{K},V)$ belongs to the kernel of $\exp^{*}$. If $T$ is an $\Ocal[G_{K}]$-submodule inside $V$, we define $\Hun_{f}(G_{K},T)$ to be the intersection of $\Hun(G_{K},T)$ with $\Hun_{f}(G_{K},V)$. Note in particular that $\Hun_{f}(G_{K},T)$ is by construction a saturated $\Ocal$-submodule, and hence is a direct summand inside $\Hun(G_{K},T)$.
\subsubsection{Co-admissible $\Lambda_\infty$-modules}\label{analytic Iwasawa}
We recall necessary some facts and definitions introduced by J.Pottharst in his studies of analytic families of Selmer complexes (\cite{JayPreprint,JayAnalytic}).
\begin{Def}\label{Iwasawa Cohomology}(\cite[Definition 3.1]{NakamuraDeRham})
Let $\Lambda_{n}$ be the $p$-adic completion of the ring $\Oiwa[\mathfrak{m}^n/p]$. The analytic Iwasawa algebra $\Lambda_\infty$ is the inverse limit $\varprojlim_n\Lambda_{n}[1/p]$. If $D$ is a $(\varphi,\Gamma_K)$-module over $B^\dag_{\mathrm{rig},K}$, the analytic Iwasawa cohomology is the $\Lambda_\infty$-module
$$H^q_{\mathrm{Iw}}(G_{K}, D){\eqdef}\varprojlim_n H^q(G_{K}, D\hat{\otimes}_K\Lambda_n^\iota)$$
where $\Gamma_{\Iw}$ acts on $\Lambda_{n}^{\iota}$ through the inverse of the natural action by multiplication on $\Lambda_{n}$.
\end{Def}
The ring $\Lambda_{\infty}$ is the ring of rigid analytic functions on the open unit disc and is a Bezout domain. If $V$ is a $p$-adic $G_{\qp}$-representation, we write $H^q_{\mathrm{Iw}}(G_{K}, V)$ for $H^q_{\mathrm{Iw}}(G_{K}, D^{\dagger}_{\rig}(V))$. Then the natural map  
$$H^q_{\mathrm{cl,Iw}}(G_{K},V)\otimes_\Lambda\Lambda_\infty\fleche H^q_{\mathrm{Iw}}(G_{K}, V)$$
is an isomorphism (\cite[Theorem 1.9]{JayAnalytic}).
\begin{Def}
A $\Lambda_\infty$-module $M$ is said to be co-admissible if there exists an inverse system $(M_n)_n$ of finitely generated $\Lambda_n[1/p]$-modules such that the map $M_{n+1}\rightarrow M_n$ induces isomorphisms $M_{n+1}\otimes_{\Lambda_{n+1}[1/p]}\Lambda_n[1/p]\simeq M_n$ for all $n\in\N$ and such that
\begin{equation}\nonumber
M=\limproj{n}\ M_{n}.
\end{equation}
\end{Def}
Let $W$ be the generic fiber of the formal spectrum $\Spf\Lambda_{\Iw}$ of $\Lambda_{\Iw}$ (the weight space). 
A $\Lambda_\infty$-module $M$ is co-admissible if it is the module of global sections of a coherent analytic sheaf on $W$. Even though the ring $\Lambda_\infty$ is not Noetherian, co-admissible $\La_{\infty}$-module admit characteristic ideals (\cite{JayPreprint}).
\begin{Prop}
\begin{enumerate}
\item The torsion submodule $M_{\mathrm{tors}}$ of a co-admissible $\Lambda_\infty$-module $M$ is also co-admissible.
\item A co-admissible $\Lambda_\infty$-module is torsion if and only if there exists a collection $\{\mathfrak{p}_\alpha\}_{\alpha\in I}$ of closed points of $\union{n}{}\Spec\La_{n}[1/p]$ such that the subset
\begin{equation}\nonumber
I_{n}=\left\{\alpha\in I|\pid_{\alpha}\in\Spec\Lambda_{n}[1/p]\right\}
\end{equation}
is finite for all $n\in\N$ and a function $n:I\fleche\N$ such that there is an isomorphism
\begin{equation}\nonumber
M\simeq\produit{\alpha\in I}{}\Lambda_{\infty}/\pid_{\alpha}^{n_{\alpha}}
\end{equation}\end{enumerate}
\end{Prop}
\begin{proof}
For all $n\in\N$, the ring $\La_{n}[1/p]$ is a principal ideal domain and the natural map $\La_{n+1}[1/p]\plonge\La_{n}[1/p]$ induces an injection $\Spec\La_{n}[1/p]\fleche\Spec\La_{n+1}[1/p]$ which sends the generic point to itself. The structure theorem for finitely generated module over principal ideal domain then yields both claims.
\end{proof}
In fact, the quotient of a co-admissible $\La_{\infty}$-module $M$ by its torsion submodule $M_{\tors}$ is a free $\La_{\infty}$-module, but we will not make use of this fact.
\begin{Def}\label{DefCharacteristicIdealAnalytic}
Let $M$ be a co-admissible, torsion $\La_{\infty}$-module such that there is an isomorphism
\begin{equation}\nonumber
M\simeq\produit{\alpha\in I}{}\Lambda_{\infty}/\pid_{\alpha}^{n_{\alpha}}.
\end{equation}
The \emph{analytic characteristic ideal} $\carac_{\La_{\infty}}M$ of $M$ is the principal ideal
\begin{equation}\nonumber
\carac_{\La_{\infty}}M\eqdef\produit{\alpha\in I}{}\pid^{n_{\alpha}}_{\alpha}.
\end{equation}
\end{Def}
Here, the fact that $\carac_{\La_{\infty}}M$ is a principal ideal follows from the fact that it is a closed ideal and the property that closed ideals of $\La_{\infty}$ are principal (\cite{Lazard}).

\subsection{Unramified Iwasawa Theory}\label{Unramified Iwasawa Theory}
In this subsection, we prove some key facts about the uniform boundedness of Bloch-Kato logarithm map along unramified field extensions. These results are of crucial importance in our later study of $\zp$-extensions of quadratic imaginary extensions of $\Q$ totally ramified at exactly one prime above $p$ and unramified at the other.

In this subsection, we assume that $(T,\rho,\Ocal)$ is a two-dimensional $G_{\qp}$-representation such that $V{\eqdef}T\otimes_{\Ocal}E$ is de Rham with Hodge-Tate weights $(r,s)$ such that $r>0$ and $s\leq 0$.
\subsubsection{Boundedness of the exponential map}
Let $v_1,v_2$ be a basis of $T$. The $E$-vector space $D_{{\dR}}(V)/\Fil^0D_{{\dR}}(V)$ is one-dimensional. Let $\omega_V$ be a fixed generator $D_{{\dR}}(V)/\Fil^0D_{{\dR}}(V)$ over $E$. Then $t^r\omega_V\in B_{{\dR}}^+\otimes T$. By density of $\Atilde^+[\frac{1}{p}]$ in its completion $B_{\dR}^{+}$, there is an element $z\in \Atilde^{+}[\frac{1}{p}]\otimes T$ such that
\begin{equation}\label{Five}
t^r\omega_V-z\in \Fil^{r+1}B_{{\dR}}\otimes T.
\end{equation}
Fix an $n$ such that $p^nz\in \Atilde^{+}\otimes_{\Ocal} T$.

Write $U^1=\{x\in\tilde{E}|v(x-1)\geq 1\}$. We define $\Uni_r^1$ to be the $p$-adic closure of the $\mathbb{Z}_p$-submodule of $\Fil^{-r}B_{\crys}^{\varphi=1}$ generated by elements of the form $\frac{a_1\cdots a_i}{t^i}$ where $a_i$ are elements in $\log(U^1)$, $i\leq r$ (See \cite[6.1.3]{FontaineOuyang} for details about the $\log$ map).
\begin{Lem}\label{LemUr}
Let ${\image}(U_r^1)$ be the image of $U^1_r$ in ${B_{{\dR}}}/{B^+_{{\dR}}}$. As in lemma \ref{2.16}, let $a$ be the largest integer such that $(p-1)a<r$. Let $m_r$ be an integer large enough so that $p^{m_r}\nmid p^a\cdot a!b!p^b/\theta(q/\pi)^j$ for all $b(p-1)\leq r$ and all $j+(p-1)b\geq r$. Suppose $m$ is greater than $r^2+(r-1)m_r+n$. Then ${\image}(U_r^1)\otimes_{\mathbb{Z}_p}T$ contains $p^mw \omega_V$ for all $w\in W(\bar{\mathbb{F}}_p)$.
\end{Lem}
\begin{proof}
We argue by descending induction. Let $w$ be in $W(\bar{\mathbb{F}}_p)$.

Choose $\tilde{v}$ such that $\theta(\tilde{v})$ is equal to $p$. As $\theta(\log(U^1))\supseteq p\mathcal{O}_{\mathbb{C}_p}$, there exists $\tilde{b}_{r-1}\in\log U^1$ such that
$$\theta(\tilde{b}_{r-1})\cdot\theta(\tilde{v})^{r-1}\cdot t^{-r}\equiv p^{r+n}w\omega_V\modulo\Fil^{1-r}B_{{\dR}}\otimes T.$$
Suppose we found $\tilde{b}_{r-1},\cdots, \tilde{b}_i\in\log (U^1)$ and $\tilde{c}_{r-1},\cdots, \tilde{c}_i\in\log (U^1)$ such that
\begin{equation}\nonumber
p^{r(r-i)+m_rr(r-i-1)+n}w\cdot \omega_V
\end{equation}
and 
\begin{equation}\nonumber
\left(\frac{\tilde{b}_{r-1}}{t}\cdot \left(\frac{\tilde{v}}{t}\right)^{r-1}+\cdots+\frac{\tilde{b}_i}{t}\left(\frac{\tilde{v}}{t}\right)^i\right)v_1+\left(\frac{\tilde{c}_{r-1}}{t}\cdot \left(\frac{\tilde{v}}{t}\right)^{r-1}+\cdots+\frac{\tilde{c}_i}{t}\left(\frac{\tilde{v}}{t}\right)^i\right)v_2
\end{equation}
are equal modulo $\Fil^{-i}B_{{\dR}}\otimes T$ (here we recall that $(v_{1},v_{2})$ is our fixed basis of $T$). Then the element 
\begin{equation}\nonumber
t^r\left(p^{r(r-i)+m_r(r(r-i)-1)+n}w\cdot \omega_V-\left(\somme{j=r-1}{i}\frac{\tilde{b}_{j}}{t}\cdot \left(\frac{\tilde{v}}{t}\right)^{j}\right)v_1-\left(\somme{j=r-1}{i}\frac{\tilde{c}_{j}}{t}\cdot \left(\frac{\tilde{v}}{t}\right)^{j}\right)v_2\right)
\end{equation}
belongs to $(\Fil^{r-i}A_{\crys}+\Fil^{r+1}B_{{\dR}}^+)\otimes T$ by (\ref{Five}). By lemma \ref{2.16}, we know that both the image under $\theta$ of the coefficient of $v_1$ in
$$t^i\left(p^{r(r-i)+m_r(r(r-i)-1)+n}w\cdot \omega_V-\left(\frac{\tilde{b}_{r-1}}{t}\cdot \left(\frac{\tilde{v}}{t}\right)^{r-1}+\cdots+\frac{\tilde{b}_i}{t}\left(\frac{\tilde{v}}{t}\right)^i\right)v_1\right)$$
and the image under $\theta$ of the coefficient of $v_2$ in
$$t^i\left(p^{r(r-i)+m_r(r(r-i)-1)+n}w\cdot \omega_V-\left(\frac{\tilde{c}_{r-1}}{t}\cdot \left(\frac{\tilde{v}}{t}\right)^{r-1}+\cdots+\frac{\tilde{c}_i}{t}\left(\frac{\tilde{v}}{t}\right)^i\right)v_2\right)$$
belong to $p^{-m_r}\mathcal{O}_{\mathbb{C}_p}$. For $i\leq j\leq r$, define $\tilde{b}'_{j}$ and $\tilde{c}'_{j}$ to be respectively $\tilde{b}_{j}$ and $\tilde{c}_{j}$ multiplied by $p^{r+m_r}$. Then there is a choice of $\tilde{b}'_{i-1}$ and $\tilde{c}'_{i-1}$ such that
\begin{equation}\nonumber
p^{r(r-i+1)+m_{r}r(r-i)+n}w\cdot \omega_V-\left(\somme{j=r-1}{i-1}\frac{\tilde{b}_{j}}{t}\cdot \left(\frac{\tilde{v}}{t}\right)^{j}\right)v_1-\left(\somme{j=r-1}{i-1}\frac{\tilde{c}_{j}}{t}\cdot \left(\frac{\tilde{v}}{t}\right)^{j}\right)v_2
\end{equation}
belongs to $\Fil^{1-i}B_{{\dR}}\otimes T$. Continuing this process we can find $\tilde{b}_{r-1},\cdots,\tilde{b}_0$ and $\tilde{c}_{r-1},\cdots,\tilde{c}_0$
such that
\begin{equation}\nonumber
p^{r^{2}+m_{r}r(r-1)+n}w\cdot \omega_V-\left(\somme{j=r-1}{0}\frac{\tilde{b}_{j}}{t}\cdot \left(\frac{\tilde{v}}{t}\right)^{j}\right)v_1-\left(\somme{j=r-1}{0}\frac{\tilde{c}_{j}}{t}\cdot \left(\frac{\tilde{v}}{t}\right)^{j}\right)v_2
\end{equation}
is in $B_{\mathrm{dR}}^+\otimes T$.
This proves the claim.
\end{proof}
 The following proposition follows immediately.
\begin{Prop}\label{boundedness}
We consider the co-boundary map
$$\exp: \frac{(B_{{\dR}}\otimes V)^{I_{p}}}{(B^+_{{\dR}}\otimes V)^{I_{p}}+(B_{\crys}\otimes V)^{I_{p},\varphi=1}}\fleche H^1(I_{p}, V).$$
Then there is an $m>0$ such that for any $w\in W(\bar{\mathbb{F}}_p)$ we have $\exp(w\cdot\omega_{V})\in p^{-m}H^1(I_{p}, T)$.
\end{Prop}
\begin{proof}
Let $C$ be in the intersection of $U_r^1$ with the kernel of the map
\begin{equation}\nonumber
B_{\crys}^{\varphi=1}\fleche {B_{{\dR}}}/{B_{{\dR}}^+}.
\end{equation}
Then $t^rC$ belongs to $\Fil^rA_{\crys}$ so $\theta(C)$ is in $p^{-m_r}\mathcal{O}_{\mathbb{C}_p}$ by lemma \ref{2.16}. Combining this with lemma \ref{LemUr} and the diagram
\begin{equation}\nonumber
\xymatrix{
0\ar[r]&\ker\ar[r]\ar[d]&U_r^1\oplus B_{{\dR}}^+\ar[r]\ar[d]&U_r^1+B_{{\dR}}^{+}\ar[d]\ar[r]&0\\
0\ar[r]&\qp\ar[r]&B_{{\crys}}^{\varphi=1}\oplus B^+_{{\dR}}\ar[r]&B_{{\dR}}\ar[r]&0
}
\end{equation}
yields the statement.
\end{proof}
\subsubsection{Yager modules}\label{SubYager}
We mainly follow \cite[Section 3.2]{LZ} to present the theory of Yager modules. Let $K/\mathbb{Q}_p$ be a finite unramified extension. Let $y_{K/\qp}$ be the map
\begin{equation}\nonumber
y_{K/\mathbb{Q}_p}:\Ocal_{K}\fleche\mathcal{O}_K[\mathrm{Gal}(K/\mathbb{Q}_p)]
\end{equation}
sending $x\in\Ocal_{K}$ to
\begin{equation}\label{EqYager}
y_{K/\mathbb{Q}_p}(x)=\sum_{\sigma\in \mathrm{Gal}(K/\mathbb{Q}_p)}x^\sigma[\sigma]\in\mathcal{O}_K[\mathrm{Gal}(K/\mathbb{Q}_p)]
\end{equation}
(note that our convention is slightly different from \cite{LZ}). Recall that $\qppur/\mathbb{Q}_p$ is the unramified $\mathbb{Z}_p$-extension and that $U=\Gal(\qppur/\qp)$. For $m\in\N$, we write $U_{m}$ for the quotient group of $U$ of cardinal $p^{m}$.

On the power-series ring $\zppurhat[[U]]$, we consider the following two actions of $U$: the multiplication action $[\cdot]$ in which $u\in U$ acts as the group-like element $u^{-1}\in\zppurhat[[U]]\croix$ and the Galois action $^{\cdot}$ of $U$ on $\zppurhat$ extended to $\left(\zppurhat\right)^{\N}$. Then the map \eqref{EqYager} induces an isomorphism of $\zp[[U]]$-modules
\begin{equation}\nonumber
y_{\qppur/\mathbb{Q}_p}:\varprojlim_{\mathbb{Q}_p\subseteq K\subseteq \qppur}\mathcal{O}_K\simeq \mathcal{S}_{\qppur/\mathbb{Q}_p}=\left\{f\in\zppurhat[[U]]\mid\forall u\in U, f^u=[u]f\right\}.
\end{equation}
The module $\mathcal{S}_{\mathbb{Q}_p^{\ur}/\mathbb{Q}_p}$ is called the Yager module in \cite{LZ}.

Denote by $\qp\subset K_{m}\subset\qppur$ the sub-extension with Galois group $U_{m}$ over $\qp$. As $K_{m}/\qp$ is unramified, the trace map $\Tr_{K_{m}/\qp}$ is surjective, $\Ocal_{K}$ is a free $\zp[U_{m}]$-module of rank 1 and $x\in\Ocal_{K_{m}}$ is a $\zp[U_{m}]$-generator if and only if $\Tr_{K_{m}/\qp}x\in\zp\croix$. Consider $\limproj{m}\ \Ocal_{K_{m}}$ the inverse system with respect to the trace maps and let $(d_{m})_{m\in\N}\in\limproj{m}\ \Ocal_{K_{m}}$ be an element with $d_{0}\in\zp\croix$. Then
\begin{equation}\nonumber
d\eqdef\limproj{m}\ d_{m}\in \limproj{\mathbb{Q}_p\subseteq K_m\subseteq \qppur}\mathcal{O}_{K_m}
\end{equation}
is a generator of $\limproj{m}\ \Ocal_{K_{m}}$ as $\zp[[U]]$-module and $y_{\qppur/\mathbb{Q}_p}(d)$ is a generator of $\mathcal{S}_{\qppur/\mathbb{Q}_p}$ as $\zp[[U]]$-module, which is thus free of rank 1 as $\zp[[U]]$-module.

Let $\rho$ on $\mathcal{F}$ be a one-dimensional representation of $U$. Mapping $u$ to $\mathrm{Aut}(\mathcal{F}\otimes \zppurhat)$ and extending by linearity defines a map
\begin{equation}\nonumber
\rho: \zppurhat[[U]]\rightarrow \mathrm{Aut}(\mathcal{F}\otimes \zppurhat).
\end{equation}
Identifying $d$ with its image in $\mathcal{S}_{\qppur/\mathbb{Q}_p}\subset \zppurhat[[U]]$, we may thus define $\rho(d)\in\mathrm{Aut}(\mathcal{F}\otimes \zppurhat)$. The defining property of $d$ ensures that image $\rho(d)\cdot x$ belongs to $(\mathcal{F}\otimes \zppurhat)^{G_{\mathbb{Q}_p}}$ for all $x\in\Fcal$. Similarly, we define $\rho(d)^\vee$ to be the image under $\rho$ of the element in $\zppurhat[[U]]$ which is the inverse of
\begin{equation}\nonumber
\limproj{m}\somme{\sigma\in U_{m}}{}d_m^\sigma[\sigma^{-1}]
\end{equation} in this group ring.
\subsubsection{Explicit description of the exponential map and Galois cohomology}
The existence of $d$ and proposition \ref{boundedness} imply the following corollary.
\begin{Cor}\label{specialization}
Let $V(\rho)$ be the twist of $V$ by an unramified $p$-adic character $\rho$ of $G_{\mathbb{Q}_p}$ such that $\rho(\mathrm{Frob}_p)=1+x\in\mathcal{O}_{\mathbb{C}_p}$ with $v_p(x)>0$. Then the map $$\exp: \left(\frac{B_{\mathrm{dR}}}{B_{\mathrm{dR}}^+}\otimes V(\rho)\right)^{G_{\mathbb{Q}_p}}\rightarrow H^1(G_{\mathbb{Q}_p}, V(\rho))$$
can be constructed as
$$\exp\left(\lim_n\sum_{\sigma\in U_n}d_n^\sigma\rho(\sigma)\cdot \omega_V\right)=\lim_n\sum_{\sigma\in U_n}\rho(\sigma)\exp(d_n^\sigma\cdot\omega_V).$$
The right-hand side is well-defined thanks to proposition \ref{boundedness}.
\end{Cor}
\begin{proof}
Consider the natural unramified rank-one Galois representation of $U=\mathrm{Gal}(\mathbb{Q}_{p}^{\mathrm{\ur}}/\mathbb{Q}_p)$ over the Iwasawa algebra $\mathcal{O}[[U]]$. We consider the map from $\mathcal{O}[[U]]$ to $\rho$ mapping $u$ to $\rho(u)$. Tensoring this map with the short exact sequence (\ref{fundamental}) and taking the long exact sequence of Galois cohomology, we get the required formula. 
\end{proof}
\begin{Cor}\label{CorDirectSummand}
Assume that for every specialization $T'$ of $T\otimes_{\Ocal}\Ocal[[U]]$, the representation $V'\eqdef T'[1/p]$ satifisfes
\begin{equation}\nonumber
D^{\varphi=1}_{{\crys}}(V')=0.
\end{equation}
Fix $m\in\N$ as proposition \ref{boundedness}. Identifying it with its image through $y_{\qppur/\mathbb{Q}_p}$, we write $d$ for our fixed generator of the $\zp[[U]]$-module $\mathcal{S}_{\qppur/\mathbb{Q}_p}$. Let 
\begin{equation}\nonumber
\exp: (\varprojlim_n\mathcal{O}_{F_n})\cdot \omega_V\rightarrow p^{-m}H^1(G_{\mathbb{Q}_p}, T\otimes\Ocal[[U]])
\end{equation}
be the inverse limit of the maps
$$\exp_{F_{n}}: \mathcal{O}_{F_n}\cdot\omega_V\rightarrow p^{-m}H^1(G_{F_n}, T)$$
where $\mathbb{Q}_p\subseteq F_n\subset\qppur$ runs over finite extensions. Then $\exp(d\cdot\omega_V)$ generates a $\Ocal[[U]]{\tenseur_{\Ocal} E}$-direct summand of
\begin{equation}\nonumber
H^1(G_{\mathbb{Q}_p}, T\otimes_{\Ocal}\Ocal[[U]])\otimes_{\Ocal}E\simeq \left(\Ocal[[U]]{\tenseur_{\Ocal} E}\right)^{2}.
\end{equation}
\end{Cor}
\begin{proof}
As $\Ocal[[U]]{\tenseur_{\Ocal} E}$ is a Bezout domain, the module generated by $\exp(d\cdot\omega_V)$ is a direct summand if the specialization of $\exp(d\cdot\omega_V)$ through any $\phi\in\Spec(\Ocal[[U]]\otimes_{\Ocal}E)$ with $u$ mapping to $1+x$ with $v_p(x)>0$ is non-zero. Since the $\exp$ map for $V(\rho)$ is injective on $\frac{D_{\mathrm{dR}}(V(\rho))}{D_{\mathrm{dR}}^+(V(\rho))}$, this follows from corollary \ref{specialization}.
\end{proof}
\paragraph{Remark:}
In our later applications, the representation $V$ will be the representation $V_{f}$ twisted by a finite order ramified character. In that case, the first assumption is satisfied.

\newcommand{\bp}{\mathsf{b}_{p}}
\newcommand{\bdp}{\mathsf{d}_{p}}
\newcommand{\dbp}{\mathsf{d}_{p}}

The previous corollary ensures that there exists $b\in E\croix$ such that $b\exp(d\cdot\omega_V)$ is an $\Ocal[[U]]$-submodule of $H^1(G_{\qp}, T\otimes\Ocal[[U]])$ with torsion-free quotient. For any such, the localization of $b\exp(d\cdot\omega_V)$ at any height-one prime $\Pcal\in\Spec\Ocal[[U]]$ is a rank one direct summand of the rank two, free module $H^1(G_{\qp}, T\otimes\Ocal[[U]])_{\Pcal}$. Equivalently, this means that if we write $b\exp(d\cdot\omega_V)$ as linear combination of an $\Ocal[[U]]$-basis, then the coefficients are co-prime to each other.

We define two non-zero elements $\bp,\bdp$ in $E$ which are closely related to Tamagawa numbers.
\begin{Def}\label{DefPseudoTamagawaBD}
Let the assumptions be the same as in the previous corollary. Choose $\bp\in E^\times$ such that $\mathsf{b}_p\exp(d\cdot\omega_V)\in H^1(G_{\qp}, T\otimes_{\Ocal}\Ocal[[U]])$ and such that the torsion submodule of the quotient
$$\frac{H^1(G_{\qp}, T\otimes_{\Ocal}\Ocal[[U]])}{\mathsf{b}_p\Ocal[[U]]\cdot\exp(d\cdot\omega_V)}$$
is zero. We define
\begin{equation}\label{EqDefHunf}
H^1_f(G_{\qp}, T\otimes_{\Ocal}\Ocal[[U]])\eqdef\mathsf{b}_p\Ocal[[U]]\cdot\exp(d\cdot\omega_V).
\end{equation}
Let $\bdp\in E\croix$ be such that the two $\Ocal$-submodules $\bdp H^1_f(G_{\qp},T)$ and $\mathsf{b}_p\exp(\omega_V)$ coincide in $\Hun_{f}(G_{\qp},V)$ (this specifies $\bdp$ up to an element of $\Ocal\croix$).
\end{Def}
When $T=T_{f}$, we may pick the dual of $\omega_{f}$ as $\omega_{V}$.

\subsection{Control Theorem of Selmer Groups}\label{Control}
In this subsection, we study the descent properties of Selmer groups for $T_{f}$ alongside the cyclotomic extension. We show in particular that for an appropriate choice of local condition at $p$, these Selmer modules have non non-trivial pseudo-null submodules and satisfy a perfect control theorem at suitable points of the classical Iwasawa algebra.
\subsubsection{Notations and definitions}
\paragraph{Galois representations}
Put $\Tfiwa\eqdef T_{f}\tenseur_{\Ocal}\Oiwa$, $\Vfiwa\eqdef\Tfiwa\tenseur_{\Oiwa}\Frac(\Oiwa)$. Let $\Afiwa$ be the Pontrjagin dual
\begin{equation}\nonumber
\Afiwa\eqdef\Hom_{\zp}\left(\Tfiwa,\qp/\zp(1)\right)
\end{equation}
of $\Tfiwa$ (here $\Hom_{\zp}$ means continuous morphisms).

If $x$ is an element in a system of parameter of $\Oiwa$ such that $\Ocal_{x}\eqdef\Oiwa/x$ is of characteristic zero, we write $\Pcal_{x}$ for the height-one prime $x\Oiwa\in\Spec\Oiwa$, $E_{x}$ for the fraction field of $\Ocal_{x}$ and $\varpi_{x}$ for a choice of uniformizing parameter of $\Ocal_{x}$. We write $(T_{x},\rho_{x},\Ocal_{x})$ (resp. $(V_{x},\rho_{x},E_{x})$, resp. $(A_{x},\rho_{x},E_{x}/\Ocal_{x})$, resp. $A^{*}_{x}(1)$) for the $G_{\Q,\Sigma}$-representation $\Tfiwa\tenseur_{\Oiwa}\Ocal_{x}$ (resp. $T_{x}\tenseur_{\Ocal_{x}}E_{x}$, resp. $T_{x}\tenseur_{\Ocal_{x}}E_{x}/\Ocal_{x}$, resp. $\Hom_{\zp}\left(T_{x},\qp/\zp(1)\right)$).

If $T$ is equal to $\Tfiwa$ or $T_{x}$, if $A^{*}(1)$ is equal to $\Afiwa$ or $A^{*}_{x}(1)$, if $\Hcal$ is $\Frac(\Ocal_{x})/\Ocal_{x}$ and if $M\subset\Hun(G_{\qp},T)$ is an $\Oiwa$ or $\Ocal_{x}$-submodule, we denote by $M^{\vee}\subset\Hun(G_{\qp},A)$ the orthogonal complement of $M$ under Tate's local duality
\begin{equation}\nonumber
\Hun(G_{\qp},T)\times\Hun(G_{\qp},A^{*}(1))\fleche\Hcal(1).
\end{equation}

We recall that $\RGamma(G_{\qp},\Tfiwa)$ and $\RGamma(G_{\qp},T_{x})$ are perfect complexes of $\Oiwa$ and $\Ocal_{x}$-modules respectively whose cohomology is concentrated in degree 1, that $\Hun(G_{\qp},\Tfiwa)$ and $\Hun(G_{\qp},T_{x})$ are free modules of rank 2 and that the natural map
\begin{equation}\nonumber
\Hun(G_{\qp},\Tfiwa)/x\Hun(G_{\qp},\Tfiwa)\fleche\Hun(G_{\qp},T_{x})
\end{equation}
is an isomorphism (all these assertions follow from lemma \ref{LemCohoLocale}). We also recall that $\Hun(G_{\Q,\Sigma},\Tfiwa)$ is an $\Oiwa$-module free of rank 1 containing the non-zero element $\z(f)_{\Iw}$ (lemma \ref{LemEuler}). According to \cite[Theorem 12.5]{KatoEuler}, there exist infinitely many $x\in\Oiwa$ such that the image of $\z(f)_{\Iw}\in\Hun(G_{\Q,\Sigma},\Tfiwa)$ in $\Hun(G_{\qp},T_{x})/\Hun_{f}(G_{\qp},T_{x})$ is non-zero. In particular, the localization map
\begin{equation}\nonumber
\loc_{p}:\Hun(G_{\Q,\Sigma},\Tfiwa)\fleche\Hun(G_{\qp},\Tfiwa)
\end{equation}
is not identically zero. As its source is a free module of rank 1 and its target is a free module, it is injective. Let $\Hun_{\glob}(G_{\qp},\Tfiwa)$ denote the image of $\loc_{p}$.

According to \cite[Theorem 12.5]{KatoEuler}, the image of $\z(f)_{\Iw}$ under the map
\begin{equation}\label{EqLocSpec}
\Hun_{\glob}(G_{\qp},\Tfiwa)/x\Hun_{\glob}(G_{\qp},\Tfiwa)\plonge\Hun(G_{\qp},T_{x})\fleche\Hun(G_{\qp},T_{x})/\Hun_{f}(G_{\qp},T_{x})
\end{equation}
is non-zero for infinitely many $x\in\Oiwa$. For such an $x$, the composition \eqref{EqLocSpec} is a map between $\Ocal_{x}$-modules free of rank 1 which is not identically zero, hence injective. We fix $x_{0}\in\Oiwa$ such an element and let $\vbf$ be a pre-image inside $\Hun(G_{\qp},\Tfiwa)$ of a generator of $\Hun_{f}(G_{\qp},T_{x_{0}})$. Then $\Hun_{\glob}(G_{\qp},\Vfiwa)$ and $\Frac\left(\Oiwa\right)\cdot\vbf$ are in direct sum inside $\Hun(G_{\qp},\Vfiwa)\simeq\Frac(\Oiwa)^{2}$. In particular, the submodule they generate together inside $\Hun(G_{\qp},\Vfiwa)$ is rank 2.

For $x\in\Oiwa$ such that $\Oiwa/(x)$ is a characteristic zero discrete valuation ring, let $\bar{\vbf}_{x}$ denote the reduction of $\vbf$ modulo $x\Oiwa$ (by definition, $\bar{\vbf}_{x_{0}}$ is thus equal to $\Hun_{f}(G_{\qp},T_{x_{0}})$). As $\Hun_{\glob}(G_{\qp},\Vfiwa)$ and $\Frac\left(\Oiwa\right)\cdot\vbf$ generate a rank 2 submodule, the composition
\begin{equation}\label{EqInjV}
\Hun(G_{\Q,\Sigma},\Tfiwa)/x\Hun(G_{\Q,\Sigma},\Tfiwa)\plonge\Hun(G_{\Q,\Sigma},T_{x})\fleche\Hun(G_{\qp},T_{x})/\Ocal_{x}\cdot\bar{\vbf}_{x}
\end{equation} 
is injective for all $x$ except possibly finitely many. Henceforth, we say that $x\in\Oiwa$ is suitable if $\Oiwa/(x)$ is a characteristic zero discrete valuation ring and if the map \eqref{EqInjV} is injective. Observe that if $x$ is suitable, then the image of the map \eqref{EqInjV} is a free $\Ocal_{x}$-module and so has trivial intersection with the torsion submodule of its target.

\paragraph{Selmer modules and Selmer complexes}Under Tate's local duality, the orthogonal complement of $\Oiwa\cdot\vbf\subset\Hun\left(G_{\qp},\Tfiwa\right)$  is a submodule $\Hun(G_{\Q,\Sigma},\Afiwa^{*}(1))$. Hence, our choice of $\vbf$ allows us to define in our setting an analogue of the Greenberg or $\pm$-Selmer conditions in the ordinary case or supersingular elliptic curves case respectively. In particular, we show that the Pontrjagin duals of $\vbfd$-Selmer modules, that is to say Selmer modules for the Selmer condition obtained by setting the usual Selmer condition at primes $\ell\nmid p$ and the orthogonal complement under Tate's local duality of $\Oiwa\cdot\vbf$ at $p$, have no non-trivial pseudo-null submodules and satisfy a control theorem. 
\begin{Def}
Let $x\in\Oiwa$ be suitable. For $A^{*}(1)$ equal to $\Afiwa$ or $A^{*}_{x}(1)$, $S$ equal to $\Oiwa$ or $\Ocal_{x}$, $\bar{\vbf}$ equal to $\vbf$ or $\bar{\vbf}_{x}$ and $\ell$ a prime, define $\Hun_{\str,\vbf}\left(G_{\ql},A^{*}(1)\right)\subset\Hun\left(G_{\ql},A^{*}(1)\right)$ by
\begin{equation}\nonumber
\Hun_{\str,\vbf}\left(G_{\ql},A\right)\eqdef\begin{cases}
0&\textrm{if $\ell\nmid p$,}\\
\left(S\cdot\bar{\vbf}\right)^{{\vee}}&\textrm{otherwise.}
\end{cases}
\end{equation}
We define the strict $\mathbf{v}^\vee$-Selmer group as follows. Let $\Sel_{\vbfd}(A^{*}(1))$ be the $\Oiwa$-module
\begin{equation}\nonumber
\Sel_{\vbfd}(A^{*}(1))\eqdef\ker\left(\Hun\left(G_{\Q,\Sigma},A^{*}(1)\right)\fleche\sommedirecte{\ell\in\Sigma}{}\Hun\left(G_{\ql},A^{*}(1)\right)/\Hun_{\str,\vbf}\left(G_{\ql},A^{*}(1)\right)\right).
\end{equation}
We write $X_{\mathbf{v}^\vee}(A^{*}(1))$ for the Pontrjagin dual of $\Sel_{\vbfd}(A^{*}(1))^{\iota}$ where $\iota$ indicates that the $\Gamma_{\Iw}$-action on $\Sel_{\vbfd}(A^{*}(1))^{\iota}$ is the inverse of the $\Gamma_{\Iw}$-action on $\Sel_{\vbfd}(A^{*}(1))$. 
\end{Def}

\paragraph{Tamagawa numbers}Let $x\in\Oiwa$ such that $\Pcal_{x}$ is a height-one prime. Then the $\Ocal_{x}$-module $\Hun(G_{\qp},T_{x})/\Hun_{f}(G_{\qp},T_{x})$ is free of rank 1. Choose $\omega_{x}$ a basis of this module. Then $\exp^{*}\omega_{x}$ belongs to $\Fil^{0}D_{\dR}(V_{x})$ and the one-dimensional vector spaces it generates inside $\Fil^{0}D_{\dR}(V_{x})$ coincides with the $E_{x}$-vector space generated by the image $\omega_{f}$ of $f$ in $\Fil^{0}D_{\dR}(V_{x})$ through the comparison isomorphism.
\begin{Def}\label{DefTamagawa}
 Let $\ell\nmid p$ be a prime. The Tamagawa number $c_{x,\ell}$ or $c_{\Pcal_{x},\ell}$ at $\ell$ is
\begin{equation}\nonumber
c_{x,\ell}\eqdef\varpi_{x}^{\length_{\Ocal_{x}}\Hun(I_{\ell},T_{x})^{\Fr(\ell)=1}}.
\end{equation}
Choose an $\Ocal_{x}$-basis $\omega_{x}$ of $\Hun(G_{\qp},T_{x})/\Hun_{f}(G_{\qp},T_{x})$. The Tamagawa number $c_{x,p}$ or $c_{\Pcal_{x},p}$ is such that
\begin{equation}\nonumber
\exp^{*}\omega_{x}=c_{x,p}\omega_{f}.
\end{equation}
\end{Def}
\paragraph{Remark:}If in definition \ref{DefPseudoTamagawaBD}, we choose $\omega_{V}$ to be the element pairing to 1 with $\omega_{f}$ under the duality between $D_{\dR}/\Fil^{0}D_{\dR}$ and $\Fil^{0}D_{\dR}$, then $c_{x,p}$ is equal up to a $p$-adic unit to the quotient $\bdp/\bp$. When $f$ is ordinary at $p$, or $f$ corresponds to a supersingular elliptic curve with $a_{p}(f)=0$, $c_{x,p}$ may be computed explicitly using local theory, see for instance \cite{XinWanIMC}. If $2\leq k<p$, $c_{x,p}$ is actually a local number, as can be seen by using the integral comparison theorem between crystalline and de Rham cohomology (see \cite[Section 14.17]{KatoEuler}.


Let $x^\iota\in \mathrm{Spec}\Lambda$ be the point $x$ composed with the involution $\gamma\mapsto \gamma^{-1}$ on $\Lambda$ for each $\gamma\in\Gamma$. We have from definition, $\langle\exp^* a,\log b\rangle=(a,b)$ ($(,)$ is the local Tate pairing).
Also the pairing of cup product
$$(,): \frac{H^1(\mathbb{Q}_p, T_x)}{H^1_f(\mathbb{Q}_p, T_x)}\times H^1_f(\mathbb{Q}_p, T_{x^\iota})\rightarrow \mathcal{O}_L$$
is surjective.
Thus for some $\mathcal{O}_{x^\iota}$ basis $\bar{\mathcal{V}}$ of $H^1_f(\mathbb{Q}_p, T_{x^\iota})$, we have
\begin{equation}\label{80}
\log\bar{\mathcal{V}}=\frac{1}{c_{x,p}}\cdot\omega_f^\vee.
\end{equation}

We consider the control theorem for $\mathbf{v}^\vee$-Selmer groups. We look at the following diagram
\[\begin{CD}0@>>>\mathrm{Sel}_{\mathbf{v}^\vee}(A^*_{x}(1))@>>>H^1(\mathbb{Q}^\Sigma/\mathbb{Q}, A^*_{x}(1))@>>>\mathcal{P}_{\mathbf{v}^\vee}(\mathbb{Q}, A^*_{x}(1)) \\
@. @VVsV @VVhV @VVgV \\
0@>>>\mathrm{Sel}_{\mathbf{v}^\vee}(A^*_{f,\mathrm{Iw}}(1))^{\mathcal{P}_x} @>>> H^1(\mathbb{Q}^\Sigma/\mathbb{Q}, A^*_{f,\mathrm{Iw}}(1))^{\mathcal{P}_x}@>>>\mathcal{P}_{\mathbf{v}^\vee}(\mathbb{Q},A^*_{f,\mathrm{Iw}}(1))^{\mathcal{P}_x}\end{CD}\]
where $\mathcal{P}_{\mathbf{v}^\vee}(\mathbb{Q}, A^*_{f,\mathrm{Iw}}(1))=\prod_{\ell\nmid p} H^1(G_\ell, A^*_{f,\mathrm{Iw}}(1))\times \frac{H^1(G_p, A^*_{f,\mathrm{Iw}}(1))}{(\Lambda\cdot\mathbf{v})^\vee}$ and $$\mathcal{P}_{\mathbf{v}^\vee}(\mathbb{Q}, A^*_{x}(1))=\prod_{\ell\nmid p} H^1(G_\ell, A^*_{x}(1))\times \frac{H^1(G_p, A^*_{x}(1))}{((\mathcal{O}_x)\cdot\mathrm{Im}(\mathbf{v}))^\vee}.$$

A standard argument using Poitou-Tate exact sequence, as in \cite[Section 3.3.2]{JetchevSkinnerWan} implies that if $\sharp\mathrm{Sel}_{\mathbf{v}^\vee}(A^*_{x}(1))<\infty$ and $H^1_{\mathbf{v}}(\mathbb{Q}^S/\mathbb{Q}, T_{x})=0$, then
\begin{equation}\label{controlTh}
\prod_{\ell\not=p}c_{x,\ell}\mathrm{Fitt}\mathrm{Sel}_{\mathbf{v}^\vee}(\mathbb{Q}, A^*_{x}(1))=\mathrm{Fitt} X_{\mathbf{v}^\vee}/\mathcal{P}_xX_{\mathbf{v}^\vee}.
\end{equation}

\subsubsection{No pseudo-null submodules}
\begin{Lem}\label{Lemma 3.3}
The cardinality $\sharp (H^2(\mathbb{Q}^\Sigma/\mathbb{Q},T_{f,\mathrm{Iw}}))[\mathcal{P}_x]<\infty$ for all but finitely many $m$'s and $x=\gamma-(1+p)^m$.
\end{Lem}
\begin{proof}
The $H^2(\mathbb{Q}^\Sigma/\mathbb{Q},T_{f,\mathrm{Iw}})$ is a finitely generated $\Lambda$-module. Then the lemma follows from the well known structure theorem of finitely generated $\Lambda$-modules.
\end{proof}
\begin{Lem}\label{Lemma 3.4}
The $\frac{H^1(\mathbb{Q}^\Sigma/\mathbb{Q}, A^*_{f,\mathrm{Iw}}(1))}{x H^1(\mathbb{Q}^\Sigma/\mathbb{Q}, A^*_{f,\mathrm{Iw}}(1))}=0$ for all but finitely many $m$'s and $x=\gamma-(1+p)^m$.
\end{Lem}
\begin{proof}
We have $$\frac{H^1(\mathbb{Q}^\Sigma/\mathbb{Q}, A^*_{f,\mathrm{Iw}}(1))}{xH^1(\mathbb{Q}^\Sigma/\mathbb{Q}, A^*_{f,\mathrm{Iw}}(1))}\hookrightarrow H^2(\mathbb{Q}^\Sigma/\mathbb{Q}, A^*_{x}(1)).$$
It is easy to see by Tate local duality that for all places $v$ and all but finitely many $x$, $$H^2(\mathbb{Q}_v, A^*_{x}(1))=0.$$
From the Global duality for these $x$, the $H^2(\mathbb{Q}^\Sigma/\mathbb{Q}, A^*_{x}(1))$ is dual to
$$\mathrm{ker}\{H^1(\mathbb{Q}^\Sigma/\mathbb{Q}, T_x)\rightarrow \prod_{v\in \Sigma} H^1(\mathbb{Q}_v, T_x)\}.$$
We claim this term is $0$ for all but finitely many $m$. Indeed $H^1(\mathbb{Q}^\Sigma/\mathbb{Q}, T_x)$ is $p$-torsion free by (Irred). Moreover we have exact sequence
$$\frac{H^1(\mathbb{Q}^\Sigma/\mathbb{Q},T_{f,\mathrm{Iw}})}{xH^1(\mathbb{Q}^\Sigma/\mathbb{Q}, T_{f,\mathrm{Iw}})}\hookrightarrow H^1(\mathbb{Q}^\Sigma/\mathbb{Q}, T_x)\rightarrow H^2(\mathbb{Q}^\Sigma/\mathbb{Q}, T_{f,\mathrm{Iw}})[x].$$
The last term is finite for all but finitely many $m$ by lemma \ref{Lemma 3.3}. The $H^1(\mathbb{Q}^\Sigma/\mathbb{Q}, T_{f,\mathrm{Iw}})$ is a torsion-free rank one $\Lambda$-module such that the localization map $H^1(\mathbb{Q}^\Sigma/\mathbb{Q}, T_{f,\mathrm{Iw}})\rightarrow H^1(\mathbb{Q}_p,T_{f,\mathrm{Iw}})$ is injective. (Because by \cite{Rohrlich} the image of $\z(f)_{\Iw}$ under this map is non-zero). Now it is easy to see that
$$\mathrm{ker}\{H^1(\mathbb{Q}^\Sigma/\mathbb{Q}, T_x)\rightarrow H^1(\mathbb{Q}_p, T_x)\}$$
is $0$ for all but finitely many $m$ (for example we take a nonzero element $a$ of $H^1(\mathbb{Q}^\Sigma/\mathbb{Q}, T_{f,\mathrm{Iw}})$, and let $x$ to avoid the points where the image of $a$ in $H^1(\mathbb{Q}_p, T_x)$ is zero). The lemma follows readily.
\end{proof}
\begin{Prop}
The $X_{\mathbf{v}^\vee}$ has no pseudo-null submodules.
\end{Prop}
\begin{proof}
Let $x=\gamma-(1+p)^m$ for some integer $m$ be a suitable point. Then we claim for all but finitely many integers $m$ we have surjection
\begin{equation}
H^1(\mathbb{Q}^\Sigma/\mathbb{Q}, A^*_{x}(1))\twoheadrightarrow \mathcal{P}_{\mathbf{v}^\vee}(\mathbb{Q}, A^*_{x}(1)).
\end{equation}
We first look at the exact sequence
\[\begin{CD}H^1(\mathbb{Q}^\Sigma/\mathbb{Q},T_{f,\mathrm{Iw}})@>\times x>>H^1(\mathbb{Q}^\Sigma/\mathbb{Q}, T_{f,\mathrm{Iw}})@>>> H^1(\mathbb{Q}^\Sigma/\mathbb{Q}, T_x)@>>> H^2(\mathbb{Q}^\Sigma/\mathbb{Q}, T_{f,\mathrm{Iw}})[x].\end{CD}\]
From lemma \ref{Lemma 3.3} the last term is torsion for all but finitely many $m$. By our assumption that $x$ is suitable, we get
$$H^1_\mathbf{v}(\mathbb{Q}^\Sigma/\mathbb{Q}, T_x)=0$$ for these $x$.
From Poitou-Tate exact sequence
$$H^1_\mathbf{v}(\mathbb{Q}^\Sigma/\mathbb{Q}, T_x)\rightarrow \mathrm{Hom}_{\mathbb{Z}_p}(\mathcal{P}_{\mathbf{v}^\vee}(\mathbb{Q}, A^*_{x}(1)),\mathbb{Q}_p/\mathbb{Z}_p)\rightarrow \mathrm{Hom}_{\mathbb{Z}_p}(H^1(\mathbb{Q}^\Sigma/\mathbb{Q}, A^*_{x}(1)),\mathbb{Q}_p/\mathbb{Z}_p).$$

It is also clear that the map $H^1(\mathbb{Q}^\Sigma/\mathbb{Q}, A^*_{x}(1))\rightarrow H^1(\mathbb{Q}^\Sigma/\mathbb{Q}, A^*_{f,\mathrm{Iw}}(1))[x]$ is an isomorphism, and that the map $$\mathcal{P}_{\mathbf{v}^\vee}(\mathbb{Q}, A^*_{x}(1))\rightarrow \mathcal{P}_{\mathbf{v}^\vee}(\mathbb{Q}, A^*_{f,\mathrm{Iw}}(1))[x]$$ is surjective. These altogether imply
$$H^1(\mathbb{Q}^\Sigma/\mathbb{Q}, A^*_{f,\mathrm{Iw}}(1))[x]\rightarrow \mathcal{P}_{\mathbf{v}^\vee}(\mathbb{Q}, A^*_{f,\mathrm{Iw}}(1))[x]$$
is surjective.

Then consider the following diagram
\[\begin{CD}0@>>>\mathrm{Sel}_{\mathbf{v}^\vee}(\mathbb{Q}, A^*_{f,\mathrm{Iw}}(1))@>>>H^1(\mathbb{Q}^\Sigma/\mathbb{Q}, A^*_{f,\mathrm{Iw}}(1))@>>>\mathcal{P}_{\mathbf{v}^\vee}(\mathbb{Q}, A^*_{f,\mathrm{Iw}}(1)) \\
@. @VVxV @VVxV @VVxV \\
0@>>>\mathrm{Sel}_{\mathbf{v}^\vee}(\mathbb{Q}, A^*_{f,\mathrm{Iw}}(1))@>>>H^1(\mathbb{Q}^\Sigma/\mathbb{Q}, A^*_{f,\mathrm{Iw}}(1))@>>>\mathcal{P}_{\mathrm{v}^\vee}(\mathbb{Q}, A^*_{f,\mathrm{Iw}}(1))\end{CD}\]
By Snake lemma and lemma \ref{Lemma 3.4} the $\frac{\mathrm{Sel}_{\mathbf{v}^\vee}(\mathbb{Q},A^*_{f,\mathrm{Iw}}(1))}{x\mathrm{Sel}_{\mathbf{v}^\vee}(\mathbb{Q},A^*_{f,\mathrm{Iw}}(1))}=0$ for all but finitely many $m$ and $x=\gamma-(1+p)^m$. By Nakayama's lemma, there is no quotient of $\mathrm{Sel}_{\mathbf{v}^\vee}(\mathbb{Q},A^*_{f,\mathrm{Iw}}(1))$ of finite cardinality. Thus $X_{\mathbf{v}^\vee}$ has no pseudo-null submodules.
\end{proof}
\begin{Def}
Write $\mathscr{F}$ for the characteristic polynomial of $X_{\mathbf{v}^\vee}$. We also write $X_{\mathbf{v}^\vee,x}$ for the dual Selmer group for $A^*_{x}(1)$.
\end{Def}
The control theorem as before implies that
\begin{equation}
\prod_{\ell\nmid p}c_{x,\ell}(f)\carac(X_{\mathbf{v}^\vee,x})=\carac(X_{\mathbf{v}^\vee}/\mathcal{P}_xX_{\mathbf{v}^\vee}).
\end{equation}
Note for any finitely generated torsion $\Oiwa$-module $M$, if $(x)$ is a prime ideal of $\Oiwa$ with $\sharp(\frac{M}{xM})<\infty$, and $\mathscr{F}$ is a generator of $\carac_{\Oiwa}(M)$, then
$$\carac_{\Oiwa}\left(\frac{M}{xM}\right)\subseteq \carac_{\Oiwa}\left(\frac{\Lambda}{(\mathscr{F},x)}\right).$$
If $M$ has no pseudo-null submodule then the above inclusion is an equality. The control theorem argument proved above thus implies that
\begin{equation}\label{ctTh1}
\prod_{\ell\nmid p}c_{x,\ell}(f)\carac_{\Ocal_x}(X_{\mathbf{v}^\vee,x})=\carac_{\Ocal_x}\left(\Ocal_{\Iw}/(\mathscr{F},x)\right).
\end{equation}
Equality (\ref{ctTh1}) is main result of this section and is used later to prove Kato's Iwasawa main conjecture.
\subsection{Beilinson-Flach Elements}\label{SubBF}
In this section, we recall constructions and results from \cite{XinWanIMC} and then use the theory of Beilinson-Flach elements (\cite{LLZ,KKLZ}) to prove conjecture \ref{ConjGreenberg} for $f$.
\subsubsection{Hida families of CM representations}\label{SubHidaCM}
 Let $\La_{\gb}$ be the completed group algebra $\zp[[\Gal(\Kcal_{\infty}/\Kcal_{\cyc})]]\simeq \zp[[\Gamma_{\bar{v}_0}]]$, which we identify through the choice of a topological generator of $\Gal(\Kcal_{\infty}/\Kcal_{\cyc})$ with the power-series ring in one variable $\zp[[Y]]$. We write $\Lcal_{\gb}$ for the fraction field of $\La_{\gb}$.
 
 There exists a unique Hida family $\gb$ parametrized by $\La_{\gb}$ of normalized CM (and hence necessarily ordinary) forms attached to characters of $\Gamma_{\Kcal}$ which passes through the CM form corresponding to the trivial character (we refer to \cite[Theorem 6.2]{HidaTilouine} for the definition and details). To $\gb$ is attached a morphism $\la_{\gb}$ of the ordinary $p$-adic Hida-Hecke algebra $\Hecke_{\infty}^{\ord}$, or equivalently a system of Hecke eigenvalues $\la_{\gb}:\Hecke_{\infty}^{\ord}\fleche\La_{\gb}$, as well as a unique maximal ideal $\mgot_{\gb}\in\Spec\Hecke_{\infty}^{\ord}$ of the Hida-Hecke algebra. To an arithmetic specialization $\phi:\Lambda_{\gb}\fleche E$ is attached a CM character $\psi_{\phi}$ and an automorphic representation $\pi(\psi_{\phi})$ of $\GL_{2}(\A_{\Q})$ which is the base-change of $\psi_{\phi}$ from $\Kcal$ to $\Q$. Up to conjugation, there is a unique character $\Psi_{\gb}:G_{\Kcal}\fleche\La_{\gb}$ which interpolates the $\psi_{\phi}$ as $\phi$ ranges over arithmetic points of $\La_{\gb}$. To $\gb$ is attached a rank two $G_{\Q}$-representation $(V(\gb),\rho_{\gb},\Lcal_{\gb})$, which may be constructed either geometrically or using induction. As we need both constructions in the following, we briefly recall them.
 
 From the geometric point of view, let $D_{\Kcal}$ be the discriminant of $\Kcal$. Let $ES_p(D_{\Kcal})$ and $GES_p(D_{\Kcal})$ be the $\zp[G_{\Q}]$-modules
\begin{equation}\nonumber
ES_p(D_{\Kcal})\eqdef\limproj{r}\ \Hun_{\et}(X_1(D_{\Kcal}p^r)\times_{\Q}\bar{\mathbb{Q}},\mathbb{Z}_p)
\end{equation}
and
\begin{equation}\nonumber
GES_p(D_{\Kcal})\eqdef\limproj{r}\ \Hun_{\et}(Y_1(D_{\Kcal}p^r)\times_{\Q}\bar{\mathbb{Q}},\mathbb{Z}_p).
\end{equation}
Let $e^{*}$ be the ordinary idempotent attached to the covariant Hecke operator $U_{p}$ acting on $\Hun_{\et}(X_1(D_{\Kcal}p^r)\times_{\Q}\bar{\mathbb{Q}},\mathbb{Z}_p)$ and $\Hun_{\et}(Y_1(D_{\Kcal}p^r)\times_{\Q}\bar{\mathbb{Q}},\mathbb{Z}_p)$ (\cite{OhtaTowers}). Then $e^{*}ES_p(D_{\Kcal})$ is a $\Hecke_{\infty}^{\ord}$-module of finite type. Let $T(\gb)$ be the quotient of $e^{*}ES_p(D_{\Kcal})_{\mgot_{\gb}}$ by its $\La_{\gb}$-torsion submodule. Then $V(\gb)$ is $T(\gb)\tenseur_{\La_{\gb}}\Lcal_{\gb}$.  We record for further use the fact that the modules of $I_{p}$-invariants $e^*ES_p(D_{\Kcal})^{I_p}$ and $e^*GES_p(D_{\Kcal})^{I_p}$ are equal (\cite[Theorem]{OhtaTowers}) and denote them both by $\mathfrak{A}_\infty^*$. Let $\mathfrak{B}_\infty^*$ (resp. $\tilde{\mathfrak{B}}_\infty^*$) be the quotient of $e^*ES_p(D_{\Kcal})$ (resp. $e^*GES_p(D_{\Kcal})$) by $\mathfrak{A}_\infty^*$. The proof of \cite[Corollary 2.3.6]{OhtaTowers} shows that there is a natural isomorphism
\begin{equation}\label{EqIsomOmegaOhta}
\Aid_{\infty}^{*}\tenseur_{\zp[[Y]]}\zppurhat[[Y]]\simeq\Hom_{\zppurhat}(S^{\ord}(D_{\Kcal}, \chi_{\Kcal},\zppurhat[[Y]]),\zppurhat[[Y]])
\end{equation}
where $S^{\ord}(D_{\Kcal}, \chi_{\Kcal},\zppurhat[[Y]])$ is the $\zppurhat[[Y]]$-module of ordinary eigenforms with level $D_{\Kcal}p^{\infty}$, central character $\chi_{\Kcal}$ and coefficients in $\zppurhat[[Y]]$.

From the automorphic point of view, there is by construction an isomorphism
\begin{equation}\label{EqIndHidafamily}
\Ind_{G_{\Kcal}}^{G_{\Q}}\left(\Psi_{\gb}\tenseur_{\La_{\gb}}\Lcal_{\gb}\right)\simeq V(\gb)
\end{equation}
 of $\Lcal[G_{\Q}]$-modules. Burungale-Skinner-Tian proved the following strengthening of \eqref{EqIndHidafamily}, which identifies $T(\gb)$ with the induction of $\Psi_{\gb}$ and which is crucial to our purpose.
\begin{Prop}[\cite{BST}]
There is an isomorphism
\begin{equation}\label{EqIndHidaFamily}
\Ind_{G_{\Kcal}}^{G_{\Q}}\Psi_{\gb}\simeq T(\gb)
\end{equation}
 of $\Lambda_{\gb}[G_{\Q}]$-modules. 
\end{Prop}
After restriction to $G_{\qp}$, $V(\gb)$ fits in a short exact sequence
\begin{equation}\nonumber
\suiteexacte{}{}{\Fcal_{\gb}^{+}}{\rho_{\gb}|G_{\qp}}{\Fcal_{\gb}^{-}}
\end{equation} 
of non-zero $\Lcal[G_{\qp}]$-modules which is split as $p$ splits in $\Kcal$. We take the convention that $\Fcal_{\gb}^{-}$ is an unramified $\Lcal[G_{\qp}]$-module. Similarly, since $p$ splits as $v_0\bar{v}_0$ in ${\Kcal}$, there is an identification
\begin{equation}\label{EqIsomIndLocal}
\left(\Ind_{G_{\Kcal}}^{G_{\Q}}\Psi_{\gb}\right)|G_{\qp}\simeq\Psi_{\gb}|G_{\Kcal_{{v}_{0}}}\oplus\Psi_{\gb}|G_{\Kcal_{\bar{v}_{0}}}.
\end{equation}
We choose the convention that the isomorphism \eqref{EqIndHidaFamily} sends $\Psi_{\gb}|G_{\Kcal_{{v}_{0}}}$ to $\Fcal_{\gb}^{-}$. We also fix a $\La_{\gb}$-basis of $\left(\Ind_{G_{\Kcal}}^{G_{\Q}}\Psi_{\gb}\right)|G_{\qp}$ of the form $(v,c\cdot v)$ for some $v$ (here we recall that $c$ denotes complex conjugation).

Let $\chi_{\gb}$ be the central character of $\gb$. Then $\mathscr{F}^+_{\gb}(\chi_{{\gb}}^{-1})$ and $\mathscr{F}^{-}_\gb$ are in addition unramified at $p$, they may be identified with quotients of $\Aid^{*}_{\infty}$ and $\Bid_{\infty}^{*}$ respectively. Let $\omega_{\gb}^\vee\in\mathscr{F}^+_{\gb}(\chi_{{\gb}}^{-1})\tenseur\Frac(\zppurhat[[Y]])$ be the functional which maps $\gb$ to $1$ in the isomorphism \eqref{EqIsomOmegaOhta}. The product of local root numbers of $g_{\phi}$ at places prime to $p$ moves $p$-adic analytically and is a unit when $\phi$ ranges over the arithmetic point of $\La_{\gb}$. Hence, the product over places prime to $p$ of the local root numbers of ${\gb}$ at places prime to $p$ product is well-defined element $\epsi_{\gb}\in\La_{\gb}$. Let $\eta_{\gb}^\vee\in\mathfrak{B}_\infty^*$ be the element which pairs with $\omega_{\gb}^\vee$ to the $\epsi_{\gb}$ under the pairing of \cite[Theorem 2.3.5]{OhtaTowers} (classes dual to $\omega_{\gb}^{\vee}$ and $\eta_{\gb}^{\vee}$ are studied in \cite[Section 10]{KKLZ}).

Suppose $(v^+,v^-)$ is a $\Frac(\Lambda_{\gb})$-basis of the $\Frac(\Lambda_{\gb})$-vector space generated by the $\zppurhat[[Y]]$-lattice $\mathscr{F}^{+}_{\gb}(\chi_{\gb}^{-1})\tenseur_{\La_{\gb}}\zppurhat[[Y]]\oplus\mathscr{F}^{-}_{\gb}\tenseur_{\La_{\gb}}\zppurhat[[Y]]$ with respect to which $\omega_{\gb}^\vee$ and $\eta_{\gb}^\vee$ are equal to $\rho(d)^\vee v^+$ and $\rho(d)v^-$ respectively (here $\rho(d)$ is as in section \ref{SubYager} applied to the unramified representations $\Fcal^{+}_{\gb}(\chi_{\gb}^{-1})$ and $\Fcal^{-}_{\gb}$). According to our choices of isomorphisms \eqref{EqIndHidaFamily} and \eqref{EqIsomIndLocal}, the family $(v^{+},c\cdot v^{+})$ generates a free $\La_{\gb}$-module which we identify with $T(\gb)$. Note the important fact that the $\La_{\gb}$ lattices $\La_{\gb}v^{+}\oplus\La_{\gb}c\cdot v^{+}$ and $\La_{\gb}v^{+}\oplus\La_{\gb}v^{-}$ are \emph{not} the same.

Finally, let $\alpha_{\gb}$ be the $U_{p}$-eigenvalue on $\gb$.
\subsubsection{Analytic families of Beilinson-Flach elements}
Let $f$ be an eigencuspform as in the introduction of this section (in particular, $\rho_{f}|G_{\qp}$ satisfies hypothesis \ref{HypCrys}). If $\phi:\Lambda_{\gb}\fleche E$ is an arithmetic specialization, then the Rankin-Selberg representation $G_{\Q}$-representation $V(f)\tenseur_{E}\rho_{g_{\phi}}$ is isomorphic to $\Ind_{G_{\Kcal}}^{G_{\Q}}(V(f)|G_{\Kcal}\tenseur_{E}\psi_{\phi})$. Let
\begin{equation}\nonumber
\eta_f\in D_{\mathrm{dR}}(V(f))/\mathrm{Fil}^0D_{\mathrm{dR}}(V(f))
\end{equation}
be the element which pairs to $1$ with $\omega_{f^*}$ under the de Rham pairing.

As in the convention of \cite{KLZ} (see also \cite[Proposition 10.1.2]{KKLZ}) we identify the $f$-component of the cohomology of the modular curve $X_0(N)$ with the quotient of the $f$-component of the cohomology of the modular curve $X_0(Np)$ on which $U_p$ acts through the eigenvalue $\alpha$ and write $(\mathrm{Pr}^\alpha)^*$ for the corresponding map. We have
\begin{equation}\nonumber
(\mathrm{Pr}^\alpha)^*\omega_f=\omega_{f_\alpha}.
\end{equation}

\begin{Def}\label{DefBF}
Let $U\simeq \mathbb{Z}_p$ be the Galois group $\Gal(\qppur/\qp)$. For $r>0$ sufficiently large, let $\mathcal{A}$ be the affinoid algebra $\Ocal\diamant{p^{-r}U}$ and let $\Lambda_{\mathcal{A},\infty}$ be $\Lambda_\infty\hat\otimes\mathcal{A}$. The \emph{analytic Beilinson-Flach class}
\begin{equation}\nonumber
\BF_\alpha=\BF_{f_\alpha,{\gb}}\in H^1_{\mathrm{cl,Iw}}(G_{\Q_{\infty}}, T_{f}\otimes T(\gb)\otimes\Lambda_{\mathcal{A},\infty})
\end{equation}
 on $\Acal$ is the class constructed in \cite[Theorem A]{LoefflerZerbesColeman}. If $\phi:\Lambda_{\mathcal{A},\infty}\fleche S$ is a point of $\Lambda_{\mathcal{A},\infty}$, then we write $\BF_{\alpha,\phi}$ for the image of $\BF_{\alpha}$ through $\phi$. In particular, the \emph{one-dimensional Beilinson-Flach class} 
 \begin{equation}\nonumber
\BF_{\alpha,\phi}\in H^1(G_{\Q,\Sigma},T_{f}\otimes\chi_{\phi}\otimes T({\gb})).
\end{equation}
attached to a classical point $\phi\in\Spec\Oiwa$ is the classe constructed in \cite[Section 6.9]{LLZ}\footnote{There, a supplementary hypothesis is put on $\alpha$ in the definition of $\BF_{\alpha,\phi}$. We thank D.Loeffler for informing us that this condition is not necessary if we do not impose local conditions on $\BF_{\alpha,\phi}$ at primes in $\Sigma$.}. Finally, we denote by $\BF_{\alpha}^{\cyc}$ the specialization of $\BF_{\alpha}$ to the cyclotomic deformation.
%
\end{Def}
Before continuing, we explain how we choose the quadratic imaginary field ${\Kcal}$.
\begin{Def}
The quadratic imaginary field ${\Kcal}$ is henceforth chosen so that it satisfies the following properties.
\begin{enumerate}
\item Let $\ell$ be in the assumption of theorem \ref{TheoIntro}. If the $\ell$ is not $2$, then we take $\mathcal{K}$ which is split at $2$ and at any prime divisor of $N$ except $\ell$, and is ramified in $\ell$. If $\ell$ is $2$, then by assumption $2||N$, and we take $\mathcal{K}$ to be ramified at $2$ and split at all other primes divisors of $N$.
\item The prime $p$ is split in $\Kcal$.
\item All other primes dividing $N$ are split in $\Kcal$.
\item The $G_{\Kcal,\Sigma}$-representation $\rhobar_{f}$ is irreducible.
\end{enumerate}
\end{Def}

\subsection{Selmer Complexes and Iwasawa Main Conjecture}\label{Pott}
\subsubsection{\Nekovar-Selmer complexes}In this subsection, we study \Nekovar-Selmer complexes of analytic families Galois representations following \cite{JayAnalytic}. In this section, $(T,\rho,\Ocal)$ and $(V,\rho,E)$ denote $G_{\Q,\Sigma}$-representations. As our actual goal is the study of $T_{f}$ and $V_{f}$, we will come to assume that $T$ and $V$ satisfy the essential properties of $T_{f}$ and $V_{f}$, that is to say that they are of rank 2, that $\rho|G_{\qp}$ is crystalline with two distinct eigenvalues and that $\rhobar|G_{\qp}$ is absolutely irreducible.

Let $A,M$ and $G$ be as in \cite[Section 1.1]{JayAnalytic} (in particular $G$ is a topological group and $M$ is a continuous $A[G]$-module). We write $\Ccont(G, M)$ for the complex of continuous cochains of the $G$-module $M$. Suppose $G$ is equal to $G_{F,\Sigma}$ for $F/\Q$ a finite extension. For $\ell\in\Sigma$ and $v|\ell$ a finite place of a finite extension $F/\ql$, a local condition at $v$ for $M$ is a pair $(U_{v}^{\bullet},i_{v})$ formed of a bounded complex of finite type $A$-modules $U_v^\bullet$ together with a morphism
\begin{equation}\nonumber
i_{v}:U_{v}^{\bullet}\fleche\Ccont(G_{F_{v}},M).
\end{equation}
This definition applies in particular to the $\Lambda_{\mathcal{A},\infty}[G_{\Kcal,\Sigma}]$-module $T\tenseur_{\Ocal}\Lacal$ and more generally to the $A[G_{\Kcal,\Sigma}]$-module $T_{\phi}\eqdef T\tenseur_{\Ocal}\Lacal\tenseur_{\Lacal,\phi}S$ when $\phi:\Lacal\fleche A$ is a map of separated, flat $\Ocal$-algebras complete with respect to a proper ideal containing $p$.
\begin{Def}[\cite{SelmerComplexes,JayAnalytic}]
Let $F/\Q$ be a finite extension. We identify $\Sigma$ with the finite set of finite places of $F$ above places in $\Sigma$. Suppose that for all $v\in\Sigma$, we are given a local condition $(U_{v}^{\bullet,i_{v}})$. The \emph{\Nekovar-Selmer complex} $\RGammaf\left(G_{F,\Sigma},T\tenseur_{\Ocal}\Lacal\right)$ of the analytic family $T\tenseur_{\Ocal}\Lacal$ is the image in the derived category of the mapping cone
\begin{equation}\label{EqDefSelmerJay}
\Cone\left(\Ccont(G_{F,\Sigma},T\tenseur_{\Ocal}\Lacal)\oplus\sommedirecte{v\in\Sigma}{}U_{v}^{\bullet}\overset{\sommedirecte{v\in\Sigma}{}d_{v}}{\fleche}\sommedirecte{v\in\Sigma}{}\Ccont(G_{F_{v}},T\tenseur_{\Ocal}\Lacal)\right)[-1]
\end{equation}
where $d_{v}$ is equal to $\loc_{v}-i_{v}$ for all $v\in\Sigma$. More generally, if there is a map $\phi:\Lacal\fleche A$, we write $\RGammaf\left(G_{F,\Sigma},T_{\phi}\right)$ for the image in the derived category of the mapping cone of \eqref{EqDefSelmerJay} but with $T\tenseur_{\Ocal}\Lacal$ replaced with $T_{\phi}=T\tenseur_{\Ocal}\Lacal\tenseur_{\Lacal,\phi}A$.
\end{Def}
As in \cite{SelmerComplexes,JayAnalytic}, if $v\in\Sigma$ does not divide $p$, we always assume that $(U_{v}^{\bullet},i_{v})$ is the unramified local condition $(\Ccont(G_{F_{v}}/I_{v},(-)^{I_{v}}),i_{v})$ where $i_{v}$ is the inflation map 
\begin{equation}\nonumber
i_{v}:\Ccont(G_{F_{v}}/I_{v},(-)^{I_{v}})\fleche\Ccont(G_{F_{v}},-)
\end{equation}
\paragraph{Local conditions at $p$}
Write $\Rcali$ for the Robba ring $B_{\mathrm{rig},\mathbb{Q}_p}^\dag$ over $\mathbb{Q}_p$, $\Rcali^+$ for $B_{\mathrm{rig},\mathbb{Q}_p}^+$ and $\Rcali_E$ for $\Rcali\otimes_{\mathbb{Q}_p}E$. We recall the notions of triangulation of a $(\p,\Gamma)$-module $D$ and of refinement of a crystalline $G_{\qp}$-representation (\cite[Definitions 2.3.2,2.4.1]{BellaicheChenevier}).
\begin{Def}
Let $D$ be a rank-two $(\varphi,\Gamma)$-modules $D$ over $\Rcali_{E}$. A \emph{triangulation} of $D$ over $\Rcali_{E}$ is a short exact sequence
\begin{equation}\nonumber
0\fleche \mathcal{F}^+D\fleche D\fleche \mathcal{F}^-D\fleche0
\end{equation}
where $\mathcal{F}^\pm D$ are $(\varphi,\Gamma)$-submodules which are free and direct summands as $\Rcali_{E}$-submodules. Let $V$ be a two-dimensional crystalline representation of $G_{\qp}$. A \emph{refinement} of $V$ is a full $\p$-stable $E$-filtration of $D_{{\crys}}(V)$
\begin{equation}\nonumber
\mathcal{F}_0=0\subsetneq\mathcal{F}_1\subsetneq\mathcal{F}_2=D_{\crys}(V).
\end{equation}
\end{Def}
When $V$ is a crystalline $G_{\qp}$-representation with distinct eigenvalues, a refinement of $V$ is equivalent to a choice of ordering of eigenvalues. More generally, according to \cite[Proposition 2.4.1]{BellaicheChenevier}, there is a one-to-one correspondence between triangulations of $D_{\crys}(V)$ and refinements of $V$, given by
$\mathcal{F}^+D=\Rcali[1/t]\mathcal{F}_1\cap D$ and $\mathcal{F}_1=\mathcal{F}^+D[1/t]^\Gamma$.

Let the Robba ring over $\mathcal{A}$ be $\Rcali_\mathcal{A}\eqdef\Rcali\hat\otimes \mathcal{A}$. We also write $\Rcali^+_\mathcal{A}=\Rcali^+\hat\otimes\mathcal{A}$. There is a natural action $U\hookrightarrow \mathcal{A}^\times$. Then we can define a $(\varphi,\Gamma)$-module $D_\mathcal{A}$ over $\Rcali_\mathcal{A}$ by pulling back the action of $\Gamma$ on $D$ but twisting the action of $\varphi$ on $D$ by the Frobenius action of $U$ as above. One can define triangulation for this family of $(\varphi,\Gamma_{\Kcal})$-modules over $\Rcali_\mathcal{A}$ in an obvious way,
and define the analytic Iwasawa cohomology for $D_\mathcal{A}$ in the same way as in section \ref{Iwasawa Cohomology}.
\begin{Lem}
The Iwasawa cohomology groups $H^i_{\Iw}(G_{\qp}, D_\mathcal{A})$ are the cohomology groups of the complex
\begin{equation}\nonumber
\RGamma_{\Iw}\left(G_{\qp},D_{\Acal}\right)\eqdef[D_{\Acal}\overset{\psi-1}{\fleche}D_{\Acal}]
\end{equation}
concentrated in degrees $1$ and $2$.
\end{Lem}
\begin{proof}
See \cite[Theorem 4.4.8]{KPX}.
\end{proof}
Let $\Rcali(\alpha^{-1})$ be the unramified rank one $(\varphi,\Gamma)$-module with Frobenius action given by the scalar $\alpha^{-1}$ where $\alpha$ is a Weil number of weight $k-1$.
\begin{Prop}\label{444}
Let $D$ be $\Rcali(\alpha^{-1})$. Then, for $N$ sufficiently large, there is an exact sequence
\begin{equation}\nonumber
0\fleche\sommedirecte{m=0}{\infty}\left(t^{m}D_{\crys}\left(D_{\Acal}\right)\right)^{\p=1}\fleche\left(\Rcali_{\Acal}^{+}\tenseur D\right)^{\psi=1}\overset{\p-1}{\fleche}\left(\Rcali_{\Acal}^{+}\tenseur D\right)^{\psi=0}\fleche\bigoplus_{m=0}^N\frac{t^m\otimes D_{{\crys}}(D_\mathcal{A})}{(1-\varphi)(t^m\otimes D_{{\crys}}(D_\mathcal{A}))}.
\end{equation}
in which the term $\sommedirecte{m=0}{\infty}\left(t^{m}D_{\crys}\left(D_{\Acal}\right)\right)^{\p=1}$ vanishes.
\end{Prop}
\begin{proof}
This is the family version of \cite[Lemma 3.17,3.18]{NakamuraDeRham} proved in \cite[Section 2]{ChenevierDensite}. The vanishing of $\left(t^{m}D_{\crys}\left(D_{\Acal}\right)\right)^{\p=1}$ follows from the fact that $k-1$ is odd.
\end{proof}
The Galois group $\Gamma_{p,\infty}^{p,\ur}$ is isomorphic to the product $\Gamma_{\Iw}\times\Gamma_{\ur}$ where $\Gamma_{\ur}$ is the Galois group of the unramified $\zp$-extension of $\qp$ and where each factor is isomorphic to $\zp$. Let $Y'+1\in\Gamma_{\mathrm{\ur}}$ be the Frobenius element. For $m\in\N$ sufficiently large, $1-\alpha p^{-m}$ is not a power of $p$ so $Y'+1-\alpha p^{-m}$ is a unit in $\Ocal[[\Gamma_{p,\infty}^{p,\ur}]]$. For $D$ the rank one $(\varphi,\Gamma)$-module $D=\Rcali(\alpha^{-1})$ above, we consider the set of height-one prime ideals $S(D)$ of $\Ocal[[\Gamma_{p,\infty}^{p,\ur}]]$ of the form $(Y'+1-\alpha p^{-m})$. Then the set $S(D)$ is finite, as we just observed. According to proposition \ref{444}, for any height-one prime $\Pcal\in\Spec\Ocal[[\Gamma_{p,\infty}^{p,\ur}]]$ not in $S(D)$, the map
\applicationsimple{\p-1}{\left(\Rcali^{+}_{\Acal}\otimes D\right)^{\psi=1}}{\left(\Rcali^{+}_{\Acal}\otimes D\right)^{\psi=0}}
localized at $\Pcal$ is an isomorphism.

Now suppose $D$ is the $(\varphi,\Gamma)$-module attached to the $G_{\qp}$-representation $V_f$. We fix a triangulation of $D$ by asking that $\mathcal{F}^-\eqdef D/\mathcal{F}^+$ be the $(\varphi,\Gamma)$-module $\Rcali(\alpha^{-1})$ and extend this triangulation to a triangulation
\begin{equation}\label{EqSuiteTriangulation}
\suiteexacte{}{}{\Fcal^{+}(D_{\Acal})}{D_{\Acal}}{\Fcal^{-}(D_{\Acal})}
\end{equation} of $D_\mathcal{A}$ in the obvious way. For simplicity, we occasionally write the corresponding modules as $\mathcal{F}^\pm_\mathcal{A}=\mathcal{F}^\pm(D_\mathcal{A})$. If $\phi:\Acal\fleche S$ is a map of flat $\Ocal$-algebra, we define $\Fcal^{+}(D_{\Acal}\tenseur_{\Acal,\phi}S)$ to be $\Fcal^{+}(D_{\Acal})\tenseur_{\Acal,\phi}S$. To \eqref{EqSuiteTriangulation} is attached a long exact sequence in cohomology which simplifies under assumption \ref{HypCrys} to a short exact sequence
\begin{equation}\label{EqSuiteTriangulationLongue}
0\fleche H^1_{\mathrm{Iw}}(G_{\mathbb{Q}_p},\mathcal{F}^+(D_\mathcal{A}))\fleche H^1_{\mathrm{Iw}}(G_{\mathbb{Q}_p}, D_\mathcal{A})\fleche H^1_{\mathrm{Iw}}(G_{\mathbb{Q}_p},\mathcal{F}^-(D_\mathcal{A}))\fleche 0
\end{equation}
where exactness on the right follows by local duality. The inclusion
\begin{equation}\nonumber
H^1_{\mathrm{Iw}}(G_{\mathbb{Q}_p},\mathcal{F}^+(D_\mathcal{A}))\fleche H^1_{\mathrm{Iw}}(G_{\mathbb{Q}_p}, D_\mathcal{A})
\end{equation}
of the previous short exact sequence may be viewed as a morphism between cohomology complexes of $(\p,\Gamma)$-modules and hence as a morphism between cohomology complexes of $G_{\qp}$-representations with coefficients in the affinoid algebra $\Acal$ through the functorial isomorphism of \cite[Theorem 2.8]{JayAnalytic}.
\begin{Def}\label{DefConditionsEnP}
Let $v|p$ be a finite place.  Let $\phi:\Lacal\fleche A$ be a map of flat $\Ocal$-algebras as above. The $\alpha$ local condition at $v$ is defined by the inclusion 
\begin{equation}\nonumber
i_{p}:H^1_{\mathrm{Iw}}(G_{F_{v}},\mathcal{F}^+(D_\mathcal{A}))\fleche H^1_{\mathrm{Iw}}(G_{F_{v}}, D_\mathcal{A})
\end{equation}
seen as a morphism in the derived category of cohomology complexes of $G_{\qp}$-representations with coefficients in $\Acal$. The relaxed local condition at $v$ is the local condition $\left(\Ccont(G_{F_{v}},D_{\Acal}),\Id_{v}\right)$. The strict local condition at $v$ is the local condition $\left(0,0\right)$. We define $\Hun_{\alpha}(G_{\qp},T_{\phi})$ to be the kernel of the map
\begin{equation}\nonumber
\Hun(G_{\qp},T_{\phi})\fleche\Hun_{\Iw}(G_{\qp},(D_{\Acal}\tenseur A)/\Fcal^{-}(D_{\Acal}\tenseur A)).
\end{equation}
\end{Def}
Note that the exactness of \eqref{EqSuiteTriangulationLongue} shows that the $\alpha$ condition is also defined by the morphism
\begin{equation}\nonumber
i_{p}:\RGamma_{\Iw}(G_{\qp},\Fcal^{+}(D_{\Acal}))\fleche\RGamma_{\Iw}(G_{\qp},D_{\Acal}).
\end{equation}
As $\Fcal^{+}(D_{\Acal}\tenseur_{\Acal,\phi}S)$ is equal to $\Fcal^{+}(D_{\Acal})\tenseur_{\Acal,\phi}S$ when $\phi:\Acal\fleche A$ is a map of flat $\Ocal$-algebra, we see that the local conditions $\alpha$, $\str$ and $\rel$ all commute with $-\Ltenseur_{\Acal,\phi}A$ in the sense that the diagram
\begin{equation}\label{EqFcalCommute}
\xymatrix{
\RGamma_{\Iw}(G_{\qp},\Fcal^{+}(D_{\Acal}))\ar[d]_{-\Ltenseur_{\Acal,\phi}A}\ar[r]^{i_{p}}&\RGamma_{\Iw}(G_{\qp},D_{\Acal})\ar[d]^{-\Ltenseur_{\Acal,\phi}A}\\
\RGamma_{\Iw}(G_{\qp},\Fcal^{+}(D_{\Acal}\tenseur_{\Acal,\phi}A))\ar[r]^(0.54){i_{p}}&\RGamma_{\Iw}(G_{\qp},D_{\Acal}\tenseur_{\Acal,\phi}A)
}
\end{equation}
and its obvious counterpart for the two other conditions are commutative.

For $\phi:\Lacal\fleche A$, let $\RGamma_{?_{1},?_{2}}(G_{\Kcal,\Sigma},T_{\phi})$ be the \Nekovar-Selmer complex attached to the analytic family of $G_{\Kcal,\Sigma}$-representations $T_{\phi}$ with the unramified condition at $v\in\Sigma$ prime to $p$, the condition $?_{1}$ at $v_{0}$ and the condition $?_{2}$ at $\bar{v}_{0}$. In particular, $\RGammalpha(G_{\Kcal,\Sigma},T\otimes\Lambda_{\Acal,\infty})$ is the \Nekovar-Selmer complex attached to the analytic family of $G_{\Kcal,\Sigma}$-representations $T\otimes\Lambda_{\Acal,\infty}$ with the unramified condition at $v\in\Sigma$ prime to $p$ and the condition of definition \ref{DefConditionsEnP} when $v|p$, $\RGammatilde_{\str,\rel}(G_{\Kcal,\Sigma},T\otimes\Lambda_{\Acal,\infty})$ (resp. $\RGammatilde_{\alpha,\rel}(G_{\Kcal,\Sigma},T\otimes\Lambda_{\Acal,\infty})$) is the \Nekovar-Selmer complex attached to the strict condition at $v_{0}$ and the relaxed condition at $\bar{v}_{0}$ (resp. the alpha condition at $v_{0}$ and the relaxed condition at $\bar{v}_{0}$). We write $\Htilde^{i}_{?_{1},?_{2}}(G_{\Kcal,\Sigma},-)$ or sometimes more simply $\Htilde^{i}_{?_{1},?_{2}}(-)$ for the $i$-th cohomology module of $\RGammalpha(G_{\Kcal,\Sigma},-)$.
\begin{Prop}\label{PropSelmerControl}
For $?_{i}\in\{\str,\rel,\alpha\}$ and $\phi:\Lacal\fleche A$ a map of flat $\Ocal$-algebras, the natural map
\begin{equation}\nonumber
\RGammatilde_{?_{1},?_{2}}(G_{\Kcal,\Sigma},T\otimes\Lambda_{\Acal,\infty})\Ltenseur_{\La_{\Acal,\infty},\phi}A\fleche\RGamma_{?_{1},?_{2}}(G_{\Kcal,\Sigma},T_{\phi})
\end{equation}
is an isomorphism. In particular, there is an isomorphism
\begin{equation}\nonumber
\RGammalpha(G_{\Kcal,\Sigma},T\otimes\Lambda_{\Acal,\infty})\Ltenseur_{\La_{\Acal,\infty}}\La_{\infty}\simeq\RGammalpha(G_{\Kcal,\Sigma},T\tenseur\La_{\infty})
\end{equation}
which induces short exact sequences
\begin{equation}\nonumber
\suiteexacte{}{}{\Htilde^{i}_{\alpha,\alpha}(T\otimes\Lambda_{\Acal,\infty})\tenseur_{\La_{\Acal},\infty}\La_{\infty}}{\Htilde^{i}_{\alpha,\alpha}(T\otimes\La_{\infty})}{\Tor_{1}^{\Lambda_{\Acal,\infty}}\left(\Htilde^{i+1}_{\alpha,\alpha}(T\otimes\Lambda_{\Acal,\infty}),\La_{\Acal,\infty}\right)}
\end{equation}
for all $i\in\Z$ and $\Htildeun_{?_{1},?_{2}}(G_{\Kcal,\Sigma},T\otimes\Lambda_{\Acal,\infty})$ is torsion-free.
\end{Prop}
\begin{proof}
As noted above, the diagram \eqref{EqFcalCommute} commutes for each local conditions $\alpha,\str$ and $\rel$ commute with $-\Ltenseur_{\Lacal}\Lambda_{\Q,\infty}$ and all $v\in\Sigma$. The results then follow formally.
\end{proof}
It follows from the definition of the relevant Selmer complexes that there are long exact sequences
\begin{equation}\label{DiagLongSelmerAlpha}
\xymatrix{
0\ar[r]&\Htilde^{1}_{\alpha,\rel}\left(G_{\Kcal,\Sigma},T\tenseur_{\Ocal}\Lacal\right)\ar[r]&\Hun\left(G_{\Kcal,\Sigma},T\tenseur_{\Ocal}\Lacal\right)\ar[r]^(0.52){A}&\Hun(G_{\Kcal_{v_{0}}},\Fcal^{-}(D_{\Acal}))\ar[d]\\
&\sommedirecte{v\in\Sigma}{}H^{2}(G_{\Kcal_{v}},T\tenseur_{\Ocal}\Lacal)&\ar[l]H^{2}\left(G_{\Kcal,\Sigma},T\tenseur_{\Ocal}\Lacal\right)&\ar[l]\Htilde^{2}_{\alpha,\rel}\left(G_{\Kcal,\Sigma},T\tenseur_{\Ocal}\Lacal\right)
}
\end{equation}
and
\begin{equation}\label{DiagLongSelmerStrict}
\xymatrix{
0\ar[r]&\Htilde^{1}_{\str,\rel}\left(G_{\Kcal,\Sigma},T\tenseur_{\Ocal}\Lacal\right)\ar[r]&\Hun\left(G_{\Kcal,\Sigma},T\tenseur_{\Ocal}\Lacal\right)\ar[r]^(0.55){B}&\Hun(G_{\Kcal_{v_{0}}},D_{\Acal})\ar[d]\\
&\sommedirecte{v\in\Sigma}{}H^{2}(G_{\Kcal_{v}},T\tenseur_{\Ocal}\Lacal)&\ar[l]H^{2}\left(G_{\Kcal,\Sigma},T\tenseur_{\Ocal}\Lacal\right)&\ar[l]\Htilde^{2}_{\str,\rel}\left(G_{\Kcal,\Sigma},T\tenseur_{\Ocal}\Lacal\right)
}.
\end{equation}
As in the diagram above, denote the third arrows of \eqref{DiagLongSelmerAlpha} and \eqref{DiagLongSelmerStrict} by $A$ and $B$ respectively. Define $\Ker$ by the exact sequence
\begin{equation}\nonumber
\suiteexacte{}{}{\Ker}{\frac{\Hun(G_{\Kcal_{v_{0}}},D_{\Acal})}{\image(B)}}{\frac{\Hun(G_{\Kcal_{v_{0}}},\Fcal^{-1}(D_{\Acal}))}{\image(A)}}
\end{equation}
Then
\begin{equation}\nonumber
\Ker=\frac{{\image}(B)+H^1(G_{\Kcal_{v_{0}}},\mathcal{F}^+_\mathcal{A})}{{\image}(B)}\simeq \frac{H^1(G_{\Kcal_{v_{0}}}, \mathcal{F}^+_\mathcal{A})}{{\image}(B)\cap H^1(G_{\Kcal_{v_{0}}},\mathcal{F}^+_\mathcal{A})}\simeq \frac{H^1(G_{\Kcal_{v_{0}}}, \mathcal{F}^+_\mathcal{A})}{{\image}(\Htildeun_{\alpha,\rel}(G_{\Kcal,\Sigma}, T\otimes\Lambda_{\mathcal{A},\infty}))}.
\end{equation}
Combining this with (\ref{DiagLongSelmerAlpha}) and (\ref{DiagLongSelmerStrict}) we obtain
\begin{equation}\label{Poitou-Tate One}
0\rightarrow \frac{H^1(G_{\Kcal_{v_{0}}},\mathcal{F}^+_\mathcal{A})}{{\image}(\Htildeun_{\alpha,\rel}(G_{\Kcal,\Sigma}, T\otimes\Lambda_{\mathcal{A},\infty}))}\rightarrow \tilde{H}^2_{\str,\rel}(G_{\Kcal,\Sigma}, T\otimes\Lambda_{\mathcal{A},\infty})\rightarrow \tilde{H}^2_{\alpha,\mathrm{rel}}(G_{\Kcal,\Sigma}, T\otimes\Lambda_{\mathcal{A},\infty})\rightarrow 0.
\end{equation}
Similarly we get
\begin{equation}\label{Poitou-Tate Two}
0\rightarrow \frac{H^1(G_{\Kcal_{\bar{v}_{0}}},\mathcal{F}^-_\mathcal{A})}{{\image}(\Htildeun_{\alpha,\rel}(G_{\Kcal,\Sigma}, T\otimes\Lambda_{\mathcal{A},\infty}))}\rightarrow \tilde{H}^2_{\alpha,\alpha}(G_{\Kcal,\Sigma}, T\otimes\Lambda_{\mathcal{A},\infty})\rightarrow \tilde{H}^2_{\alpha,\mathrm{rel}}(G_{\Kcal,\Sigma}, T\otimes\Lambda_{\mathcal{A},\infty})\rightarrow 0.
\end{equation}
\begin{Lem}
Let $\chi$ be a finite order character of $\Gamma_{\Iw}$. Then
\begin{equation}\label{EqBellaicheKatoLocal}
\Hun_{\alpha}(G_{\qp},T\otimes\chi)=\Hunf(G_{\qp},T\otimes\chi)
\end{equation}
and 
\begin{equation}\label{EqBellaicheKato}
\Htildeun_{\alpha,\alpha}(G_{\Kcal,\Sigma},T\otimes\chi)=\Hunf(G_{\Kcal,\Sigma},T\otimes\chi).
\end{equation}
In particular $\Htildeun_{\alpha,\alpha}(G_{\Kcal,\Sigma},T\otimes\Lacal)$ is a torsion $\Lambda_{\mathcal{A},\infty}$-module and $\Htildeun_{\alpha,\rel}(G_{\Kcal,\Sigma}, T\otimes\Lambda_{\mathcal{A},\infty})$ is a torsion-free $\Lambda_{\mathcal{A},\infty}$-module of rank one.
\end{Lem}
\begin{proof}
By \cite[Proposition 5]{BellaicheRank}, there is an equality $\Hun_{\alpha}(G_{\qp},V\otimes\chi)=\Hunf(G_{\qp},V\otimes\chi)$. As $\Hun_{\alpha}(G_{\qp},T\otimes\chi)$ and $\Hunf(G_{\qp},T\otimes\chi)$ are free $\zp$-modules by lemma \ref{LemCohoLocale}, \eqref{EqBellaicheKatoLocal} holds. Equation \eqref{EqBellaicheKato} then follows as the local condition at $v\nmid p$ are the same for  $\Htildeun_{\alpha,\alpha}(G_{\Kcal,\Sigma},T\otimes\chi)$ and $\Hunf(G_{\Kcal,\Sigma},T\otimes\chi)$. It follows from \cite[Theorem 12.5]{KatoEuler} that the $\Ocal$-modules $\Hunf(G_{\Q,\Sigma},T\otimes\chi)$ and $\Hunf(G_{\Q,\Sigma},T\otimes\chi\otimes\chi_{\Kcal})$ have rank zero, and hence vanish, for all finite order characters $\chi\in\hat{\Gamma}_{\Iw}$ except possibly finitely many (here $\chi_{\Kcal}$ is the quadratic character attached to the extension $\Kcal/\Q$). As
\begin{equation}\nonumber
\Htildeun_{\alpha,\alpha}(G_{\Kcal,\Sigma},T\tenseur\chi)\simeq\Hunf(G_{\Q,\Sigma},T\otimes\chi)\oplus\Hunf(G_{\Q,\Sigma},T\otimes\chi\otimes\chi_{\Kcal})
\end{equation}
and as $\Htildeun_{\alpha,\alpha}(G_{\Kcal,\Sigma},T\otimes\Lacal)\tenseur_{\Lacal,\chi}E$ maps injectively into $\Htildeun_{\alpha,\alpha}(G_{\Kcal,\Sigma},V\tenseur\chi)$ for all $\chi\in\hat{\Gamma}_{\Iw}$ with values in $\Ocal$ by proposition \ref{PropSelmerControl}, we obtain that $\Htildeun_{\alpha,\alpha}(G_{\Kcal,\Sigma},T\otimes\Lacal)$ is a torsion $\Lambda_{\mathcal{A},\infty}$-module. As the cokernel of 
\begin{equation}\nonumber
\Htildeun_{\alpha,\alpha}(G_{\Kcal,\Sigma},T\otimes\Lacal)\fleche\Htildeun_{\alpha,\rel}(G_{\Kcal,\Sigma},T\otimes\Lacal)
\end{equation} 
is included inside $\Hun(G_{\Kcal_{\bar{v}_{0}}},T\otimes\Lacal)/\Hun(G_{\Kcal_{\bar{v}_{0}}},\Fcal^{-}D_{\Acal})$, it is of rank at most 1. It follows that $\Htildeun_{\alpha,\rel}(G_{\Kcal,\Sigma},T\otimes\Lacal)$ is of rank at most 1. According to \cite[Section 7]{KLZ}, the class $\loc_{v_{0}}\BF_{\alpha}$ belongs to $\Hun(G_{\qp},\Fcal^{+}(D_{\Acal}))$ and is not torsion. This entails that the family of Beilinson-Flach elements over $\Acal$ is a non-torsion element of $\Htildeun_{\alpha,\rel}(G_{\Kcal,\Sigma}, T\otimes\Lambda_{\mathcal{A},\infty})$, hence that this module is of rank at least 1. Hence, it is of rank exactly 1.
\end{proof}
\paragraph{The regulator map and duality}For any $(\varphi,\Gamma)$-module of the form $\Rcali_\mathcal{A}(\alpha^{-1})$ for some $\alpha\in\mathcal{A}^*$, we define a regulator map as in \cite[(6.2.1)]{KLZ}
\begin{align}\nonumber
\Reg_{\Rcali_{\Acal}(\alpha^{-1})}:\Hun_{\Iw}(G_{\qp},\Rcali_{\Acal}(\alpha^{-1}))&\isom\Rcali_{\Acal}(\alpha^{-1})^{\psi=1}\isom\Rcali_{\Acal}^{+}(\alpha^{-1})^{\psi=1}\\\nonumber
&\overset{\p-1}{\plonge}\Rcali_{\Acal}^{+}(\alpha^{-1})^{\psi=0}\isom\Acal\hat{\tenseur}\La_{\infty}.
\end{align}
Here, the last map is the Mellin transform. By construction, this applies in particular to $\Fcal^{-}(D_{\Acal})$. As we observed, the map $\Reg_{\Fcal^{-}}$ is an isomorphism after localization at a height-one prime ideal not in $S(\Fcal^{-}_{\Acal})$. The $(\p,\Gamma)$-module $\Fcal^{+}_{\Acal}$ being a twist of $\Rcali(\beta^{-1})$ by a character of $\Gamma_{p,\infty}^{p,\ur}$ which factors through $\Gamma_{\Iw}$, after a reparametrization of the weight space, we may by the same method define also a regulator map
\begin{equation}\nonumber
\Reg_{\Fcal^{+}}:\Hun_{\Iw}(G_{\qp},\Fcal^{+}(D_{\Acal}))\plonge\Acal\hat{\tenseur}\La_{\infty}.
\end{equation}
Then $\Reg_{\Fcal^{+}}$ is an isomorphism after localization outside the finite set $S(\Rcali_\mathcal{A}(\beta^{-1}))$, which we denote $S(\Fcal^{+})$ in a slight abuse of notation.

For $(*,?)\in\{+,-\}^{2}$, we define 
\begin{equation}\nonumber
\Fcal^{*,?}\left(D\hat{\tenseur}D_{\Acal}(\gb)\right)\eqdef\Fcal^{*}(D_)\hat{\tenseur}\Fcal^{?}(D_{\Acal}(V(\gb))).
\end{equation}
Then the $G_{\qp}$-representation $\Fcal^{-,+}\left(D_{\Acal}\hat{\tenseur}D_{\Acal}(\gb)(\chi_{\gb}^{-1})\right)$ is unramified. There is then a map\begin{equation}\nonumber
\Hun\left(G_{\qp},\Fcal^{-,+}\left(D\hat{\tenseur}D_{\Acal}(\gb)(\chi_{\gb}^{-1})\right)\right)\fleche\left(\Fcal^{-,+}\left(D_{\Acal}\hat{\tenseur}D(\gb)(\chi_{\gb}^{-1})\right)\right)^{G_{\qp}}\hat{\tenseur}\left(\Acal\hat{\tenseur}\La_{\infty}\right)
\end{equation}
by \cite[Theorem 8.2.3]{KLZ}. We denote the previous map by $\Reg_{\Fcal^{-,+}}$.
\subsubsection{Rankin-Selberg $p$-adic $L$-functions}
\paragraph{Definitions}In this subsection, we introduce three $p$-adic $L$-functions attached to $f$ and state their fundamental Iwasawa-theoretic properties. 
\begin{Def}
The \emph{cyclotomic $p$-adic $L$-function} $\Lcal_{\alpha}(f)$ of $f_{\alpha}$ is the unique locally analytic function on $\Hom_{\operatorname{cont}}(\Gal(\Q(\zeta_{p^{\infty}})/\Q),\Cp\croix)$  satisfying 
\begin{equation}\nonumber
\Lcal_{\alpha}(f)(\chi_{\cyc}^r\chi^{-1})=\frac{(r+\frac{k-2}{2}-1)!
p^{n(r+\frac{k-2}{2})}\alpha^{-n}G(\chi)^{-1}}{(2\pi i)^{\frac{k-2}{2}+r}\Omega_f^{(-1)^{\frac{k}{2}-r}}}L_{\{p\}}\left(f,\chi,r+\frac{k-2}{2}\right)
\end{equation}
for all character $\chi\in\hat{\Gamma}_{\Iw}$ of finite order $p^{n}$ and all integer $r$ such that $1\leq r+\frac{k-2}{2}\leq k-1$.
\end{Def}
As $f_{\alpha}$ has non-critical slope, existence and unicity of $\Lcal_{\alpha}(f)$ follow from \cite{ManinPadic,Vishik,MazurTateTeitelbaum}.

Recall that $\gb$ is the Hida family of section \ref{SubHidaCM} and write $\pi(\psi_{\phi})$ an arithmetic specialization thereof. In \cite{LoefflerNoteRankin}), it is shown that that the special values of the $L$-function of $f\tenseur\pi(\psi_{\phi})$ admit a $p$-adic interpolation when the weight of $\pi(\psi_{\phi})$ is lower than the weight of $f$. According to \cite[Theorem 7.5.1]{LoefflerZerbesColeman}, these special values define a 3-variables $p$-adic $L$-function. 
\begin{Def}
 The \emph{Rankin-Selberg $p$-adic $L$-function} $\Lcal_{\alpha,\alpha}(f\tenseur\gb)$ attached to $f$ and to the Hida family of CM forms $\gb$ is the $p$-adic $L$-function of \cite[Theorem 7.5.1]{LoefflerZerbesColeman}.
 \end{Def}
 The last $p$-adic $L$-function attached to $f$ is the Rankin-Selberg $p$-adic $L$-function constructed as constant terms of $p$-adic families of Klingen Eisenstein series for $\GU(2,0)$ in \cite{EischenWan}.
 \begin{Def}\label{DefGreenbergRankinSelberg}
The \emph{Greenberg Rankin-Selberg $p$-adic $L$-function} $\Lcal^{\Gr}_{\Kcal}(f)\in\La_{\Kcal}$ is the $p$-adic $L$-function of \cite[Theorem 1.2]{EischenWan}.
\end{Def}
The construction and properties of $\Lcal^{\Gr}_{\Kcal}(f)$ are recalled in appendix \ref{AppendixIwasawaGreenberg} below.
\paragraph{The regulator map and duality}Crucially for our purpose, the Rankin-Selberg $p$-adic $L$-function $\Lcal_{\alpha,\alpha}(f\tenseur\gb)$ is related to the image of Beilinson-Flach elements through the regulator map by the explicit reciprocity law (\cite[Theorem 6.5.9]{KLZ2},\cite[Theorem 7.5.1]{LoefflerZerbesColeman}). In particular, choosing a suitable specialization of $\gb$ shows that $\Lcal_{\alpha,\alpha}(f\tenseur\gb)$ specializes to a $p$-adic $L$-function closely related $\Lcal_{\alpha}(f)$ while exchanging the role of $f$ and $\gb$ in $\Lcal_{\alpha,\alpha}$, one may also obtain a $p$-adic $L$-function which is closely related to $\Lcal^{\Gr}_{\Kcal}(f)$.  There are however two difficulties worth noting. The first is that different construction of the $p$-adic $L$-function corresponds to different choices of $p$-adic periods. The second is that exchanging the role of $f$ and $\gb$ replaces them with their dual while the explicit reciprocity law involves Poincarduality, which does not induce the self-duality of the $\Ocal[G_{\Q,\Sigma}]$-module $T_{f}$. In this paragraph, we deal with these subtleties

Let $\pi: E\rightarrow Y(N)$ be the universal elliptic curve over the open modular curve $Y(N)$ and let $\mathbb{L}_{k-2}$ be the local system $\mathrm{Sym}^{k-2}R^1\pi_*\Ocal$, which is of rank $k-1$ over $\Ocal$. Consider the Poincarduality pairing
\begin{equation}\label{EqPoincareDuality}
\diamant{\cdot,\cdot}:H^1_{{\et}}(Y(N)\times_{\Q}\Qbar, \mathbb{L}_{k-2})_{\mathfrak{m}_f}\times H^1_{c}(Y(N)\times_{\Q}\Qbar, \mathbb{L}_{k-2})_{\mathfrak{m}_f}\fleche E.
\end{equation}
Note that as $\bar{\rho}_f$ is absolutely irreducible, we have
$$H^1_{{\et}}(Y(N)\times_{\Q}\Qbar, \mathbb{L}_{k-2})_{\mathfrak{m}_f}=H^1_{c}(Y(N)\times_{\Q}\Qbar, \mathbb{L}_{k-2})_{\mathfrak{m}_f}=H^1_{{\et},!}(Y(N)\times_{\Q}\Qbar, \mathbb{L}_{k-2})_{\mathfrak{m}_f}$$
(the last cohomology group being interior cohomology). Inverting $p$ and applying Faltings comparison map, Poincarduality induces Serre duality 
\begin{equation}\nonumber
H^1_{\dR}(X(N), \omega_{-k}\otimes\Omega^1_{X(N)})_{\mathfrak{m}_f}\times H^0_{\dR}(X, \omega_k)_{\mathfrak{m}_f}\fleche E
\end{equation}
on the (graded piece of) algebraic de Rham cohomology. As $\Ocal[G_{\Q,\Sigma}]$-module, $T_{f}$ may be identified with $H^1_{{\et}}(Y(N)\times_{\Q}\Qbar, \mathbb{L}_{k-2})_{\mathfrak{m}_f}[\lambda_f](k/2)$ and this latter $\Ocal[G_{\Q,\Sigma}]$-module naturally comes with the pairing induced from \eqref{EqPoincareDuality}. However, this pairing does not induce the self-duality of $T_{f}$. To be more precise, let $\{x,y\}$ be an $\Ocal$-basis of the free $\Ocal$-module $H^1_{{\et}}(Y, \mathbb{L}_{k-2}\otimes_{\mathbb{Z}_p}\Ocal)_{\mathfrak{m}_f}[\lambda_f]$ and let $\diamant{f,f}_{X(N)}$ be the Petersson inner product. Up to a $p$-adic unit, we then have
\begin{equation}\label{EqCompPoincare}
\diamant{x,y}=\frac{\langle f,f\rangle_{X(N)}}{\Omega_f^+\Omega_f^-}
\end{equation}
under \eqref{EqPoincareDuality} by \cite[Lemma 4.17 to Theorem 4.20]{DarmonDiamondTaylor}. Denoting by $\Cid_{f}\in\Cp$ the right-hand side of \eqref{EqCompPoincare}, the self-duality pairing $(\cdot,\cdot)$ on $T_{f}$ then satisfies
\begin{equation}\label{EqCompPoincare2}
(\cdot,\cdot)=\frac{1}{\Cid_{f}}\diamant{\cdot,\cdot}.
\end{equation}
In particular, the vector which pairs to $1$ with $\omega_f$ under \eqref{EqPoincareDuality} is paired to $\frac{1}{\Cid_{f}}$ under the perfect self-duality pairing of $T_f$.

The following proposition, which is due to \cite{KLZ2,LoefflerZerbesColeman}, provides the precise relation between our three $p$-adic $L$-functions and Beilinson-Flach elements.
\begin{Prop}\label{PropReciprocitylaw}
Let $\Lcal^{\mathrm{Katz}}_{\Kcal}$ be Katz $p$-adic $L$-function as in \cite[(8.2)]{HidaTilouine}. Write $\Ecal(f)$ and $\Ecal^*(f)$ for $1-\frac{\beta}{p\alpha}$ and $1-\frac{\beta}{\alpha}$ respectively and denote the Atkin-Lehner pseudo-eigenvalue of $f$ by $\lambda_N$ (see \cite[Corollary 6.4.3]{KLZ}). Finally, recall that $\omega_f$ and $\eta_f$ are identified with $(\mathrm{Pr}^\alpha)^*(\omega_f)$ and $(\mathrm{Pr}^\alpha)^*(\eta_f)$ respectively. Under the convention at the end of section \ref{SubYager}, the following equalities then hold
\begin{equation}\label{Reci1}
\diamant{ \mathrm{Reg}_{\bar{v}_0,\mathcal{F}^-}(\loc_{\bar{v}_{0}}\BF_\alpha),\eta_f}=\Ecal(f)\Ecal^*(f)\lambda_N\Lcal_{\alpha,\alpha}(f\tenseur\gb),
\end{equation}
\begin{equation}\label{EqReciTrivial}
(\eta_f, \omega_f)=\frac{1}{\mathfrak{C}_f}
\end{equation}
and
\begin{equation}\label{Reci2}
\langle \mathrm{Reg}_{v_0,\mathcal{F}^+}(\loc_{v_{0}}\BF_\alpha),\omega_f\rangle=\Lcal^{\Gr}_{\Kcal}(f)\cdot
\frac{1}{h_{\Kcal}
\Lcal^{\mathrm{Katz}}_{\Kcal}}\cdot\frac{v^-}{c\cdot v^+}
\end{equation}
up to $p$-adic units.
\end{Prop}
\begin{proof}
Equation \eqref{Reci1} follows from \cite[Theorem 6.5.9]{KLZ2} and the discussion immediately before the statement of the proposition. After exchanging the role $f$ and $g$ using \cite[Theorem 7.5.1]{LoefflerZerbesColeman}, \cite[Theorem 6.5.9]{KLZ2} describes the interpolation property of the left-hand side of equation \eqref{Reci2} and thus shows that $\langle \mathrm{Reg}_{v_0,\mathcal{F}^+}(\loc_{v_{0}}\BF_\alpha),\omega_f\rangle$ is equal to the Rankin-Selberg $p$-adic $L$-function of \cite[Theorem I]{HidaFourier} multiplied by the scalar ratio comparing $v^{-}$ and $c\cdot v^{+}$. After multiplying by $h_{\Kcal}
\Lcal^{\mathrm{Katz}}_{\Kcal}$, \cite[Equation (7.5)]{XinWanRankinSelberg} then shows that up to a $p$-adic unit, the coefficient $x$ such that  
\begin{equation}\nonumber
\langle \mathrm{Reg}_{v_0,\mathcal{F}^+}(\loc_{v_{0}}\BF_\alpha),\omega_f\rangle\left(h_{\Kcal}
\Lcal^{\mathrm{Katz}}_{\Kcal}
\right)=x\frac{v^-}{c\cdot v^+}
\end{equation}
satisfies the same interpolation properties as $\Lcal^{\Gr}_{\Kcal}(f)$. Equation \eqref{Reci2} follows. Finally, equation \eqref{EqReciTrivial} is a restatement of the definition of $\eta_{f}$ taking into account \eqref{EqCompPoincare2}.
\end{proof}
The following integrality proposition is needed later.
\begin{Lem}\label{LemIntegrality}
Recall that we have identified the coefficient ring $\La_{\gb}$ of $\gb$ with the power-series ring $\zp[[Y]]$ by setting $Y=\gamma_{\bar{v}_{0}-1}$ for $\gamma_{\bar{v}_{0}}$ a topological generator of $\Gamma_{\bar{v}_{0}}$. Then
\begin{equation}\label{EqKatz}
\frac{h_{\Kcal}
\Lcal^{\mathrm{Katz}}_{\Kcal}c\cdot v^+}{v^-}
\end{equation}
is in $\mathbb{Z}_p^{{\ur}}[[Y]]$ up to some powers of $Y$.
\end{Lem}
\begin{proof}
See \cite[Proposition 8.3]{XinWanRankinSelberg}.
\end{proof}
The following lemma relates $\Lcal_\alpha(f)$ to the restriction $\Lcal^{\cyc}_{\alpha,\alpha}(f)$ of $\Lcal_{\alpha,\alpha}(f)$ to the cyclotomic line.
\begin{Lem}\label{LemTwoVariablesOneVariable}
There is an equality
\begin{equation}\nonumber
\Lcal_{\alpha,\alpha}^{\cyc}(f)=\Lcal_\alpha(f)\cdot\Lcal_\alpha(f\tenseur{\chi_{{\Kcal}}})
\end{equation}
up to a non-zero constant. In particular, $\Lcal_{\alpha,\alpha}^{\cyc}(f)$ is not identically zero.
\end{Lem}
\begin{proof}
We view $(f,\alpha)$ as a non-critical point $x_{f}$ on the eigencurve $E(1,M)$ of \cite{ColemanMazur,EmertonInterpolationEigenvalues} (of some tame level $M$). According to \cite[Theorem 4.5.7]{EmertonInterpolationEigenvalues} and \cite[Theorem 7.5.1]{LoefflerZerbesColeman}, there exists an affinoid neighborhood $U$ of $x$ and $p$-adic $L$-functions $\Lcal^{U}_{\alpha}(f),\Lcal^{U}_{\alpha}(f\tenseur\chi_{\Kcal})$ and $\Lcal^{U}(f\tenseur\gb)$ on $\Hom_{\operatorname{cont}}(\Gal(\Q(\zeta_{p^{\infty}})/\Q),\Cp\croix)\times U$ and $\Hom_{\operatorname{cont}}(\Gal(\Q(\zeta_{p^{\infty}})/\Q),\Cp\croix)\times U\times\Spec\zp[[Y]]$ respectively. Denote by $U^{\crys}$ the subset of points such that the underlying $G_{\qp}$-representation is crystalline. Then $U^{\crys}$ is a Zariski-dense subset. All three $p$-adic $L$-functions satisfy an interpolation property at the points of the subset $W\subset\Hom_{\operatorname{cont}}(\Gal(\Q(\zeta_{p^{\infty}})/\Q),\Cp\croix)\times U\times\Spec\zp[[Y]]$ of points of the form $(\chi,x,y)$ with $x\in U^{\crys}$ and $y\in\Spec\zp[[Y]]$ equal to the augmentation ideal. Let us denote by $\Lcal^{\{x\}}(f_{x}\tenseur\gb_{y})$ (resp. $\Lcal_{\alpha}^{\{x\}}(f_{x})$, $\Lcal^{\{x\}}_{\alpha}(f_{x}\tenseur\chi_{\Kcal})$) the specialization of $\Lcal^{U}(f\tenseur\gb)$ at $(x,y)\in U\times\{y\}$ (resp. $x\in U$). In particular, $\Lcal^{\{x_{f}\}}(f_{x_{f}}\tenseur\gb_{y})$ is $\Lcal^{\cyc}_{\alpha,\alpha}(f)$ and $\Lcal_{\alpha}^{\{x\}}(f_{x})\cdot\Lcal^{\{x\}}_{\alpha}(f_{x}\tenseur\chi_{\Kcal})$ is $\Lcal_\alpha(f)\cdot\Lcal_\alpha(f\tenseur{\chi_{{\Kcal}}})$ and we may compare the special valeurs interpolated by $\Lcal^{\{x\}}(f_{x}\tenseur\gb_{y})$ and $\Lcal_{\alpha}^{\{x\}}(f_{x})\cdot\Lcal^{\{x\}}_{\alpha}(f_{x}\tenseur\chi_{\Kcal})$ if $x\in U^{\crys}$.

As $\Lcal^{\{x\}}(f_{x}\tenseur\gb_{y})$ is non-zero for $x\in U^{\crys}$, the quotient 
\begin{equation}\label{EqFunRigid}
x\longmapsto\frac{\Lcal_{\alpha}^{\{x\}}(f_{x})\cdot\Lcal^{\{x\}}_{\alpha}(f_{x}\tenseur\chi_{\Kcal})}{\Lcal^{\{x\}}(f_{x}\tenseur\gb_{y})}
\end{equation}
is defined on $U^{\crys}$ and extends to a unique rigid analytic function defined outside the closed locus of vanishing of $\Lcal^{U}(f\tenseur\gb_{y})$. Comparing the interpolation properties of $\Lcal_{\alpha}^{\{x\}}(f_{x})\cdot\Lcal^{\{x\}}_{\alpha}(f_{x}\tenseur\chi_{\Kcal})$ and $\Lcal^{\{x\}}(f_{x}\tenseur\gb_{y})$, we see that \eqref{EqFunRigid} coincides with
\begin{equation}\nonumber
x\longmapsto\frac{\diamant{f_{x},f_{x}}}{\Omega^{+}_{f_{x}}\Omega^{+}_{f_{x}\tenseur\chi_{\Kcal}}}
\end{equation}
on $U^{\crys}$. Specializing to $x_{f}$, we get the statement of the lemma (with the non-zero constant $\diamant{f,f}/\Omega^{+}_{f}\Omega^{+}_{f\tenseur\chi_{\Kcal}}$).\end{proof}

\subsection{The Iwasawa Main Conjecture in the crystalline case}\label{SubCrysGreenberg}
In this subsection, we state the two-variable Greenberg-Iwasawa Rankin-Selberg Main Conjecture for the eigencuspform $f$ and show that it implies the Iwasawa Main Conjecture $f$. In particular, we prove conjecture \ref{ConjIMCweak} up to powers of $p$ when our ongoing assumptions on $f$ and the assumptions of \ref{TheoGreenberg} of appendix \ref{AppendixIwasawaGreenberg} are satisfied. More precisely, we prove the following theorem.
\begin{Theo}\label{TheoSansP}
Let $f\in S_{k}(\Gamma_{0}(N))$ be a normalized eigencuspform of weight $k\geq2$ satisfying the following hypotheses.
\begin{enumerate}
\item The $G_{\qp}$-representation $\rho_{f}|G_{\qp}$ is crystalline and short with irreducible residual representation.
\item There exists $\ell||N$, $\ell\not=p$ such that $\rhobar_{f}|G_{\ql}$ is a ramified extension
\begin{equation}\nonumber
\suiteexacte{}{}{\mu\chi_{\cyc}^{1-k/2}}{\rhobar_{f}|G_{\ql}}{\mu\chi_{\cyc}^{-{k/2}}}
\end{equation}
where $\mu:G_{\ql}\fleche\{\pm1\}$ is the non-trivial unramified character.
\end{enumerate}
Then the Iwasawa Main Conjecture (conjecture \ref{ConjIMCweak}) holds for $f$.
\end{Theo}

\subsubsection{Statement of the {Greenberg-Iwasawa Rankin-Selberg Main Conjecture}}
Let $\Kcal\subset_{f} F\subset\Kcal_{\infty}$ be a finite subextension and let $v$ be a finite place of $\Ocal_{F}$. The Greenberg local condition $\Hun(G_{F_{v}},A)\subset\Hun(G_{F_{v}},A)$ at $v$ is defined to be
\begin{equation}\nonumber
\image\left(\Hun(G_{F_{v}}/I_{v},V^{I_{v}})\fleche\Hun(G_{F_{v}}/I_{v},A^{I_{v}})\right)\subset\Hun(G_{F},A)
\end{equation}
if $v\nmid p$, to be $\Hun_{\Gr}(G_{F_{v}},A)$ if $v|v_{0}$ and to be 0 if $v|\bar{v}_{0}$.
\begin{Def}
The \emph{Greenberg Rankin-Selberg Selmer group} of $f$ is the $\La_{\Kcal}$-module
\begin{equation}\nonumber
\Sel^{\Gr}_{\Kcal}(f)\eqdef\liminj{\Kcal\subset_{f}F\subset\Kcal_{\infty}}\ker\left(\Hun(G_{F,\Sigma},A)\fleche\sommedirecte{v\in\Sigma}{}\Hun(G_{F_{v}},A)/\Hun_{\Gr}(G_{F_{v}},A)\right).
\end{equation}
We write
\begin{equation}\nonumber
X^{\Gr}_{\Kcal}(f)\eqdef\Hom\left(\Sel^{\Gr}_{\Kcal}(f),\qp/\zp\right)
\end{equation}
for the Pontrjagin dual of $\Sel^{\Gr}_{\Kcal}(f)$. 
\end{Def}
We recall the statement of the two-variable Greenberg-Iwasawa Rankin-Selberg Main Conjecture.
\begin{Conj}(Greenberg Main Conjecture)\label{ConjGreenberg}
The $\La_{\Kcal}$-module $X^{\Gr}_{\Kcal}(f)$ is torsion and there is an equality of ideals of $\La_{\Kcal}$
\begin{equation}\nonumber
\carac_{\La_{\Kcal}}\left(X^{\Gr}_{\Kcal}(f)\right)=\left(\Lcal^{\Gr}_{\Kcal}(f)\right).
\end{equation}
\end{Conj}
In appendix \ref{AppendixIwasawaGreenberg}, we prove the following results towards conjecture \ref{ConjGreenberg}.
\begin{Theo}\label{TheoGreenberg}
Let $f\in S_{k}(\Gamma_{0}(N))$ be an eigencuspform of even weight $k$ satisfying the following assumptions.
\begin{enumerate}
\item $\rhobar_{f}|G_{\Kcal}$ is absolutely irreducible.
\item $\rho_{f}|G_{\qp}$ is a crystalline representation. Moreover either $\bar{\rho}_f|_{G_{\qp}}$ is absolutely irreducible, or $f$ is ordinary at $p$.
\item There exists $q||N$ (in particular $q\nmid p$) which is not split in $\Kcal$.
\item If $\ell|N$ is not split in $\mathcal{K}$, then $\ell||N$. Moreover if $2$ is non-split in $\mathcal{K}$, then $2||N$.
\suspend{enumerate}
Then the following inclusion of ideals of $\Ocal^{\ur}[[\Gamma_{\Kcal}]]$ holds
$$\carac_{\mathcal{O}^{{\ur}}[[\Gamma_\mathcal{K}]]}(X_{\Kcal}^{\Gr}(f)
\otimes_{\mathcal{O}}\mathcal{O}^{{\ur}})\subseteq(\Lcal_{\Kcal}^{\Gr}(f))$$
up to height-one primes which are pullbacks of primes in $\mathcal{O}[[\Gamma^+]]$. Assume in addition that the following assumption holds.
\resume{enumerate}
\item If $\ell|N$ is not split in $\mathcal{K}$, then $\ell$ is ramified in $\Kcal$ and $\pi(f)_\ell$ is a special Steinberg representation twisted by $\chi_{{\ur}}$ for $\chi_{{\ur}}$ the unramified character sending $\ell$ to $(-1)\ell^{\frac{k}{2}-1}$.
\end{enumerate}
Then
$$\carac_{\mathcal{O}^{{\ur}}[[\Gamma_\mathcal{K}]]}(X_{\Kcal}^{\Gr}(f)
\otimes_{\mathcal{O}}\mathcal{O}^{{\ur}})\subseteq(\Lcal_{\Kcal}^{\Gr}(f))$$
holds.
\end{Theo}

\subsubsection{Proof of the Iwasawa Main Conjecture up to powers of $p$}
In this subsection, we prove the following proposition.
\begin{Prop}\label{Prop}
The following equalities hold.
\begin{equation}\label{EqTheoSansPanalytic}
\carac_{\La_{\infty}}\Htilde^{2}_{\alpha,\alpha}(G_{\Q,\Sigma},T\tenseur\La_{\infty})=\left(\Lcal_{\alpha}(f)\right)\La_{\infty}
\end{equation}
\begin{equation}\label{EqTheoSansP}
\carac_{\La_{\Iw}[1/p]}H^{2}_{\et}(\Z[1/p],T\tenseur\La_{\Iw})=\carac_{\La_{\Iw}[1/p]}H^{1}_{\et}(\Z[1/p],T\tenseur\La_{\Iw})/\La_{\Iw}\cdot\z(f)_{\Iw}.
\end{equation}
\end{Prop}
Extending slightly the notations of definition \ref{Iwasawa Cohomology}, we let $\La_{\Kcal,n}$ be the $p$-adic completion of $\La_{\Kcal}[\mathfrak{m}^n/p]$ and $\La_{\Kcal,\infty}$ be the inverse limit on $n$ of the $\La_{\Kcal,n}$. We let $\La_{\Acal}$ (resp. $\La_{\Acal,n}$) be the analytic ring $\La_{\infty}\hat{\tenseur}\Acal$ (resp. $\La_{n}\hat{\tenseur}\Acal$) and similarly for $\La_{\Kcal,\Acal}$ and $\La_{\Kcal,\Acal,n}$.

By \cite[5.1.6]{SelmerComplexes} and since $\La_{\Kcal,\Acal,\infty}$ is a flat $\La_{\Kcal}$-module, there is an isomorphism of co-admissible $\La_{\Kcal,\Acal,\infty}$-modules
\begin{equation}\nonumber
\Htilde^{2}_{\str,\rel}(G_{\Kcal,\Sigma},T\otimes\La_{\Kcal,\Acal,\infty})\simeq X_{\Kcal}^{\Gr}(f)\hat{\tenseur}_{\La_{\Kcal}}\La_{\Kcal,\Acal,\infty}
\end{equation}
and hence an equality of characteristic ideals in the sense of \ref{DefCharacteristicIdealAnalytic} if they are torsion co-admissible $\La_{\Kcal,\Acal,\infty}$-modules.
\begin{Prop}\label{PropAppendix}
For $*\in\{v_{0},\bar{v}_{0}\}$, let $S_{*}(\Fcal^{\pm}(D_{\Acal}))$ be the set of height-one primes $S(\Fcal^{\pm}(D_{\Acal}))$, in which we identify $\Kcal_{*}$ with $\qp$. Let $\Pcal\in\Spec\Lambda_{\Kcal,\Acal,\infty}$ be a height-one prime which is neither in $S_{v_0}(\Fcal^+(D_{\Acal}))$ nor in $S_{\bar{v}_0}(\Fcal^-(D_{\Acal}))$. Let $n\in\N$ be an integer and let $\Pcal_{n}$ be the intersection of $\Pcal$ with $\La_{\Kcal,\Acal,n}$. Assume the hypotheses of theorem \ref{TheoGreenberg}. Then 
\begin{equation}\label{EqDivFundamental}
\ord_{\Pcal_{n}}\carac_{\La_{\Kcal,\Acal,n}}\left(\Htilde^{2}_{\alpha,\alpha}(G_{\Kcal,\Sigma},T\otimes\La_{\Kcal,\Acal,n})\right)\geq\ord_{\Pcal_{n}}\left(\Lcal_{\alpha,\alpha}(f)\right).
\end{equation}
\end{Prop}
\begin{proof}
According to \eqref{Poitou-Tate Two} and the multiplicativity of characteristic ideals in short exact sequences, there is an equality
\begin{align}\label{EqDivPleinDeTermesUn}
\carac_{\La_{\Kcal,\Acal,n}}\left(\Htilde^{2}_{\alpha,\alpha}(G_{\Kcal,\Sigma},T\otimes\La_{\Kcal,\Acal,n})\right)&=\carac_{\La_{\Kcal,\Acal,n}}\left(\Htilde^{2}_{\alpha,\rel}(G_{\Kcal,\Sigma},T\otimes\La_{\Kcal,\Acal,n})\right)\\\nonumber
&\times\carac_{\La_{\Kcal,\Acal,n}}\left(\frac{H^1(G_{\Kcal_{\bar{v}_{0}}},\mathcal{F}^-_\mathcal{A})}{{\image}(\Htildeun_{\alpha,\rel}(G_{\Kcal,\Sigma},T\otimes\Lambda_{\Kcal,\Acal,n}))}\right).
\end{align}
Similarly, there is an equality
\begin{align}\label{EqDivPleinDeTermes}
\carac_{\La_{\Kcal,\Acal,n}}\left(\Htilde^{2}_{\alpha,\rel}(G_{\Kcal,\Sigma},T\otimes\La_{\Kcal,\Acal,n})\right)&=\carac_{\La_{\Kcal,\Acal,n}}\left(\Htilde^{2}_{\str,\rel}(G_{\Kcal,\Sigma},T\otimes\La_{\Kcal,\Acal,n})\right)\\\nonumber
&\times\left(\carac_{\La_{\Kcal,\Acal,n}}\left(\frac{H^1(G_{\Kcal_{{v}_{0}}},\Fcal^{+}_{\Acal})}{{\image}(\Htildeun_{\alpha,\rel}(G_{\Kcal,\Sigma}, T\otimes\Lambda_{\Kcal,\Acal,n}))}\right)\right)^{-1}
\end{align}
according to \eqref{Poitou-Tate One}. Under the hypothesis of the proposition, theorem \ref{TheoGreenberg} of appendix \ref{AppendixIwasawaGreenberg} applies so
\begin{equation}\nonumber
\left(\Lcal^{\Gr}_{\Kcal}(f)\right)|\carac_{\La_{\Kcal,\Acal,n}}\left(\Htilde^{2}_{\str,\rel}(G_{\Kcal,\Sigma},T\otimes\La_{\Kcal,\Acal,n})\right).
\end{equation}
Let $\Pcal_{n}\in\Spec\La_{\Kcal,\Acal,n}$ be as in the statement of the proposition. Assume in addition that $Y\notin\Pcal_{n}$. Then
\begin{equation}\nonumber
\ord_{\Pcal_{n}}\left(\Lcal^{\Gr}_{\Kcal}(f)\right)=\ord_{\Pcal_{n}}\carac_{\La_{\Kcal,\Acal,n}}\left(\frac{H^1(G_{\Kcal_{{v}_{0}}},\Fcal^{+}_{\Acal})}{{\image}(\BF_{\alpha})}\right)
\end{equation}
by equation \eqref{Reci2} in proposition \ref{PropReciprocitylaw} and lemma \ref{LemIntegrality}. Then \eqref{EqDivPleinDeTermes} implies that the $\Pcal_{n}$-adic valuation of
\begin{equation}\label{EqDeuxGrosTrucs}
\left(\carac_{\La_{\Kcal,\Acal,n}}\left(\frac{H^1(G_{\Kcal_{{v}_{0}}},\Fcal^{+}_{\Acal})}{{\image}(\BF_{\alpha})}\right)\right)\left(\carac_{\La_{\Kcal,\Acal,n}}\left(\frac{H^1(G_{\Kcal_{{v}_{0}}},\Fcal^{+}_{\Acal})}{{\image}(\Htildeun_{\alpha,\rel}(G_{\Kcal,\Sigma}, T\otimes\Lambda_{\Kcal,\Acal,n}))}\right)\right)^{-1}
\end{equation}
is less than the $\Pcal_{n}$-adic valuation of 
\begin{equation}\nonumber
\carac_{\La_{\Kcal,\Acal,n}}\left(\Htilde^{2}_{\alpha,\rel}(G_{\Kcal,\Sigma},T\otimes\La_{\Kcal,\Acal,n})\right).
\end{equation}
Then, \eqref{EqDivPleinDeTermesUn} further implies that the $\Pcal_{n}$-adic valuation of the product of \eqref{EqDeuxGrosTrucs} with
\begin{align}\nonumber
\carac_{\La_{\Kcal,\Acal,n}}\left(\frac{H^1(G_{\Kcal_{\bar{v}_{0}}},\mathcal{F}^-_\mathcal{A})}{{\image}(\Htildeun_{\alpha,\rel}(G_{\Kcal,\Sigma}, T\otimes\Lambda_{\Kcal,\Acal,n}))}\right)
\end{align}
is less than the $\Pcal_{n}$-adic valuation of
\begin{equation}\nonumber
\carac_{\La_{\Kcal,\Acal,n}}\left(\Htilde^{2}_{\alpha,\alpha}(G_{\Kcal,\Sigma},T\otimes\La_{\Kcal,\Acal,n})\right).
\end{equation}
Using again the multiplicativity of characteristic ideals in short exact sequences, we find that
\begin{equation}\nonumber
\ord_{\Pcal_{n}}\carac_{\La_{\Kcal,\Acal,n}}\left(\frac{H^1(G_{\Kcal_{\bar{v}_{0}}},\mathcal{F}^-_\mathcal{A})}{{\image}(\BF_{\alpha})}\right)\leq\ord_{\Pcal_{n}}\carac_{\La_{\Kcal,\Acal,n}}\left(\Htilde^{2}_{\alpha,\alpha}(G_{\Kcal,\Sigma},T\otimes\La_{\Kcal,\Acal,n})\right).
\end{equation}
Equation \eqref{Reci1} of proposition \ref{PropReciprocitylaw} and the fact that $\Cid_{f}$ is a unit in $\La_{\Kcal,\Acal,n}$ then shows that 
\begin{equation}\nonumber
\ord_{\Pcal_{n}}\carac_{\La_{\Kcal,\Acal,n}}\left(\Htilde^{2}_{\alpha,\alpha}(G_{\Kcal,\Sigma},T\otimes\La_{\Kcal,\Acal,n})\right)\geq\ord_{\Pcal_{n}}\left(\Lcal_{\alpha,\alpha}(f)\right)
\end{equation}
as desired. If now $\Pcal=(Y)$, then $\Lcal_{\alpha,\alpha}(f)\modulo\Pcal$ is $\Lcal_{\alpha,\alpha}^{\cyc}(f)$, which is not identically zero by lemma \ref{LemTwoVariablesOneVariable}. As the right-hand side is equal to zero, \eqref{EqDivFundamental} also holds in that case.
\end{proof}
Recall that $\Lcal^{\cyc}_{\alpha,\alpha}(f)$ is the specialization of $\Lcal_{\alpha,\alpha}(f)$ to the cyclotomic line. \begin{Cor}\label{CorWan}
Let $n\in\N$ be an integer. Under the assumptions of proposition \ref{PropAppendix}, 
\begin{equation}\label{EqDivLan}
\ord_{\Pcal}\carac_{\La_{n}}\left(\Htilde^{2}_{\alpha,\alpha}(G_{\Kcal,\Sigma},T\otimes\La_{n})\right)\geq\ord_{\Pcal}\left(\Lcal_{\alpha,\alpha}^{\cyc}(f)\right)\La_{n}
\end{equation}
for all height-one prime $\Pcal\in\Spec\La_{n}$.
\end{Cor}
\begin{proof}
By proposition \ref{PropSelmerControl} (and its proof), there is an isomorphism
\begin{equation}\nonumber
\RGammalpha(G_{\Kcal,\Sigma},T\otimes\Lambda_{\Kcal,\Acal,n})\Ltenseur_{\La_{\Kcal,\Acal,\infty}}\La_{n}\simeq\RGammalpha(G_{\Kcal,\Sigma},T\otimes\La_{n})
\end{equation}
and hence an isomorphism
\begin{equation}\label{EqIsomSpec}
\Htilde^{2}_{\alpha,\alpha}(G_{\Kcal,\Sigma},T\otimes\La_{\Kcal,\Acal,\infty})\tenseur_{\La_{\Kcal,\Acal,\infty}}\La_{n}\simeq\Htilde^{2}_{\alpha,\alpha}(G_{\Kcal,\Sigma},T\otimes\La_{n})
\end{equation}
as $\Htilde^{3}_{\alpha,\alpha}(G_{\Kcal,\Sigma},T\otimes\La_{\Kcal,\Acal,\infty})$ vanishes. Let $\Pcal\notin S_{v_0}(\Fcal^+(D_{\Acal}))\cup S_{\bar{v}_0}(\Fcal^-(D_{\Acal}))$ be a height-one prime of $\La_{\Kcal,\Acal,\infty}$ and denote by $\Pcal_{\cyc,n}$ its image through the natural map to $\La_{n}$. Let $\Lcal^{\cyc,n}_{\alpha,\alpha}(f)$ the image of $\Lcal^{\cyc}_{\alpha,\alpha}(f)$ in $\La_{n}$. Then
\begin{equation}\nonumber
\ord_{\Pcal_{\cyc,n}}\left(\Lcal^{\cyc}_{\alpha,\alpha}(f)\right)\leq\ord_{\Pcal_{\cyc,n}}\left(\left(\carac_{\La_{\Kcal,\Acal,\infty}}\Htilde^{2}_{\alpha,\alpha}(G_{\Kcal,\Sigma},T\otimes\La_{\Kcal,\Acal,\infty})\right)\tenseur_{\La_{\Kcal,\Acal,\infty}}\La_{n}\right)
\end{equation}
by \eqref{EqDivFundamental} and so
\begin{equation}\nonumber
\ord_{\Pcal_{\cyc,n}}\left(\Lcal^{\cyc}_{\alpha,\alpha}(f)\right)\leq\ord_{\Pcal_{\cyc,n}}\carac_{\La_{n}}\left(\Htilde^{2}_{\alpha,\alpha}(G_{\Kcal,\Sigma},T\otimes\La_{\Kcal,\Acal,\infty})\tenseur_{\La_{\Kcal,\Acal,\infty}}\La_{n}\right)
\end{equation}
as characteristic ideals can only shrink by specialization. Combining this with \eqref{EqIsomSpec}, we obtain
\begin{equation}\nonumber
\ord_{\Pcal_{\cyc,n}}\left(\Lcal^{\cyc}_{\alpha,\alpha}(f)\right)\leq\ord_{\Pcal_{\cyc,n}}\carac_{\La_{n}}\left(\Htilde^{2}_{\alpha,\alpha}(G_{\Kcal,\Sigma},T\otimes\La_{n})\right)
\end{equation}
Suppose that $\Pcal$ now belongs to $S_{v_0}(\Fcal^+(D_{\Acal}))\cup S_{\bar{v}_0}(\Fcal^-(D_{\Acal}))$. Then its intersection $\Pcal_{\cyc,n}$ with $\La_{\infty}$ is the unit ideal and so the inequality \eqref{EqDivLan} still holds.
\end{proof}
Once corollary \ref{CorWan} is known, theorem \ref{TheoSansP} follows from \cite[Theorem 5.4]{JayPreprint}. As this reference might not be widely available, we sum up the argument below.
\begin{proof}[Proof of theorem \ref{TheoSansP}]
Letting $n$ go to infinity in corollary \ref{CorWan} yields the divisibility 
\begin{equation}\nonumber
\left(\Lcal_{\alpha}(f)\right)\La_{\infty}
|\carac_{\La_{\infty}}\Htilde^{2}_{\alpha,\alpha}(G_{\Q,\Sigma},T\tenseur\La_{\infty}).
\end{equation}
According to theorem \ref{TheoKatoIntro}, there is on the other hand a divisibility
\begin{equation}\nonumber
\carac_{\La_{\Iw}[1/p]}H^{2}_{\et}(\Z[1/p],T\tenseur\La_{\Iw}[1/p])|\carac_{\La_{\Iw}[1/p]}H^{1}_{\et}(\Z[1/p],T\tenseur\La_{\Iw}[1/p])/\La_{\Iw}[1/p]\cdot\z(f)_{\Iw}.
\end{equation}
After extension of coefficients to $\La_{\infty}$, localization at $p$ and projection to $\Fcal^{-}$, the regulator map sends Kato's zeta element to $\Lcal_{\alpha}$, thus just as before \cite[Theorem 16.6]{KatoEuler} yields
\begin{equation}\nonumber
\carac_{\La_{\infty}}\Htilde^{2}_{\alpha,\alpha}(G_{\Q,\Sigma},T\tenseur\La_{\infty})|\left(\Lcal_{\alpha}(f)\right)\La_{\infty}.
\end{equation}
Equations \eqref{EqTheoSansPanalytic} and \eqref{EqTheoSansP} follow. 
\end{proof}
\subsubsection{Powers of $p$}\label{Section 4.4}
In this subsection, we provide the ingredient which remains missing after theorem \ref{TheoSansP} to establish conjecture \ref{ConjIMCweak} for the modular motive attached to $f$, namely we show that the power of $p$ dividing $\carac_{\La_{\Iw}}H^{2}_{\et}\left(\Z[1/p],T(f)_{\Iw}\right)$ is at least as large as predicted by the Iwasawa Main Conjecture. The main idea is as follows: as we already know that $\carac_{\La_{\Iw}}H^{2}_{\et}\left(\Z[1/p],T(f)_{\Iw}\right)$ divides $\carac_{\La_{\Iw}}H^{1}_{\et}\left(\Z[1/p],T(f)_{\Iw}\right)/\La_{\Iw}\cdot\z(f)_{\Iw}$, it is enough to show (as we have seen at the end of the control theorem section) that that there is an arithmetic point of $\Spec\La_{\Iw}$ where the orders of a certain Selmer group is as predicted by Kato's Iwasawa main conjecture.

\begin{Def}
A \emph{cyclotomic arithmetic point} $\tilde{\phi}\in\Spec\Oiwa$ of conductor $p^n$ is a $\bar{\mathbb{Q}}_p$-point mapping $(1+X)$ to $\zeta(1+p)^{1+j}$ for $\zeta$ a root of unity of order $p^{n}$ and $0\leq j\leq k-2$ an integer. If $\ptilde$ has values in a discrete valuation ring $\Ocal_{\ptilde}$, we write $T_{\ptilde}$ for $\Tfiwa\tenseur_{\Oiwa,\ptilde}\Ocal_{\ptilde}$.
\end{Def}
If $\ptilde$ is an arithmetic point of conductor $p^{n}$ and if $\chi_{\ptilde}:\left(\Z/p^{n}\Z\right)\croix\fleche\C\croix$ is a character, then the specialization to $\tilde{\phi}$ of the $p$-adic $L$-function $\Lcal_{\alpha}(f)$ interpolates the $L$-value $L(f,\chi_{\tilde{\phi}}^{-1},1+j)$ up to normalization factors. We also write $\tilde{\phi}^{-1}$ for the arithmetic point which is symmetric to $\tilde{\phi}$ with respect to the central point of the functional equation. At $\tilde{\phi}^{-1}$, the $L$-values interpolated are at the critical integers $k-1-j$.

We recall that $\Hun(G_{\Kcal_{\bar{v}_{0}}},T_{\ptilde})$ is a free $\Ocal_{\ptilde}[G_{n}]$-module of rank 2 and write $(v_1, v_2)$ for an $\Ocal_{\ptilde}[G_{n}]$-basis of $\Hun(G_{\Kcal_{\bar{v}_{0}}},T_{\ptilde})$ such that $v_1$ is a generator of $\Hunf(G_{\Kcal_{\bar{v}_{0}}},T_{\ptilde})$. We sometimes write them as $v_{1,\bar{v}_0}, v_{2,\bar{v}_0}$. 
Write $\Gamma_\mathcal{K}=\Gamma_{\cyc}\times\Gamma_{\bar{v}_0}$. Let $\Ycal,\Xcal$ and $\pr_{\cyc}$ be respectively $\Spec\La_{\Kcal},\Spec\La_{\Iw}$ and the natural projection
\begin{equation}\nonumber
\pr_{\cyc}:\Ycal\fleche\Xcal.
\end{equation}
Define the fiber $\Ycal_{\ptilde}$ at $\ptilde$ to be $\Ycal\times_{\Xcal,\ptilde}\Spec\Ocal_{\ptilde}$ and denote by $\BF_{\alpha}^{\Ycal_{\ptilde}}$ the one variable family $\BF_{\alpha,\ptilde}$ of Beilinson-Flach element in Definition \ref{DefBF}.
\begin{Def}
An arithmetic point $\tilde{\phi}\in\mathcal{X}$ is \emph{generic} if $L(f,\chi_{\tilde{\phi}},1+j)\not=0$ and if the restriction of $\mathcal{L}^{\mathrm{Gr}}_{\mathcal{K}}(f)$ to ${\mathcal{Y}_{\tilde{\phi}}}$ is not identically $0$.
\end{Def}
It is clear that all but finite many arithmetic points are generic. We fix an arithmetic point $\ptilde$ of conductor $p^r$ and equal to $\chi_{\cyc}^{1+j-\frac{k}{2}}$ times a finite order $\chi_{\tilde{\phi}}$ and such that $\tilde{\phi}$ and $\tilde{\phi}^{-1}$ are both generic. We freely consider $(v_{1,\bar{v}_0},v_{2,\bar{v}_0})$ as a basis of the free rank-two $\Ocal_{\ptilde}[[\Gamma_{\bar{v}_0}]]$-module $H^1(G_{\qp}, T_{\tilde{\phi}}\otimes_{\Ocal_{\ptilde}}\mathcal{O}_{\mathcal{Y}_{\tilde{\phi}}})$.

\paragraph{Remark:}The family of $p$-adic Galois representation deforming $T$ alongside $\Ycal_{\ptilde}$ and $\Ycal_{\ptilde^{-1}}$ is a family of twists of $T_{f}$ which are unramified as $G_{\qp}$-representations (where we identify $G_{\qp}$ with $G_{\Kcal_{v_{0}}}$).

Define an element $\mathcal{L}^1_{\tilde{\phi}}\in\Ocal_{\ptilde}[[\Gamma_{\bar{v}_0}]]$ such that along $\mathcal{Y}_{\tilde{\phi}}$,
\begin{equation}\label{EqLalongYcal}
\BF^{\Ycal_{\ptilde}}_{\alpha}\equiv (k-2-j)!\left(\mathcal{L}^1_{\tilde{\phi}}\right)\left(G(\chi^{-1}_{\tilde{\phi}})
\left(\frac{\beta_f}{p^{1+j}}\right)^r\right) v_{2,\bar{v}_0}\modulo v_{1,\bar{v}_0}.
\end{equation}
Then $\mathcal{L}^1_{\tilde{\phi}}(0)\not=0$ and
\begin{equation}\nonumber
\alpha^{2r}_f\cdot\frac{\ptilde\left({\Lcal_{\alpha}^{\cyc}(f)}{\Lcal_{\alpha}^{\cyc}(f\tenseur\chi_{\Kcal})}\right)}
{G(\chi_{\ptilde})^2p^{2rj}}=\frac{(j!)^2L_\mathcal{K}(f,\chi^{-1}_{\ptilde},1+j)}{(2\pi i)^{2+2j}\Omega_f^+\Omega_f^-}.
\end{equation}
From our assumption on $\ptilde$, we deduce
\begin{equation}\nonumber
\alpha^{2r}_f\cdot\frac{\ptilde({\Lcal_{\alpha}^{\cyc}(f)}{\Lcal_{\alpha}^{\cyc}(f\tenseur\chi_{\Kcal})})}
{G(\chi_{\ptilde})^2p^{2rj}}=\frac{(j!)^2L_\mathcal{K}(f,\chi^{-1}_{\ptilde},1+j)}{(2\pi i)^{2+2j}\Omega_f^+\Omega_f^-}.
\end{equation}
Let
\begin{equation}\nonumber
\Bcal\eqdef\frac{h_{\Kcal}
\Lcal^{\mathrm{Katz}}_{\Kcal}c\cdot v^+}{v^-}
\end{equation}
be the displayed element \eqref{EqKatz} in lemma \ref{LemIntegrality}. Let $\phi\in\Spec\zp[[Y]]$ be a generic point. Then 
\begin{equation}\label{EqRecPhi}
\mathrm{Reg}_{v_0,\mathcal{F}^+}(\loc_{v_{0}}\BF_{\alpha,\phi})=\frac{1}{\mathfrak{C}_f}\phi(\mathcal{B})(k-2-j)!\phi\left(\Lcal^{\Gr}_{\Kcal}(f)\right)
G(\chi^{-1}_{\ptilde})\cdot\left(\frac{\phi(\alpha_{\mathbf{g}})\beta_f}{p^{1+j}}\right)^r
\omega_f^\vee
\end{equation}
according to \cite[Theorem 7.1.4]{KLZ}. Note that $\beta_f, j, r$ only depend on $\ptilde$ and that $\alpha_{\mathbf{g}}$ is an element in $\mathbb{Z}_p[[\Gamma_{v_0}]]$ such that $\alpha_{\mathbf{g}}(0)=1$. If $\Cid_{f}$ is the number defined after \eqref{EqCompPoincare}, we have
\begin{equation}\nonumber
\exp^*_{\bar{v}_0}\BF^{}_{\alpha,\ptilde}=\frac{1}{j!\mathfrak{C}_f}\ptilde
({\Lcal_{\alpha}^{\cyc}(f)}\Lcal_{\alpha}^{\cyc}(f\tenseur\chi_{\Kcal}))\cdot G(\chi^{-1}_{\ptilde})\cdot\left(\frac{\alpha_f\beta_g}{p^{1+j}}\right)^r\cdot\omega_f.
\end{equation}
Hence
\begin{align*}
\exp^*_{\bar{v}_0}\BF^{}_{\alpha,\ptilde}&=\frac{L_\mathcal{K}(f,\chi^{-1}_{\ptilde},1+j)\cdot G(\chi^{-1}_{\ptilde})}{j!\mathfrak{C}_f(2\pi i)^{2+2j}\Omega_f^+\Omega_f^-}G^2(\chi_{\ptilde})\left(\frac{\alpha_f}{p^{1+j}}\right)^rp^{2rj}\alpha_f^{-2r}\cdot\omega_f&\\
&=\frac{L_\mathcal{K}(f,\chi^{-1}_{\ptilde},1+j)}{j!\mathfrak{C}_f(2\pi i)^{2+2j}\Omega_f^+\Omega_f^-}G(\chi_{\ptilde})\alpha_f^{-r}{p^{rj}}\cdot\omega_f.&\\&&
\end{align*}
Note here that the factor $\mathcal{E}(f)\mathcal{E}(f^*)$ in the interpolation formula for Rankin-Selberg $p$-adic $L$-function is cancelled by the factor $(1-\frac{\beta}{\alpha})(1-\frac{\beta}{p\alpha})$ in \cite[Corollary 6.4.3]{KLZ}. Thus
\begin{equation}\label{EqExpDualAtPhiTilde}
\exp^*_{\bar{v}_0}\ptilde(\mathrm{BF}_{\alpha})=\frac{L_{\mathcal{K}}(f,\chi^{-1}_{\ptilde},1+j)}{j!\mathfrak{C}_f\Omega_f^+\Omega_f^-(2\pi i)^{2+2j}}G(\chi_{\ptilde})\left(\frac{\beta_f}{p^{1+j}}\right)^r\cdot p^{r(2j+2-k)}\cdot\omega_f.
\end{equation}
Now we repeat these constructions for the arithmetic point $\tilde{\phi}^{-1}$. We define a basis $(v_{1,\bar{v}_0}, v_{2,\bar{v}_0})$ at $\tilde{\phi}^{-1}$ as at $\tilde{\phi}$. To emphasis the dependence on the point we use $\tilde{\phi}(v_{i,\bar{v}_0})$ or $\tilde{\phi}^{-1}(v_{i,\bar{v}_0})$ to denote them. Similarly as for $\tilde{\phi}$, we may compute the image through the regulator map of the specialization of the Beilinson-Flach class at $\ptilde^{-1}$ and their images though arithmetic points $\phi$ and $\phi^{-1}$. As above, we obtain
\begin{equation}\label{EqRegulatorPhiTildeInverse}
\mathrm{Reg}_{v_0,\mathcal{F}^+}(\loc_{v_{0}}\BF_{\alpha,\phi^{-1}})=\frac{1}{\mathfrak{C}_f}\phi^{-1}(\mathcal{B})j!\phi^{-1}(\Lcal^{\Gr}_{\Kcal}(f)) G(\chi_{\tilde{\phi}})\left(\frac{\phi^{-1}(\alpha_\mathbf{g})\beta_f}{p^{1+j}}\right)^r \omega_f^\vee
\end{equation}
and
\begin{equation}\label{EqExpPhiTildeInverse}
\exp^*_{\bar{v}_0}\tilde{\phi}^{-1}(\BF_{\alpha})=\frac{L_{\mathcal{K}}(f,\chi^{-1}_{\tilde{\phi}^{-1}},k-1-j)}{(k-2-j)!\mathfrak{C}_f(2\pi i)^{2k-2-2j}\Omega_f^+\Omega_f^-}
\left(G(\chi_{\tilde{\phi}})\left(\frac{\beta_f}{p^{1+j}}\right)^r\right)(p^{2r(k-2-j)+r}/p^{(k-1)r})\omega_f.
\end{equation}
Putting everything together, we obtain the equation
\begin{equation}\label{L^1}\mathcal{L}^1_{\tilde{\phi}}(0)=\frac{j!}{(k-2-j)!}\frac{G(\chi_{\tilde{\phi}})}{G(\chi^{-1}_{\tilde{\phi}})}\frac{L_{\mathcal{K}}(f,1+j,\chi^{-1}_{\tilde{\phi}})(p^{2jr+r}/p^{(k-1)r})}
{\mathfrak{C}_f(2\pi i)^{2+2j}\Omega_f^+\Omega_f^-\exp^*\tilde{\phi}(v_{2,\bar{v}_0})\cdot}.
\end{equation}

In definition \ref{DefPseudoTamagawaBD}, we defined
\begin{equation}\nonumber
H^1_f(G_{\Kcal_{v_{0}}}, T_{\tilde{\phi}}\otimes\Ocal_{\ptilde}[[\Gamma_{\bar{v}_0}]])
\end{equation}
such that for any height one prime $\Pcal\in\Spec\Ocal_{\ptilde}[[\Gamma_{\bar{v}_0}]]$, the $H^1_f(G_{\Kcal_{v_{0}}}, T_{\tilde{\phi}}\otimes\Ocal_{\ptilde}[[\Gamma_{\bar{v}_0}]])_{\Pcal}$
is an $\Ocal_{\ptilde}[[\Gamma_{\bar{v}_0}]]_{\Pcal}$-direct summand of
$$H^1(G_{\mathcal{K}_{v_0}}, T_{\tilde{\phi}}\otimes\Ocal[[\Gamma_{\bar{v}_0}]])_{P'}\simeq \Ocal_{\ptilde}[[\Gamma_{\bar{v}_0}]]^2_{\Pcal}.$$
\begin{Def}
The \emph{unramified local Selmer condition} at $v_0$ is the subspace
\begin{equation}\nonumber
\Hunf(G_{\Kcal_{v_{0}}}, T_{\tilde{\phi}}\otimes_{\Ocal_{\ptilde}}\Ocal_{\ptilde}[[\Gamma_{\bar{v}_0}]])
\end{equation}
along the one-variable family over $\Ocal_{\ptilde}[[\Gamma_{\bar{v}_0}]]$.
\end{Def}
Let $\Kcal\subset_{f} F\subset\Kcal_{\infty}$ be a finite subextension and let $v$ be a finite place of $\Ocal_{F}$. The Greenberg local condition $\Hun(G_{F_{v}},A)\subset\Hun(G_{F_{v}},A)$ at $v$ is defined to be
\begin{equation}\nonumber
\image\left(\Hun(G_{F_{v}}/I_{v},V^{I_{v}})\fleche\Hun(G_{F_{v}}/I_{v},A^{I_{v}})\right)\subset\Hun(G_{F},A)
\end{equation}
if $v\nmid p$, to be $\Hun_{\Gr}(G_{F_{v}},A)$ if $v|v_{0}$ and to be 0 if $v|\bar{v}_{0}$. Recall that $\Gamma_{\bar{v}_0}=\mathrm{Gal}({\Kcal}^{\bar{v}_0}/{\Kcal})$ is the Galois group of the maximal subextension of $\Kcal_{\infty}$ unramified outside $\bar{v}_{0}$. Using the identification of $\Ocal[[Y]]$ with $\Ocal[[\Gamma_{\bar{v}_0}]]$ through the map $Y\longmapsto\gamma_{\bar{v}_0}-1$, we identify the $\Ocal_{\ptilde}[[Y]]$-modules
\begin{equation}\nonumber
H^1(G_{\Kcal,\Sigma},T_{\ptilde}\tenseur\Ocal_{\Ycal_{\ptilde}}),\limproj{\Kcal\subset F\subset\Kcal^{\bar{v}_{0}}}\Hun(G_{F,\Sigma},T_{\ptilde})
\end{equation}

The following proposition is the payoff of our working alongside the unramified family $\Ycal_{\ptilde}$ passing through the classical point $\ptilde$ of the cyclotomic Iwasawa algebra.
\begin{Prop}\label{PropContainmentYcal}
\begin{equation}\nonumber
\prod_v c_{\mathcal{K},\tilde{\phi},v}(f)\prod_v c_{\mathcal{K},\tilde{\phi}^{-1},v}(f)\carac_{\Ocal_{\ptilde}}(X_{\mathrm{BK},v_1,\tilde{\phi}})
\carac_{\Ocal_{\ptilde}}(X_{\mathrm{BK},v_1,\tilde{\phi}^{-1}})
\end{equation}
is contained in
\begin{equation}\nonumber
\left(\frac{L_\mathcal{K}(f,1+j,\chi^{-1}_{\tilde{\phi}})}{(2\pi i)^{2+2j}\Omega_f^+\Omega_f^-}\right)\left(\frac{L_{\mathcal{K}}(f,k-1-j,
\chi^{-1}_{\tilde{\phi}^{-1}})}{(2\pi i)^{2k-2-2j}\Omega_f^+\Omega_f^-}\right).
\end{equation}
\end{Prop}
\begin{proof}

Recall that $H^1(G_{\mathcal{K}_{\bar{v}_0}}, T_{\tilde{\phi}}\otimes \mathcal{O}_{\mathcal{Y}_{\tilde{\phi}}})$ is a free rank two module over $\mathcal{O}_{\mathcal{Y}_{\tilde{\phi}}}$, and that we lifted $v_1$ to a generator of a rank one direct summand of this module, which we still denote as $v_1$.
Now we look at the following exact sequences of $\Ocal_{\ptilde}[[U]]$-modules
\begin{equation}
\label{Poitou-Tate one}0\rightarrow H^1_{{\ur},\mathrm{rel}}(G_{\Kcal,\Sigma}, T_{\tilde{\phi}}\otimes\mathcal{O}_{\mathcal{Y}_{\tilde{\phi}}})\rightarrow \Hun_{\ur}(G_{\Kcal_{v_0}}, T_{\tilde{\phi}}\otimes\mathcal{O}_{\mathcal{Y}_{\tilde{\phi}}})\rightarrow X_{\mathrm{rel},\mathrm{str}}\rightarrow X_{{\ur},\mathrm{str}}\rightarrow 0
\end{equation}

\begin{equation}\label{Poitou-Tate two}
0\rightarrow H^1_{{\ur},\mathrm{rel}}(G_{\Kcal,\Sigma}, T_{\tilde{\phi}}\otimes\mathcal{O}_{\mathcal{Y}_{\tilde{\phi}}})\rightarrow \frac{H^1(\mathcal{K}_{\bar{v}_0}, T_{\tilde{\phi}}\otimes\mathcal{O}_{\mathcal{Y}_{\tilde{\phi}}})}{\mathcal{O}_{\mathcal{Y}_{\tilde{\phi}}}v_1}\rightarrow X_{{\ur},v_1}\rightarrow X_{{\ur},\mathrm{str}}\rightarrow 0.
\end{equation}
Recall we wrote $\Ocal[[\Gamma_{\bar{v}_0}]]=\Ocal[[Y]]$ for the variable $Y=\gamma_{\bar{v}_0}-1$ and definition \ref{DefPseudoTamagawaBD}. By Corollaries \ref{specialization} and \ref{CorDirectSummand} there is an isomorphism
$$\Hun_{\ur}(G_{\Kcal_{v_0}}, T\otimes\Ocal_{\Ycal_{\ptilde}})\simeq \Ocal_{\ptilde}[[Y]]$$
which interpolates $1/\mathsf{b}_{\tilde{\phi},p}$ times the regulator map, divided by the specialization of $\rho(d)$ there. Theorem \ref{GBMC} establishes an inclusion in the Greenberg-Iwasawa main conjecture. Combined with the definition of Tamagawa numbers in definition \ref{DefTamagawa}, equation \eqref{EqLalongYcal} describing the interpolation properties of the Beilinson-Flach classes, equations \eqref{EqRecPhi},\eqref{EqExpDualAtPhiTilde},\eqref{EqRegulatorPhiTildeInverse}, \eqref{EqExpPhiTildeInverse} computing the image of $\BF_{\alpha,\ptilde}$ through the regulator and dual exponential maps, the integrality of $\mathcal{B}$ in lemma \ref{LemIntegrality} and the short exact sequence \eqref{Poitou-Tate one}, it therefore implies that
\begin{equation}\nonumber
\frac{1}{\mathsf{b}_{\tilde{\phi},p}}\cdot\frac{1}{\mathsf{b}_{\tilde{\phi}^{-1},p}}\carac_{\Ocal_{\ptilde}[[Y]]}\left(X_{
{\ur},\mathrm{str},\tilde{\phi}}\right)\carac_{\Ocal_{\ptilde}[[Y]]}\left(X_{{\ur},\mathrm{str},\tilde{\phi}^{-1}}\right)
\end{equation}
is contained in
\begin{equation}\nonumber
\carac_{\Ocal_{\ptilde}[[Y]]}\left(\frac{H^1_{{\ur},\mathrm{rel}}(G_{\Kcal,\Sigma}, T_{\tilde{\phi}}\otimes\mathcal{O}_{\mathcal{Y}_{\tilde{\phi}}})}{\frac{\mathfrak{C}_f}{(k-2-j)!(G(\chi^{-1}_{\tilde{\phi}})
(\frac{\beta_f}{p^{1+j}})^r)}\Ocal_{\ptilde}\cdot[[Y]]\BF^{\Ycal_{\ptilde}}_{\alpha}}\right)
\cdot \carac_{\Ocal_{\ptilde}[[Y]]}\left(\frac{H^1_{{\ur},\mathrm{rel}}(G_{\Kcal,\Sigma}, T_{\tilde{\phi}^{-1}}\otimes\mathcal{O}_{\mathcal{Y}_{\tilde{\phi}^{-1}}})}{\frac{\mathfrak{C}_f}{j!(G(\chi_{\tilde{\phi}})
(\frac{\beta_f}{p^{k-1-j}})^r)}\Ocal_{\ptilde}\cdot[[Y]]\BF^{\Ycal_{\ptilde^{-1}}}_{\alpha}}\right)
\end{equation}
up to powers of $Y$ (as fractional ideals for the right hand side). Here, the containments
\begin{equation}\nonumber
\frac{1}{\mathsf{b}_{\tilde{\phi},p}}\carac_{\Ocal_{\ptilde}[[Y]]}\left(X_{{\ur},\mathrm{str},\tilde{\phi}}\right)\subseteq\carac_{\Ocal_{\ptilde}[[Y]]}\left(\frac{H^1_{{\ur},\mathrm{rel}}(G_{\Kcal,\Sigma}, T_{\tilde{\phi}}\otimes\mathcal{O}_{\mathcal{Y}_{\tilde{\phi}}})}{\frac{\mathfrak{C}_f}{(k-2-j)!(G(\chi^{-1}_{\tilde{\phi}})
(\frac{\beta_f}{p^{1+j}})^r)}\Ocal_{\ptilde}\cdot[[Y]]\BF^{\Ycal_{\ptilde}}_{\alpha}}\right)
\end{equation}
and
\begin{equation}\nonumber
\frac{1}{\mathsf{b}_{\tilde{\phi}^{-1},p}}\carac_{\Ocal_{\ptilde}[[Y]]}\left(X_{{\ur},\mathrm{str},\tilde{\phi}^{-1}}\right)\subseteq\carac_{\Ocal_{\ptilde}[[Y]]}\left(\frac{H^1_{{\ur},\mathrm{rel}}(G_{\Kcal,\Sigma}, T_{\tilde{\phi}^{-1}}\otimes\mathcal{O}_{\mathcal{Y}_{\tilde{\phi}^{-1}}})}{\frac{\mathfrak{C}_f}{(k-2-j)!(G(\chi_{\tilde{\phi}})
(\frac{\beta_f}{p^{1+j}})^r)}\Ocal_{\ptilde}\cdot[[Y]]\BF^{\Ycal_{\ptilde}}_{\alpha}}\right)
\end{equation}
are both deduced from the containment
\begin{equation}\nonumber
\carac(X^{\mathrm{Gr}}_{\mathcal{K}}(f))\subseteq (\Lcal^{\Gr}_{\Kcal}(f))
\end{equation}
of theorem \ref{GBMC} by specializing to $\Ycal_{\ptilde}$ and $\Ycal_{\ptilde^{-1}}$ and using the fact that characteristic ideals can only shrink by specialization (in fact, any prime $v$ of $\mathcal{K}$ not dividing $p$ is finitely decomposed in the $\mathbb{Z}_p$-extension $\mathcal{K}^{\bar{v}_0}/\mathcal{K}$ and is completely split in $\mathcal{K}_\infty/\mathcal{K}^{\bar{v}_0}$ so the characteristic ideal of $X^{\mathrm{Gr}}_{\mathcal{K}}(f)$ specialized to $\Ocal_{\ptilde}[[Y]]$ and $\Ocal_{\ptilde^{-1}}[[Y]]$ is exactly the characteristic ideal of $X_{{\ur},\mathrm{str},\tilde{\phi}}$ and $X_{{\ur},\mathrm{str},\tilde{\phi}^{-1}}$ respectively).

Just as before applying the Poitou-Tate exact sequences (\ref{Poitou-Tate two}), we get the lower bound for the $({\ur},v_1)$ Selmer groups over $\Ocal[[\Gamma_{\bar{v}_0}]]$ by $\mathcal{L}^1_{\tilde{\phi}}$ up to powers of $Y$. Recall also that $Y$ is not a divisor of $\mathcal{L}^1_{\tilde{\phi}}$ and $\mathcal{L}^1_{\tilde{\phi}^{-1}}$ since $\mathcal{L}^1_{\tilde{\phi}}(0)$ and $\mathcal{L}^1_{\tilde{\phi}^{-1}}(0)$ are nonzero. Specializing $Y$ to zero, we obtain that
\begin{equation}\nonumber
\prod_v c_{\mathcal{K},\tilde{\phi},v}(f)\prod_v c_{\mathcal{K},\tilde{\phi}^{-1},v}(f)\carac_{\Ocal_{\ptilde}}(X_{\mathrm{BK},v_1,\tilde{\phi}})
\carac_{\Ocal_{\ptilde}}(X_{\mathrm{BK},v_1,\tilde{\phi}^{-1}})
\end{equation}
is contained in
\begin{equation}\nonumber
\left(\frac{L_\mathcal{K}(f,1+j,\chi^{-1}_{\tilde{\phi}})}{(2\pi i)^{2+2j}\Omega_f^+\Omega_f^-}\right)\left(\frac{L_{\mathcal{K}}(f,k-1-j,
\chi^{-1}_{\tilde{\phi}^{-1}})}{(2\pi i)^{2k-2-2j}\Omega_f^+\Omega_f^-}\right).
\end{equation}
Above, the $v$ may or may not divide $p$. Pay attention that in the formula (\ref{L^1}), the $p^{2jr+r}/p^{(k-1)r}$ and $\frac{j!}{(k-2-j)!}$ at $\tilde{\phi}$ and $\tilde{\phi}^{-1}$ cancel out. Here the subscript $\mathcal{K}$ in $c_\mathcal{K}$ means the local Tamagawa numbers over $\mathcal{K}$. It is known that $c_{\mathcal{K},\ell,\phi}(f)=c_{\ell,\phi}(f)\cdot c_{\ell,\phi}(f\tenseur{\chi_\mathcal{K}})$. The subscript $\mathrm{BK}$ stands for the Bloch-Kato Selmer condition at $v_0$. 
\end{proof}
Finally, we record the following.
\begin{Lem}\label{Periods}
The period $\Omega_{f\tenseur{\chi_\mathcal{K}}}^\mp$ is an $\Ocal$-multiple of $\Omega_f^\pm$.
\end{Lem}
\begin{proof}
This a much easier variant of \cite[Lemma 9.6]{SkinnerZhang}. As our running assumptions are different from those of \cite{SkinnerZhang} (and our result correlatively much weaker), we recall the proof.

The eigencuspforms $f$ and $f\tenseur{\chi_\mathcal{K}}$ are new of level $N$ and of level $M$ for some $M$ dividing $ND_{\Kcal}^2$. Consider the map
$$H^1(\Gamma_0(N),L^{k-2}_{/\Ocal})_{\mathfrak{m}_f}\fleche H^1(\Gamma_0(ND_{\Kcal}^2), L^{k-2}_{/\Ocal})_{\mathfrak{m}_f}$$
constructed in \cite[Proof of Lemma 9.6]{SkinnerZhang}. As $\mgot_{f}$ is non-Eisenstein, this is a map between free modules which sends the canonical differentials  $\omega_f^\pm$ are mapped to $\omega^{\mp}_{f\tenseur{\chi_\mathcal{K}}}$ multiplied by the Gauss sums attached to $\chi_{\Kcal}$ and the elements $\gamma_f^+,\gamma_{f}^{-}$ forming a basis of the left-hand side to some $\Ocal$-multiple of $\gamma^-_{f\tenseur{\chi_\mathcal{K}}},\gamma^+_{f\tenseur{\chi_\mathcal{K}}}$. Hence $\Omega_{f\tenseur{\chi_\mathcal{K}}}^\mp$ is an $\Ocal$-multiple of $\Omega_f^\pm$ divided by the Gauss sum attached to $\chi_{\Kcal}$, which is a $p$-adic unit as $p$ is prime to $D_{\Kcal}$.
\end{proof}
Let $\mathscr{F}_1$ be the image of $\z(f)_{\Iw}$ in
$\frac{H^1(\mathbb{Q}^S/\mathbb{Q},\mathbf{T})}{\Lambda\mathbf{v}}$, and $\mathscr{F}$ be the characteristic polynomial of the $\mathbf{v}^\vee$-Selmer group $X_{\mathbf{v}^\vee}(f)$ as in Section \ref{Control}. We have the following
\begin{Lem}\label{LemIMCwithpadic}
Conjecture \ref{ConjIMCweak} for $M(f)$ is equivalent to the equality
$$(\mathscr{F}_1)=(\mathscr{F})$$
as principal ideals of $\Lambda_{\Iw}$.
\end{Lem}
\begin{proof}
This is proved using Poitou-Tate duality in the same way as in the proof of \cite[Theorem 7.4]{KobayashiIMC}, once it is observed that the role $\mathscr{L}^\pm_p$ in \emph{loc.cit} is played here by $\mathscr{F}_1$.
\end{proof}
We finally obtain the following theorem.
\begin{Theo}\label{TheoCrys}
Let $f\in S_{k}(\Gamma_{0}(N))$ be a normalized eigencuspform of weight $k\geq2$ satisfying the following hypotheses. 
\begin{enumerate}
\item The $G_{\qp}$-representation $\rho_{f}|G_{\qp}$ is crystalline (equivalently $p\nmid N$) and short.
\item The local residual representation $\rhobar_{f}|G_{\qp}$ is absolutely irreducible.
\item There exists $\ell||N$ such that $\dim_{\Fp}\rhobar^{I_{\ell}}=1$ and $\dim_{\Fp}\rhobar^{G_{\ql}}=0$. 
\end{enumerate}
Then conjecture \ref{ConjIMCweak} holds for $M(f)$.
\end{Theo}
\begin{proof}
Similarly to the proof of \cite[Theorem 17.4]{KatoEuler} and \cite[1.3]{KobayashiIMC}, theorem \ref{TheoKatoIntro} and lemma \ref{LemIMCwithpadic} shows that $\mathscr{F}_1$ is divisible by $\mathscr{F}$. Hence, it is enough to check that $\tilde\phi(\mathscr{F})$ is divisible by $\tilde\phi(\mathscr{F}_1)$ for the arithmetic point $\tilde{\phi}$. Let $P$ be the height one prime corresponding to $\tilde{\phi}$.

We have $\tilde\phi(\mathscr{F}_1)=\exp^*(\tilde\phi(\z(f)_{\Iw}))/c'_{P,p}$.
On the other hand Kato proved that
$$\exp^*\tilde{\phi}(\z(f)_{\Iw})=\frac{L(f,1+j,\chi^{-1}_{\tilde{\phi}})}{(2\pi i)^{1+j}\Omega_f^{(-1)^j}}\omega_f.$$
The theorem then follows from proposition \ref{PropContainmentYcal}, equality (\ref{ctTh1}) and lemma \ref{Periods}.
\end{proof}

\subsection{The $p$-irregular case}\label{irregular}
Previously, we have assumed the Satake parameters at $p$ satisfy $\alpha\not=\beta$. Conjecturally this is always true, however only known when $k=2$. For completeness of the result, we treat the case $\alpha=\beta$ as well\footnote{We thank David Loeffler and Chris Williams for discussions on this part, and thank Betina-Williams for writing up \cite{BetinaWilliams} for us.}. Note that since $k$ is even, this can happen only when $k\geq 4$.
In this section we establish (\ref{needed}), which is the only missing ingredient for proving the Iwasawa main conjecture. (Loeffler-Zerbes assumed $\alpha\not=\beta$ in \cite{KLZ}.) This replaces equation \ref{Reci1} in proposition \ref{PropReciprocitylaw} (equation \ref{Reci2} remains unchanged ).
\subsubsection{Geometry of the eigencurve}
We start by briefly discussing the local geometry of the eigencurve near the $p$-irregular point, using freely the notions in \cite{AshStevens} (more details are given in \cite{BetinaWilliams}). Let $U$ be a small affinoid neighborhood in the weight space for the point $k$.
We write $D_k$ ($D_U$) for the coefficient sheave of weight $k$ (over $U$) on the modular curve $Y_0(N)$ used to define modular symbols (following the notation of Ash-Stevens \cite{AshStevens}). We consider the long exact sequence
$$\cdots H^0(Y_0(N), D_k)\rightarrow H^1 (Y_0(N), D_U)\rightarrow H^1(Y_0(N), D_U)\rightarrow H^1(Y_0(N), D_k)\rightarrow H^2(Y_0(N), D_U)\cdots$$
induced from
$$0\rightarrow D_U\rightarrow D_U\rightarrow D_k\rightarrow 0,$$
where the second arrow is multiplying by a uniformizer $m_k$ at $k$.
Ash-Stevens proved that one can do slope decomposition for the cohomology for slope $h$ of $f_\alpha$. This slope is non-critical, \textit{i.e.} the $p$-adic valuation of $\alpha$ is less than $k-1$. According to Coleman inequality (\cite{AshStevens}), the cohomology of slope $\leq h$ and weight $k$ is then classical. On the other hand the Hecke eigensystem for $f$ does not appear in $H^0(Y_0(N), V_k)$, so the localization $H^0(Y_0(N), D_k)^{\leq h}_f$ of $H^0(Y_0(N), D_k)^{\leq h}$ at $f$ vanishes.  Thus the multiplying by $m_k$ map on $H^1(Y_0(N), D_U)^{\leq h}_f$ is injective. So $H^1(Y_0(N), D_U)^{\leq h}_f$ is a torsion-free module over the discrete valuation ring $\mathcal{O}(U)_k$, thus free. Now by \cite[Lemma 2.15]{BarreraDimitrovJorza}, one can lift this to obtain freeness of an appropriate localization of $H^1(Y_0(N), D_U)^{\leq h}$ over a neighborhood of $k$ in $U$. From Coleman inequality and multiplicity one for $\mathrm{GL}_2$, there is a Hecke operator $t$ outside $p$ over $U$, such that the $tH^1(Y_0(N), D_U)^{\leq h}$ is free of rank two over $\mathcal{O}(U)$, which specializes to the $2$-dimensional fixed space by $\Gamma_0(p)$ of $\pi_{f,p}$ (we simply choose it to kill all forms at weight $k$ outside the eigensystem of $f$).  Thus locally the $U_p$ is acting on this rank two module, and satisfies a quadratic equation over the weight ring. It is easy to show that the local coefficient ring of the eigencurve is obtained from joining the $U_p$ operator (on this rank two module), $U_q$ operators for all $q|N$ and the $T_q$ operators for $q\nmid N$ to the weight space (\textit{i.e.} a neighborhood of $k$). It is of degree $2$ over the weight ring. We denote the normalization of the local eigencurve as $\mathcal{C}$. The $\mathcal{C}$ is a smooth curve with one or two irreducible components.  There is a rigid analytic function $Z$ on $\mathcal{C}$ giving the $U_p$ eigenvalue of the corresponding eigenform. Now we take base change extension of the coefficient ring to $\mathcal{C}$. We have the following easy lemma.
\begin{Lem}
The subspace of $tH^1(Y_0(N), D_\mathcal{C})^{\leq h}$ killed by $U_p-Z$ is free of rank one over $\mathcal{C}$.
\end{Lem}
\begin{proof}
This submodule is clearly generically of rank one. We just need to note that the $\alpha$-eigenspace for the $U_p$ operator in $\pi_{f,p}$ is $1$-dimensional to get the lemma.
\end{proof}
\subsubsection{The argument}
Now we consider a Coleman family $\mathcal{C}$ as above in a small neighborhood of $f_\alpha$, and denote the form as $\mathbf{f}$ with coefficient ring a normal domain $\mathbb{I}$, and let $\mathbf{g}$ be a Hida family of CM forms with respect to the field $\mathcal{K}$. In \cite{BetinaWilliams} it is constructed a $3$-variable $p$-adic $L$-function interpolating critical values of Rankin-Selberg $L$-values for $\mathbf{f}$ and $\mathrm{g}$ using modular symbols, which we denote as $\mathcal{L}_{\mathbf{f},\mathbf{g},\mathrm{Symb}}$. There also constructed the two-variable $p$-adic $L$-function of $\mathbf{f}$ over $\mathbb{Q}$ using modular symbols. We refer the details of the notion of Coleman family $\mathbf{f}$ and interpolation properties of $\mathcal{L}_{\mathbf{f},\mathbf{g},\mathrm{Symb}}$ to \emph{loc.cit.}.
\paragraph{Remark:}
Strictly speaking the construction in \cite{BetinaWilliams} gives $L$-values for $\mathbf{f}$ over $\mathcal{K}$ twisted by CM characters corresponding to forms in $\mathbf{g}$, which are nothing but Rankin-Selberg $L$-values for $\mathbf{f}$ and $\mathbf{g}$ as in \cite[Theorem 7.1.5]{KLZ}, except for periods used.

Note that when $\alpha=\beta$ we are not in the ``noble eigenform'' case as defined in \cite[Definition 4.6.3]{KLZ2}. As a consequence, the $\eta_\mathbf{f}$ (our $\mathbf{f}$ is $\mathcal{F}$ there) defined in \cite[Corollary 6.4.3]{KLZ2} is only a meromorphic section instead of a basis of the free rank one module there.
\begin{Prop}\label{Prop.4.37}
We use the same notations as in \cite[Theorem 6.5.9]{KLZ2}. Then the right hand side of identity of \emph{loc.cit.} equals
$$(-1)^{k'-j+1}(k')!\begin{pmatrix}k\\j \end{pmatrix}L(f,g,1+j).$$
\end{Prop}
Note that if $\alpha\not=\beta$, and thus $\mathcal{E}^*(f)$ is nonzero, then the statement is the same as \cite[Theorem 6.5.9]{KLZ2}. When $\alpha=\beta$ the same argument as the proof of \cite[Theorem 6.5.9]{KLZ2}, using the computations in \cite[Proposition 5.3.5]{LLZ} actually gives the above proposition.\\

We denote $\mathcal{L}_{\mathbf{f},\mathbf{g}}$ for the $p$-adic $L$-function of \cite[Theorem 7.1.5]{KLZ}. The $\mathcal{L}_{\mathbf{f},\mathbf{g},\mathrm{Symb}}$ has the same interpolation formula as $\mathcal{L}_{\mathbf{f},\mathbf{g}}$, except that the period is a nonzero constant $\Omega_{\Pi_{f,\mathcal{K}}}$ depending on $f$, while the period factor for $\mathcal{L}_{\mathbf{f},\mathbf{g}}$ is $\langle f,f\rangle\mathcal{E}(f_\alpha)\mathcal{E}^*(f_\alpha)$ (here $f$ is the normalized eigenform of level $N$ whose $\alpha$-stabilization is $f_\alpha$). Note that the latter period takes $0$ in the $\alpha=\beta$ case, which is a key difficulty. Note also that the $\mathcal{E}(f_\alpha)$ and $\mathcal{E}^*(f_\alpha)$ are not $p$-adically rigid analytic. We write $\mathcal{G}$ for the meromorphic function
$\frac{\mathcal{L}_{\mathbf{f},\mathbf{g},\mathrm{Symb}}}{\mathcal{L}_{\mathbf{f},\mathbf{g}}}$.
\paragraph{Remark:}
An automorphic construction of the meromorphic function $\mathcal{G}$ is done in \cite{KimEigencurve}. It is also explained there that it annihilates the congruence module up to a fixed power of $p$. We do not need these facts here.

From the interpolation formulas, we see that $\mathcal{G}$ involves only the variable corresponding to the Coleman family $\mathbf{f}$, namely it is a nonzero element in $\mathrm{Frac}(\mathbb{I})$.
Now we consider an arithmetic point $\phi'$ where $f_{\phi'}$ takes the form $f_\alpha$ we start with, and choose $\mathbf{g}_{\phi'}$ is a CM form with trivial character whose conductor is prime to $p$ and whose weight is lower than that of $f$, and $L(f_{\phi'},g_{\phi'},1+j_{\phi'})$ is a critical value which is nonzero. Note that by our assumption that the weight of $f$ is at least $4$, this is always possible (\textit{e.g.} by taking the critical value to be non-central).

We also consider another arithmetic point $\phi'_0$, where again $f_{\phi'_0}$ is our $f_\alpha$, but the $g_{\phi'_0}$ corresponds to a finite order character of $\Gamma_\mathbb{Q}$. We require the $L(f_{\phi'_0}, g_{\phi'_0},1+j_{\phi'_0})$ to correspond to the nonzero $L$-value $L_\mathcal{K}(f,\chi_\phi,\frac{k}{2})$, where $\chi_\phi$ is the finite order character of $\Gamma_\mathbb{Q}$ in the previous section. (Our convention on the notation is the $j_{\phi'_0}$ denotes a character corresponding to an arithmetic point -- it incorporates an integer together with a finite order Hecke character). Note also that we allow that the $\phi'_0$ and $\phi'$ are not in the same irreducible component of CM families.
\paragraph{Remark:}
If we can prove Proposition \ref{Prop.4.37} at these points $\phi'_0$ then we are done. However the $g_{\phi'_0}$ has weight one and is not covered in the work of Loeffler-Zerbes \cite{KLZ}. Nevertheless we can still achieve the formula (\ref{needed}) by the argument below.

Let $\langle,\rangle$ and $\mathcal{L}$ be the pairing and regulator map as in \cite[Theorem 7.1.5]{KLZ}. We also use the notations $\eta_\mathbf{f}$ and $\omega_\mathbf{g}$ as in \emph{loc.cit.}. Let $\mathrm{BF}_{\mathbf{f},\mathbf{g}}$ be the Beilinson-Flach element constructed for $\mathbf{f}$ and $\mathbf{g}$ as in \cite{KLZ}. Its construction, and its relation to the Rankin-Selberg $p$-adic $L$-function can be done in the same way as in \emph{loc.cit.}, thanks to our previous discussion on the geometry of the eigencurve (especially that it is locally of degree $2$ over the weight space). In particular we have
$$\langle \mathcal{L}(\mathrm{BF}_{\mathbf{f},\mathbf{g}}),\eta_\mathbf{f}\mathcal{G}\otimes\omega_{\mathbf{g}}\rangle=
\mathcal{L}_{\mathbf{f},\mathbf{g},\mathrm{Symb}}.$$
However we need some additional work to prove the following (\ref{need}) and (\ref{needed}). Let $\eta^\alpha_f$ be as in \cite[Theorem 6.5.9]{KLZ2}.
Using proposition \ref{Prop.4.37}, the argument to prove \cite[Theorem 7.1.5]{KLZ} (note that the weight of $\mathbf{g}_{\phi'}$ is at least $2$) gives
\begin{equation}\label{A}\langle\phi'(\mathcal{L}(\mathrm{BF}_{\mathbf{f},\mathbf{g}})),\eta^\alpha_f\otimes\phi'(\omega_{\mathbf{g}})\rangle=
\phi'(\mathcal{L}_{\mathbf{f},\mathbf{g},\mathrm{Symb}})
\frac{\Omega_{\Pi_{f,\mathcal{K}}}}{\langle f,f\rangle}.\end{equation}
By taking specializations, we have
\begin{equation}\label{B}\langle \phi'(\mathcal{L}(\mathrm{BF}_{\mathbf{f},\mathbf{g}})),\phi'(\eta_\mathbf{f}\mathcal{G}\otimes\omega_{\mathbf{g}})\rangle=
\phi'(\mathcal{L}_{\mathbf{f},\mathbf{g},\mathrm{Symb}}).\end{equation}
Similarly we have
\begin{equation}\label{C}\langle \phi'_0(\mathcal{L}(\mathrm{BF}_{\mathbf{f},\mathbf{g}})),\phi'_0(\eta_\mathbf{f}\mathcal{G})\otimes
\phi'_0(\omega_{\mathbf{g}})\rangle=
\phi'_0(\mathcal{L}_{\mathbf{f},\mathbf{g},\mathrm{Symb}}).\end{equation}
By comparing the equations (\ref{A}), (\ref{B}) and (\ref{C}), we see that $\phi'(\eta_\mathbf{f}\mathcal{G})$ (or equivalently $\phi'_0(\eta_\mathbf{f}\mathcal{G})$) is a finite nonzero multiple of $\eta^\alpha_f$, and that
\begin{equation}\label{need}\langle\phi'_0(\mathcal{L}(\mathrm{BF}_{\mathbf{f},\mathbf{g}})),
\eta^\alpha_f\otimes\phi'_0(\omega_{\mathbf{g}})\rangle=\phi'_0(\mathcal{L}_{\mathbf{f},\mathbf{g},\mathrm{Symb}})
\frac{\Omega_{\Pi_{f,\mathcal{K}}}}{\langle f,f\rangle},
\end{equation}
By the same argument we can also get that
\begin{equation}\label{needed}\langle\mathcal{L}(\mathrm{BF}_{\mathbf{f},\mathbf{g}})|_{f_\alpha}, \eta^\alpha_f\otimes\omega_{\mathbf{g}}\rangle=\mathcal{L}_{\mathbf{f},\mathbf{g},\mathrm{Symb}}|_{f_\alpha}\cdot\frac{\Omega_{\Pi_{f,\mathcal{K}}}}{\langle f,f\rangle}.
\end{equation}
These are precisely the results we need on the explicit reciprocity laws for Beilinson-Flach elements in the previous section to prove the Iwasawa main conjecture.
\paragraph{Remark:}
To summarize, the key in the argument above is the comparison between the specialization of $\eta_\mathbf{f}\mathcal{G}$ to $f_\alpha$ and the $\eta^\alpha_f$. This is made possible thanks to the construction in \cite{BetinaWilliams}.

\subsection{The ordinary case}\label{SubOrdGreenberg}
As mentioned in the introduction, the available literature might not contain a full argument for the proof of theorem \ref{TheoSUintro}. Here, we briefly indicate how the methods of section \ref{SubCrysGreenberg} carry out to this case.
\begin{proof}[Proof of theorem \ref{TheoSUintro}]
The argument is similar as that of section \ref{SubCrysGreenberg} but much easier. We first choose an auxiliary quadratic field $\mathcal{K}$ in the same way and prove the two-variable Iwasawa-Greenberg Rankin-Selberg Main Conjecture over $\mathcal{K}$ involving the $p$-adic $L$-function, by reducing it to the Greenberg type main conjecture (theorem \ref{GBMC}) via explicit reciprocity law for Beilinson-Flach elements and Poitou-Tate exact sequence. The details of this argument are given in \cite[Theorem 3.8]{CastellaWan} (there the result is proved for weight $2$. However the argument goes through in general by replacing the Greenberg's main conjecture there by Theorem \ref{GBMC} here). Then completely as above, we combine it with Kato's result to get full equality for the Iwasawa main conjecture over $\mathbb{Q}$.
\end{proof}
\section{The main theorems}
\subsection{Statements}
\begin{Theo}\label{TheoMain}
Let $p\geq3$ be a prime. Let $f\in S_{k}(\Gamma_{0}(Np^{r}))$ be a normalized eigencuspform with $k\geq2$. Assume that $\rhobar_{f}$ satisfies the following properties.
\begin{enumerate}
\item The $G_{\Q,\Sigma}$-representation $\rhobar_{f}$ is absolutely irreducible.
\item The semisimplification of $\rhobar_{f}|G_{\qp}$ is not equal to $\chibar\oplus\chibar_{\cyc}\chibar$.
\item There exists $\ell\nmid p$ such that $\rhobar_{f}|G_{\ql}$ is a ramified extension
\begin{equation}\nonumber
\suiteexacte{}{}{\mu\chi_{\cyc}^{1-k/2}}{\rhobar|G_{\ql}}{\mu\chi_{\cyc}^{-{k/2}}}
\end{equation}
where $\mu:G_{\ql}\fleche\{\pm1\}$ is an unramified character. If moreover $\rhobar_{f}|G_{\qp}$ is irreducible, then $\mu$ is not trivial.
\end{enumerate}
If $\ell||N$, then the zeta morphism is an isomorphism 
\begin{equation}\nonumber
\triv_{\z(f)_{\Iw}}:\Delta_{\Oiwa}(T(f)_{\Iw})\isocan\Oiwa.
\end{equation}
Equivalently, conjecture \ref{ConjIMC} is true.
\end{Theo}
\paragraph{Remark:}We note that under the hypothesis of the theorem, the determinant of $\rho_{f}$ is an odd power of ${\chi}_{\cyc}$. So the assumption that the semisimplification of $\rhobar_{f}|G_{\qp}$ be different from $\chibar\oplus\chibar$ is automatically satisfied.
\begin{Theo}\label{TheoUnivMain}
Under the same assumption as theorem \ref{TheoMain}, the trivialization morphism
\begin{equation}\nonumber
\triv_{\zs}:\Ds(\Ts)\plonge\Frac(\Lambdaf)
\end{equation}
given by the zeta morphism
\begin{equation}\nonumber
\zs:\Ts(-1)^{+}\fleche\Hun_{\et}(\Z[1/\Sigma],\Ts)
\end{equation}
induces an isomorphism
\begin{equation}\nonumber
\triv_{\zs}:\Ds(\Ts)\isocan\Lambdaf.
\end{equation}
Equivalently, conjecture \ref{ConjUnivWeak} is true.
\end{Theo}
\paragraph{}
We state and prove the following corollaries.
\begin{Cor}\label{CorSelmerMain}
Let $f\in S_{k}(\Gamma_{0}(N))$ be a normalized eigencuspform satisfying all the hypotheses of theorem \ref{TheoMain}. Then $\Sel_{\Q}(f)$ is a finite group if and only if $L(f,k/2)\neq0$.
\end{Cor}
\begin{proof}
If $L(f,k/2)\neq0$, then $\Sel_{\Q}(f)$ is finite by \cite[Theorem]{KatoEuler} so we may and do assume that $L(f,k/2)$ is equal to zero. Let $\gamma$ be a topological generator of $\Gamma_{\Iw}$ and write $\psi$ for the point in $\Spec\Oiwa$ corresponding to $\gamma-1\longmapsto 0$. If the specialization of $\z(f)_{\Iw}$ through $\psi$ is $0$, then $\gamma-1$ belongs to $\carac_{\Oiwa}\left(\Hun_{\et}(\Z[1/p],T(f)_{\Iw})/\Oiwa\cdot\z(f)_{\Iw}\right)$. As conjecture \ref{ConjIMC} is true for $M(f)$, 
\begin{equation}\nonumber
\carac_{\Oiwa}\left(H^{2}_{\et}(\Z[1/p],T(f)_{\Iw})\right)=\carac_{\Oiwa}\left(\Hun_{\et}(\Z[1/p],T(f)_{\Iw})/\Oiwa\cdot\z(f)_{\Iw}\right).
\end{equation}
Hence $\gamma-1$ also belongs to $\carac_{\Oiwa}\left(H^{2}_{\et}(\Z[1/p],T(f)_{\Iw})\right)$ and the $\Ocal$-corank of $\Sel_{\Q}(f)$ is strictly positive. If the specialization of $\z(f)_{\Iw}$ through $\psi$ does not vanish, then our assumption that $L(f,k/2)$ vanishes implies that the image of $\z(f)$ generates a non-zero line in $\Hunf(G_{\Q,\Sigma},V(f)(k/2))$. As before, this means that $\mathrm{corank}\Sel_{\Q}(f)\geq 1$.
\end{proof}
The following corollary contributes to the study of the Birch and Swinnerton-Dyer Conjecture for abelian varieties of $\GL_{2}$-type whose $L$-function does not vanish at 1. 
\begin{Cor}\label{CorEllMain}
Let $A/\Q$ be an abelian variety of $\GL_{2}$-type of conductor $N$. According to \cite{RibetGL2type} and the proof of Serre's conjecture \cite{WintenbergerSerreI,WintenbergerSerreII}, $A$ is modular. Let $f$ be the associated weight two cusp form. Assume that $L(A,1)\neq0$ and that $f$ satisfies the assumptions of theorem \ref{TheoMain}. Then
\begin{equation}\nonumber
v_{p}\left(L(A,1)/\Omega_f\right)=v_{p}\left(\cardinal{\Sha(A/\Q)[p^{\infty}]}\produit{q|N}{}\Tam_{q}(A/\Q)\right).
\end{equation}
Equivalently, the $p$-part of the Birch and Swinnerton-Dyer Conjecture for $A$ holds.
\end{Cor}
\begin{proof}
According to theorem \ref{TheoMain} conjecture \ref{ConjIMC} holds and the zeta morphism
\begin{equation}\nonumber
\triv_{\z(A)_{\Iw}}:\Delta_{\Oiwa}(T_{p}A\tenseur\Oiwa)\plonge\Frac(\Oiwa)
\end{equation}
is an isomorphism
\begin{equation}\nonumber
\triv_{\z(A)_{\Iw}}:\Delta_{\Oiwa}(T_{p}A\tenseur\Oiwa)\isocan\Oiwa.
\end{equation}
is an isomorphism. As $L(A,1)\neq0$, $\z(A)$ does not vanish. If $\psi:\Oiwa\fleche\Ocal$ is the quotient map by the augmentation ideal, the diagram 
\begin{equation}\nonumber
\xymatrix{
\Delta_{\Oiwa}(T_{p}A\tenseur\Oiwa)\ar[rr]^(0.6){\triv_{\z(A)_{\Iw}}}\ar[d]_{-\tenseur_{\Oiwa,\psi}\Ocal}&&\Oiwa\ar[d]^{-\tenseur_{\Oiwa,\psi}\Ocal}\\
\Delta_{\Ocal}(T_{p}A)\ar[rr]^{\triv_{\z(A)}}&&\Ocal
}
\end{equation}
is thus commutative and its horizontal arrows are isomorphisms. Under our hypotheses, $\Sha(A/\Q)[p^{\infty}]$ is finite by \cite[Theorem]{KatoEuler}. Consequently, the isomorphism
\begin{equation}\nonumber
\triv_{\z(A)_{\Iw}}:\Delta_{\Oiwa}(T_{p}A\tenseur\Oiwa)\isocan\Oiwa
\end{equation}
implies the Birch and Swinnerton-Dyer Conjecture at $p$ by \cite[Proposition 1.55]{BurnsFlachMotivic}.
\end{proof}

\subsection{Proofs of theorems \ref{TheoMain} and \ref{TheoUnivMain}}
\begin{Prop}\label{PropFamille}
The trivialization map $\triv_{\zs}:\Ds(\Ts)\plonge\Frac(\Lambdaf)$ of conjecture \ref{ConjUnivWeak} sends $\Ds(\Ts)^{-1}$ inside $\Lambdaf$. 
\end{Prop}
\begin{proof}
By way of contradiction, we assume that $\triv_{\zs}:\Ds(\Ts)\plonge\Frac(\Lambdaf)$ does not send $\Ds(\Ts)^{-1}$ inside $\Lambdaf$. Let $\lambda:\Hsmr\fleche\Ocal$ be an Iwasawa-suitable specialization and let $T_{\lambda,\Iw}$ be $\Ts\tenseur_{\Lambdaf,\lambda}\Oiwa$ (note that $T_{\lambda,\Iw}$ is a $\Oiwa$-module of rank $2d$ where $d$ is the rank of $\Hsmr$ as $\Lambdaf$-module). In order to distinguish them from their counterparts with coefficients in $\Hsmr$, we write $\zsl$ and $\zsl(\lambda)$ for the morphisms of theorem \ref{TheoNakamura} but where all objects are regarded as having coefficients in $\Lambdaf$ (in particular $\lambda$ is then seen as a morphism $\lambda:\Lambdaf\fleche\Ocal$ by restriction). According to proposition \ref{PropDetComSigma}, there is then a commutative diagram
\begin{equation}\nonumber
\xymatrix{
\Delta_{\Sigma}(\Ts)^{-1}\ar[d]_{-\tenseur_{\Lambdaf}\Oiwa}\ar[rr]^{\triv_{\zsl}}&&\frac{x}{y}\Lambdaf\ar[d]^{\lambda}\\
\Delta_{\Sigma}(T_{\lambda,\Iw})^{-1}\ar[rr]^(0.55){\triv_{\zsl(\lambda_{\Iw})}}&&\frac{x'}{y'}\Oiwa
}
\end{equation}
where $\Ts$ is viewed as a $\Lambdaf$-module. Our hypothesis is that $x/y$ does not belong to $\Lambdaf$. As $\Lambdaf$ is a factorial ring, there exists $\pid\in\Spec\Lambdaf$ a height-one prime containing $y$ but not $x$. As $\Lambdaf$ is regular, $\pid$ is principal. Let $y_{0}$ be one of its generator. Let $n$ be a sufficiently large integer. According to propositions \ref{PropDense} and \ref{PropZetaDense}, the set of specialization $\lambda':\Hsmr\fleche\Ocal$ such that one of the specialization in the fiber above the point of $\Lambdaf$ below $\lambda'$ is not \Iwagood or such that $\Spec\Hsmr$ is not étale over $\Lambdaf$ at $\lambda'$ is of large codimension. Hence, it does not contain the set of specialization $\lambda:\Hsmr\fleche\Oiwa$ such that $\lambda(y_{0})$ belongs to $\mgot_{\Oiwa}^{n}$. For $n$ large enough and such a $\lambda$, $\lambda(x)$ does not belong to $\mgot_{\Oiwa}^{n}$. Hence, there exists a specialization $\lambda_{\Iw}:\Hsmr\fleche\Oiwa$ such that $\Spec\Hsmr$ is étale over $\Lambdaf$ at $\lambda_{\Iw}$, such that all specializations in the fiber above the point of $\Lambdaf$ below are \Iwagood and such that  $\triv_{\zsl(\lambda_{\Iw})}\left(\Ds(T_{\lambda,\Iw})^{-1}\right)$ does not belong to $\Oiwa$. 

Because $\Spec\Hsmr$ is étale over $\Lambdaf$ at $\lambda_{\Iw}$, enlarging $\Ocal$ if necessary, we may assume that $T_{\lambda,\Iw}$ is a lattice inside a direct sum \begin{equation}\nonumber
\sommedirecte{r=1}{d}V_{\lambda,r,\Iw}
\end{equation} 
of $G_{\Q,\Sigma}$-representations $V_{\lambda,r,\Iw}$ which are free modules of rank 2 over $\Frac(\Oiwa)$. There is thus a short exact sequence
\begin{equation}\label{EqSuiteCourteLattice}
\suiteexacte{}{}{T_{\lambda,\Iw}}{\sommedirecte{r=1}{d}T_{\lambda,r,\Iw}}{C}
\end{equation}
where $T_{\lambda,r,\Iw}$ are specializations of $\Ts$ attached to specializations $\lambda_{r}:\Hsmr\fleche\Ocal$ and where $C$ a torsion $\Oiwa$-module. The short exact sequence \eqref{EqSuiteCourteLattice} induces a canonical isomorphism
\begin{equation}\label{EqIsomSigmaC}
\produittenseur{r=1}{d}\Delta_{\Sigma}(T_{\lambda,r,\Iw})\isocan\Delta_{\Sigma}(T_{\lambda,\Iw})\tenseur_{\Oiwa}\Delta_{\Sigma}(C).
\end{equation}
According to \cite[Proposition 1.20]{BurnsFlachMotivic} (taking into account the duality between $\RGamma_{c}$ and $\RGamma_{\et}$, see for instance \cite[Appendix]{VenjakobETNC}), $\Delta_{\Sigma}(C)$ is the unit object in the category of graded invertible modules. Hence 
\begin{equation}\nonumber
\triv_{\zsl(\lambda_{\Iw})}\left(\produittenseur{r=1}{d}\Ds(T_{\lambda,r,\Iw})^{-1}\right)\not\subset\Oiwa.
\end{equation}
By construction,
\begin{equation}\nonumber
\triv_{\zsl(\lambda_{\Iw})}=\produittenseur{r=1}{d}\triv_{\zs(\lambda_{r,\Iw})}
\end{equation}
so there must exists an $r$ such that 
\begin{equation}\nonumber
\triv_{\zs(\lambda_{r,\Iw})}\left(\Ds(T_{\lambda,r,\Iw})^{-1}\right)\not\subset\Oiwa.
\end{equation}
By our choice of $\lambda$, $\lambda_{r,\Iw}$ is \Iwagood. Consequently, proposition \ref{PropDivETNC} yields the inclusion 
\begin{equation}\nonumber
\triv_{\zs(\lambda_{r,\Iw})}\left(\Ds(T_{\lambda,r,\Iw})^{-1}\right)\subset\Oiwa.
\end{equation}
This is a contradiction.
\end{proof}
\begin{Prop}\label{PropCrys}
Suppose that $\psi(f):\Hsmr\fleche\Ocal$ is a classical point which is either crystalline and short, or crystalline and short up to the quadratic twist $\omega^{\frac{p-1}{2}}$ (recall $\omega$ is the Teichmuller character), or good ordinary, or good ordinary up to a quadratic twist. Then the trivialization morphism
\begin{equation}\nonumber
\triv_{\z(f)_{\Iw}}:\Delta_{\Oiwa}(T(f)_{\Iw})\plonge\Frac(\Oiwa)
\end{equation}
given by the zeta morphism
\begin{equation}\nonumber
\z(f)_{\Iw}:V(f)_{\Iw}(-1)^{+}\fleche\Hun_{\et}(\Z[1/p],V(f)_{\Iw})
\end{equation}
induces an isomorphism
\begin{equation}\nonumber
\triv_{\z(f)_{\Iw}}:\Delta_{\Oiwa}(T(f)_{\Iw})\isocan\Oiwa.
\end{equation}
Equivalently, conjecture \ref{ConjIMC} is true for $\psi(f)$.
\end{Prop}
\begin{proof}
If $\rho_{f}|G_{\qp}$ is crystalline and short or crystalline and short, this is theorem \ref{TheoCrys}. If $\rho_{f}|G_{\qp}$ is ordinary or ordinary up to a quadratic twist, this is theorem \ref{TheoSUintro} (proved in the end of the last section). The cases involving quadratic twists are proved in the completely same way.
\end{proof}
\begin{Prop}\label{PropCrysRel}
Suppose that there exists a point $x:\Lambdaf\fleche\Ocal$ at which
\begin{equation}\nonumber
\xymatrix{
\Spec\Hsmr\ar[d]\\
\Spec\Lambdaf
}
\end{equation}
is étale and such that all points $\psi:\Hsmr\fleche\Ocal$ above $x$ satisfy conjecture \ref{ConjIMC}. Then the trivialization map
$$\triv_{\zs}:\Ds(\Ts)\plonge\Frac(\Lambdaf)$$
of conjecture \ref{ConjUnivWeak} is an isomorphism
\begin{equation}\nonumber
\Ds(\Ts)\simeq\Lambdaf.
\end{equation}
\end{Prop}
\begin{proof}
According to proposition \ref{PropFamille}, there exists $\alpha\in\Lambdaf$ such that the image of $\Delta_{\Sigma}(\Ts)^{-1}$ through $\triv_{\zs}$ is $\alpha\Lambdaf$. Let $S_{x}$ be the set of points $\lambda:\Hsmr\fleche\Ocal$ above $x$. Then $\Spec\Hsmr\fleche\Spec\Lambdaf$ is étale at $x$ and all $\lambda_{\Iw}$ are \Iwagood specializations. According to proposition \ref{PropDetComSigma}, there thus exists $\alpha\in\Lambdaf$ and a commutative diagram
\begin{equation}\nonumber
\xymatrix{
\Delta_{\Sigma}(\Ts)^{-1}\ar[dd]_{-\tenseur_{\Lambdaf,x_{\Iw}}\Oiwa}\ar[rr]^{{\zsl}}&&\alpha\Lambdaf\ar[dd]^{\lambda}\\
\\
\produittenseur{\lambda\in S_{x}}{}\Delta_{\Sigma}(T_{\lambda,\Iw})^{-1}\ar[rr]^(0.65){{\produittenseur{\lambda\in S_{x}}{}\zs(\lambda_{\Iw})}}&&\frac{x'}{y'}\Oiwa
}
\end{equation}
According to our hypothesis, for all $\lambda\in S_{x}$, $\triv_{\zs(\lambda_{\Iw})}\left(\Delta_{\Sigma}(T_{\lambda,\Iw})^{-1}\right)=\Oiwa$. Hence, the commutative diagram above may be written
\begin{equation}\nonumber
\xymatrix{
\Delta_{\Sigma}(\Ts)^{-1}\ar[dd]_{-\tenseur_{\Lambdaf,x_{\Iw}}\Oiwa}\ar[rr]^{{\zsl}}&&\alpha\Lambdaf\ar[dd]^{\lambda}\\
\\
\produittenseur{\lambda\in S_{x}}{}\Delta_{\Sigma}(T_{\lambda,\Iw})^{-1}\ar[rr]^(0.55){{\produittenseur{\lambda\in S_{x}}{}\zs(\lambda_{\Iw})}}&&\Oiwa
}
\end{equation} 
This implies that $\alpha$ is a unit and thus that the trivialization map
$$\triv_{\zs}:\Ds(\Ts)\plonge\Frac(\Lambdaf)$$
of conjecture \ref{ConjUnivWeak} is an isomorphism
\begin{equation}\nonumber
\Ds(\Ts)^{-1}\simeq\Lambdaf.
\end{equation}
\end{proof}
\begin{proof}[Proof of theorem \ref{TheoUnivMain}]
According to lemma \ref{LemFontaineLaffailleTwist}, there exists a point of $\Spec\Hsmr$ which is crystalline and short, or crystalline and short up to a quadratic twist, or good ordinary, or good ordinary up to a quadratic twist. By proposition \ref{PropCrys}, there then exists a point $x:\Lambdaf\fleche\Ocal$ at which
\begin{equation}\nonumber
\xymatrix{
\Spec\Hsmr\ar[d]\\
\Spec\Lambdaf
}
\end{equation}
is étale and such that all points $\psi:\Hsmr\fleche\Ocal$ above $x$ satisfy conjecture \ref{ConjIMC}. By proposition \ref{PropCrysRel}, the trivialization map
$$\triv_{\zs}:\Ds(\Ts)\plonge\Frac(\Lambdaf)$$
of conjecture \ref{ConjUnivWeak} is an isomorphism
\begin{equation}\nonumber
\triv_{\zs}:\Ds(\Ts)\simeq\Lambdaf
\end{equation}
and conjecture \ref{ConjUnivWeak} is thus true.
\end{proof}

\begin{Prop}\label{PropPointClassique}
Suppose that the trivialization map
$$\triv_{\zs}:\Ds(\Ts)\plonge\Frac(\Lambdaf)$$
of conjecture \ref{ConjUnivWeak} is an isomorphism
\begin{equation}\nonumber
\Ds(\Ts)\simeq\Lambdaf.
\end{equation}
Let $\lambda(f):\Hsm\fleche\Ocal$ be a classical point. Then there is an isomorphism
\begin{equation}\nonumber
\triv_{\z(f)_{\Iw}}:\Ds(T(f)_{\Iw})\simeq\Oiwa.
\end{equation}
Equivalently, conjecture \ref{ConjIMC} is true for $\lambda(f)$.
\end{Prop}
\begin{proof}
We know that there is an inclusion $\triv_{\z(f)_{\Iw}}\left(\Ds(T(f)_{\Iw})^{-1}\right)\subset\Oiwa$. Hence, there exists $\gamma\in\Oiwa$ such that 
\begin{equation}\nonumber
\triv_{\z(f)_{\Iw}}:\Ds(T(f)_{\Iw})\simeq \frac{1}{\gamma}\Oiwa.
\end{equation}
Moreover, there is a commutative diagram
\begin{equation}\nonumber
\xymatrix{\Ds(\Ts)\ar[dd]_{-\tenseur_{\Hsmr,\lambda(f)_{\Iw}}\Oiwa}\ar[rr]^{{\zs}}&&\frac{\alpha}{\beta}\Hsmr\ar[dd]^{\lambda(f)_{\Iw}}\\
\\
\Delta_{\Sigma}(T(f)_{\Iw})\ar[rr]^(0.55){{\zs(f)_{\Iw}}}&&\frac{1}{\gamma}\Oiwa
}
\end{equation}
According to proposition \ref{PropDetComSigma}, we may choose $\beta$ such that $\lambda(f)_{\Iw}(\beta)\neq0$.

Let $Z$ be the set of specialization sending $\beta$ to zero. According to propositions \ref{PropDense} and \ref{PropZetaDense}, there exists a specialization $\lambda:\Hsmr\fleche\Ocal$ which is not in $Z$ above a point $x$ of $\Lambdaf$ such that $\Spec\Hsmr\fleche\Spec\Lambdaf$ is étale over $x$ and such that all specializations $\psi:\Hsmr\fleche\Ocal$ over $x$ are \Iwagood\!\!. By our assumption, there are then elements $\gamma_{\psi}\in\Oiwa$ for all $\psi$ over $x$ and a commutative diagram
\begin{equation}\nonumber
\xymatrix{
\Delta_{\Sigma}(\Ts)\ar[dd]_{-\tenseur_{\Lambdaf,x_{\Iw}}\Oiwa}\ar[rr]^{{\zsl}}&&\Lambdaf\ar[dd]^{x}\\
\\
\produittenseur{\psi\in S_{x}}{}\Delta_{\Sigma}(T_{\psi,\Iw})\ar[rr]^(0.55){{\produittenseur{}{}\zs(\psi_{\Iw})}}&&\left(\produit{\psi\in S_{x}}{}\frac{1}{\gamma_{\psi}}\right)\Oiwa
}
\end{equation} 
This shows that all the $\gamma_{\psi}$ are units and that the diagram above may be written
\begin{equation}\nonumber
\xymatrix{
\Delta_{\Sigma}(\Ts)\ar[dd]_{-\tenseur_{\Lambdaf,x_{\Iw}}\Oiwa}\ar[rr]^{{\zsl}}&&\Lambdaf\ar[dd]^{x}\\
\\
\produittenseur{\psi\in S_{x}}{}\Delta_{\Sigma}(T_{\psi,\Iw})\ar[rr]^(0.55){{\produittenseur{}{}\zs(\psi_{\Iw})}}&&\Oiwa
}
\end{equation} 
Because this holds in particular for $\lambda$ and because $\lambda\notin Z$, we have a commutative diagram
\begin{equation}\nonumber
\xymatrix{\Ds(\Ts)\ar[dd]_{-\tenseur_{\Hsmr,\lambda_{\Iw}}\Oiwa}\ar[rr]^{{\zs}}&&\frac{\alpha}{\beta}\Hsmr\ar[dd]^{\lambda_{\Iw}}\\
\\
\Delta_{\Sigma}(T_{\lambda,\Iw})\ar[rr]^(0.55){{\zs(\lambda_{\Iw})}}&&\Oiwa.
}
\end{equation}
Hence, $\lambda_{\Iw}(\alpha/\beta)$ is a unit while $\lambda(f)_{\Iw}(\alpha/\beta)$ is $1/\gamma$. Taking $\lambda$ close to $\lambda(f)$, we see that $1/\gamma$ is a unit.
\end{proof}
\begin{proof}[Proof of theorem \ref{TheoMain}]
Let $\rhobar$ be the residual representation attached to $M(f)(-\frac{k-2}{2})$. Then by proposition \ref{PropFamille}, the trivialization map $\triv_{\zs}:\Ds(\Ts)\plonge\Frac(\Lambdaf)$ of conjecture \ref{ConjUnivWeak} sends $\Ds(\Ts)^{-1}$ inside $\Lambdaf$. According to lemma \ref{LemFontaineLaffailleTwist}, there exists a point of $\Spec\Hsmr$ which is crystalline and short, or crystalline and short up to a quadratic twist, or good ordinary, or good ordinary up to a quadratic twist. By proposition \ref{PropCrys}, there then exists a point $x:\Lambdaf\fleche\Ocal$ at which
\begin{equation}\nonumber
\xymatrix{
\Spec\Hsmr\ar[d]\\
\Spec\Lambdaf
}
\end{equation}
is étale and such that all points $\psi:\Hsmr\fleche\Ocal$ above $x$ satisfy conjecture \ref{ConjIMC}. By proposition \ref{PropCrysRel}, the trivialization map
$$\triv_{\zs}:\Ds(\Ts)\plonge\Frac(\Lambdaf)$$
of conjecture \ref{ConjUnivWeak} is an isomorphism
\begin{equation}\nonumber
\Ds(\Ts)\simeq\Lambdaf.
\end{equation}
By proposition \ref{PropPointClassique}, the classical point of $\Hsmr$ corresponding to $M(f)(-\frac{k-2}{2})$ then satisfies conjecture \ref{ConjIMC}.
\end{proof}
%
%
\section{Appendix A: completed cohomology and essential vectors for automorphic representations}\label{AppendixCompleted}
\subsection{Generalities}
\subsubsection{Notations}
In this section, $A$ is a complete, local, reduced, noetherian ring with residue field. Let $\ell\nmid p$ be a rational prime. Let $F$ be a finite extension of $\ql$ with ring of integers $\Ocal_{F}$, uniformizing parameter $\varpi$ and residual field $\Fp$ of cardinal $q$. For $n\geq0$, we denote by $\G_{n}$ the group $\GL_{n}(F)$ (so that $\G_{0}$ is the trivial group). We consider $\G_{n}$ as a subgroup of $\G_{n+1}$ through the embedding
\begin{equation}\nonumber
g\mapsto\matrice{g}{0}{0}{1}
\end{equation}
where the two zeroes 0 indicate that the last line and column of $\matrice{g}{0}{0}{1}$ have all their coefficients equal to zero.

Let $K_{n}$ be the maximal compact subgroup $\GL_{n}(\Ocal_{F})$ of $\G_{n}$. If $m\geq0$, let $K_{n}(m)$ be the compact subgroup 
\begin{equation}\nonumber
K_{n}(m)=\left\{\matricetype\in K_{n}|\matricetype\equiv\matrice{1}{0}{0}{1}\modulo\varpi^{m}\right\}.
\end{equation}
For $n\geq2$, we denote by $U_{n},P_{n}$ and $N_{n}$ the following subgroups of $\Gd$ :
\begin{equation}\nonumber
U_{n}=\left\{\matrice{I_{n-1}}{v}{0}{1}|v\in F^{n-1}\right\}
\end{equation}
The mirabolic subgroup $P_{n}$ of $\G_{n}$ is the subgroup $\G_{n-1}U_{n}$ and $N_{n}$ is the subgroup of upper-triangular, unipotent matrices. The abstract group $U_{n}$ is isomorphic to the additive group $F^{n-1}$. If $n=1$, we set all these subgroups equal to the trivial subgroup of $F\croix$.

For $H$ equal to any of the group $\G_{n},K_{n}(m),U_{n},P_{n}$ or $N_{n}$, we denote by $\Rep_{A}H$ the category of smooth $A[H]$-modules.

If $H$ is a closed subgroup of $\G_{n}$ and if $(\s,W)$ is an $A[H]$-module, we denote by $\Ind_{H}^{\G_{n}}\s$ the unnormalized induction, that is to say the $\G_{n}$-representation $(\pi,V)$ where $V$ is the $A$-module of functions $f:\G_{n}\fleche W$ satisfying 
\begin{enumerate}
\item For all $(g,h)\in \Gd\times H$, $f(hg)=\s(h)f(g)$.
\item There is a compact open subgroup $K$ of $\G_{n}$ (which depends possibly on $f$) such that $f(gx)=f(g)$ for all $(g,x)\in\G_{n}\times K$.
\end{enumerate}
and where $\pi(g)f=f(\cdot g)$. The sub-representation $\cInd_{H}^{\G_{n}}$ is the sub-$A$-module of functions which satisfy in addition the property that their support is included in $HK$ for some compact subset $K$ of $\G_{n}$.

Let $\ktilde$ be a Galois extension of $\Fp$ containing all $\ell$-power roots of unity and let $\Atilde$ be $A\tenseur_{W(\Fp)}W(\ktilde)$. We fix $\psi:F\fleche W(\ktilde)\croix$ an additive character of $F$ whose kernel is equal to $\varpi\Ocal_{F}$ and extend it to a character of $N_{n}$ by setting $\psi(n)=\psi(n_{1,2}+n_{2,3}+\cdots+n_{n-1,n})$ for $n\in N_{n}$.

If $H\subset\G_{n}$ is a subgroup, if $(\psi,H)$ is a $W(\ktilde)[H]$-character and if $V$ is an $\Atilde[H]$-module, we write $V(H,\psi)$ for the $\Atilde$-submodule of $V$ generated by elements of the form $hv-\psi(h)v$ and we write $V_{H,\psi}$ for the quotient $V/V(H,\psi)$, which is the largest quotient of $V$ on which $H$ acts through $\psi$. In particular, if $M$ is a Levi subgroup of $\G_{n}$ with unipotent radical $U$ and if $\indicatrice$ is the constant character $1$, then $J_{U}:V\mapsto V_{U,\indicatrice}$ is the Jacquet functor (\cite{BushnellHenniart}).

\subsubsection{Bernstein-Zelevinsky functors}
In \cite{EmertonHelm}, the following functors first introduced in \cite{BernsteinZelevinsky} are extended to our setting.
\begin{enumerate}
\item
\application{\Phi^{+}}{\Rep_{\Atilde}P_{n-1}}{\Rep_{\Atilde}P_{n}}{V}{\cInd_{P_{n-1}U_{n}}^{P_{n}}V}
Here $U_{n}$ acts on $\Phi^{+}(V)$ via $\psi$.
\item
\application{\Phi^{-}}{\Rep_{\Atilde}P_{n}}{\Rep_{\Atilde}P_{n-1}}{V}{V_{U_{n},\psi}}

\item
\application{\Psi^{+}}{\Rep_{\Atilde}\GL_{n-1}}{\Rep_{\Atilde}P_{n}}{V}{V}
Here $U_{n}$ acts on $\Phi^{+}(V)$ trivially.
\item
\application{\Psi^{-}}{\Rep_{\Atilde}P_{n}}{\Rep_{\Atilde}\GL_{n-1}}{V}{V_{U_{n},\indicatrice}}
\item
\application{(-)^{(r)}}{\Rep_{\Atilde}\GL_{n}}{\Rep_{\Atilde}\GL_{n-r}}{V}{\Psi^{-}(\Phi^{-})^{r-1}\Res_{\G_{n}}^{P_{n}}V}
\end{enumerate}
\begin{Prop}\label{PropDerivativeDescends}
If
\begin{equation}\nonumber
\Thetatilde:\Rep_{\Atilde}G\fleche\Rep_{\Atilde}H
\end{equation} is any of the functor $\Phi^{+},\Phi^{-},\Psi^{+},\Psi^{-}$ or $(-)^{(r)}$, then $\Thetatilde$ descends to a functor
\begin{equation}\nonumber
\Theta:\Rep_{A}G\fleche\Rep_{A}H
\end{equation}
in the sense that for all $V\in\Rep_{A}G$
\begin{equation}\nonumber
\Theta(V)\tenseur_{A}\Atilde=\Thetatilde\left(V\tenseur_{A}\Atilde\right).
\end{equation}
\end{Prop}
\begin{proof}
See \cite[Proposition 3.1.4]{EmertonHelm}.
\end{proof}
We use the same notations $\Phi^{+},\Phi^{-},\Psi^{+},\Psi^{-}$ and $(-)^{(r)}$ for the functors on $\Rep_{A}H$ whose existence is asserted in proposition \ref{PropDerivativeDescends}. All these functors are exact and commute with arbitrary base-change of ring of coefficients. The functor $\Phi^{+}$ is left-adjoint to $\Phi^{-}$. The functor $\Psi^{-}$ is left-adjoint to $\Psi^{+}$. The functors $\Psi^{-}\Psi^{+}$ and $\Phi^{-}\Phi^{+}$ are isomorphic to the identity functor. 
\begin{Def}
The \BerZ derivative of order $r$ is the functor $(-)^{(r)}$.
\end{Def}
By construction, the \BerZ derivative of order $n$ of $V\in\Rep_{A}\G_{n}$ is simply an $A$-module and the first \BerZ derivative of a character is a free $A$-module of rank 1. We record the following fact.
\begin{Prop}\label{PropDerivIrr}
The \BerZ derivative of order $n$ is multiplicative with respect to parabolic induction in the following sense. If $V$ and $W$ are respectively an admissible $A[\G_{n}]$-module and an admissible $A[\G_{m}]$-module, then the parabolic induction $\Ind_{P}^{\G_{n+m}}V\tenseur W$ of the representation $V\tenseur W$ or the parabolic subgroup
\begin{equation}\nonumber
P\eqdef\matrice{\G_{n}}{*}{0}{\G_{m}}
\end{equation}
of $\G_{n+m}$ satisfies $\left(\Ind_{P}^{\G_{n+m}}V\tenseur W\right)^{(n+m)}\simeq V^{(n)}\tenseur W^{(m)}$. Suppose that $A$ is a domain and that $V\in\Rep_{A}\G_{n}$ is an absolutely irreducible, admissible $A[\G_{n}]$-module. Then $\rank_{A}V^{(n)}\leq 1$ and $\rank_{A}V^{(n)}=1$ if $V$ is cuspidal. 
\end{Prop}
\begin{proof}
The first assertion is \cite[Lemme 1.10]{VignerasLivre}. After extension of scalars to the field of fractions of $A$, the second assertion is \cite[III.5.10]{VignerasLivre} as slightly extended in \cite[Theorem 3.1.7]{EmertonHelm}.
\end{proof}
By definition 
\begin{align}\nonumber
(\Phi^{+})^{n-1}\Psi^{+}\left(V^{(n)}\right)&=(\Phi^{+})^{n-1}\Psi^{+}\Psi^{-}(\Phi^{-})^{n-1}(V)\\\nonumber
&=(\Phi^{+})^{n-1}(\Phi^{-})^{n-1}(V)
\end{align}
As the functor $\Phi^{+}$ is left adjoint to $\Phi^{-}$, there is thus a natural injective map of $A[P_{n}]$-modules 
\begin{equation}\label{EqLeftAdjoint}
(\Phi^{+})^{n-1}\Psi^{+}\left(V^{(n)}\right)\plonge V.
\end{equation}
The image of this injection is a sub-$A[P_{n}]$-module denoted by $\Icali(V)$ and called the space of Schwartz function. Notice that by construction, the sub-$A[\G_{n}]$-module $A[\G_{n}]\Icali(V)$ generated by $\Icali(V)$ satisfies $\left(A[\G_{n}]\Icali(V)\right)^{(n)}=V^{(n)}$.
\subsubsection{Co-Whittaker modules}
In this subsection, we collect the results we need on the theory of co-Whittaker modules and the integral Bernstein center. All the results are due to D.Helm (\cite{HelmBernstein,HelmWhittaker}) and D.Helm and M.Emerton (\cite{EmertonHelm}).

Let $(\pi,V)$ be an $A[\G_{n}]$-module. 
\begin{Prop}\label{PropAIG}
Suppose that $A$ is a field and that $V$ is smooth, admissible. Let $C$ be the cosocle of $V$, that is to say the largest semisimple quotient of $V$. Then the following assertions are equivalent.
\begin{enumerate}
\item\label{ItemUn} The $A[\G_{n}]$-module $V$ has finite length, $C$ is an absolutely irreducible $A[\G_{n}]$-module and
\begin{equation}\nonumber
V^{(n)}=C^{(n)}\neq0.
\end{equation}
\item\label{ItemDeux} $V^{(n)}\simeq A$ and for all quotient $W$ of $V$, $W^{(n)}=0$ if and only if $W=0$.
\item\label{ItemTrois} $V^{(n)}\simeq A$ and $\Icali(V)$ generates $V$ as $A[\G_{n}]$-module.
\end{enumerate}
\end{Prop}
\begin{proof}
This is a reformulation of comparable results in \cite[Section 6.3]{EmertonHelm}.
\end{proof}
\begin{Def}\label{DefAIG}
When $A$ is a field, we say that a smooth, admissible $A[\G_{n}]$-module $V$ has \essAIG dual if it satisfies one of the equivalent properties of proposition \ref{PropAIG}.
\end{Def}

We return to the general case where $A$ is a complete, local, reduced noetherian $W(\Fp)$-algebra.
\begin{Def}\label{DefCoWhitt}
An $A[\G_{n}]$-module $(\pi,V)$ is co-Whittaker if it satisfies the following conditions.
\begin{enumerate}
\item The $A[\G_{n}]$-module $(\pi,V)$ is smooth, admissible.
\item The top \BerZ derivative $V^{(n)}$ is a free $A$-module of rank 1.
\item If $\pid$ is a prime ideal in $\Spec A$ with residual field $\kg(\pid)$, then $V\tenseur_{A}\kg(\pid)$ has \essAIG dual in the sense of definition \ref{DefAIG}.
\end{enumerate}
\end{Def}
We recall the following unicity result (\cite[Theorem 6.2.1]{EmertonHelm})
\begin{Prop}\label{PropUniqueAIG}
Let $A$ be a reduced, complete, Noetherian, local $p$-torsion free $W(\Fp)$-algebra, with residue field $\Fp$, and let $(V,\rho,A)$ be an $A[G_{F}]$-module. Then there is, up to isomorphism, at most one admissible $A[G]$-module $\rho)$ such that:
\begin{enumerate}
\item The $A[\G_{n}]$-module $\pi(\rho)$ is co-Whittaker and torsion-free as $A$-module.
\item For each minimal prime $\aid\in\Spec A$, the $\kg(\aid)[\G_{n}]$-module attached to $\rho\tenseur_{A}\kg(\aid)$ through the generic Local Langlands Correspondence of \cite{BreuilSchneider} is isomorphic to $\pi(\rho)\tenseur_{A}\kg(\aid)$.\end{enumerate}
\end{Prop}
\begin{proof}
Suppose $\pi(\rho)$ is such an $A[\G_{n}]$-module. Then $\pi(\rho)\tenseur_{A}A/\mgot$ has \essAIG dual so $\pi(\rho)$ satisfies all the conditions of \cite[Theorem 6.2.1]{EmertonHelm}.
 \end{proof}
\subsection{Description of co-Whittaker torsion-modules for $\G_{2}$}
Let $A$ be a complete, local, noetherian domain of residual characteristic zero with fraction field $K$. Let 
\begin{equation}\nonumber
\rho:G_{F}\fleche\GL_{2}(A)
\end{equation}
be a continuous Galois representation such that the attached Weil-Deligne representation is semisimple. Let $\Ncal$ be the monodromy operator acting on the Weil-Deligne representation attached to $\rho\tenseur_{A}K$. Let $\la:A\fleche\Ocal$ be a local morphism of rings and let $\rho_{\la}$ be the $G_{F}$-representation $\lambda\circ\rho$. Let $\Ncal_{\la}$ be the monodromy operator acting on the Weil-Deligne representation attached to $\rho_{\la}\tenseur_{\Ocal}E$. We assume that the monodromy filtration on $\rho_{\la}\tenseur_{\Ocal}E$ is equal to the monodromy filtration on $\rho\tenseur_{A}K$ (this means that $\rho$ is a minimal lift of $\rho_{\la}$ in the sense of \cite[Definition 4.5.9]{EmertonHelm}).
\begin{Prop}\label{PropClassificationGL2}
There exist co-Whittaker $A[\G_{2}]$modules $\pi(\rho)$ and $\pi(\rho_{\la})$ attached to $\rho$ and $\rho_{\la}$ through the Local Langlands Correspondance in the sense of \cite[Theorem 1.2.1]{EmertonHelm}. Moreover, $\pi(\rho)$ and $\pi(\rho_{\la})$ admit explicit descriptions directly parallel to the case of characteristic zero field coefficients. Namely, if $S$ is either $A$ or $\Ocal$ and $\pi$ is either $\pi(\rho)$ or $\pi(\rho_{\la})$, then:
\begin{enumerate}
\item Either $\pi$ is in the principal series: there exist two characters $\chi_{i}:F\croix\fleche S\croix$ such that 
\begin{equation}\nonumber
\pi\simeq\Ind_{B}^{\Gd}(\chi_{1}\tenseur\chi_{2}).
\end{equation}
\item Or it is generically special Steinberg: there exists a character $\mu:F\croix\fleche S\croix$ and a short exact sequence of $S[\Gd]$-modules
\begin{equation}\nonumber
\suiteexacte{}{}{\mu\circ\det}{\Ind_{B}^{\Gd}(\mu\tenseur\mu)}{\pi}.
\end{equation}
\item Or it is supercuspidal: there exists a subgroup $J$ of $\G$ which is compact modulo center, an $S[J]$-module $(\s,\Lambda)$ and a character $(\phi,S)$ of $F\croix$ such that
\begin{equation}\nonumber
\pi\simeq\left(\cInd_{J}^{\Gd}\Lambda\right)\tenseur(\phi\circ\det).
\end{equation}
Moreover, $\pi(\rho)\tenseur_{A,\la}\Ocal\simeq\pi(\rho_{\la})$.
\end{enumerate}
\end{Prop}
\begin{proof}
This follows in much greater generality from the Local Langlands Correspondence in $p$-adic families conjectured in \cite{EmertonHelm} and proven by D.Helm, and Helm-Moss in \cite{HelmBernstein,HelmWhittaker,HelmMoss} (note that under our hypotheses, the automorphic representations attached by the Local Langlands Correspondence in $p$-adic families must be irreducible). However, under our specific setting, the proof is easy and presumably well-known. We briefly recall it.

Let $L$ denote $\Frac(S)$ and let $\pi\tenseur L$ be the $L[\Gd]$-module attached by the Local Langlands Correspondance to $\rho\tenseur_{A}K$ if $S=A$ and to $\rho_{\la}\tenseur_{\Ocal}E$ if $S=\Ocal$. Then $\pi\tenseur L$ falls in exactly one of the three categories of the proposition and our hypothesis on the monodromy action on $\rho_{\la}$ ensures that the category of $\pi\tenseur L$ does not depend on whether $S$ is equal to $A$ or to $\Ocal$. Each of the $S[\Gd]$-module $(\pi,V)$ defined in the statement of the proposition is then by construction a smooth, admissible, $S[\Gd]$-module which is torsion-free as $S$-module and such that $\pi\tenseur_{S}L$ is irreducible and isomorphic to $\pi\tenseur L$. As $V\tenseur_{S}L$ is an irreducible $S[\Gd]$-module and since $\Icali(V\tenseur_{S}L)\simeq (V\tenseur_{S}L)^{(2)}\simeq V^{(2)}\tenseur_{S}L$, $\Icali(V\tenseur_{S}L)$ generates $V$ as $S[\Gd]$-module if $V^{(2)}$ has positive $S$-rank. We compute $V^{(2)}$.

In the principal series case, $V^{(2)}\simeq\chi_{1}^{(1)}\tenseur\chi_{2}^{(1)}\simeq S$ by proposition \ref{PropDerivIrr}. In the generically special Steinberg case, $V^{(2)}$ is equal to $\Psi^{-}\circ\Phi^{-}(V)$. As $\psi$ is a non-trivial character while the action of $N_{2}$ on $\mu\circ\det$ is trivial, $\Phi^{-}(\mu\circ\det)$ vanishes. As $\Phi^{-}$ is exact, $\Phi^{-}(\St\mu)$ is equal to $\Phi^{-}\left(\Ind_{B}^{\Gd}\mu\tenseur\mu\right)$. Consequently, we may use again proposition \ref{PropDerivIrr} again to obtain
\begin{equation}\nonumber
\left(\St\mu\right)^{(2)}\simeq\left(\Ind_{B}^{\G}\mu\tenseur\mu\right)^{(2)}\simeq\mu^{(1)}\tenseur_{S}\mu^{(1)}\simeq S.
\end{equation}
Finally, if $(\pi,V)$ is an $S[\Gd]$-module as in the supercuspidal case, then $V\tenseur_{A}\Fp$ is cuspidal, in particular absolutely irreducible. Then proposition \ref{PropDerivIrr} shows that $(V\tenseur_{S}L)^{(2)}$ and $(V\tenseur_{S}\Fp)^{(2)}$ are vector spaces of dimension 1 over $L$ and $\Fp$ respectively. As $(-)^{(r)}$ commutes with arbitrary change of coefficients, $V^{(2)}$ is of dimension 1 after scalar extension to $L$ and $\Fp$. It then follows from Nakayama's lemma that $V^{(2)}$ is a free $A$-module of rank 1.

In all three cases then, $(\pi,V)$ is a co-Whittaker $S[\Gd]$-module compatible with the Local Langlands Correspondence after extension of scalars to $L$. According to proposition \ref{PropUniqueAIG}, it must be the unique such $S[\Gd]$-module. As $\pi(\rho)\tenseur_{A,\la}\Ocal$ also satisfies all these properties, there is an isomorphism $\pi(\rho)\tenseur_{A,\la}\Ocal\simeq\pi(\rho_{\la})$.

\end{proof}
Let $N(\rho)$ be the Artin conductor of $\rho$. By our hypothesis, this is also the Artin conductor of $\rho_{\la}$.

\begin{Prop}\label{PropEssentialVectors}
Let $m$ be the $\varpi$-adic valuation of $N(\rho)$ and let $U\subset\Gd$ be the compact open subgroup
\begin{equation}\nonumber
U=\left\{\matricetype\in\GL_{2}(\Ocal_{F})|\matricetype\equiv\matrice{*}{*}{0}{1}\modulo \varpi^{m}\right\}.
\end{equation}
The $A$-module $\pi(\rho)_{U}$ and the $\Ocal$-module $\pi(\rho_{\la})_{U}$ are free of rank 1 and the natural map
\begin{equation}\nonumber
\pi(\rho)_{U}\tenseur_{A,\la}\Ocal\fleche\pi(\rho_{\la})_{U}
\end{equation}
is an isomorphism.
\end{Prop}
\begin{proof}
As the representation $\pi(\rho_{\la})$ is co-Whittaker and has coefficients in $\Ocal$, the $\Ocal[\Gd]$-module $\pi(\rho^{*}_{\la})\eqdef\Hom_{\Ocal}(\pi(\rho_{\la}),\Ocal)$ is smooth, admissible and essentially absolutely irreducible and generic. By the theory of essential vectors (\cite{CasselmanAtkin}), $(\pi(\rho^{*}_{\la})\tenseur_{\Ocal}E)^{U}$ is free of dimension 1. As $\pi(\rho^{*}_{\la})^{U}$ is $\Ocal$-torsion free, it is a free $\Ocal$-module of rank 1. By duality, we get that $\pi(\rho_{\la})_{U}$ is $\Ocal$-free of rank 1. As
\begin{equation}\nonumber
\pi(\rho)_{U}\tenseur_{A,\la}\Ocal\simeq(\pi(\rho)\tenseur_{A,\la}\Ocal)_{U}\simeq\pi(\rho_{\la})_{U},
\end{equation}
and $\pi(\rho)_{U}\tenseur_{A}\Frac(A)$ is of rank 1, the $A$-module $\pi(\rho)_{U}$ is $A$-free of rank 1.
\end{proof}

\subsection{Completed cohomology}
We return to the setting of the main text and in particular assume that $\mgot_{\rhobar}$ satisfies assumption \ref{HypTW}.

In the following, the letter $A$ stands for a ring which is equal either to $\Ocal/\varpi^{s}$ for some $s\geq1$, or to $\Ocal$ or to $E$. Let $U$ be an allowable subgroup attached to the maximal ideal $\mgot_{\rhobar}$. We assume that  $U=U_{p}U_{\Sigma}U^{\Sigma}$ with $U_{p}$ a compact open subgroup of $\GL_{2}(\qp)$, $U_{\Sigma}$ a compact open subgroup of $\produit{\ell\in\Sp}{}\GL_{2}(\ql)$ and $U^{\Sigma}$ a maximal compact open subgroup of $\GL_{2}(\Afinis{\Sigma\infty}{\Q})$. Then $H^{i}_{\et}\left(U,A\right)_{\mgot_{\rhobar}}$ is an $\Hsmr(U)\tenseur_{\Ocal}A$-module which we denote by $H^{i}_{\et}(U,A)_{\mgot_{\rhobar}}$ for simplicity.

For $s\geq1$, the completed cohomology $\Hun_{c}\left(U^{p},\Ocal/\varpi^{s}\right)$ with compact support, tame level $U_{\Sigma}$ and finite coefficients is defined as
\begin{equation}\nonumber
\Htildeun_{c}\left(U_{\Sigma},\Ocal/\varpi^{s}\right)_{\mgot_{\rhobar}}\eqdef\liminj{U_{p}}\ \Hun_{c}\left(U^{\Sigma}U_{\Sigma}U_{p},\Ocal/\varpi^{s}\right)_{\mgot_{\rhobar}}.
\end{equation}
The completed cohomology with compact support, tame level $U_{\Sigma}$ and coefficients in $\Ocal$ is defined as
\begin{equation}\nonumber
\Htildeun_{c}\left(U_{\Sigma},\Ocal\right)_{\mgot_{\rhobar}}\eqdef\limproj{s}\ \liminj{U_{p}}\ \Hun_{c}\left(U^{\Sigma}U_{\Sigma}U_{p},\Ocal/\varpi^{s}\right)_{\mgot_{\rhobar}}
\end{equation}
and the completed cohomology with compact support, tame level $U_{\Sigma}$ and coefficients in $E$ is defined as
\begin{equation}\nonumber
\Htildeun_{c}\left(U_{\Sigma},E\right)_{\mgot_{\rhobar}}\eqdef\left(\limproj{s}\ \liminj{U_{p}}\ \Hun_{c}\left(U^{\Sigma}U_{\Sigma}U_{p},\Ocal/\varpi^{s}\right)_{\mgot_{\rhobar}}\right)\tenseur_{\Ocal}E.
\end{equation}
In all three cases, $\Htilde_{c}(U_{\Sigma},A)_{\mgot_{\rhobar}}$ is naturally a faithful $\Hsmr\tenseur_{\Ocal}A$-modules endowed with an action of $G_{\Q,\Sigma}$, with an admissible action of $\GL_{2}(\qp)$ and with an action of $\GL_{2}(\R)$ through the natural action of this latter group on $\C-\R$ (this action factors through the sign of the determinant). 

The completed cohomology with compact support and coefficients in $A$ is then the direct limit on all sufficiently small compact open subgroups 
\begin{equation}\nonumber
\Htildeun_{c}\left(A\right)_{\mgot_{\rhobar}}\eqdef\liminj{K\subset U_{\Sigma}}\ \Htildeun_{c}\left(K,A\right)_{\mgot_{\rhobar}}.
\end{equation}
Then $\Htildeun_{c}\left(A\right)_{\mgot_{\rhobar}}$ is an admissible $\Hr[\produit{\ell\in\Sigma^{p}}{}\GL_{2}(\ql)]$-module and
\begin{equation}\nonumber
\Htildeun_{c}\left(A\right)_{\mgot_{\rhobar}}^{K}=\Htildeun_{c}\left(K,A\right)_{\mgot_{\rhobar}}
\end{equation} 
for all compact open subgroup $K\subset\produit{\ell\in\Sigma^{p}}{}\GL_{2}(\ql)$ contained in $U_{\Sigma}$ (\cite[Theorem 2.2.16]{EmertonInterpolationEigenvalues}). 

For all such $K\subset U_{\Sigma}$, we also consider 
\begin{equation}\nonumber
\Htildeun_{\et}(\Ocal)_{\mgot_{\rhobar}}=\limproj{U'\subset U}{}\Hun_{\et}(U',\Ocal)_{\mgot_{\rhobar}},\ \Htildeun_{\et}(K,\Ocal)_{\mgot_{\rhobar}}=\limproj{U'\subset U_{p}}{}\Hun_{\et}(U'KU^{\Sigma},\Ocal)_{\mgot_{\rhobar}}
\end{equation}
where the inverse limits are respectively over all compact open subgroup of $\GL_{2}(\Afiniq)$ which are maximal at all primes not in $\Sigma$ and other all compact open subgroups of $U_{p}$. The $\Ocal$-modules $\Htildeun_{\et}(\Ocal)$ and $\Htildeun_{c}(\Ocal)$ are then related by Poincarduality
\begin{equation}\label{EqPoincareDualityComp}
\Htildeun_{\et}(\Ocal)_{\mgot_{\rhobar}}\times\Htildeun_{c}(\Ocal)_{\mgot_{\rhobar}}\isom\Ocal(-1).
\end{equation}
In particular, $\Htildeun_{\et}(K,\Ocal)_{\mgot_{\rhobar}}$ is the module of $K$-coinvariants of $\Htildeun_{\et}(\Ocal)_{\mgot_{\rhobar}}$.

The following combines the proof of the Local Langlands Correspondance in $p$-adic families (\cite{HelmBernstein,HelmWhittaker,HelmMoss}) and \cite[Theorem 6.2.13]{EmertonLocalGlobal}.
\begin{Theo}\label{TheoEmerton}[M.Emerton]
Let $\pi_{\ell}(\rho_{\Sigma}|G_{\ql})$ be the continuous $\Ocal$-dual of the co-Whittaker module $\Hr[\GL_{2}(\ql)]$-module $\pitilde_{\ell}(\rho_{\Sigma}|G_{\ql})$ attached to $\rho_{\Sigma}|\G_{\ql}$ through the Local Langlands Correspondence. There is an isomorphism of $\Hr[G_{\Q,\Sigma}\times\GL_{2}(\qp)\times\produit{\ell\in\Sp}{}\GL_{2}(\ql)]$-modules
\begin{equation}\nonumber
\Htildeun_{c}(\Ocal)_{\mgot_{\rhobar}}\simeq\rho_{\Sigma}^{*}(1)\tenseur\pi_{p}(\rho_{\Sigma}|G_{\qp})\hat{\tenseur}\produittenseur{\ell\in\Sp}{}\pi_{\ell}(\rho_{\Sigma}|G_{\ql}).
\end{equation}
\end{Theo}
\section{Appendix B: Iwasawa-Greenberg Main Conjecture}\label{AppendixIwasawaGreenberg}
In this appendix, we prove theorem \ref{TheoGreenberg}. We strongly advise the reader to look at \cite[Introduction]{XinWanRankinSelberg} for a concise outline of the argument proving this Greenberg main conjecture.
\subsection{The Scalar Weight Case: Review}
The idea in proving \cite[Theorem 8.2.1]{CLW} when $f$ has weight two is roughly summarized as follows. We first construct families of Klingen Eisenstein series $E_{\mathrm{Kling}}$ on the unitary group $\Uni(3,1)$ using \cite{EischenWan}. This corresponds to Eisenstein series induced from the Klingen parabolic subgroup of $\Uni(3,1)$. (The motivation for Klingen Eisenstein family is as an automorphic object corresponding to reducible family of Galois representations containing $\rho_\pi$ as a direct summand.) The Hida theory developed in \cite[Sections 2-4]{CLW} (especially the fundamental exact sequence there) for semi-ordinary forms enables us to construct a family of cusp forms, which is congruent to the Klingen Eisenstein family modulo $\Lcal_{\Kcal}^{\Gr}(f)$. Then we proved there is a functional (constructed via Fourier-Jacobi expansion map) acting on the space of families of semi-ordinary forms on $\Uni(3,1)$, which maps $E_{\mathrm{Kling}}$ to an element which is a unit of the coefficient ring $\mathcal{O}^{\ur}[[\Gamma_\mathcal{K}]]$, up to multiplying by an element in $\bar{\mathbb{Q}}_p^\times$. (\textit{i.e.} Proposition 7.11.3 of \emph{loc.cit}. This is the hard part of the whole argument). With this in hand, this functional and the cuspidal family we mentioned above gives a map from the cuspidal Hecke algebra to $\mathcal{O}^{\ur}[[\Gamma_\mathcal{K}]]$ which, modulo $\Lcal_{\Kcal}^{\Gr}(f)$ gives the Hecke eigenvalues acting on the Klingen Eisenstein family (\textit{i.e.} a congruence of Hecke eigenvalues between Eisenstein family and cusp forms). Passing to the Galois side, such congruence enables us to construct enough elements in the Selmer groups from the ``lattice construction'', proving the lower bound of the Selmer group. (See the proof of \cite[Theorem 8.2.1]{CLW} and \cite[Section 9.3]{XinWanRankinSelberg}).

\subsection{Vector Valued Cases}

Now we return to the situation in this paper (\textit{i.e.} general weight). All ingredients are available, except that we need some new idea in the vector valued case to construct and study the corresponding functional on semi-ordinary forms on $\Uni(3,1)$ using Fourier-Jacobi expansion map, so that its value on the Klingen Eisenstein family is an element in $\mathcal{O}^{{\ur}}[[\Gamma_\mathcal{K}]]^\times$, up to multiplying by a non-zero constant in $\bar{\mathbb{Q}}^\times_p$ (\textit{i.e.} primitivity of the Klingen Eisenstein family constructed). In this section we present the whole argument for the entirety of the logic, be brief and refer to the specific part of \cite{CLW} or \cite{XinWanRankinSelberg} for parts which are completely the same as \emph{loc.cit.}, and explain full details for the new ingredients (Archimedean argument involving Ikeda's theory), which are mainly in Subsection \ref{Study}.

\subsection{Unitary groups and Hida Theory for Semi-Ordinary Forms}\label{3}
In \cite{CLW} we developed Hida theory assuming the weight of $f$ is two for ease of presentation (as only scalar weight forms are needed there and this was enough for the application there). See Proposition 2.9.1 and Remark 2.9.2 of \emph{op.cit.}. Here for completeness we briefly develop the Hida theory needed here for general weight. 
\subsubsection{Unitary Groups}
Define $G_n=\mathrm{GU}(n,n)$ for the unitary similitude group for the skew-Hermitian matrix $\begin{pmatrix}&1_n\\-1_n&\end{pmatrix}$ and $\mathrm{U}(n,n)$ for the corresponding unitary group.

Let $\delta\in\mathcal{K}$ be a totally imaginary element such that $-i\delta$ is positive, and $d=\mathrm{Nm}(\delta)$ is a $p$-adic unit. Let $\mathrm{U}(2)=\mathrm{U}(2,0)$ (resp. $\mathrm{GU}(2)=\mathrm{GU}(2,0)$) be the unitary group (resp. unitary similitude group) associated to the skew-Hermitian matrix $\zeta=\begin{pmatrix}\mathfrak{s}\delta&\\&\delta\end{pmatrix}$ for some $\mathfrak{s}\in\mathbb{Z}_+$ prime to $p$. More precisely $\mathrm{GU}(2)$ is the group scheme over $\mathbb{Z}$ defined by: for any $\mathbb{Z}$-algebra $A$,
$$\mathrm{GU}(2)(A)=\{g\in\mathrm{GL}_2(A\otimes_\mathbb{Z}\mathcal{O}_\mathcal{K})|{}^t\!\bar{g}\zeta g=\lambda(g)\zeta,\ \lambda(g)\in A^\times.\}$$
In application below we are going to choose the groups as in \cite[Section 5.4]{CLW}. The map $\lambda: \mathrm{GU}(2)\rightarrow \mathbb{G}_m$, $g\mapsto\lambda(g)$ is called the similitude character and $\mathrm{U}(2)\subseteq \mathrm{GU}(2)$ is the kernel of $\mu$. Let $W$ be the corresponding Hermitian space over $\mathcal{K}$ and fix a lattice $L\subset W$ over $\mathcal{O}_\mathcal{K}$ such that $\mathrm{Tr}_{\mathcal{K}/\mathbb{Q}}\langle L, L\rangle\subset \mathbb{Z}$.
Let $G=\mathrm{GU}(3,1)$ (resp. $\mathrm{U}(3,1)$) be the similarly defined unitary similitude group (resp. unitary group) over $\mathbb{Z}$ associated to the skew-Hermitian matrix $\begin{pmatrix}&&1\\&\zeta&\\-1&&\end{pmatrix}$. Let $P\subseteq G$ be the parabolic subgroup of $\mathrm{GU}(3,1)$ consisting of those matrices in $G$ of the form $\begin{pmatrix}\times&\times&\times&\times\\&\times&\times&\times\\&\times&\times&\times\\&&&\times\end{pmatrix}$. Let $N_P$ be the unipotent radical of $P$. Then
$$M_P{\eqdef}\mathrm{GL}(X_\mathcal{K})\times \mathrm{GU}(2)\hookrightarrow \mathrm{GU}(V),\ (a,g_1)\mapsto \mathrm{diag}(a,g_1,\mu(g_1)\bar{a}^{-1})$$
is the Levi subgroup. Let $G_P{\eqdef}\mathrm{GU}(2)(\subseteq M_P)$ be $\mathrm{diag}(1,g_1,\lambda(g_1))$. Let $\delta_P$ be the modulus character for $P$. We usually use a more convenient character $\delta'$ such that $\delta^{'3}=\delta_P$.\\

\noindent Since $p$ splits as $v_0\bar{v}_0$ in $\mathcal{K}$, $\mathrm{GL}_4(\mathcal{O}_\mathcal{K}\otimes\mathbb{Z}_p)\stackrel{\sim}{\rightarrow}\mathrm{GL}_4
(\mathcal{O}_{\mathcal{K}_{v_0}})\times\mathrm{GL}_4(\mathcal{O}_{\mathcal{K}_{\bar{v}_0}})$. Here $\mathrm{U}(3,1)(\mathbb{Z}_p)\stackrel{\sim}{\rightarrow}\mathrm{GL}_4(\mathcal{O}_{\mathcal{K}_{v_0}})
=\mathrm{GL}_4(\mathbb{Z}_p)$ with the projection onto the first factor. Let $B$ and $N$ be the upper triangular Borel subgroup of $\mathrm{GL}_4$ and its unipotent radical, respectively. Let $K_p=\mathrm{GU}(3,1)(\mathbb{Z}_p)\simeq \mathrm{GL}_4(\mathbb{Z}_p)$, and for any $n\geq 1$ let $K_0^n$  be the subgroup of $K$ consisting of matrices which are upper-triangular modulo $p^n$. Let $K_1^n\subset K_0^n$ be the subgroup of matrices whose diagonal
elements are $1$ modulo $p^n$.

\begin{Def}
A weight $\underline{k}$ \index{$\underline{k}$} is defined to be an $(r+s)$-tuple
$$\underline{k}=(a_1,\cdots, a_r;b_1,\cdots, b_s)\in\mathbb{Z}^{r+s}$$ with $a_1\geq \cdots \geq a_r\geq -b_1\geq \cdots -b_s$.
\end{Def}
We refer to \cite[Section 3.1]{HsiehCM} for the definition of the algebraic representation $L_{\underline{k}}$ of $H$ with the action denoted by $\rho_{\underline{k}}$ and define a model $L^{\underline{k}}$ of the representation $H$ with the highest weight $\underline{k}$ as follows. The underlying space of $L^{\underline{k}}$ is $L_{\underline{k}}$ and the group action is defined by
$$\rho^{\underline{k}}(h)=\rho_{\underline{k}}({}^t\!h^{-1}),h\in H.$$
In \cite[Section 3.1]{HsiehCM} also defined a distinguished functional $l_{\underline{k}}:L_{\underline{k}}(R)\rightarrow R$ for any ring $R$.
We refer to \cite[Section 3.4]{HsiehCM} for the notion of holomorphic automorphic forms of weight $\underline{k}$, and to \cite[Definition 3.2]{HsiehCM} for the automorphic sheaf $\omega_{\underline{k}}$ of weight $\underline{k}$.

\subsubsection{Shimura varieties for Unitary Similitude Groups}
In the following we follow closely \cite[Section 2, 3]{HsiehCM} and refer to some of the details there. We consider the group $\mathrm{GU}(3,1)$. For any open compact subgroup $K=K_pK^p$ of $\mathrm{GU}(3,1)(\mathbb{A}_f)$ whose $p$-component is $K_p=\mathrm{GU}(3,1)(\mathbb{Z}_p)$ and whose prime to $p$ component is $K^p$, we refer to \cite[Section 2.1]{HsiehCM} for the definition and arithmetic models of the associated Shimura variety, which we denote as $S_G(K)_{/\mathcal{O}_{\mathcal{K},(v_0)}}$. The scheme $S_G(K)$ represents the following functor $\underline{A}$: for any $\mathcal{O}_{\mathcal{K},(v_0)}$-algebra $R$, $\underline{A}(R)=\{(A,\bar{\lambda},\iota,\bar{\eta}^p)\}$ where $A$ is an abelian scheme over $R$ with CM by $\mathcal{O}_\mathcal{K}$ given by $\iota$, $\bar{\lambda}$ is an orbit of prime-to-$p$ polarizations and $\bar{\eta}^p$ is an orbit of prime-to-$p$ level structures. We denote $\bar{S}_G(K)$ a smooth toroidal compactification and $S^*_G(K)$ the minimal compactification. We refer to \cite[Section 2.7]{HsiehCM} for details. The boundary components of $S^*_G(K)$ are in one-to-one correspondence with the set of cusp labels defined below. For $K=K_pK^p$ as above we define the set of cusp labels to be:
$$C(K){\eqdef}(\mathrm{GL}(X_\mathcal{K})\times G_P(\mathbb{A}_f))N_P(\mathbb{A}_f)\backslash G(\mathbb{A}_f)/K.$$
This is a finite set. We denote by $[g]$ the class represented by $g\in G(\mathbb{A}_f)$. For each such $g$ whose $p$-component is $1$ we define $K_P^g=G_P(\mathbb{A}_f)\cap gKg^{-1}$ and denote $S_{[g]}{\eqdef}S_{G_P}(K_P^g)$ the corresponding Shimura variety for the group $G_P$ with level group $K_P^g$. By strong approximation we can choose a set $\underline{C}(K)$ of representatives of ${C}(K)$ consisting of elements $g=pk^0$ for $p\in P(\mathbb{A}_f^{\Sigma})$ and $k^0\in K^0$ for $K^0$ the maximal compact subgroup of $G(\mathbb{A}_f)$ defined in \cite[Section 1.10]{HsiehCM}.
\subsubsection{Igusa varieties and $p$-adic automorphic forms}\label{2.4}
Now we recall briefly the notion of Igusa varieties in \cite[Section 2.3]{HsiehCM}. Let $M$ be the standard lattice of $V$ and $M_p=M\otimes_{\mathbb{Z}}\mathbb{Z}_p$. Let $\mathrm{Pol}_p=\{N^{-1}, N^0\}$ be a polarization of $M_p$. Recall this means that if $N^{-1}$ and $N^0$ are maximal isotropic $\mathcal{O}_\mathcal{K}\otimes\mathbb{Z}_p$-submodules in $M_p$, that they are dual to each other with respect to the Hermitian metric on $V$, and also that:
$$\mathrm{rank}_{\mathbb{Z}_p} N_{v_0}^{-1}=\mathrm{rank}_{\mathbb{Z}_p}N_{\bar{v}_o}^0=3,\ \mathrm{rank}_{\mathbb{Z}_p}N_{\bar{v}_0}^{-1}=\mathrm{rank}_{\mathbb{Z}_p}N_{v_0}^0=1.$$

We mainly follow \cite[Section 2.3]{HsiehCM} in this subsection. The Igusa variety of level $p^n$ and tame level $K$ is the scheme over $\mathcal{O}_{\mathcal{K},(v_0)}$ representing the quadruple $\underline{A}(R)=\{(A,\bar{\lambda},\iota,\bar{\eta}^p)\}$ for the Shimura variety of $\mathrm{GU}(3,1)$ as above, together with an injection of group schemes
$$j:\mu_{p^n}\otimes_\mathbb{Z}N^0\hookrightarrow A[p^n]$$
over $R$ which is compatible with the $\mathcal{O}_\mathcal{K}$-action on both sides.
Note that the existence of $j$ implies that $A$ must be ordinary along the special fiber. There is also a theory of Igusa varieties over $\bar{S}_G(K)$. As in \emph{loc.cit.} let $\bar{H}_{p-1}\in H^0(S_G(K)_{/\bar{\mathbb{F}}},\mathrm{det}(\underline{\omega})^{p-1})$ be the Hasse invariant. Over the minimal compactification some power (say the $t$-th) of the Hasse invariant can be lifted to $\mathcal{O}_{v_0}$. We denote such a lift by $E$. By the Koecher principle we can regard $E$ as in $H^0(\bar{S}_G(K),\mathrm{det}(\underline{\omega}^{t(p-1)}))$. Let $\mathcal{O}_m{\eqdef}\mathcal{O}_{\mathcal{K},v_0}/p^m\mathcal{O}_{\mathcal{K},v_0}$. Set $T_{0,m}{\eqdef}\bar{S}_G(K)[1/E]_{/\mathcal{O}_m}$. For any positive integer $n$ define $T_{n,m}{\eqdef}I_G(K^n)_{/\mathcal{O}_m}$ and $T_{\infty,m}=\varprojlim_n T_{n,m}$. Then $T_{\infty,m}$ is a Galois cover over $T_{0,m}$ with Galois group $\mathbf{H}\simeq \mathrm{GL}_3(\mathbb{Z}_p)\times \mathrm{GL}_1(\mathbb{Z}_p)$. Let $\mathbf{N}\subset \mathbf{H}$ be the upper triangular unipotent radical. Define:
$$V_{n,m}=H^0(T_{n,m},\mathcal{O}_{T_{n,m}}).$$
Let $V_{\infty,m}=\varinjlim_n V_{n,m}$ and $V_{\infty,\infty}=\varprojlim_m V_{\infty,m}$ be the space of $p$-adic automorphic forms on $\mathrm{GU}(3,1)$ with level group $K$. We also define $W_{n,m}=V_{n,m}^\mathbf{N}$, $W_{\infty,m}=V_{\infty,m}^\mathbf{N}$ and $\mathcal{W}=\varinjlim_m\varinjlim_n W_{n,m}$. We define $V_{n,m}^0$, etc, to be the cuspidal part of the corresponding spaces.\\

\noindent We can make similar definitions for the definite unitary similitude groups $G_P$ as well and define $V_{n,m,P}$,$V_{\infty,m,P}$, $V_{\infty,\infty,P}$, $V_{n,m,P}^{\mathbf{N}}$, $\mathcal{W}_P$, etc.\\

\noindent Let $K_0^n$ and $K_1^n$ be the subgroup of $\mathbf{H}$ consisting of matrices which are in $B_3\times {}^t\!B_1$ or $N_3\times {}^t\!N_1$ modulo $p^n$. (These notations are already used for level groups of automorphic forms. The reason for using the same notation here is that automorphic forms with level group $K_\bullet^n$ are $p$-adic automorphic forms of level group $K_\bullet^n$). We sometimes denote $I_G(K_1^n)=I_G(K^n)^{K_1^n}$ and $I_G(K_0^n)=I_G(K^n)^{K_0^n}$. We define
$$M_{\underline{k}}(K^n_\bullet, R){\eqdef}H^0(I_G(K^n_\bullet)_{/R}, \omega_{\underline{k}}).$$

\noindent We can define Igusa varieties for $G_P$ as well. For $\bullet=0,1$ we let $K_{P,\bullet}^{g,n}{\eqdef} gK^n_\bullet g^{-1}\cap G_P(\mathbb{A}_f)$ and let $I_{[g]}(K_\bullet^n){\eqdef}I_{G_P}(K_{P,\bullet}^{g,n})$ be the corresponding Igusa variety over $S_{[g]}$. We denote the coordinate ring of $I_{[g]}(K_1^n)$ over $\mathcal{O}_m$ by $A_{[g],m}^n$ . Let $A_{[g],m}^{\infty}=\varinjlim_n A_{[g],m}^n$ and let $\hat{A}_{[g]}^\infty$ be the $p$-adic completion of $A^\infty_{[g],m}$. This is the space of $p$-adic automorphic forms for the group $\mathrm{GU}(2,0)$ of level group $gKg^{-1}\cap G_P(\mathbb{A}_f)$.\\

\noindent\underline{For Unitary Groups}\\
\noindent Assume the tame level group $K$ is neat. For any $c$ an element in $\mathbb{Q}_+\backslash \mathbb{A}_{\mathbb{Q},f}^\times /\mu(K)$, we refer to \cite[2.5]{HsiehCM} for the notion of $c$-Igusa schemes $I_{\mathrm{U}(2)}^0(K,c)$ for the unitary groups $\mathrm{U}(2,0)$ (not the similitude group). It parameterizes quintuples $(A,\lambda,\iota,\bar{\eta}^{(p)},j)_{/S}$ similar to the Igusa schemes for unitary similitude groups but requires $\lambda$ to be a prime to $p$ $c$-polarization of $A$ such that $(A,\bar{\lambda},\iota,\bar{\eta}^{(p)},j)$ is a quintuple as in the definition of Shimura varieties for $\mathrm{GU}(2)$. Let $g_c$ be such that $\mu(g_c)\in\mathbb{A}_\mathbb{Q}^\times$ is in the class of $c$. Let ${}^c\!K=g_cKg_c^{-1}\cap U(2)(\mathbb{A}_{\mathbb{Q},f})$. Then the space $I_{\mathrm{U}(2)}^0(K,c)$ is isomorphic to the space of forms on $I_{\mathrm{U}(2)}^0({}^c\!K,1)$ (see \emph{loc.cit.}).\\

\noindent\underline{Fourier-Jacobi Expansions}\\
Define $N_H^1{\eqdef}\{\begin{pmatrix}1& 0\\ * & 1_2\end{pmatrix}\}\times\{1\}\subset H$. For an automorphic form or $p$-adic automorphic form $F$ on $\mathrm{GU}(3,1)$ we refer to \cite[Section 2.8]{EischenWan} for the notion of analytic Fourier-Jacobi expansions
$$FJ_P(g,f)=a_0(g,f)+\sum_\beta a_\beta(y,g,f)q^\beta$$
at $g\in\mathrm{GU}(3,1)(\mathbb{A}_\mathbb{Q})$ for $a_\beta(-,g,f):\mathbb{C}^2\rightarrow L_{\underline{k}}(\mathbb{C})$ being theta functions with complex multiplication.
Also there is an algebraic Fourier-Jacobi expansion
$$FJ_{[g]}^h(f)_{N_H^1}=\sum_\beta a_{[g]}^h(\beta,f)q^\beta,$$
at a $p$-adic cusp $([g],h)$, and $a_{[g]}^h(\beta,f)\in L_{\underline{k}}(A_{[g]}^\infty)_{N_H^1}\otimes_{A_{[g]}}H^0(\mathcal{Z}_{[g]}^\circ,\mathcal{L}(\beta))$ (see \cite[(3.9)]{HsiehCM}. Note that the subscript $N^1_H$ is important to take it out of the $H^0$). We define the Siegel operator to be taking the $0$-th Fourier-Jacobi coefficient as in \emph{loc.cit.}.
Over $\mathbb{C}$ the analytic Fourier-Jacobi expansion for a holomorphic automorphic form $f$ is given by:
$$FJ_\beta(f,g)=a_\beta(y,g,f)=\int_{\mathbb{Q}\backslash \mathbb{A}} f(\begin{pmatrix}1&&n\\&1_2&\\&&1\end{pmatrix}g)e_\mathbb{A}(-\beta n)dn.$$

\subsubsection{Semi-Ordinary Forms}\label{semi-ordinary}

In this subsection we develop a theory for families of ``semi-ordinary'' forms over a two dimensional weight space (the whole weight space for $\mathrm{U}(3,1)$ is three dimensional). The idea goes back to the work of Hida (also Tilouine-Urban for $\mathrm{GSp}(4)$) who they defined the concept of being ordinary with respect to general parabolic subgroups (the usual definition of ordinary is with respect to the Borel subgroup), except that we are working with coherent cohomology while Hida and Tilouine-Urban use cohomology of arithmetic groups. In our case it means being ordinary with respect to the parabolic subgroup of $\mathrm{GL}_4(\mathbb{Q}_p)\cong\mathrm{U}(3,1)(\mathbb{Q}_p)$ consisting of matrices of the form $\begin{pmatrix}*&*&*&*\\ *&*&*&*\\&&*&*\\&&&*\end{pmatrix}$. The crucial point is that our families are defined over the two dimensional \emph{Iwasawa algebra}, as Hida theory for ordinary forms instead of Coleman-Mazur theory for finite slope forms, which is over some affinoid domain. Our argument here will sometimes be an adaption of the argument in the ordinary case in \cite{HsiehCM} and we will sometimes be brief and refer to \emph{loc.cit.} for some computations so as not to introduce too many notations.\\

\noindent We always use the identification $\mathrm{U}(3,1)(\mathbb{Q}_p)\simeq \mathrm{\mathrm{GL}}_4(\mathbb{Q}_p)$. Define $\alpha_i=\mathrm{diag}(1_{4-i}, p\cdot 1_i)$. We let $\alpha=\begin{pmatrix}1&&&\\&1&&\\&&p&\\&&&p^2\end{pmatrix}$ and refer to \cite[3.7, 3.8]{HsiehCM} for the notion of Hida's $U_{\alpha}$ and $U_{\alpha_i}$ operators associated to $\alpha$ or $\alpha_i$. We define $e_\alpha=\lim_{n\rightarrow \infty}U_\alpha^{n!}$. (That this is well-defined follows as in \cite[Section 4.3]{HsiehCM}.) We are going to study forms and families invariant under $e_\alpha$ and call them ``semi-ordinary'' forms. Suppose $\pi$ is an irreducible automorphic representation on $\mathrm{U}(3,1)$ with weight $\underline{k}$ and suppose that $\pi_p$ is an unramified principal series representation. If we write $\kappa_1=b_1$ and $\kappa_{i}=-a_{5-i}+5-i$ for $2\leq i \leq 4$, then there is a semi-ordinary vector in $\pi$ if and only if we can re-order the Satake parameters as $\lambda_1,\lambda_2,\lambda_3,\lambda_4$ such that
$$\mathrm{val}_p(\lambda_3)=\kappa_3-\frac{3}{2},\ \mathrm{val}_p(\lambda_4)=\kappa_4-\frac{3}{2}.$$
\\

\subsubsection{Control Theorems}
\noindent We define $K_0(p,p^n)=\prod_{\ell\not=p}K_\ell\times K_0(p,p^n)_p$, for $K_0(p,p^n)_p$ consisting of matrices which are of the form $\begin{pmatrix}*&*&*&*\\&*&*&*\\&&*&*\\&&&*\end{pmatrix}$ modulo $p$ and are of the form $\begin{pmatrix}*&*&*&*\\ *&*&*&*\\&&*&*\\&&&*\end{pmatrix}$ modulo $p^n$. We are going to prove some control theorems for the level group $K_0(p,p^n)$. We also define a $\mathrm{GL}_2$ level group $K'_0(p)\subset\mathrm{GL}_2(\mathbb{Z}_p)$ to be the set of matrices congruent to $\begin{pmatrix}*&*\\0&*\end{pmatrix}$ modulo $p$. Let $N'$ be the set of matrices $\begin{pmatrix}1&\mathbb{Z}_p\\&1\end{pmatrix}$. For the definition of the automorphic sheaves $\omega_{\underline{k}}$ of weight $\underline{k}$ we refer to \cite[section 3.2]{HsiehCM}. There also defined a subsheaf $\omega_{\underline{k}}^\flat$ in Section 4.1 of \emph{loc.cit.} as follows. Let $\mathcal{D}=\bar{S}_G(K)-S_G(K)$ be the boundary of the toroidal compactification and $\underline{\omega}$ the pullback to identity of the relative differential of the Raynaud extension of the universal Abelian variety. Let $\underline{k}''=(a_1-a_3,a_2-a_3)$. Let $\mathcal{B}$ be the abelian part of the Mumford family of the boundary. Its relative differential is identified with a subsheaf of $\underline{\omega}|_{\mathcal{D}}$. The $\omega_{\underline{k}}^\flat\subset \omega_{\underline{k}}$ is defined to be $\{s\in\omega_{\underline{k}}, s|_{\mathcal{D}}\in\mathscr{F}_\mathcal{D}\}$ for $\mathscr{F}_\mathcal{D}{\eqdef}\det (\underline{\omega}|_{\mathcal{D}})^{a_3}\otimes\underline{\omega}_{\mathcal{B}}^{\underline{k}''}$, where the last term means the automorphic sheaf of weight $\underline{k}''$ for $\mathrm{GU}(2,0)$ (see \cite[Section 4.1]{HsiehCM}).\\

\noindent\underline{Weight Space}\\
\noindent Let $H=\mathrm{GL}_3\times \mathrm{GL}_1$ and $T$ be the diagonal torus. Then $\mathbf{H}=H(\mathbb{Z}_p)$. We let $\Lambda_{3,1}=\Lambda$ be the completed group algebra $\mathbb{Z}_p[[T(1+p\mathbb{Z}_p)]]$. This is a formal power series ring with four variables. There is an action of ${T}(\mathbb{Z}_p)$ given by the action on the $j:\mu_{p^n}\otimes_\mathbb{Z} N^0\hookrightarrow  A[p^n]$. (see \cite[3.4]{HsiehCM}) This gives the space of $p$-adic modular forms a structure of $\Lambda$-algebra. A $\bar{\mathbb{Q}}_p$-point $\phi$ of $\mathrm{Spec}\Lambda$ is called arithmetic if it is determined by a character $[\underline{k}].[\zeta]$ of $T(1+p\mathbb{Z}_p)$ where $\underline{k}$ is a weight and $\zeta=(\zeta_1,\zeta_2,\zeta_3;\zeta_4)$ for $\zeta_i\in \mu_{p^\infty}$. Here $[\underline{k}]$ is the character of $T(1+\mathbb{Z}_p)$ by ${[\underline{k}]}(t_1,t_2,t_3,t_4)=(t_1^{a_1}t_2^{a_2}t_3^{a_3}t_4^{-b_1})$ and $[\zeta]$ is the finite order character given by mapping $(1+p\mathbb{Z}_p)$ to $\zeta_i$ at the corresponding entry $t_i$ of $T(\mathbb{Z}_p)$. We often write this point $\underline{k}_\zeta$. We also define $\omega^{[\underline{k}]}$ a character of the torsion part of $T(\mathbb{Z}_p)$ (isomorphic to $(\mathbb{F}_p^\times)^4$) given by $\omega^{[\underline{k}]}(t_1,t_2,t_3,t_4)=\omega(t_1^{a_1}t_2^{a_2}t_3^{a_3}t_4^{-b_1})$.
\begin{Def}
We fix $\underline{k}'=(a_1,a_2)$ and $\rho=L_{\underline{k}'}$. Let $\mathcal{X}_\rho$ be the set of arithmetic points $\phi\in \mathrm{Spec}{\Lambda}_{3,1}$ corresponding to the weight $(a_1,a_2,a_3;b_1)$ such that $a_1\geq a_2\geq a_3\geq -b_1+4$. (The $\zeta$-part is trivial). Let $\mathrm{Spec}\tilde{\Lambda}=\mathrm{Spec}\tilde{\Lambda}_{(a_1,a_2)}$ be the Zariski closure of $\mathcal{X}_\rho$.
\end{Def}
We define for $q=0,\flat$
$$V_{\underline{k}}^q(K_0(p,p^n),\mathcal{O}_m){\eqdef}H^0(T_{n,m},\omega_{\underline{k}}^q)^{K_0(p,p^n)\cap \mathrm{GL}_3(\mathbb{Z}_p)\times\mathrm{GL}_1(\mathbb{Z}_p)}.$$

\noindent As in \cite[3.3]{HsiehCM} we have a canonical isomorphism given by taking the ``$p$-adic avartar''
$$H^0(T_{n,m},\omega_{\underline{k}})\simeq V_{n,m}\otimes L_{\underline{k}}, f\mapsto \hat{f}$$
and $\beta_{\underline{k}}: V_{\underline{k}}(K_1^n,\mathcal{O}_m)\rightarrow V_{n,m}^\mathbf{N}$
by $f\mapsto \beta_{\underline{k}}(f){\eqdef}l_{\underline{k}}(\hat{f})$.
The following lemma is \cite[lemma 4.2]{HsiehCM}.
\begin{Lem}\label{3.1}
Let $q\in\{0,\flat\}$ and let $V_{\underline{k}}^q(K_0(p,p^n),\mathcal{O}_m){\eqdef}H^0(T_{n,m},\omega_{\underline{k}}^q)^{K_0(p,p^n)}$. Then we have
$$H^0(I_G(K_0(p,p^n))[1/E],\omega_{\underline{k}}^q)\otimes\mathcal{O}_m=V_{\underline{k}}^q(K_0(p,p^n),\mathcal{O}_m).$$
\end{Lem}
We record a contraction property for the operator $U_\alpha$.
\begin{Lem}
If $n>1$, then we have
$$U_\alpha\cdot V_{\underline{k}}(K_0(p,p^n), \mathcal{O}_m)\subset V_{\underline{k}}(K_0(p,p^{n-1}),\mathcal{O}_m).$$
\end{Lem}
The proof is the same as \cite[Proposition 4.4]{HsiehCM}.
The following proposition follows from the contraction property for $e_\alpha$:
\begin{Prop}\label{3.2}
$$e_\alpha V_{\underline{k}}^q(K_0(p,p^n),\mathcal{O}_m)=e_\alpha V_{\underline{k}}(K_0(p),\mathcal{O}_m).$$
\end{Prop}
The following lemma tells us that to study semi-ordinary forms one only needs to look at the sheaf $\omega_{\underline{k}}^\flat$.
\begin{Lem}
Let $n\geq m>0$, then
$$e_\alpha.V_{\underline{k}}^\flat(K_0(p,p^n),\mathcal{O}_m)=e_\alpha\cdot V_{\underline{k}}(K_0(p,p^n),\mathcal{O}_m).$$
\end{Lem}
\begin{proof}
Same as \cite[lemma 4.10]{HsiehCM}.
\end{proof}
\noindent Similar to the $\beta_{\underline{k}}$ we define a more general $\beta_{\underline{k},\rho}$ as follows: Let $\rho$ be the algebraic representation $L_\rho=L_{\underline{k}'}$ of $\mathrm{\mathrm{GL}}_2$ with lowest weight $-\underline{k}'=(-a_1,-a_2)$. We identify $L_{\underline{k}}$ with the algebraically induced representation $\mathrm{Ind}_{\mathrm{\mathrm{GL}}_2\times \mathrm{\mathrm{GL}}_1\times \mathrm{\mathrm{GL}}_1}^{\mathrm{\mathrm{GL}}_3\times \mathrm{\mathrm{GL}}_1}\rho\otimes\chi_{a_3}\otimes\chi_{b_1}$ ($\chi_a$ means the algebraic character defined by taking the ($-a$)-th power). We define the functional $l_{\underline{k},\rho}: L_{\underline{k}}\rightarrow L_{\underline{k}'}$ by evaluating at identity (similar to the definition of $l_{\underline{k}}$). We define $\beta_{\underline{k},\rho}$ similar to $\beta_{\underline{k}}$ but replacing $l_{\underline{k}}$ by $l_{\underline{k},\rho}$.
\begin{Prop}
If $n\geq m>0$, then the morphism
$$\beta_{\underline{k},\rho}:V_{\underline{k}}(K_1(p^n),\mathcal{O}_m)\rightarrow (V_{n,m}\otimes L_\rho)^{N'}$$
is $U_\alpha$-equivariant (here $N'\subset\mathrm{GL}_2(\mathbb{Z}_p)$ is embedded into $\mathrm{GL}_3(\mathrm{Z}_p)\times\mathrm{GL}_1(\mathbb{Z}_p)$ as $\mathrm{diag}(N',1,1)$), and there is a Hecke-equivariant homomorphism $s_{\underline{k},\rho}: (V_{n,m}\otimes L_\rho)^{N'}\rightarrow V_{\underline{k}}(K_1(p^n),\mathcal{O}_m)$ such that $\beta_{\underline{k},\rho}\circ s_{\underline{k},\rho}=U_\alpha^m$ and $s_{\underline{k},\rho}\circ\beta_{\underline{k},\rho}=U_\alpha^m$. So the kernel and the cokernel of $\beta_{\underline{k},\rho}$ are annihilated by $U_\alpha^m$.
\end{Prop}
\begin{proof}
We follow \cite[Proposition 4.7]{HsiehCM}. Our $s_{\underline{k},\rho}$ is defined as follows: for $(\underline{A},\bar{j})$ over a $\mathcal{O}_m$-algebra $R$,
$$s_{\underline{k},\rho}(\alpha^m)f(\underline{A},\bar{j}){\eqdef}\sum_{v_{\chi'}\in \rho\otimes\chi_{a_3}\otimes\chi_{b_1}}\sum_u\frac{1}{\chi_{r,1}(\alpha^m)}\cdot \mathrm{Tr}_{R_0^{\alpha^mu}/R}(f(\underline{A}_{\alpha^m u}.j_{\alpha^m u}))\rho_{\underline{k}}(u)v_{\chi'}.$$
Here the character $\chi_{r,1}$ is defined by
$$\chi_{r,1}(\mathrm{diag}(a_1,a_2,a_3;d)){\eqdef}(a_1a_2a_3)^{-1}d.$$
The $v_{\chi'}$'s form a basis of the representation $\rho\otimes\chi_{a_3}\otimes\chi_{b_1}$ which are eigenvectors for the diagonal torus action with eigenvalues $\chi'$'s  (the eigenvalues appear with multiplicity one so we use the subscript $\chi'$ to denote the corresponding vector). The $u$ runs over a set of representatives of $$\alpha^{-m}N_H(\mathbb{Z}_p)\alpha^m\cap N_H(\mathbb{Z}_p)\backslash N_H(\mathbb{Z}_p).$$ The $(\underline{A}_{\alpha u}, j_{\alpha u})$ is a certain pair with $\underline{A}_{\alpha u}$ an abelian variety admitting an isogeny to $\underline{A}$ of type $\alpha$ (see \cite[3.7.1]{HsiehCM} for details) and $R_0^{\alpha u}/R$ being the coordinate ring for $(\underline{A}_{\alpha u}, j_{\alpha u})$ (see 3.8.1 of \emph{loc.cit.}). Note that the twisted action of
$$\tilde{\rho}_{\underline{k}}(\alpha^{-1})v_{\chi'}{\eqdef}p^{-\langle\mu,\underline{k}+\chi'\rangle}v_{\chi'}$$
satisfies $\tilde{\rho}_{\underline{k}}(\alpha^{-1})v_{\chi'}=1$ for all the $\chi'$ above. Write $\chi$ for $\chi_{a_3}\boxtimes \chi_{b_1}$. Note also that for any eigenvector $v_{\chi'}\in\mathrm{Ind}_{\mathrm{\mathrm{GL}}_2\times\mathrm{\mathrm{GL}}_1\times\mathrm{\mathrm{GL}}_1}
^{\mathrm{\mathrm{GL}}_3\times\mathrm{\mathrm{GL}}_1}\rho\otimes\chi$ for the torus action such that $v_{\chi'}\not\in \rho\otimes \chi$, and $\mu\in X_*(T)$ (the co-character group) with $\mu(p)=\alpha$, we have $\langle\mu,\underline{k}+\chi'\rangle <0$. By the definition of $U_\alpha^m=U_{\alpha^m}$, if $f=\sum_\chi g_\chi\otimes v_\chi$, then
\begin{align}\nonumber
U_{\alpha^m}\cdot f(\underline{A},j)=&\sum_{v_{\chi'}\in\rho\otimes\chi}s_{\underline{k},\rho}(\alpha^m)g_{\chi'}(\underline{A},j)\\\nonumber
&+\sum_{v_{\chi'}
\not\in\rho\otimes\chi}p^{-\langle m\mu, \underline{k}\chi'\rangle}\frac{1}{\chi_{r,1}(\alpha^m)}\mathrm{Tr}_{R_0^{\alpha^mu}/R}(g_{\chi'}(\underline{A}_
{\alpha^mu}, j))\otimes\rho_{\underline{k}}(u)v_{\chi'}.
\end{align}
For the notation $R^{\alpha^mu}_0$ see \cite[3.8.1]{HsiehCM} for an explanation. So $\beta_{\underline{k},\rho}\circ s_{\underline{k},\rho}(\alpha^m)=U_{\alpha^m}$ and $s_{\underline{k},\rho}(\alpha^m)\circ\beta_{\underline{k},\rho}=U_{\alpha^m}$. Taking $s_{\underline{k},\rho}{\eqdef}s_{\underline{k},\rho}(\alpha^m)$ , then we proved the proposition.
\end{proof}
The next proposition follows from the above one as \cite[Proposition 4.9]{HsiehCM}. Let $\underline{k}$ and $\rho$ be as before.
\begin{Prop}\label{3.5}
If $n\geq m>0$, then there is an isomorphism
$$\beta_{\underline{k},\rho}:e_\alpha\cdot V_{\underline{k}}(K_0(p,p^n),\mathcal{O}_m)\simeq e_\alpha (V_{n,m}\otimes L_\rho)^{K'_0(p)}[\underline{k}].$$
Here the $k$ in $[k]$ is regarded as a character of $\mathrm{diag}(1,1,\mathbb{Z}^\times_p, \mathbb{Z}^\times_p)$.
\end{Prop}

We are going to prove some control theorems and fundamental exact sequence for semi-ordinary forms along this two-dimensional weight space $\Spec\tilde{\Lambda}$. The following proposition follows from Lemma \ref{3.1} and Proposition \ref{3.2} in the same way as \cite[Lemma 4.10, Proposition 4.11]{HsiehCM}, noting that the level group is actually in $K_0(p)$ by the contraction property.
\begin{Prop}\label{3.7}
Let $e_\alpha.\mathcal{V}_{\underline{k}}(K_0(p,p^n)){\eqdef}\varinjlim_m e_\alpha\cdot V_{\underline{k}}(K_0(p,p^n),\mathcal{O}_m)$. Then $e_\alpha.\mathcal{V}(K_0(p,p^n))$ is $p$-divisible and
$$e_\alpha\cdot \mathcal{V}_{\underline{k}}(K_0(p,p^n))[p^m]=e_\alpha\cdot V_{\underline{k}}(K_0(p,p^n),\mathcal{O}_m)=e_\alpha\cdot H^0(I_G(K^n_1)[1/E],\omega_{\underline{k}})\otimes\mathcal{O}_m.$$
\end{Prop}
The next proposition is crucial to prove control theorems for semi-ordinary forms along the weight space $\mathrm{Spec}\tilde{\Lambda}$.
\begin{Prop}\label{3.8}
The dimensions of the spaces $e_\alpha M_{\underline{k}}(K_0(p,p^n),\mathbb{C})$ are uniformly bounded for all $\underline{k}\in\mathcal{X}_\rho$.
\end{Prop}
\begin{proof}
This is \cite[Proposition 3.0.5]{CLW}.
\end{proof}
The following theorem says that all semi-ordinary forms of sufficiently regular weights are classical, and can be proved in the same way as \cite[Theorem 4.19]{HsiehCM} using Proposition \ref{3.8}.
\begin{Theo}\label{3.9}
For each weight $\underline{k}=(a_1,a_2,a_3;b_1)\in\mathcal{X}_\rho$, there is a positive integer $A(\underline{a})$ depending on $\underline{a}=(a_1,a_2,a_3)$ such that if $b_1>A(\underline{a},n)$ then the natural restriction map
$$e_\alpha M_{\underline{k}}(K_0(p),\mathcal{O})\otimes\mathbb{Q}_p/\mathbb{Z}_p\simeq e_\alpha\cdot\mathcal{V}_{\underline{k}}(K_0(p))$$
is an isomorphism.
\end{Theo}
For $q=0$ define $\mathcal{W}^q$ as $\mathcal{W}$ but with the structure sheaf replaced by its cuspidal part. For $q=0$ or $\emptyset$ define the space of $\tilde{\Lambda}$-adic semi-ordinary forms
$$V_{\mathrm{so}}^q{\eqdef}\mathrm{Hom}(e_\alpha\cdot(\mathcal{W}^q\otimes L_\rho)^{K'_0(p)},\mathbb{Q}_p/\mathbb{Z}_p)\otimes_{\Lambda_{3,1}}\tilde{\Lambda}$$
$$\mathcal{M}_{\mathrm{so}}^q(K,\tilde{\Lambda}){\eqdef}\mathrm{Hom}_{\tilde{\Lambda}}(V_{\mathrm{so}}^q,\tilde{\Lambda}).$$

Thus from the finiteness results above Proposition \ref{3.7}, we get a form of Hida's control theorem.
\begin{Theo}\label{3.10}
Let $q=0$ or $\emptyset$. Then
\begin{itemize}
\item[(1)] $V^q_{\mathrm{so}}$ is a free $\tilde{\Lambda}$-module of finite rank.
\item[(2)] For any $k\in \mathcal{X}_\rho$ satisfying the assumption of Theorem \ref{3.9} we have $\mathcal{M}^q_{\mathrm{so}}(K,\tilde{\Lambda})\otimes\tilde{\Lambda}/P_{\underline{k}}\simeq e_\alpha\cdot M_{\underline{k}}^q(K,\mathcal{O}).$
\end{itemize}
\end{Theo}
\begin{proof}
Same as \cite[Theorem 4.21]{HsiehCM} using propositions \ref{3.2}, \ref{3.5}, theorem \ref{3.9} and proposition \ref{3.7}.
\end{proof}

\noindent\underline{Descent to Prime to $p$-Level}\\
The following proposition is needed to apply the description of local Galois representations at $p$ of semi-ordinary forms with prime to $p$ level above.
\begin{Prop}\label{3.11}
Suppose $\underline{k}$ is such that
$$a_1=a_2=0, a_3\equiv b_1\equiv 0\modulo p-1, a_2-a_3\gg 0, a_3+b_1>>0.$$
Suppose $F\in e_\alpha M_{\underline{k}}^0(K_0(p),\mathbb{C})$ is an eigenform with trivial nebentypus at $p$ such that $\rhobar_{F}$ can be written $\rhobar_{f}\oplus \psi_{1}\oplus\psi_{2}$ for $\psi_{i}$ a character of $G_{\Kcal}$. Let $\pi_F$ be the associated automorphic representation. Then $\pi_{F,p}$ is an unramified principal series representation.
\end{Prop}
\begin{proof}
Similar to \cite[proposition 4.17]{HsiehCM}. As $F$ is semi-ordinary at $p$, it has a $\pi_{F,p}$ has a fixed vector for $K_0(p)$. By the classification of admissible representations with $K_0(p)$-fixed vector (see e.g. \cite[Theorem 3.7]{CartierPadic}) we know $\pi_{F,p}$ has to be a subquotient of $\mathrm{Ind}_B^{\mathrm{GL}_4}\chi$ for $\chi$ an unramified character of $T_n(\mathbb{Q}_p)$. If this induced representation is irreducible then we are done. Suppose this is not the case. As $a_2-a_3>>0$, $a_3+b_1>>0$ and as $F$ is semi-ordinary, the character $\chi$ may be written as $\chi=\chi_1\otimes\chi_2\otimes\chi_3\otimes\chi_4$ with $\chi_1=\chi_2|\cdot|$ and $\chi_3,\chi_4$ having $p$-adic weights $\kappa_1=b_1$ and $\kappa_2=3-a_3$ respectively. This implies that $F$ is in fact ordinary and so that the local representation $\rhobar_{F}$ at $p$ is reducible. By our assumption, $\bar{\rho}_F^{ss}$ is the direct sum of $\bar{\rho}_f$ with two characters. So $\rhobar_{f}|G_{\qp}$ must itself be reducible. This contradicts our ongoing hypothesis that $\bar{\rho}_{f}|{G_{\mathbb{Q}_p}}$ is irreducible. Thus $\pi_{F,p}$ must be unramified.
\end{proof}

\noindent\underline{A Definition Using Fourier-Jacobi Expansion}\\
We can define a $\tilde{\Lambda}$-adic Fourier-Jacobi expansion map for families of semi-ordinary families as in \cite[4.6.1]{HsiehCM} by taking the $\tilde{\Lambda}$-dual of the Pontryagin dual of the usual Fourier-Jacobi expansion map (replacing the $e$'s in \emph{loc.cit.} by $e_\alpha$'s). We also define the $\Lambda$-adic Siegel operators $\Phi_{[g]}^h$'s by taking the $0$-th Fourier-Jacobi coefficient.
\begin{Def}
Let $A$ be a finite torsion free $\Lambda$-algebra. Let $\mathcal{N}_{\mathrm{so}}(K,A)$ be the set of formal Fourier-Jacobi expansions:
$$F=\{\sum_{\beta\in\mathscr{S}_{[g]}}a(\beta, F)q^\beta, a(\beta,F)\in (A\hat{\otimes}\hat{A}_{[g]}^\infty)^\Lambda \otimes H^0(\mathcal{Z}_{[g]}^\circ,\mathcal{L}(\beta))\}_{g\in X(K)}$$
such that for a Zariski dense set $\mathcal{X}_F\subseteq \mathcal{X}_\rho$ of points $\phi\in\mathrm{Spec}A$ where the induced point in $\mathrm{Spec}\Lambda$ is some arithmetic weight $\underline{k}_\zeta$, the specialization $F_\phi$ of $F$ is the highest weight vector of the Fourier-Jacobi expansion of a semi-ordinary modular form with tame level $K^{(p)}$, weight $\underline{k}$ and nebentype at $p$ given by $[\underline{k}][\underline{\zeta}]\omega^{-[\underline{k}]}$ as a character of $K_0(p)$.
\end{Def}
Then we have the following
\begin{Theo}\label{Theorem 3.14}
$$\mathcal{M}_{\mathrm{so}}(K,A)=\mathcal{N}_{\mathrm{so}}(K,A).$$
\end{Theo}
The proof is the same as \cite[Theorem 4.25]{HsiehCM}.\\

\noindent \underline{Fundamental Exact Sequence}

\noindent Now we prove a fundamental exact sequence for semi-ordinary forms. Let $w_3'=\begin{pmatrix}&&1&\\1&&&\\&1&&\\&&&1\end{pmatrix}$.
\begin{Lem}
Let $\underline{k}\in \mathcal{X}_\rho$ and $F\in e_\alpha M_{\underline{k}}(K_0(p,p^n), R)$ and $R\subset \mathbb{C}$. Let $W_2=\begin{pmatrix}1&&&\\&&1&\\&1&&\\&&&1\end{pmatrix}\cup \mathrm{Id}$ be the Weyl group for $G_P(\mathbb{Q}_p)$. There is a constant $A$ such that for any $\underline{k}\in \mathcal{X}_\rho$ such that $a_2-a_3>A, a_3+b_1>A$, for each $g\in G(\mathbb{A}_f^{(p)})$,
$\Phi_{P,wg}(F)=0$ for any $w\not\in W_2w_3'$.
\end{Lem}
The lemma can be proved using the computations in the proof of \cite[lemma 4.14]{HsiehCM}. Note that by semi-ordinarity and the contraction property the level group at $p$ for $F$ is actually $K_0(p)$.\\

\noindent The following is a semi-ordinary version of \cite[Theorem 4.16]{HsiehCM}, noting that $e_\alpha$ induces identity after the Siegel operator $\hat{\Phi}^{w_3'}$ by \cite[Proposition 4.5.2]{CLW}. The proof is also similar (even easier since the level group at $p$ is in fact in $K_0(p)$ by the contraction property).
\begin{Theo}\label{3.13}
For $\underline{k}\in\mathcal{X}_\rho$, we have
$$0\rightarrow e_\alpha\mathcal{M}^0_{\underline{k}}(K,A)\rightarrow e_\alpha\mathcal{M}_{\underline{k}}(K,A)\xrightarrow{\hat{\Phi}^{w_3'}=\oplus \hat{\Phi}_{[g]}^{w_3'}}\oplus_{g\in C(K)}\mathcal{M}_{\underline{k}'}(K_{P,0}^g(p),A)$$
is exact.
\end{Theo}
We need a family version of the fundamental exact sequence\begin{Theo}\label{TheoFundamental}
The following short sequence is exact :
\begin{equation}\nonumber
0\fleche e_\alpha\mathcal{M}^0(K,A)\fleche e_\alpha\mathcal{M}(K,A)\xrightarrow{\hat{\Phi}^{w_3'}=\oplus \hat{\Phi}_{[g]}^{w_3'}}\oplus_{g\in C(K)}\mathcal{M}(K_{P,0}^g(p),A)\fleche 0.
\end{equation}
\end{Theo}
\begin{proof}
The minimal compactification of the Igusa variety is affine, by \cite[Theorem 4.16]{HsiehCM}. After choosing a a weight $\underline{k}$ and specializing $\underline{k}$, the result thus follows from theorems \ref{3.9}, \ref{3.10} and \ref{3.13}.
\end{proof}
\subsection{Eisenstein Family}
\subsubsection{Klingen Eisenstein Family}
In this section we recall the construction in \cite{EischenWan} of the Klingen Eisenstein series, using the pullback formula for $$\Uni(3,1)\times\Uni(0,2)\hookrightarrow \Uni(3,3)$$
from some nearly holomorphic Siegel Eisenstein series $E_{\mathrm{sieg}}$ on $\Uni(3,3)$. We refer to \cite[Section 3.2-3.3]{EischenWan} for backgrounds on Siegel Eisenstein series and pullback formula. We suppose $\xi_0$ is a Hecke character of Archimedean type $(\kappa/2,-\kappa/2)$ for $\kappa\equiv 0(\mathrm{mod}\ 2(p-1))$. Suppose also the $p$-adic avatar of $\xi'_0{\eqdef}\xi_0\cdot(\epsilon^{-1}\circ\mathrm{Nm})$ factors through $\Gamma_\mathcal{K}$.
\begin{Prop}\label{EischenWan}\cite[Theorem 1.2]{EischenWan}
Let $\pi=\pi_f$ be the unitary automorphic representation generated by the weight $k$ form $f$. Let $\tilde{\pi}$ be the dual representation of $\pi$. Let $\Sigma$ be a finite set of primes containing all the bad primes
\begin{itemize}
\item[(i)] There is an element $\mathcal{L}^\Sigma_{f,\mathcal{K},\xi_0}\in\Lambda_{\mathcal{K},\mathcal{O}^{\ur}}\otimes_{\mathbb{Z}_p}\mathbb{Q}_p$ such that for any character $\xi_\phi$ of $\Gamma_\mathcal{K}$, which is the avatar of a Hecke character of conductor $p$, infinite type $(\frac{\kappa_\phi}{2}+m_\phi,-\frac{\kappa_\phi}{2}-m_\phi)$ with $\kappa_\phi$ an even integer which is at least $6$, $m_\phi\geq \frac{k-2}{2}$, we have
$$\phi(\mathcal{L}^\Sigma_{f,\mathcal{K},\xi_0})=\frac{L^\Sigma(\tilde{\pi},\xi_\phi,\frac{\kappa_\phi-1}{2})
\Omega_p^{4m_\phi+2\kappa_\phi}}{\Omega_\infty^{4m_\phi+2\kappa_\phi}}
c_{\phi}'.p^{\kappa_\phi -3}\mathfrak{g}(\xi_{\phi,2})^2\prod_{i=1}^2(\chi_{i}^{-1}\xi_{\phi,2}^{-1})(p)$$
where $c_{\phi}'$ is a constant coming from an Archimedean integral.
\item[(ii)] There is a set of formal $q$-expansions $\mathbf{E}_{f,\xi_0}{\eqdef}\{\sum_\beta a_{[g]}^t(\beta)q^\beta\}_{([g],t)}$ for $\sum_\beta a_{[g]}^t(\beta)q^\beta\in\Lambda_{\mathcal{K},\mathcal{O}^{\ur}}\otimes_{\mathbb{Z}_p}\mathcal{R}_{[g],\infty}$ , $([g],t)$ are $p$-adic cusp labels (we refer to \cite[Section 2.8.2]{EischenWan} for the notation and the notion of $p$-adic cusps), such that for a Zariski dense set of arithmetic points $\phi\in\mathrm{Spec}_{\mathcal{K},\mathcal{O}}$, $\phi(\mathbf{E}_{f,\xi_0})$ is the Fourier-Jacobi expansion of the highest weight vector of the holomorphic Klingen Eisenstein series constructed by pullback formula which is an eigenvector for $U_{t^+}$ with non-zero eigenvalue. The weight $L_{\underline{k}}$ for $\phi(\mathbf{E}_{f,\xi_0})$ is $(m_\phi+\frac{k-2}{2}, m_\phi-\frac{k-2}{2},0;\kappa_\phi)$.
\item[(iii)] The $a_{[g]}^t(0)$'s are divisible by $\mathcal{L}_{f,\mathcal{K},\xi_0}^\Sigma.\mathcal{L}_{\bar{\tau}'}^\Sigma$ where $\mathcal{L}_{\bar{\tau}'}^\Sigma$ is the $p$-adic $L$-function of a Dirichlet character as in \cite{EischenWan}.
\end{itemize}
\end{Prop}
We assumed in \cite{EischenWan} that the $\pi_{f,p}$ has distinct Satake parameters, which turns out to be unnecessary in our $\Uni(2)$ case. We also refer to \cite[Section 2.4]{EischenWan} for the convention of weights of automorphic forms on $\Uni(2)$ and $\Uni(3,1)$. This is just a translation of the main theorem of \cite{EischenWan} to the situation here.

To be compatible with arithmetic applications to our main results, we re-parameterize the family by a translation given below.
We define an automorphism of the Iwasawa algebra $\Lambda_\mathcal{K}$ as the following composition
$$t^{-c}_{\xi'_0}: \Lambda_\mathcal{K}\rightarrow\Lambda_\mathcal{K}\rightarrow\Lambda_\mathcal{K}$$
where the first map is given by $\gamma^\prime \mapsto \gamma^\prime\xi'_0(\gamma^{\prime,-1})$ for each $\gamma^\prime\in\Gamma_\mathcal{K}$, and the second map is determined by the map $\gamma^\prime\rightarrow \gamma^{\prime,-c}$ ($-c$ means inverse composed with complex conjugation) for each $\gamma^\prime\in\Gamma_\mathcal{K}$. For each $\phi\in\mathrm{Spec}\Lambda_\mathcal{K}$ corresponding to a character $\chi_\phi$ of $\Gamma_\mathcal{K}$, the composition of $\phi$ with $t^{-c}_{\xi^\prime_0}$ corresponds to the character $(\chi_\phi\xi^{\prime,-1}_0)^{-c}$. We write $\mathcal{L}^{\mathrm{Gr},\Sigma}_{f,\mathcal{K}}\in\mathrm{Frac}(W(\bar{\mathbb{F}}_p)) [[\Gamma_\mathcal{K}]])$ for the $p$-adic $L$-function $t^{-c}_{\xi'_0}(\mathcal{L}_{f,\mathcal{K},\xi_0})$ we constructed above (here $W(R)$ means the ring of Witt vectors of $R$ and we dropped the subscript $\xi$ as this $p$-adic $L$-function is distinguished corresponding to the trivial character). We can recover the full $p$-adic $L$-function $\mathcal{L}_{f,\mathcal{K}}$ by putting back the Euler factors at primes in $\Sigma$. We write this $p$-adic $L$-function as
\begin{equation}\label{GrpL}
\Lcal_{\Kcal}^{\Gr}(f).
\end{equation}

In proposition \ref{EischenWan}, an Archimedean constant $c_{\phi}'$ appears. It is perhaps not easy to compute this constant directly. Nevertheless, comparing the above $p$-adic $L$-function and Hida's Rankin-Selberg $p$-adic $L$-function, we can make the constant $c_{\phi}'$ precise.
\begin{Lem}\label{constant c}
The constant $c_\phi'$ of proposition \ref{EischenWan} is given by $$\Gamma(\kappa_\phi+m_\phi-\frac{k}{2})\Gamma(\kappa_\phi+m_\phi+\frac{k}{2}-1)2^{-3\kappa_\phi-4m_\phi+1}\pi^{1-2\kappa_\phi-2m_\phi}
i^{k-\kappa_\phi-2m_\phi-1}.$$
\end{Lem}
\begin{proof}
Let $\mathbf{g}$ be the Hida family corresponding to the family of characters of $\Gamma_\mathcal{K}$ (see \cite[Definition 7.8]{XinWanRankinSelberg}). We pick an auxiliary Hida family of ordinary forms $\mathbf{f}'$ (for example using CM forms) and compare
\begin{itemize}
\item The product $\mathcal{L}^{\mathrm{Hida}}_{\mathbf{f}'\otimes \mathbf{g}}\cdot\mathcal{L}^{\mathrm{Katz}}_\mathcal{K}h_\mathcal{K}$, where $\mathcal{L}^{\mathrm{Katz}}_\mathcal{K}h_\mathcal{K}$ is the class number $h_\mathcal{K}$ of $\mathcal{K}$ times the Katz $p$-adic $L$-function, which interpolates the Petersson inner product of specializations of $\mathbf{g}_\phi$ (we refer to \cite[Section 7.5]{XinWanRankinSelberg} for details). The $\mathcal{L}^{\mathrm{Hida}}_{\mathbf{f}'\otimes\mathbf{g}}$ is the Rankin-Selberg $p$-adic $L$-function constructed by Hida in \cite{HidaFourier} interpolating algebraic part of the critical values of Rankin-Selberg $L$-functions for specializations of $\mathbf{f}'$ and $\mathbf{g}$, where the specializations of $\mathbf{g}$ has higher weight.
\item The $p$-adic $L$-function $\mathcal{L}_{\mathbf{f}',\mathcal{K}}$ constructed using the doubling method as above.
\end{itemize}
We first look at the interpolation formulas at arithmetic points where the Siegel Eisenstein series are of scalar weight. The computations are done in \cite{XinWanANT} (although the ramifications in \emph{loc.cit} is slightly different, however those assumptions are put for constructing the family of Klingen Eisenstein series. The computations in the doubling method construction of the $p$-adic $L$-function carries out in the same way in the situation here).  We see that the above two items have the same value at these points. As these arithmetic points are Zariski dense, the $\mathcal{L}^{\mathrm{Hida}}_{\mathbf{f}'\otimes \mathbf{g}}\cdot\mathcal{L}^{\mathrm{Katz}}_\mathcal{K}h_\mathcal{K}$ and $\mathcal{L}_{\mathbf{f}',\mathcal{K}}$ should be equal. Then we look at the arithmetic points considered in the above proposition. Comparing the interpolation formulas here and in \cite[Theorem I]{HidaFourier}, we get the formulas for $c_\phi'$ (note that the critical $L$-value is not zero since it is away from center).
\end{proof}
For convenience of application we normalize our Klingen Eisenstein series by twisting by an anticyclotomic character so that the specializations are of the weight $L_{\underline{k}}=(\frac{k-2}{2}, -\frac{k-2}{2}, -m_\phi; \kappa_\phi+m_\phi)$. Now we still write $\mathcal{L}^{\mathrm{Hida}}_{f\otimes\mathbf{g}}$ for the Rankin-Selberg Hida $p$-adic $L$-function interpolating critical values of the Rankin-Selberg $L$-function for $f$ and specializations of $g$ whose weight is higher than $f$ (see \cite[Definition 7.8]{XinWanRankinSelberg}). Since the higher weight form $\mathrm{g}$ is ordinary, Hida's construction works in the same way even though $f$ is not ordinary. We have the following corollary by comparing interpolation formulas (see \cite[(7-2), (7-5)]{XinWanRankinSelberg} for details).
\begin{Cor}
$$\mathcal{L}^{\mathrm{Hida}}_{f\otimes\mathbf{g}}\cdot\mathcal{L}^{\mathrm{Katz}}_\mathcal{K}h_\mathcal{K}
=\Lcal_{\Kcal}^{\Gr}(f).$$
\end{Cor}
The corollary follows from above lemma and the interpolation formulas on both hand sides.\\
We prove the following
\begin{Lem}
The $p$-adic $L$-function $\Lcal_{\Kcal}^{\Gr}(f)$ belongs to $\mathcal{O}^{{\ur}}[[\Gamma_\mathcal{K}]]$.
\end{Lem}
\begin{proof}
From the construction the denominator of $\Lcal_{\Kcal}^{\Gr}(f)$ can only be powers of $p$ times the product of the Euler factors of a finite number of primes of $\Lcal_{\Kcal}^{\Gr}(f)$. But by the argument in \cite[Proposition 8.3]{XinWanIMC} we know the denominator can at most be powers of $Y$ if we take $\mathbb{Z}_p[[Y]]$ as the coefficient ring of $\mathbf{g}$. These two sets are disjoint. Thus the denominator must be a unit.
\end{proof}

\subsubsection{Study of the Fourier-Jacobi functional}\label{Study}
It remains to study the $p$-adic property of the Klingen Eisenstein family. Before continuing we refer to \cite[Introduction]{XinWanRankinSelberg} for an outline of the strategy to construct and study the Fourier-Jacobi functional (this is very helpful for understanding the argument below). It uses the pullback formula of $E_{\mathrm{Sieg}}$ on $\Uni(3,3)$ via
$$\Uni(3,1)\times\Uni(0,2)\hookrightarrow \Uni(3,3)$$
with
$$E_{\mathrm{Sieg}}|_{\Uni(3,1)\times\Uni(2)}=E_{\mathrm{Kling}}\boxtimes f.$$
We first compute the Fourier-Jacobi coefficient for $E_{\mathrm{Sieg}}$ and pair it with the form $f$ on $\Uni(0,2)$, and thus obtain the Fourier-Jacobi coefficient of $E_{\mathrm{Kling}}$ (detailed below).
We refer to \cite[Section 6D]{XinWanRankinSelberg} for the computation of Fourier-Jacobi expansion for $E_{\mathrm{Sieg}}$ and the notion of local Fourier-Jacobi integrals.

\noindent\underline{Idea Here}\\
The main difficulty here is that the explicit local Fourier-Jacobi computation at Archimedean place is complicated if the weight is not scalar. Our idea is not to compute such local Fourier-Jacobi integrals at $\infty$. Instead we fix the weight $(k,\kappa_\phi,m_\phi)$ (notation as before) and vary its nebentypus of the Klingen Eisenstein series at $p$. Such arithmetic points are Zariski dense in $\Spec\mathcal{O}[[\Gamma_\mathcal{K}]]$. We construct (as in previous works) a Fourier-Jacobi functional $\mathrm{FJ}^{\mathbf{h}}_{\beta,\theta_1}$ (in Definition \ref{definition 3.7}) on families on $\Uni(3,1)$, and show that we can factor out a number $C_\infty$ depending only on the Archimedean data (and is thus the same number for all arithmetic points), which can be proved to be non-zero, and an element $\mathcal{L}\in\mathcal{O}^{{\ur}}[[\Gamma_\mathcal{K}]]^\times\cdot\bar{\mathbb{Q}}_p^\times$, such that at each of such arithmetic point we have
$$\mathrm{FJ}^{\mathbf{h}}_{\beta,\theta^\star}(\phi(E_{\mathrm{Kling}})))=C_\infty\cdot\phi(\mathcal{L}).$$
The key idea is to employ Ikeda's theory on Fourier-Jacobi coefficients to get lemma \ref{LemCoreIkeda}, which says that the Fourier-Jacobi coefficient of a nearly holomorphic Siegel Eisenstein series on $\Uni(3,3)$ is a \emph{finite sum} of products of Siegel Eisenstein series on $\Uni(2,2)$ and theta functions.\footnote{We are thankful to Ikeda for showing us the argument.}\\

\noindent\underline{Interpolating Inner Products}\\
We first recall a construction of Hsieh \cite[Definition 4.3]{HsiehHidaFamilies} which is important in \cite[Section 7]{XinWanRankinSelberg}. Let $\Lambda_{\Uni(2)}=\mathcal{O}_L[[T_1, T_2]]$ be the weight algebra for $\Uni(2)$ (\textit{i.e.} the complete group ring of the diagonal torus of $\Uni(2)(1+p\mathbb{Z}_p)$). Write $\mathcal{M}_{\mathrm{ord}}(K,\Lambda_{\Uni(2)})$ for the space of ordinary Hida families on $\Uni(2)$ with respect to a tame level group $K\subset\Uni(2)(\mathbb{A}_f)$. We similarly write $\breve{\mathcal{M}}_{\mathrm{ord}}(K,\Lambda_{\Uni(2)})$ for the space of such Hida families with the weight map being the inverse of that of $\mathcal{M}_{\mathrm{ord}}(K,\Lambda_{\Uni(2)})$.
\begin{Def}\label{define measure}
For a \emph{neat} tame level group $K\subset \Uni(2)(\mathbb{A}^{(p\infty)})$ we define a $\Lambda_{\Uni(2)}$-pairing $\mathbf{B}_K\langle-,-\rangle$
$$\mathbf{B}_K: \mathcal{M}_{\mathrm{ord}}(K,\Lambda_{\Uni(2)})\times\breve{\mathcal{M}}_{\mathrm{ord}}(K,\Lambda_{\Uni(2)})\rightarrow \Lambda_{\Uni(2)}.$$
For any $\mathbf{f}\in\mathcal{M}_{\mathrm{ord}}(K, \Lambda_{\Uni(2)})$, $\mathbf{g}\in\breve{\mathcal{M}}_{\mathrm{ord}}(K, \Lambda_{\Uni(2)})$, $\phi\in\Spec\Lambda_{\Uni(2)}(\mathbb{C}_p)$ a weight two point (in other words, the specializations of $\mathbf{f}_\phi$ and $\mathbf{g}_\phi$ have scalar weights), and any $n>0$ we define
\begin{align*}
&\mathbf{B}_{K,n}\langle\mathbf{g},\mathbf{f}\rangle&&{\eqdef}\sum_{[x_i]\in\Uni(2)(\mathbb{Q})\backslash \Uni(2)/KU_0(p^n)}U^{-n}_p\mathbf{f}(x_i)\mathbf{g}(x_i\begin{pmatrix}&1\\p^n&\end{pmatrix})&\\
&&&(\mathrm{mod}(1+T_1)^{p^n}-1, (1+T_2)^{p^n}-1).&
\end{align*}
Then Hsieh proves that
$$\mathbf{B}_{K,n+1}\equiv\mathbf{B}_{K,n}(\mathrm{mod}(1+T_1)^{p^n}-1,(1+T_2)^{p^n}-1).$$
We define
$$\mathbf{B}_K\langle\mathbf{g},\mathbf{f}\rangle=\lim_n\mathbf{B}_{K,n}\langle\mathbf{g},\mathbf{f}\rangle.$$
By definition we have
$$\phi(\mathbf{B}_{K}\langle\mathbf{g},\mathbf{f}\rangle)=\sum_{[x_i]\in\Uni(2)(\mathbb{Q})\backslash \Uni(2)/KU_0(p^n)}U^{-n}_p\mathbf{f}_\phi(x_i)\mathbf{g}_\phi(x_i\begin{pmatrix}&1\\p^n&\end{pmatrix})$$
and hence
\begin{align*}&\phi(\mathbf{B}_{K}\langle\mathbf{g},\mathbf{f})\rangle&&=\mathrm{vol}(KU_0(p^n))^{-1}\int_{[\Uni(2)]}U^{-n}_p\mathbf{f}_\phi(h)\mathbf{g}_\phi
(h\begin{pmatrix}&1\\p^n&\end{pmatrix})dh&\\
&&&=\mathrm{vol}(KU_0(p^n))^{-1}\int_{[\Uni(2)]}\mathbf{f}_\phi(h\begin{pmatrix}&1\\1&\end{pmatrix}_p)\mathbf{g}_\phi
(h)dh&
\end{align*}
if $\phi$ corresponds to an ordinary form whose $p$-part conductor is $p^n$.
\index{$\mathbf{B}\langle,\rangle$}
\end{Def}

In the below we fix and suppress the tame level group $K$ and write the pairing as $\mathbf{B}\langle,\rangle$. Hsieh also proved (\cite[Lemma 4.4]{HsiehHidaFamilies}) when specializing the pairing to arithmetic points $Q$ of conductor $p^\alpha$ (allowed to have vector valued weight), this gives
$$\langle U^{-\alpha}_p\mathbf{f}_Q,\mathbf{f}_Q\rangle.$$
(Petersson inner product of vector valued forms in \emph{loc.cit}).

We recall some explicit constructions in \cite[Section 7.7]{CLW} (also in \cite[Section 8B]{XinWanRankinSelberg}). These are important in our study of the Fourier-Jacobi functional.

\begin{Def}\label{define notation}
\begin{itemize}
\item Hida families $\boldsymbol{\theta}$, $\tilde{\boldsymbol{\theta}}_3$, $\mathbf{h}$, $\tilde{\mathbf{h}}_3$ of CM forms on $\Uni(2)$ constructed in \cite[Section 8B]{XinWanIMC}, using the Hecke characters chosen in \emph{loc.cit.}.
\item Elements $g_1$, $g_2$, $g_3$ and $g_4$ defined in \cite[Definition 8.20]{XinWanIMC} in $\mathrm{GU}_2(\mathbb{A}_\mathbb{Q})$ whose $p$ components are $1$.
\item Some $p$-adic $L$-functions $\mathcal{L}_1, \mathcal{L}_2, \cdots, \mathcal{L}_6$, where
$\mathcal{L}_1$, $\mathcal{L}_5$ and $\mathcal{L}_6$ are units in the Iwasawa algebra, while $\mathcal{L}_2$ is a nonzero constant fixed throughout the family.
\end{itemize}
\end{Def}

Note that the families $\mathbf{h}$ and $\tilde{\mathbf{h}}_3$ are constructed via interpolating at weight two arithmetic points.
By standard facts about Hida control theorems, we know the families $\mathbf{h}$ and $\tilde{\mathbf{h}}_3$ also interpolate the highest weight vectors of CM forms at arithmetic points of weight $(\frac{k-2}{2},-\frac{k-2}{2})$ and $(-\frac{k-2}{2},\frac{k-2}{2})$ respectively. From the computations in \cite[Section 8D]{XinWanRankinSelberg}, we know $\mathbf{B}\langle\mathbf{h},\tilde{\mathbf{h}}_3\rangle$ and $\mathbf{B}\langle\boldsymbol{\theta},\tilde{\boldsymbol{\theta}}_3\rangle$ (Hsieh's $\Lambda$-adic pairing) are interpolated by $p$-adic $L$-functions $\mathcal{L}_3$ and $\mathcal{L}_4$. So Hsieh's theory in \cite{HsiehHidaFamilies} also implies the Petersson inner products at weight $k$ arithmetic points are also interpolated by $\mathcal{L}_3$ and $\mathcal{L}_4$ respectively.\\

\noindent\underline{Geometry and Fourier-Jacobi Functional}\\
Now we look at the Fourier-Jacobi expansion theory for forms on $\mathrm{GU}(3,1)$.
Recall \cite[Section 3.6]{XinWanRankinSelberg} that for $\beta\in\mathbb{Q}^\times$ there is a line bundle $\mathcal{L}(\beta)$ on the boundary component $\mathcal{Z}_{[g]}$ ($[g]$ is some cusp label) of the Shimura variety for $\mathrm{GU}(3,1)$. We refer to \emph{loc.cit.} and \cite[Section 3.6]{HsiehCM} for the theory of Fourier-Jacobi coefficients for forms on $\mathrm{GU}(3,1)$. Some points are worth pointing out when working with vector valued forms. We refer to \cite{LanComparison} for the comparison between algebraic and analytic Fourier-Jacobi coefficients.  As noted in \cite[Section 3.6]{HsiehCM}, the algebraic Fourier-Jacobi coefficient takes values only in the $N^1_H$-coinvariants (see 3.6.2 of \emph{loc.cit} for the notation) of the representation $L_{\underline{k}}$, due to the ambiguity in choosing a basis for the differentials of the Mumford family. Thus this $\beta$-th Fourier-Jacobi expansion takes values in the space of forms on $\Uni(2)$ tensored with the space of global sections of $\mathcal{L}(\beta)$. We only look at the $L_{\underline{k}'}{\eqdef}L_{(\frac{k-2}{2},-\frac{k-2}{2})}$-components (regarded as a sub-representation of the restriction of the representation $L_{\underline{k}}$ to $\mathrm{GL}_2$, which clearly appears with multiplicity one in this restriction). In fact the ambiguity on the choice above does not make any differences when looking at this $L_{\underline{k}'}$-component, by the description of the $N^1_H$-coinvariants in the proof of \cite[Lemma 3.12]{HsiehCM} (this corresponds to the $L_{\underline{k}'}$-component there). According to the description in \cite[Section 5.3]{LanComparison} this corresponds to looking at a quotient of the $\mathcal{E}^{\Phi^{(g)}_{\mathcal{H}},\delta^{(g)}_{\mathcal{H}}}_{M,\mathrm{an}}(W)$ in \emph{loc.cit.} which is the pullback to $\mathcal{Z}_{[g]}$ of an automorphic vector bundle of weight $\underline{k}'$ on the Igusa variety for the definite unitary group $\Uni(2)$, tensoring with the line bundle $\mathcal{L}(\beta)$. We look at theta functions $\theta^\star$ on the Klingen parabolic subgroup $P(\mathbb{A})\subseteq \mathrm{GU}(3,1)(\mathbb{A})$, which are exactly the ones considered in \cite[Section 6.10]{XinWanRankinSelberg} and \cite[Section 5]{XinWanIMC}, \textit{i.e.} corresponds to the dual of the space of global sections of $\mathcal{L}(\beta)$. It is defined as in \cite[Definition 6.45]{XinWanRankinSelberg} corresponding to the Schwartz function $\phi_1=\phi_{1,\infty}\times\prod_{v<\infty}\phi_{1,v}$ there, and define the functional $l_{\theta^\star}$ as in \cite[Section 8E]{XinWanRankinSelberg} (this corresponds to inner product integral over $N_P(\mathbb{Q}\backslash\mathbb{A}_{\mathbb{Q}})$ with the theta function $\theta_{\phi_1}$ corresponding to $\phi_1$ where $N_P$ is the unipotent radical of the Klingen parabolic subgroup $P$), such that $l_{\theta^\star}\in \mathrm{Hom}(H^0(\mathcal{Z}_{[g]}, \mathcal{L}(\beta)), \mathcal{O})$.

\begin{Def}\label{definition 3.7}
We define the Fourier-Jacobi functional on $\Uni(3,1)$ families by
$$\mathrm{FJ}^{\mathbf{h}}_{\beta,\theta^\star}: F\mapsto \mathbf{B}\langle e^{\mathrm{ord}}l_{\theta^\star}\mathrm{FJ}_{\beta}(F),\pi(g_2)\mathbf{h}\rangle.$$
We note that since the $\mathrm{FJ}_\beta$ map takes the $N^1_H$ covariant quotient, so it takes values in the $L_{(\frac{k}{2},-\frac{k}{2})}$-component.
\end{Def}
\paragraph{Remark:}
To see that
$$F\mapsto \mathbf{B}\langle e^{ord}l_{\theta^\star}(\mathrm{FJ}_\beta(F)),\pi(g_2)\mathbf{h}\rangle$$
indeed gives a functional on the space of semi-ordinary families (over the two-dimensional weight space) of forms on $\Uni(3,1)$, we note that the Archimedean weight is fixed throughout the two-dimensional family. In the theory of $p$-adic semi-ordinary families (as developed in \cite{XinWanIMC}) we are interpolating the highest weight vector of the automorphic forms. Thus each component of the $L_{k'}$-projection of $\mathrm{FJ}_\beta(F)$ is interpolated $p$-adic analytically.  Note also that in the ordinary case when $k>2$ the family of functionals constructed here is different from the two-variable specializations of the construction in \cite{XinWanRankinSelberg}. (The highest weight vector is not in the $(\frac{k-2}{2}, -\frac{k-2}{2})$-component above).

We now compute this functional on Klingen Eisenstein series at arithmetic points. We first need some preparations.\\
\noindent\underline{Ikeda Theory}\\
We refer to \cite[Section 6.1]{XinWanRankinSelberg} for the background of Siegel Eisenstein series and write $I_n(\tau)$ for the corresponding space of (local or global) Siegel sections defined using a Hecke or local character $\tau$. Write $E=E_{\mathrm{Sieg}}$ for the Siegel Eisenstein series on $\Uni(3,3)$ that we use in the pullback formula. Now we put ourselves in the context of \cite{IkedaFourierJacobi}. We are in the $m=1$ and $n=2$ of \cite[Section 2, case 2]{IkedaFourierJacobi} (see the definitions of $X,Y,Z,V$ there). Let $\psi$ be an additive character of $\mathbb{A}$ and $\omega_\psi$ be the Weil representation there.

We first make some observations. For $\beta\in\mathbb{Q}^\times$, let $\psi=\psi_\beta$ be the additive character corresponding to $\beta$. Let $\pi$ be the representation of the Jacobi group $N'\Uni(2,2)$ (see \cite[Introduction]{XinWanRankinSelberg} for the notations)  generated by the $\beta$-th Fourier-Jacobi coefficient of a nearly holomorphic Siegel Eisenstein series on $\Uni(3,3)$. Then by \cite[Proposition 1.3]{IkedaFourierJacobi},
there is a map
$$\omega_\psi\otimes I\hookrightarrow\pi$$
with dense image, where $\omega_\psi$ is the Weil representation of the Jacobi group defined there (determined by a splitting character $\lambda$ we fix throughout), and $I$ is the representation generated by the integral of \emph{loc.cit}. By \cite[Theorem 3.2]{IkedaFourierJacobi} this is the sub-representation of the automorphic representation corresponding to Siegel Eisenstein series on $\Uni(2,2)$. Let $\mathfrak{g}$ be the Lie algebra for $\Uni(2,2)(\mathbb{R})$ and $K$ be a maximal compact subgroup of it. Let $\mathfrak{k}$ be its Lie algebra. Write the Harish-Chandra decomposition of the complex Lie algebra by
$$\mathfrak{g}^\mathbb{C}=\mathfrak{k}^\mathbb{C}\oplus\mathfrak{p}^+\oplus\mathfrak{p}^-.$$
By a \emph{lowest weight representation} we mean a $(\mathfrak{g}, K)$-module generated by elements which are killed by some finite power of $\mathfrak{p}^-$. Then as $E_{\mathrm{Sieg}}$ is nearly holomorphic by construction, we know $\pi$ is a lowest weight representation. Also as $\omega_\psi$ is a Weil representation, it is also a lowest weight representation (see \cite{AdamsTheta}). By the Leibnitz rule
$$X^-(v_1\otimes v_2)=X^-v_1\otimes v_2+v_1\otimes X^-v_2$$
for $v_1\in\omega_\psi$, $v_2\in I$ and $X^-\in\mathfrak{p}^-$, we see $I$ is also a lowest weight representation. Note that under the Archimedean theta correspondence of $\Uni(1)(\mathbb{R})\times\Uni(2,2)(\mathbb{R})$, any component of the Weil representation of $\Uni(2,2)(\mathbb{R})$ with given central character must be irreducible (since $\Uni(1)(\mathbb{R})$ is compact, see \textit{e.g.} \cite{AdamsTheta}). Recall that the space $I$ is contained in the space of Siegel Eisenstein series on $\Uni(2,2)$ with fixed character (\cite[Theorem 3.2]{IkedaFourierJacobi}). We also note the well known fact that the representation of $\Uni(2,2)(\mathbb{R})$ admissibly induced from a character of the Siegel parabolic subgroup is of finite length. Now we claim that for a fixed $K$-type $\sigma$ of $\pi$, there are only finitely many $K$-types of $\omega_\psi$ and $I$, whose tensor product contains $\sigma$. This can be seen, for example by noting that in any lowest weight representation generated by a single vector, the $\Uni(1)(\mathbb{R})$-weights are monotonically going up with finite dimensional eigenspaces for each eigencharacter. It follows that the map
\begin{equation}\label{embe}
\omega_\psi\otimes I\hookrightarrow \pi
\end{equation}
is surjective for $K$-finite vectors, from the density of the image of the above $\hookrightarrow$.\\

\noindent\underline{Notation}: Write $w_1$, $w_2$, ... for a basis of $L_{(\frac{k-2}{2},-\frac{k-2}{2})}$. To save notation, in the following we also write $E_{\mathrm{sieg}}$ for the $w_1$-component adelic Siegel Eisenstein series in the $L_{(\frac{k-2}{2},-\frac{k-2}{2})}$-component of $E_{\mathrm{sieg}}$ written out using this basis.

We refer to \cite[Theorem 3.2]{IkedaFourierJacobi} for the definition of an integral operator $R(f,\phi)$, where $f$ is a Siegel Eisenstein section on $\Uni(3,3)$ and $\phi$ is a Schwartz function. It also makes sense to talk about local versions of these integrals.
From the above discussion we get the following key lemma from \cite[Lemma 1.1, Proposition 1.3]{IkedaFourierJacobi}, which states that the $\mathrm{FJ}_\beta(E_{\mathrm{sieg}})$ is a finite sum of products of theta functions and Siegel Eisenstein series on $\Uni(2,2)$.
\begin{Lem}\label{LemCoreIkeda}
There are a finite number of Archimedean Schwartz functions $\phi'_{4,i,\infty}$, $\phi'_{2,i,\infty}$ of $X(\mathbb{R})$ for $X$ in the Introduction of \cite{IkedaFourierJacobi}, such that
$\mathrm{FJ}_\beta(E_{\mathrm{sieg}})$ can be written as a finite sum of expressions as $$\Theta_{\phi'_{4,i,\infty}\otimes\prod_{v<\infty}
\phi_{i,v}}(nh)E(R(f_\infty, \phi'_{2,i,\infty})\cdot\prod_{v<\infty}f_{v,i}')$$ for some Siegel sections $\prod_{v<\infty}f_{v,i}'\in I_2(\tau)$ and Schwartz functions $\prod_{v<\infty}\phi_{i,v}$.
\end{Lem}
Note with the formula on page 620 of \cite{IkedaFourierJacobi} that the inner product of theta functions is equal to the inner product of the kernel Schwartz function. We choose $\phi_{4,i,\infty}$ to be a basis contributing to the surjection (\ref{embe}) on the $K$-type of $\omega_\psi$, and the $\phi_{2,i,\infty}$ a set of dual basis (finite dimensional) with respect to the pairing in \emph{loc.cit.}. These $\phi_{4,i,\infty}$ and $\phi_{2,i,\infty}$ depend only on $\underline{k}$ ($i=1,2,\cdots$).

Define the local $\beta$-th Fourier-Jacobi integral for $f_v\in I_3(\tau_v)$ as in \cite[Definition 6.4]{XinWanRankinSelberg} by
$$\mathrm{FJ}_\beta(f_v, nh)=\int_{\mathbb{Q}_v}f_v(\begin{pmatrix}&1_3\\-1_3&\end{pmatrix}\begin{pmatrix}1_3&\begin{matrix}S_v&0\\0&0_2
\end{matrix}\\&1_3\end{pmatrix}nh)e(-\beta S_v)dS_v.$$
Note the relation between this integral and the integral in \cite[Theorem 3.2]{IkedaFourierJacobi}.

Note also the computations in \cite[Section 6]{XinWanRankinSelberg} imply for each $v<\infty$ and $f_v$ the Siegel section on $\Uni(3,3)$ that we use to define the Siegel Eisenstein series $E_{\mathrm{sieg}}$, the local Fourier-Jacobi integral can be written in the form
$$\mathrm{FJ}_{\beta}(f_v)=\sum_{j_v=1}^{n_v}f_{j_v}\phi_{j_v}$$
for $f_{j_v}\in I_2(\tau_v)$ Siegel sections on $\Uni(2,2)$ and $\phi_{j_v}$ local Schwartz functions.
Then from lemma \ref{LemCoreIkeda} and the computation in \cite[Page 628]{IkedaFourierJacobi} on
$$\langle \mathrm{FJ}_\beta(E_{\mathrm{sieg}}),\Theta_\phi\rangle$$
and choosing the test Schwartz function $\phi$ properly, noting the duality between the $\phi_{2,i,\infty}$ and $\phi_{4,i,\infty}$'s (and take appropriate Schwartz functions at non-Archimedean places as well), we obtain that
\begin{equation}\label{FJJ}
\mathrm{FJ}_\beta(E_{\mathrm{sieg}})=\sum_i
(\prod_v\sum_{j_v})E(R(f_\infty,
\phi_{2,i,\infty})\cdot\prod_v f_{j_v}, -)\cdot \Theta_{\phi_{4,i,\infty}\cdot\prod_{v<\infty}\phi_{j_v}},
\end{equation}

\noindent\underline{A choice for $\beta$}\\
We can ensure that the $L_{(\frac{k-2}{2},-\frac{k-2}{2})}$-component of $l_{\theta^\star}(\mathrm{FJ}_\beta(E_{\mathrm{Kling}}))$ is non-zero for some $\theta^\star$ and $\beta\in\mathbb{Q}^\times\cap \mathbb{Z}_p^\times$. If for all nonzero $\beta$ this component is zero, then the $L_{(\frac{k-2}{2},-\frac{k-2}{2})}$-component of $E_{\mathrm{Kling}}$ is a constant function on the Shimura variety of $\mathrm{GU}(3,1)$. This contradicts the description of the boundary restriction of $E_{\mathrm{Kling}}$, namely the $L_{(\frac{k-2}{2},-\frac{k-2}{2})}$-component is non-zero at some cusp (namely at $w_3$ in \cite[Section 6.8]{XinWanRankinSelberg}) while is zero at other cusps.

Thus there must be a $\beta'\not=0$ such that $\mathrm{proj}_{(\frac{k-2}{2},-\frac{k-2}{2})}(\mathrm{FJ}_{\beta'}(E_{\mathrm{Kling}}))\not=0$. Let $\beta'=p^n\beta''$ for $\beta''\in\mathbb{Z}_p^\times$ and $n\in\mathbb{Z}$. Let $y$ be an element in $\mathcal{K}^\times$ which is very close to $(p,1)$ in the $p$-adic topology of $\mathcal{K}_p$. Then
$\mathrm{diag}(y\bar{y}, y, y, 1)^n\in\Uni(3,1)(\mathbb{Q})$. Set $\beta=\beta'(y\bar{y})^{-n}\in\mathbb{Z}_p^\times\cap \mathbb{Q}$ then $$\mathrm{proj}_{(\frac{k-2}{2},-\frac{k-2}{2})}(\mathrm{FJ}_\beta\rho(\mathrm{diag}(y\bar{y}, y, y, 1)^n_p) E_{\mathrm{Kling}})$$
is not the zero function. So there must be a choice of $\theta^\star$ such that the weight $(\frac{k-2}{2},-\frac{k-2}{2})$-component of $l_{\theta^\star}(\mathrm{FJ}_\beta(E_{\mathrm{Kling}}))$ above is non-zero. (The reason of making sure that $\beta\in\mathbb{Z}_p^\times$ is that only for those $\beta$ we did the Fourier-Jacobi coefficient computation at $p$ for the Klingen Eisenstein series in \cite[Section 6H]{XinWanRankinSelberg}.) Without loss of generality we assume that the $w_1$-component of it is nonzero using the basis $w_1$, $w_2$,... we fixed before.

\noindent\underline{Compute the Fourier-Jacobi Functional}\\
Suppose we are doing our computations at an arithmetic point $\mathsf{z}$ of conductor $p^t$ (notation as in \cite{XinWanRankinSelberg}) and write $\mathbf{h}_\mathsf{z}$ for the specialization of $\mathbf{h}$ (the CM family in Definition \ref{define notation}) at $\mathsf{z}$. By the doubling method for $\mathbf{h}_\mathsf{z}$ under $\iota: \Uni(2)\times\Uni(2)\hookrightarrow \Uni(2,2)$ with the Siegel Eisenstein series on $\Uni(2,2)$
$$E(R(f_\infty,
\phi_{2,i,\infty})\cdot\prod_v f_{j_v}, -)$$ above (see \cite[Proposition 6.1]{XinWanRankinSelberg} for details), we know by the pullback formula (see \cite[Proposition 3.4]{EischenWan}) expressing the Klingen Eisenstein series as pullback of Siegel Eisenstein series,
and (\ref{FJJ}) (for detail of this argument, see \cite[Proposition 8.24]{XinWanRankinSelberg}, which uses lemma 6.46 and Corollary 6.47 there), there is a constant $C_{i,\infty}$ such that the Fourier-Jacobi functional in Definition \ref{definition 3.7} computed as (recall we noted that the $L_{(\frac{k-2}{2},-\frac{k-2}{2})}$ appears with multiplicity one in the restriction of representations of $\mathrm{GL}_3$ to $\mathrm{GL}_2$),
\begin{align}\label{Fourier-Jacobi}
&&&\int_{[N_P/Z(N_P)]}\int_{[\Uni(2)]}\mathrm{FJ}_\beta(E_{\mathrm{Kling}})(ng)\mathbf{h}_\mathsf{z}(g)\theta_{\phi_1}(n)dndg&\notag\\
 &=&&p^t\mathsf{z}(\mathcal{L}_5\mathcal{L}_6)C_{i,\infty}\int_{[\Uni(2)]}\mathbf{h}_{\mathsf{z}}(g)\cdot\boldsymbol{\theta}^{\mathrm{low}}_\mathsf{z}(g)f(g)dg&
\notag \\&=&&
p^t\mathsf{z}(\mathcal{L}_5\mathcal{L}_6)C_{i,\infty}\int_{[\Uni(2)]}\mathbf{h}_\mathsf{z}(g)
\boldsymbol{\theta}^{\mathrm{low}}_\mathsf{z}(g)f(g)dg.&
\end{align}
Here we write $[G]$ to denote $G(\mathbb{Q})\backslash G(\mathbb{A}_\mathbb{Q})$ for a group $G$, and denote $Z(G)$ as the center of $G$. The $N_P$ is the unipotent radical of the Klingen parabolic subgroup $P$ of $\Uni(3,1)$. The $f(g)$ is the form on $\Uni(2)$ given by the $\mathrm{GL}_2$ modular form we study extended using the trivial character via $\mathrm{GU}(2)=\mathrm{GL}_2\times_{\mathbb{Q}^\times}\mathcal{K}^\times$. The $\mathcal{L}_5$ and $\mathcal{L}_6$ are as in Definition \ref{define notation}, and come from the pullback formula for $\mathbf{h}_\mathsf{z}$ above. The $\boldsymbol{\theta}^{\mathrm{low}}_\mathsf{z}(g){\eqdef}\boldsymbol{\theta}_\mathsf{z}(g\begin{pmatrix}&1\\1&\end{pmatrix}_p)$ as in \cite{XinWanRankinSelberg}.  The $C_{i,\infty}$ comes from the product of
\begin{itemize}
\item The local pullback integral at $\infty$ of the Siegel section $R(f_\infty,\phi_{2,i,\infty})$ with respect to specializations $\mathbf{h}_\mathsf{z}$ under $\Uni(2)\times\Uni(2)\hookrightarrow\Uni(2,2)$. (See \textit{e.g.} \cite[Section 6E]{XinWanRankinSelberg} the pullback integral. This is fixed since we fix the Archimedean datum).
\item We pair the $\Theta_{\phi_{4,i,\infty}
\prod_{v<\infty}\phi_{j_v}}$ with the $\theta^\star$ as in \cite[Lemma 6.46, Corollary 6.47]{XinWanRankinSelberg} and obtain a theta function on $1\times\Uni(2)(\mathbb{A}_\mathbb{Q})$. The corresponding Schwartz function at $\infty$ is with respect to the theta correspondence between compact unitary groups $\Uni(1)(\mathbb{R})$ and $\Uni(2)(\mathbb{R})$. By considering the central character in the above triple product integral (\ref{Fourier-Jacobi}),
the ``$\boldsymbol{\theta}_\mathsf{z}$'' part contributing non-trivially to the triple product is the eigen-component of trivial central character (the Archimedean weights of $f$ and of $\mathbf{h}_\mathsf{z}$ are dual to each other). This is a multiple of the Archimedean kernel function of the $\phi_{2,\infty}$ in \emph{loc.cit.}. This multiple is the second factor contributing to $C_{i,\infty}$.
\end{itemize}

We write $$C_\infty=\sum_i C_{i,\infty}.$$

Note that the $C_\infty$ only depends on our Archimedean datum as the $C_{i,\infty}$'s are. We see
$$C_\infty\not=0$$ as by our choices that the $w_1$-component of $l_{\theta^\star}(\mathrm{FJ}_\beta(E_{\mathrm{Kling}}))$ is nonzero for appropriate local choices at finite places for $E_{\mathrm{Kling}}$.

Now we turn to evaluate the integral (\ref{Fourier-Jacobi}) using Ichino's triple product formula. We refer to \cite[Section 8D]{XinWanRankinSelberg} for the details of Ichino's formula. As in \cite[Section 8E]{XinWanRankinSelberg} we need to pair (\ref{Fourier-Jacobi}) with the following
\begin{equation}\label{match}
\mathsf{z}(\mathbf{B}\langle\tilde{\mathbf{h}_3}\tilde{f},
\tilde{\boldsymbol{\theta}}_3\rangle),
\end{equation}
(recall the specializations of $\tilde{\mathbf{h}_3}$, $\tilde{f}$, $\tilde{\boldsymbol{\theta}}_3$ are in the dual automorphic representation space for the corresponding specializations of $\mathbf{h}$, $f$ and $\boldsymbol{\theta}$, respectively). Here by the product $\tilde{\mathbf{h}_3}\tilde{f}$ we mean the scalar valued form obtained using the natural pairing between the coefficient representations $(L_{(\frac{k-2}{2},-\frac{k-2}{2})}$ and $(L_{(-\frac{k-2}{2},\frac{k-2}{2})}$ of $\tilde{\mathbf{h}_3}$ and $\tilde{f}$ respectively. As in \cite[Sections 8D, 8E]{XinWanRankinSelberg} we appeal to Ichino's formula to evaluate the product of the two triple product integrals above (namely the product of (\ref{Fourier-Jacobi}) and (\ref{match})).
The local triple product computations are already done in \emph{loc.cit.}, except the following two cases which are different from their.\\
\noindent\underline{Archimedean triple product integral}\\
Note that at the Archimedean place, the representation $L^{k-2}$ for $\pi_f$ has dimension $k-1$. We note also the following for the vector-valued case here: for a basis $(v_1, v_2, \cdots, v_n)$ of $L^{k-2}$ and a basis $(v^\vee_1,v^\vee_2,\cdots, v^\vee_n)$ of its dual representation. If we write $L^{k-2}$-valued and $(L^{k-2})^\vee$-valued forms
$$h=h_1v_1+h_2v_2+\cdots h_nv_n$$
$$f=f_1v^\vee_1+f_2v^\vee_2+\cdots f_nv^\vee_n,$$
then the vector valued pairing $\langle h, f\rangle$ can be expressed as
$$(k-1)\int_{[\Uni(2)]}h_i((g)f_i(g)dg$$
for any $i=1,2,\cdots, k-1$. So we can choose one $i$ and do the Ichino triple product computation.
It follows easily from Schur orthogonality that the local triple product integral at the Archimedean place (see \cite[Section 8.4]{XinWanRankinSelberg}) is equal to $\frac{1}{k-1}$ (a fixed number). Thus the product of (\ref{Fourier-Jacobi}) and (\ref{match})) turns out to be some constant $C\in\bar{\mathbb{Q}}_p^\times$ times a product of several $p$-adic $L$-functions (see \cite[between Definition 8.25 to Lemma 8.26]{XinWanRankinSelberg}), which are units in $\mathcal{O}^{{\ur}}[[\Gamma_\mathcal{K}]]$ (see \cite[Proposition 8.27]{XinWanRankinSelberg}) by our choices for $\boldsymbol{\theta}$ and $\mathbf{h}$ (similar as in \cite[Section 8.2]{XinWanRankinSelberg}). Here to remove the square-free conductor assumption in \cite{XinWanRankinSelberg} for prime divisors of $N$ split in $\mathcal{K}$, we compute here a local triple product integral below for supercuspidals which is not studied in \cite[Section 8D]{XinWanRankinSelberg}.\\

\noindent\underline{Supercuspidal Triple Product Integral}
\begin{Prop}
Suppose $\pi_\ell$ is supercuspidal representation with trivial character and conductor $p^t$, $t\geq 2$ and $\varphi_\ell\in\pi_\ell$ is a new vector. Let $\tilde{\pi}_\ell$ be the contragradient representation of $\pi_\ell$ and $\tilde{\varphi}_\ell\in\tilde{\pi}_\ell$ be the new vector. Consider the matrix coefficient
\begin{equation}\nonumber
\Phi=\Phi_{\varphi_\ell, \tilde{\varphi}_\ell}(g){\eqdef}\langle\pi(g)\varphi_\ell,\tilde{\varphi}_\ell\rangle
\end{equation}
normalized such that $\Phi_{\varphi_\ell,\tilde{\varphi}_\ell}(1)=1$. Then for $g\in \mathrm{diag}(\ell^n, 1)\begin{pmatrix}1&\ell^{-n}\mathbb{Z}_\ell\\&1\end{pmatrix}K_t$, $\Phi(g)\not=0$ only when $n=0$. In that case $\Phi(g)=1$. For $g\in \mathrm{diag}(1,\ell^n)\begin{pmatrix}1&\\ \ell^{t-n}\mathbb{Z}_\ell&1\end{pmatrix}K_t$, $\Phi(g)\not=0$ only when $n=0$. In this case $\Phi(g)=1$.
\end{Prop}
This is an easy consequence of \cite[Proposition 3.1]{HuTriple}. The following corollary follows immediately from the above proposition and the computations in \cite[Section 8D, especially Lemma 8.13]{XinWanRankinSelberg}.
\begin{Cor}
Let $K_t$ be the level group defined before \cite[Lemma 8.13]{XinWanRankinSelberg}. Let $\pi_\ell$ be a supercuspidal representation of $\mathrm{GL}_2(\mathbb{Q}_\ell)$ with trivial character and conductor $\ell^t$. Let $\varphi_\ell\in\pi_\ell$ and $\tilde{\varphi}_\ell\in\tilde{\varphi}_\ell$ be as above. Let $\pi_h=\pi(\chi_{h,1},\chi_{h,2})$, $\pi_\theta(\chi_{\theta,1}, \chi_{\theta,2})$ and $\chi_{h,1}$, $\chi_{h,2}$, $\chi_{\theta,1}$, $\chi_{\theta,2}$ be characters of $\mathbb{Q}_\ell^\times$ with conductors $\ell^{t_1}$ and $t_1>t$ such that $\chi_{h,1}\chi_{\theta,1}$ and $\chi_{h,2}\chi_{\theta,2}$ are both unramified. Suppose $f_{\chi_\theta}$, $f_{\chi_h}$, $\tilde{f}_{\tilde{\chi}_\theta}$ and $\tilde{f}_{\tilde{\chi}_h}$ are as in \cite[Section 8D]{XinWanRankinSelberg}. Then Ichino's local triple production integral
$$I_\ell(\varphi_\ell\otimes f_{\chi_\theta}\otimes f_{\chi_h},\tilde{\varphi}_\ell\otimes\tilde{f}_{\tilde{\chi}_\theta}\otimes \tilde{f}_{\tilde{\chi}_h})=\mathrm{Vol}(K_{t_1}).$$
\end{Cor}

Now we arrived at the
\begin{Prop}\label{Proposition 3.12}
The
$$\mathrm{FJ}^\mathbf{h}_{\beta,\theta^\star}(E_{\mathrm{Kling}})=\langle e^{ord}l_{\theta^\star}(\mathrm{FJ}_\beta(E_{\mathrm{Kling}})),\mathbf{h}\rangle$$
is a product of an element in $\mathcal{O}^{\ur}[[\Gamma_\mathcal{K}]]^\times$ and an element in $\bar{\mathbb{Q}}_p^\times$.
\end{Prop}
\begin{proof}
We compute the ratio of the product of (\ref{Fourier-Jacobi}) and (\ref{match}) over the $p$-adic $L$-functions $\mathcal{L}_1\mathcal{L}_2\mathcal{L}_5\mathcal{L}_6$. From (\ref{Fourier-Jacobi})  and the computation of the Ichino's triple product integral (and that the $\mathbf{B}\langle\mathbf{h},\tilde{\mathbf{h}}_3\rangle$ and $\mathbf{B}\langle\mathbf{h},\tilde{\mathbf{h}}_3\rangle$ are interpolated by $\mathcal{L}_3$ and $\mathcal{L}_4$ are used here), we see that the ratio is a fixed nonzero constant times the above $C_\infty\not=0$ (see \cite[Proposition 8.29]{XinWanRankinSelberg}). Note also that $\mathcal{L}_1$, $\mathcal{L}_5$ and $\mathcal{L}_6$ are units, while $\mathcal{L}_2$ is a fixed nonzero constant. These altogether imply the $\mathrm{FJ}^\mathbf{h}_{\beta,\theta^\star}(E_{\mathrm{Kling}})$ is a divisor of a product of an element in $\mathcal{O}^{\ur}[[\Gamma_\mathcal{K}]]^\times$ and an element in $\bar{\mathbb{Q}}_p^\times$, and thus itself also is.
\end{proof}

Now we state the main theorem of this section. We start with the following definition.
\begin{Def}
We let $\mathbb{T}_\mathcal{D}$  be the reduced Hecke algebra generated by the Hecke operators at unramified primes acting on the space of the two variable family of semi-ordinary cusp forms with level group $K_\mathcal{D}$, the $U_\alpha$ operator at $p$, and then take the reduced quotient. Let the Eisenstein ideal $I_\mathcal{D}$ of $\mathbb{T}_\mathcal{D}$ be generated by $\{t-\lambda(t)\}_t$ for $t$ in the abstract Hecke algebra where $\lambda(t)$ is the Hecke eigenvalue of $t$ acting on $E_{Kling}$. Let the Eisenstein ideal $\mathcal{E}_\mathcal{D}$ be the inverse image of $I_\mathcal{D}$ in $\mathcal{O}[[\Gamma_\mathcal{K}]]\subset \mathbb{T}_\mathcal{D}$.
\end{Def}
In \cite[Lemma 9.1]{XinWanRankinSelberg}, it is proved that if $\Pcal$ is a height 1 prime ideal of $\mathcal{O}^{{\ur}}[[\Gamma_\mathcal{K}]]$ different from $(p)$, then the order of divisibility of $\mathcal{L}^{\mathrm{Gr},\Sigma}_{\mathbf{f},\mathcal{K}}$ by $\Pcal$ is less or equal to the order of divisibility of $\mathcal{E}_\mathcal{D}$ by $\Pcal$. This is generalized to our setting as follows.
\begin{Lem}
Let $\Pcal$ be a height $1$ prime of $\mathcal{O}^{{\ur}}[[\Gamma_\mathcal{K}]]$ which is not $(p)$. Then
$$\mathrm{ord}_\Pcal(\mathcal{L}^{\mathrm{Gr},\Sigma}_{\mathbf{f},\mathcal{K}})\leq \mathrm{ord}_\Pcal(\mathcal{E}_\mathcal{D}).$$
\end{Lem}
\begin{proof}
The proof is entirely similar to that \cite[Lemma 9.1]{XinWanRankinSelberg} using proposition \ref{Proposition 3.12} and theorem \ref{TheoFundamental} above in place of Theorem 3.6 and the functional constructed in section 8.5 of \cite{XinWanRankinSelberg}.
\end{proof}
\begin{Theo}\label{GBMC}
Let $f$ be a normalized cuspidal eigenform of even weight $k$, trivial character and conductor $N$ prime to $p$, such that either $\bar{\rho}_f|_{G_{\qp}}$ is absolutely irreducible, or $f$ is ordinary at $p$. Suppose for each prime divisor $q$ of $N$ non-split in $\mathcal{K}$, we have $q||N$, and that there is at least one such prime. Suppose $2$ is non-split in $\mathcal{K}$, then $2||N$. Suppose also that $\bar{\rho}_f$ is absolutely irreducible over $G_\mathcal{K}$. Then we have
$$\carac_{\mathcal{O}^{{\ur}}[[\Gamma_\mathcal{K}]]}(X_{\Kcal}^{\Gr}(f)
\otimes_{\mathcal{O}}\mathcal{O}^{{\ur}})\subseteq(\Lcal_{\Kcal}^{\Gr}(f))$$
up to height one primes which are pullbacks of primes in $\mathcal{O}[[\Gamma^+]]$. Moreover, if for each $\ell|N$ non-split in $\mathcal{K}$, we have $\ell$ is ramified in $\mathcal{K}$ and $\pi_\ell$ is the Steinberg representation twisted by $\chi_{{\ur}}$ for $\chi_{{\ur}}$ being the unramified character sending $\ell$ to $(-1)\ell^{\frac{k}{2}-1}$, then the whole containment above is true.
\end{Theo}
\begin{proof}
As in \cite[Theorem 8.2.1]{CLW}, the main theorem can be proven from Proposition \ref{Proposition 3.12} in almost the same way as the proof at the end of \cite[Section 9.2]{XinWanRankinSelberg} and \cite[Theorem 7.5]{SkinnerUrban}. First prove the $\Sigma$-primitive version: we use the lattice construction (Proposition 9.2 in \cite[Section 9.2]{XinWanRankinSelberg}) to show that $\mathcal{E}_\mathcal{D}$ contains the characteristic ideal of the dual Selmer group. The only difference is to check the condition (9) in the ``Set up'' of \emph{loc.cit.}. As in \cite[Theorem 7.5]{SkinnerUrban}, let $R$ be the four dimensional $\tilde{\Lambda}$-valued pseudo-character of $G_\mathcal{K}$ corresponding to the space of semi-ordinary cuspidal families on $\Uni(3,1)$. We suppose for contradiction our pseudo-character $R=R_1+R_2+R_3$ where $R_1$ and $R_2$ are $1$-dimensional and $R_3$ is $2$-dimensional. Then by residual irreducibility of $\bar{\rho}_f$ we can associate to $R_3$ a $2$-dimensional $\mathbb{T}_\mathcal{D}$-coefficient Galois representation. Take an arithmetic point $x$ in the absolute convergence region for Klingen Eisenstein series of sufficiently regular weight (in the sense that $a_2-a_3\gg0$ and $a_3+b_1\gg0$) and consider the specialization of the Galois representation to $x$ (the specialization of $R$ to $x$ corresponds to a classical cuspidal automorphic representation of $\Uni(3,1)$ unramified at $p$). First of all as in \cite[Theorem 7.5]{SkinnerUrban} a twist of this specialization of $R_3$ descends to a Galois representation of $G_\mathbb{Q}$ which we denote as $R_{3,x}$. We know that $R_{3,x}$ has Hodge-Tate weight $0,k-1$ and is deRham (by the corresponding property for the semi-ordinary representation $R_x=R_1+R_2+R_3$). The $R_{3,x}$ is modular by the modularity lifting result of Pan in \cite[Theorem 1.0.4]{PanFontaineMazur} which states that any odd Galois representation
$$\rho: G_\mathbb{Q}\rightarrow \mathrm{GL}_2(\mathcal{O})$$
unramified outside a finite set of primes such that $\bar{\rho}|_{G_{\mathbb{Q}_p}}$ is absolutely irreducible and de Rham with distinct Hodge-Tate weights is modular. These imply that $R_x$ is CAP (see \cite[proof of Theorem 7.5]{SkinnerUrban}), that is, it has the same system of Hecke eigenvalues as a Klingen-type Eisenstein series), contradicting the result of \cite[Theorem 2.5.6]{HarrisEisenstein}.

We first get the divisibility for $\mathcal{L}^\Sigma_{f,\mathcal{K}}$, up to height one primes which are pullbacks of height one primes of $\hat{\mathcal{O}}^{\ur}[[\Gamma^+]]$, the corresponding result for $\mathcal{L}_{f,\mathcal{K}}$ follows by putting back local Euler factors at $\Sigma$ using \cite[Proposition 2.4]{GreenbergVatsal} (note that $\mathcal{K}_\infty$ contains the cyclotomic $\mathbb{Z}_p$-extension).
Finally we use the last assumption of the theorem to apply \cite{HsiehAnticyclotomic} on the vanishing of the anticyclotomic $\mu$-invariant of the $p$-adic $L$-function in our theorem and see that it is not contained in any height one prime of $\mathcal{O}[[\Gamma^+]]$ (and thus is co-prime to the Dirichlet $p$-adic $L$-function showing up in the constant terms of the Klingen Eisenstein family).
\end{proof}

\bibliographystyle{plain}
\bibliography{olivier}
\bigskip{\footnotesize%
  \textrm{Olivier Fouquet} \par
  \textsc{Départment de Mathématiques}
  \par
  \textsc{16, route de Gray 25000 Besan\c con}
  \par
  \textsc{France} \par  
  \textit{E-mail address}: \texttt{olivier.fouquet@univ-fcomte.fr} \par
  \textit{Telephone number}: \texttt{+33622645039}\par
\par
\bigskip
  \par
   \textrm{Xin Wan}\par  
  \textsc{Morningside Center of Mathematics}\par
  \textsc{Academy of Mathematics and Systems Science, Chinese Academy of Science}
  \par
  \textsc{Beijing, 100190}
  \par
  \textsc{China}
  \par
\textit{E-mail Address}: \texttt{xwan@math.ac.cn}}

\end{document}